\documentclass[twocolumn]{macros/autart} 
\usepackage{amsmath} 
\usepackage{amssymb}  
\usepackage{graphicx}
\usepackage{latexsym}
\usepackage{verbatim}
\usepackage{ifthen}
\usepackage{psfrag}
\usepackage{macros/subfigure}

\newtheorem{nntheorem}{ Theorem}[section]
\newtheorem{nnlemma}[nntheorem]{ Lemma}
\newtheorem{nndefinition}[nntheorem]{ Definition}
\newtheorem{nncorollary}[nntheorem]{ Corollary}
\newtheorem{nnconjecture}[nntheorem]{ Conjecture}
\newtheorem{nnproposition}[nntheorem]{ Proposition}
\newtheorem{nnassumption}[nntheorem]{ Assumption}
\newtheorem{nexample}[nntheorem]{ Example}
\newtheorem{nnremark}[nntheorem]{ Remark}\newtheorem{nnapplication}[nntheorem]{ Application}
\newtheorem{nnproblem}[nntheorem]{ Exercise}

\newtheorem{nntheorem2}{ Theorem}
\newtheorem{app}[nntheorem2]{ Application}

\newenvironment{ttheorem}[1]
{\begin{nntheorem}{\rm\textrm{(#1)}}\sl}
{\end{nntheorem}}

\newenvironment{proposition}[1]
{\begin{nnproposition}{\rm\textrm{(#1)}}\sl}
{\end{nnproposition}}

\newenvironment{lemma}[1]
{\begin{nnlemma}{\rm\textrm{(#1)}}\sl}
{\end{nnlemma}}

\newenvironment{conjecture}[1]
{\begin{nnconjecture}{\rm\textrm{(#1)}}\sl}
{\end{nnconjecture}}

\newenvironment{corollary}[1]
{\begin{nncorollary}{\rm\textrm{(#1)}}\sl}
{\end{nncorollary}}

\newenvironment{definition}[1]
{\begin{nndefinition}{\rm\textrm{
\ifx\hfuzz#1\hfuzz
\else
  (#1)%
\fi
}}\sl}
{\end{nndefinition}}

\newenvironment{assumption}[1]
{\begin{nnassumption}{\rm\textrm{
\ifx\hfuzz#1\hfuzz
\else
  (#1)
\fi
}}\sl}
{\end{nnassumption}}

{\begin{app}{\rm\textrm{(#1)}}\sl}
{\end{app}}

\newenvironment{remark}[1]
{\begin{nnremark}{\rm\textrm{(#1)}}\sl}
{\end{nnremark}}

\newcommand{\eoe}
           {\hspace*{\fill}{$\vcenter{\hrule height1pt 
                     \hbox{\vrule width1pt height3pt 
            \kern3pt \vrule width1pt} \hrule height1pt}$} }

\newenvironment{example}[1]
{\begin{nexample}{\rm\textrm{(#1)}}\rm}{\eoe
\end{nexample}}

\newcommand{\eop}
           {\hspace*{\fill}{$\vcenter{\hrule height1pt 
                     \hbox{\vrule width1pt height5pt 
            \kern5pt \vrule width1pt} \hrule height1pt}$} }

\newenvironment{proofs}
{\par\noindent\textbf{Proof Sketch.}}{\eop\smallskip\vskip 3 pt}

\newcommand{\dom}{\mathop{\rm dom}\nolimits}

\newcommand{\reals}{{\mathbb R}}

\newcommand{\tto}{\;{\lower 1pt \hbox{$\rightarrow$}}\kern -12pt
           \hbox{\raise 2pt \hbox{$\rightarrow$}}\;}

\newcommand{\A}{\mathcal{A}}

\newcommand{\HS}{\mathcal{H}}
\newcommand{\R}{{\mathcal{R}}}
\newcommand{\realsgeq}{{\reals_{\geq 0}}}
\newcommand{\nats}{\mathbb{N}}      

\newcommand{\J}{{\mathcal{J}}}

\usepackage{svg}
\usepackage{array}
\usepackage{mathtools}
\usepackage{eqparbox}
\usepackage{dsfont}
\usepackage{pstool}
\usepackage{rotating,graphics,psfrag,epsfig}
\usepackage{macros/scalerel}
\usepackage{balance}
\usepackage[normalem]{ulem}
\allowdisplaybreaks

\usepackage{amsmath}
\usepackage[colorinlistoftodos,prependcaption,textsize=tiny]{todonotes}
\usepackage{varwidth}
\usepackage{parskip}
\usepackage[inline, shortlabels]{enumitem} 
\usepackage{float}

\usepackage{multicol}
\usepackage{lipsum}
\usepackage{nccmath}
\usepackage{soul}
\usepackage{cancel}

\usepackage{booktabs,arydshln,xcolor}
\newcolumntype{L}{>{\centering\arraybackslash}m{2.1cm}}
\newcolumntype{V}{>{\centering\arraybackslash}m{1.6cm}}
   
\usepackage{amsfonts}              
\usepackage{mathrsfs}
\usepackage{indentfirst}   
\usepackage{setspace}
\usepackage{wrapfig}
\newcommand{\EatDot}[1]{}

\DeclareMathOperator*{\argmin}{arg\,min}
\newcommand{\U}{\mathcal{U}}

\newcommand{\pn}[1]{{\color{black} #1}}
\newcommand{\pnn}[1]{{\color{black} #1}}
\newcommand{\pno}[1]{{\color{black} #1}}
\newcommand{\sj}[1]{{\color{black} #1}}

\newenvironment{proof}
{\par\noindent\textbf{Proof.}}{\eop\smallskip\vskip 3 pt}

\usepackage{ifthen}
\newboolean{Conference}
\setboolean{Conference}{true}
\newcommand{\NotConf}[1]{\ifthenelse{\boolean{Conference}}{}{{\color{black}#1}}}
\newcommand{\IfConf}[2]{\ifthenelse{\boolean{Conference}}{#1}{{\color{black}#2}}}   

\newboolean{Inclusion}
\setboolean{Inclusion}{false}
\newcommand{\IfInc}[2]{\ifthenelse{\boolean{Inclusion}}{ \color{cyan}#1\color{black}#2}{#2}}
\newcommand{\IfIncd}[2]{\ifthenelse{\boolean{Inclusion}}{{\color{purple}#1}}{{#2}}}

\newboolean{Personal}
\setboolean{Personal}{false}
\newcommand{\IfPers}[1]{\ifthenelse{\boolean{Personal}}{{\color{purple}#1} }{}}

\newboolean{Infinite-horizon} 
\setboolean{Infinite-horizon}{true} 
\newcommand{\IfIh}[2]{\ifthenelse{\boolean{Infinite-horizon}}{#1}{{\color{purple}#2}}}

\newboolean{Two-Players} 
\setboolean{Two-Players}{true} 
\newcommand{\IfTp}[2]{\ifthenelse{\boolean{Infinite-horizon}}{#1}{{\color{purple}#2}}}

\newboolean{Journal}
\setboolean{Journal}{true}
\newcommand{\IfJ}[1]{\ifthenelse{\boolean{Journal}}{{\color{purple}#1} }{}}

\newboolean{Automatica}
\setboolean{Automatica}{true}
\newcommand{\NotAutom}[1]{\ifthenelse{\boolean{Automatica}}{}{{\color{blue}#1}}}
\newcommand{\IfAutom}[2]{\ifthenelse{\boolean{Automatica}}{#1}{{\color{blue}#2}}}   

\newboolean{Automatica2s}
\setboolean{Automatica2s}{false}
\newcommand{\NotAutomss}[1]{\ifthenelse{\boolean{Automatica2s}}{}{\color{black}#1\color{black}}}
\newcommand{\IfAutomss}[2]{\ifthenelse{\boolean{Automatica2s}}{{\color{black}#1}}{\color{black}#2\color{black}}}  

\newboolean{Comments}
\setboolean{Comments}{false}

\begin{document}
\begin{frontmatter} 

\title{Two-Player 
Zero-Sum Hybrid Games
} 

\author{Santiago J. Leudo}\ead{sjimen28@ucsc.edu} \>and    
\author{Ricardo G. Sanfelice}\ead{ricardo@ucsc.edu}               
\address{Electrical and Computer Engineering Department, University of California, Santa Cruz, CA 95064, USA.}  
          
                                          
\IfIncd{
\IfConf{
 \title{\bf Two-Player Infinite-Horizon Zero-Sum Set-Valued Hybrid Games}
\author{Santiago J. Leudo and Ricardo G. Sanfelice
\thanks{S. Jimenez Leudo and R. G. Sanfelice are with Electrical and Computer Engineering Department, University of California, Santa Cruz, CA 95064.
     Email: {\tt\small \{sjimen28, ricardo\}@ucsc.edu}}
  }      
\maketitle
\thispagestyle{plain}
\pagestyle{plain}
}
{ 
\ititle{
 \title{\bf Two-Player Infinite-Horizon Zero-Sum Set-Valued Hybrid Games}
\iauthor{
  Santiago Jimenez Leudo\\
   {\normalsize sjimen28@ucsc.edu}\\
  Ricardo Sanfelice\\
  {\normalsize ricardo@ucsc.edu}}
\idate{\today{}} 
\iyear{2020}
\irefnr{006. {\bf Status}: NOT PUBLISHED. {\bf Readers of this material have the responsibility to inform all the authors promptly if they wish to reuse, modify, correct, publish, or distribute any portion of this report.}}
\makeititle

 \setcounter{page}{1}
\title{Noncooperative Infinite Horizon Hybrid Games}
\author{Santiago J. Leudo\\
Ricardo Sanfelice}
\maketitle
}
}
}{}
\begin{abstract}
In this paper, we formulate a two-player zero-sum game under dynamic constraints defined by hybrid dynamical \IfIncd{inclusions}{equations}. 
The game consists of a min-max problem involving a cost functional that depends on the actions and resulting solutions {to the hybrid system}, defined as functions of hybrid time and, hence, can flow or jump.
A terminal set conveniently defined allows to recast both finite and infinite horizon problems.
We present sufficient conditions given in terms of Hamilton–Jacobi–Bellman-Isaacs-like equations to \IfIncd{evaluate and upper bound the largest cost over the set of adversarial strategies. 
For a specific type of games, we provide pointwise conditions with similar structure to }{}guarantee to attain a solution to the game. It is shown that when the players select the optimal strategy, the \IfIncd{cost}{value function} can be evaluated without computing solutions {to the hybrid system}. Under additional conditions, we show that the optimal state-feedback laws render a set of interest asymptotically stable for the resulting hybrid closed-loop system. 
\IfIncd{}{\sj{Applications of these games, 
presented here as robust control problems, include disturbance rejection 
and security problems.}}
%
%
\end{abstract}
\end{frontmatter}
\NotAutom{
\IfConf{
\maketitle
}{}
} 
\NotConf
{
\newpage
\tableofcontents
\newpage
}
\pnn{
\section{Introduction}
}
\NotConf{\subsection{Background}}
\setlength{\parindent}{3ex}

\subsection{Background}

A game is an optimization problem with multiple decision makers (players), a set of constraints (potentially dynamic) that enforces the ``rules'' of the game, and a set of payoff functions to be optimized by selecting decision variables. Constraints on the {state and} decision variables formulated as dynamic relationships lead to {dynamic games}; see \cite{basar1999dynamic} and the references therein. 
In this setting, an interesting scenario arises when the players have different incentives, e.g., to minimize or maximize \sj{their own} cost function. Dynamic noncooperative games focus on the case in which the players select their actions \sj{with no coalition formation, such that an individual benefit 
does not necessarily imply a benefit
to the other players} \NotAutomss{\cite{varaiya1969existence,Basar2017},}\cite{fudenberg1991game,Owen,isaacs1999differential,Hespanha09}.
This type of dynamic game has been thoroughly studied in the literature when the dynamic constraints are given in terms of difference equations or differential equations, known as differential games.

Challenges arise when the players' dynamics exhibit both continuous and discrete behaviors, {for instance}, due to 
intermittent information availability, resets, timers in the control algorithms that expire, or non-smooth mechanical behaviors exhibiting impacts, among others. Hybrid {dynamical systems}  conveniently capture this kind of behavior \cite{65,220}.
{Under certain assumptions, differential algebraic equation (DAEs) \cite{kunkel2006differential} -- also known as descriptor systems -- can be recast as hybrid equations, see \cite[Lemma 2]{nanez2017invariance}. \sj{Specifically}, when the initial condition to a DAE is \emph{consistent} and the data pair of the system is \emph{regular} (for each subsystem, in the case of switched DAEs \cite{LIBERZON2012954}), a solution to the DAE is also a solution to the equivalent hybrid system defined as in \cite{65}. }
However, when designing 
{algorithms that make }optimal choices of the decision variables under constraints given by hybrid dynamics, relying {only on continuous-time or discrete-time} approaches potentially results in suboptimal solutions. Unfortunately, tools for the design of algorithms for games with such hybrid dynamic constraints, which we refer to as hybrid games, are not fully developed.

Particular classes of dynamic games involving hybrid dynamic constraints have been recently studied in the literature.
A game-theory-based control design approach is presented for timed automata 
in \cite{asarin1995symbolic,asarin1998controller}, for hybrid finite-state automata in \cite{Henzinger1999RectangularHG,tomlin2000game}, {and for o-minimal hybrid systems in \cite{bouyer2006control}. In these articles,} the specifications to be guaranteed by the system are defined in terms of temporal logic formulae. When the payoff is defined in terms of a terminal cost, {such} approach allows designing reachability-based controllers through the satisfaction of Hamilton-Jacobi conditions to certify safety of hybrid finite-state automata \cite{ding2011toward}.
{Following an approach that allows for richer dynamics, \cite{vladimerou2009stormed} studies} a class of reachability games between a controller and the environment, {under constraints defined by} hybrid automata {(STORMED games)} {for which} at each decision step, the players can choose either to have their variables evolve continuously or discretely, {following predefined rules.}
For {continuous-time systems with state resets,} tools for the computation of the {region of attraction} for hybrid limit cycles under the presence of disturbances are provided in \cite{choi2022computation}, {where the inputs only affect the flow.}
%
%

%
%
Efforts pertaining to differential games with impulsive elements include
\cite{platzer2017differential,platzer2015differential}, 
{where} the interaction between the players is modeled similarly to switched systems,
\cite{dharmatti2006zero}{, which} establishes continuity of bounds on value functions and viscosity solutions, \cite{gromov2017class}{, which} formulates necessary and sufficient conditions for optimality in bimodal linear-quadratic differential games, and \cite{cacace2020stochastic}{, which} studies a class of stochastic two-player differential games in match race problems.
\IfPers{
Games involving multiple players with potentially different interests emerge in multiagent systems, both in benign (or cooperative) and contested (or noncooperative) settings.  A list of examples includes and is not limited to route selection in a road network \cite{Paccagnan}, heavy-duty vehicle platooning \cite{gattami2011establishing}, control of smart grids \cite{tatarenko2014game}, trading modeling in the stock market \cite{firoozi2017mean}, and control of large populations of systems \cite{grammatico2015decentralized}.  Generally speaking, a game is an optimization problem with multiple players, constraints that enforce the ``rules'' of the game, and payoff functions to be optimized through the selection of decision variables. Constraints on the actions played by the players formulated as dynamic relationships (i.e., involving time) lead to 
{dynamic games}. Differential games pertain to the case when these constraints are given in terms of differential equations; see, e.g., \cite{basar1999dynamic} and the references therein.  Of particular interest is the contested setting, which occurs when the players have independent objectives, such as when one player aims at minimizing a cost function and another player aims at maximizing it under dynamic constraints. If the players select their actions seeking their own benefit, a dynamic noncooperative game emerges. This type of dynamic games have been thoroughly studied in the literature, when the dynamic constraints are given in terms of difference equations or differential equations --  in general, referred to as differential games -- including, to just list a few, \cite{varaiya1969existence,fudenberg1991game,Owen,isaacs1999differential,Basar2017,Hespanha09}. 

In recent years, significant progress has been made in the understanding of dynamic games with players sharing information over networks; see, e.g., \cite{marden2015game}. 
Interestingly, the combination of physics, computing, and networks 
leads to dynamic constraints that exhibit both continuous and discrete behavior.
In particular, intermittent information availability, resets of variables, such as expiring timers, and other nonsmooth and instantaneous changes lead to dynamic constraints that can be conveniently captured using hybrid \IfIncd{inclusions}{equations}.
The design of algorithms that guarantee optimality under such hybrid dynamic constraints
requires new tools, since using tools from the differential games' literature would most likely lead to  suboptimal solutions.  Unfortunately, tools for the design of algorithms for games with such hybrid dynamic constraints, which we refer to as \IfIncd{set-valued}{} 
hybrid games, are not as developed as those for differential games, as cited above.  
In \cite{Henzinger1999RectangularHG,tomlin2000game}, a control design approach 
using game theory that is applicable to a class of hybrid automata models is presented.  Specifically, the models considered therein are based on finite-state automata, the specifications are defined in terms of  temporal logic formulae, and the payoff is solely given by a terminal cost. 
Decidability for hybrid automata given a winning condition are studied in \cite{vladimerou2009stormed}.
Applications of the approach in \cite{tomlin2000game} include reachability-based controller design \cite{ding2011toward,fisac2015reach}. 
The work in \cite{platzer2017differential,platzer2015differential} pertains to a class of dynamic games 
in which the evolution of the variables associated to each of the players is modeled using differential equations, while the interactions between the players is modeled as switches that occur at isolated time instances, similar to switched systems.
Conditional viability for  impulsive systems with two competing input actions was considered in \cite{aubin2005conditional} and treated as an evolutionary game.
Other efforts pertaining to differential games with impulsive (or discontinuous) elements include establishing continuity of bounds on value functions and (viscosity) solutions \cite{dharmatti2006zero}, formulating necessary and sufficient conditions for optimal strategies for the special case of bimodal linear-quadratic differential games \cite{gromov2017class}, and a class of stochastic two-player differential games in the context of sail boat competitions \cite{cacace2020stochastic}.}

%
\subsection{Contribution and Applications}
{Motivated by the lack of tools for the design of algorithms for {general} hybrid games\IfPers{ with dynamic constraints that are richer than those allowed by finite-state automata and switched systems}, we formulate a framework for the study of two-player zero-sum games with hybrid dynamic constraints.}
{Specifically, we optimize a cost functional, 
which includes}
\begin{itemize}
  \item a stage cost that penalizes the evolution of the state and the input {during} flow, 
  \item a stage cost that penalizes the evolution {of the state and the input at} jumps, and
  \item a terminal cost to penalize the final value of the variables.
\end{itemize}
Following the framework in \cite{65,220}, we model the hybrid dynamic constraints as a hybrid {dynamical} system, 
{denoted} $\HS\IfIncd{_s}{}$ {and} given in terms of the hybrid \IfIncd{inclusion}{equation}
\begin{equation}
\HS\IfIncd{_s}{}  \left\{ 
\begin{matrix}
\dot{x} &\IfIncd{\in}{=}& F(x,u_{C1},u_{C2}) \quad  (x,u_{C1},u_{C2}) \in C \\ 
 x^+ &\IfIncd{\in}{=}& G(x,u_{D1},u_{D2}) \quad (x,u_{D1},u_{D2}) \in D 
\end{matrix}
\right. 
\IfIncd{
\label{Hinc}
}
{\label{Heq}}
\end{equation}
where $x \in \mathbb{R}^n$ is the state, $(u_{C1},u_{D1}) \in \mathbb{R}^{m_{C1}} \times \mathbb{R}^{m_{D1}}$ is the input chosen by player $P_1$, $(u_{C2},u_{D2}) \in \mathbb{R}^{m_{C2}} \times \mathbb{R}^{m_{D2}}$ is the input chosen by player $P_2$, and the data of $(C,F,D,G)$ is given as follows\footnote{Here, {$m_C=m_{C1}+m_{C2}$ and $m_D=m_{D1}+m_{D2}$.}}: 
\begin{itemize}
  \item The $\mathit{flow}$ $\mathit{map}$ $F: \mathbb{R}^n \times \mathbb{R}^{m_C}  \IfIncd{\rightrightarrows}{\rightarrow} \mathbb{R}^n$ captures the continuous evolution of the system when the state {and the input are} in the $\mathit{flow}$ $\mathit{set}$ $C {\subset \reals^n \times \reals^{m_{C1}+m_{C2}}}$. 
  \item The $\mathit{jump}$ $\mathit{map}$ ${G: \mathbb{R}^n  \times \mathbb{R}^{m_D} \IfIncd{\rightrightarrows}{\rightarrow}  \mathbb{R}^n}$ describes the discrete evolution of the system when the state {and the input are} in the $\mathit{jump}$ $\mathit{set}$ $D{\subset \reals^n \times \reals^{m_{D1}+m_{D2}}}$. 
\end{itemize}
 For such broad class of systems, 
 we consider a \sj{Bolza-form} cost functional {$\mathcal{J}$} 
 associated to the solution to $\HS\IfIncd{_s}{}$ from $\xi$ and study {a zero-sum two-player hybrid game that, informally, is given as}
 \begin{equation}
{ \min_{(u_{C1},u_{D1})} \max_{(u_{C2},u_{D2})}}\> \mathcal{J}(\xi, u_{C1},u_{C2},u_{D1},u_{D2})
\label{eq:minmaxJinformal}
\end{equation}

\sj{This game captures the dynamics of systems operating in {contested} scenarios with hybrid dynamics, such as continuous-time dynamics with logical modes,
multiple modes of operation, and dynamics or control signals that change abruptly or impulsively. Such dynamics can be represented by}
switching systems, hybrid automata, or impulsive differential equations, \sj{all of which can be modeled as in \eqref{Heq}; see \cite{65,34}.}
  
  \IfIncd{}{Several applications \sj{lead to the game in \eqref{eq:minmaxJinformal}. The following robust control scenario is of particular interest and concrete instances are considered in this paper:}

  \sj{\textbf{$(\star)$}} Given the system $\HS\IfIncd{_s}{}$ as in (\IfIncd{\ref{Hinc}}{\ref{Heq}}) with state $x$, the {robust control} 
  problem consists of establishing conditions such that player $P_1$ selects a control input $(u_{C1},u_{D1})$ that {minimizes} {a cost functional $\J$ until the game ends, which occurs when the state enters a set $X$,} in the presence of a disturbance $(u_{C2},u_{D2})$ chosen by $P_2$ \sj{(modeled as maximizing the cost)}. 

 
The solution of the game formulated in this paper, known as a {\emph{saddle-point equilibrium}}, is given in terms of the {actions 
 of} the players. Informally, {when a player unilaterally deviates from the equilibrium action, it does not improve its individual outcome.} Thus, by formulating the applications above as two-player zero-sum hybrid games, we can synthesize the saddle-point equilibrium and determine the control action that minimizes the cost $\J$ for the maximizing adversarial action. 
}
The main contributions of this paper are summarized as follows.

\begin{itemize}
  \item In Section III, we present a framework for the study of two-player zero-sum games with
   hybrid dynamic constraints.
  \IfIncd{\item In Theorem \ref{thHJBszihinc}, we present sufficient conditions to evaluate and upper bound the largest cost over the set of adversarial strategies when player $P_1$ plays optimally.
  }
  {\item We present in Theorem \ref{thHJBszih} sufficient conditions based on Hamilton–Jacobi–Bellman-Isaacs-like equations to {design} a saddle-point equilibrium and evaluate the game value function without {computing} solutions {to the hybrid system}.}
  \IfIncd{}{
  \item Connections between optimality and asymptotic stability of a 
  {set} are {proposed} in Section V 
  and framed in the game theoretical approach employed.
  \item We present in Section VI 
  applications to \sj{robust 
  control scenarios} by formulating and solving them as two-player zero-sum hybrid dynamic games.} 
  \end{itemize}

This work extends our preliminary conference paper \cite{leudohygames} 
 \IfIncd{to the case of dynamic \pnn{constraints} given in terms of hybrid \pnn{systems with} set-valued \pnn{flow and jump} maps. 
We present a less conservative approach by relaxing the assumption on uniqueness of solutions to the hybrid system that define the dynamics of the game. 
Although, in general, \pnn{it might not be possible to construct} a saddle point equilibrium solution to the game when \pnn{the dynamics admit nonunique solutions} (due to the game being ill-defined), an upper \pnn{value function is provided.}}
since a more general problem is considered. By including a terminal set,
we formulate problems with variable terminal time and with infinite horizon.
Sufficient conditions to solve the general problem are provided in Theorem 4.1. with the corresponding proof, that uses as a preliminary step Proposition
4.2. Remarks 4.5 and 4.6 provide a discussion on the existence \sj{and computation} of a value function. 
Theorem 5.4 provides connections between optimality and stability for the general problem plus considering special cases that
allow to guarantee pre-asymptotic stability. Proposition~6.1  allows to solve a robust LQR game with aperiodic jumps. Complete proofs of the results are provided.
\sj{We provide examples of robust control applications, namely, 
a disturbance rejection problem with periodic jumps and a security problem 
for a bouncing ball system under attacks.
}

\sj{Related to \eqref{eq:minmaxJinformal} are the} zero-sum games for DAEs studied {in the literature \cite{Gardner1977}. A min-max principle built upon {Pontryagin's Maximum Principle} is provided in \cite{wu1992class}. {Linear dynamics and quadratic costs} result in coupled Riccati differential equations, and conditions {for their solvability} are provided in \cite{xu1994isaacs} and  \cite{engwerda2012feedback}. In \cite{tanwani2019feedback}, noncooperative games for Markov switching DAEs are studied and Hamilton–Jacobi–Bellman-Isaacs equations are derived. {When the initial condition
  to a switching DAE is consistent and the data pair of the system
  is regular for each subsystem,} a deterministic version of the problem solved in \cite{tanwani2019feedback} can be recast as \eqref{eq:minmaxJinformal}.}

In recent {works},
 optimality {for} hybrid systems {modeled as in (\IfIncd{\ref{Hinc}}{\ref{Heq}})} is certified via Lyapunov-like conditions \cite{ferrante2019certifying}, providing cost evaluation results for the case in which the data is given in terms of set-valued maps.
 The work in \cite{GOEBEL2019153} provides sufficient conditions to guarantee the existence of optimal solutions. 
 A receding-horizon algorithm to implement these ideas is presented in \cite{193c}.
 Cost evaluation results and conditions to guarantee asymptotic stability of a set of interest are established for a discrete-time system under 
 adversarial scenarios in \cite{232}. 
 A fixed finite-horizon case of the hybrid game formulated in this paper is studied in \cite{leudo2022FH}. The conditions on the optimization problem formulated therein 
 are similar to their counterparts in the differential/dynamic game theory literature. Nevertheless, in contrast to this work, 
 the end of the game therein is attained when the time of solutions to $\HS$ reach
  a terminal set $\mathcal{T}$. 
 To account for hybrid time domains, which are introduced in Section~\ref{sec:preliminaries}, a hybrid time domain-like geometry is assumed for $\mathcal{T}$ as in \cite{193c}.
 This results in optimality conditions in terms of PDEs, 
 and the optimal feedback laws are not stationary.

\IfAutomss{
Due to space constrains, a technical report is provided in \cite{TRHGArxiv}, with proofs and details omitted \sj{here}.
}{}
\subsection{Notation} Given two vectors $x,y$, we use the equivalent notation $(x,y)=[x^\top y^\top]^\top$. The symbol $\mathbb{N}$ denotes the set of natural numbers including zero. The symbol $\mathbb{R}$ denotes the set of real numbers and $\mathbb{R}_{\geq 0}$ denotes the set of nonnegative reals. Given a vector $x$ and a nonempty set $\mathcal{A}$, the distance from $x$ to $\mathcal{A}$ is defined as $\left| x \right|_\mathcal{A}=\inf_{y \in \mathcal{A}}\left | x-y \right |$. 
{We} denote with $\mathbb{S}_+^n$ the set of real positive definite matrices of dimension $n$, \IfPers{ i.e., $A \in \mathbb{S}_+^n$ if  $A \in \mathbb{R}^{n \times n}$ $A=A^\top$ and there exists a real nonsingular matrix $M$ such that $A=M M^\top$}
and with $\mathbb{S}_{0+}^n$ the set of real positive semidefinite matrices of dimension $n$. Given a nonempty set $C$, we denote by  $\textup{int} \> C$ its interior and by $\overline{C}$ its closure. 
The $n$-dimensional identity matrix is denoted by $I_n$. 
\NotAutomss{Given a symmetric matrix $A\in \reals^{n\times n}$, the scalars $\underline{\lambda}(A)$ and $\overline{\lambda}(A)$ denote the minimum and maximum eigenvalue of $A$, respectively.}

\section{Preliminaries}
\label{sec:preliminaries}
\IfTp{}{\subsection{Multiagent network}
The network in which the interactions between players occur is modeled by the directed graph $\Gamma=(\IfTp{\{1,2\}}{\mathcal{V}},\mathcal{E},\mathcal{G})$, with  node set $\IfTp{\{1,2\}}{\mathcal{V}}=\{1,2,\dots,N\}$. The set of edges $\mathcal{E} \subset \IfTp{\{1,2\}}{\mathcal{V}}\times \IfTp{\{1,2\}}{\mathcal{V}}$ models the links between nodes, i.e., every pair $(i,k) \in \mathcal{E}$ implies there exists information sharing from the player $i$ to the player $k$. The $i,k$ entry of the adjacency matrix, $\mathcal{G}$, namely $\mathcal{G}_
{ik}$, is equal to 1 if $(i, k) \in \mathcal{E}$, and takes a zero value otherwise.  The set of indices corresponding to the neighbors that can send
information to the $i$-th player is denoted by $\mathcal{N}(i) := \{k \in \IfTp{\{1,2\}}{\mathcal{V}} : (k, i) \in \mathcal{E}\}$. The number of players is then $N$ and the number of communication links is $|\mathcal{E}|$. 
We consider the case in which there is a limited amount of communication channels. }

\subsection{Hybrid Systems with Inputs}

Since solutions to the dynamical system {$\HS\IfIncd{_s}{} $ as in \IfIncd{\eqref{Hinc}}{\eqref{Heq}}}  can exhibit both continuous and discrete behavior, we use ordinary time $t$ to determine the amount of flow, and a counter $j \in \mathbb{N}$ that counts the number of jumps. 
Based on this \sj{parametrization}, the concept of a hybrid time domain, in which solutions are \sj{defined}, \IfAutomss{\sj{
of a hybrid arc, and of a hybrid input are defined as in \cite[Section 2.2]{65}.}}{is introduced.
\begin{definition}{Hybrid time domain}
A set $E \subset \mathbb{R}_{\geq0} \times \mathbb{N}$ is a hybrid time domain if, for each $(T,J)\in E$, the set
$    E\cap ([0,T]\times \{0,1,\dots,J\}) $
is a compact hybrid time domain, i.e., it can be written in the form
   $ \bigcup_{j=0}^J ([t_j,t_{j+1}]\times\{j\})$
for some finite nondecreasing {sequence} 
$\{t_j\}^{J+1}_{j=0}$ with $t_{J+1}=T$.  
Each element $(t,j)\in E$ denotes the elapsed hybrid time, which indicates that $t$ seconds of flow time and $j$ jumps have occurred. 
\end{definition}
}
\sj{A hybrid signal is a function defined on a hybrid time domain. Given a hybrid signal $\phi$ and $j \in \mathbb{N}$, we define $I^j_{\phi}\pnn{:}=\{t:(t,j)\in \dom \phi\}$.}
\NotAutomss{\begin{definition}{Hybrid arc}
A hybrid signal $\phi:\dom \phi \rightarrow \mathbb{R}^n$ is called 
a hybrid arc if, for each $j \in \mathbb{N}$, the function $t \mapsto \phi(t,j)$ is locally absolutely continuous {on} 
$I^j_{\phi}$. A hybrid arc $\phi$ is compact if $\dom \phi$ is compact. 
\label{htd}
\end{definition}

\begin{definition}{Hybrid Input} A hybrid signal $u$ is a hybrid input if for each $j \in \mathbb{N}$, the function $t \mapsto u(t,j)$ is Lebesgue measurable and locally essentially bounded on the interval $I^j_u$.  
  \label{Def:HyInput}
\end{definition}
}
Let $\mathcal{X}$  be the set of hybrid arcs $\phi: \textup{dom}\> \phi \rightarrow \mathbb{R}^n$ and  $\mathcal{U}=\mathcal{U}_C \times \mathcal{U}_D$ the set of hybrid inputs $u=(u_{C},u_{D}): \textup{dom}\>u \rightarrow \mathbb{R}^{m_C} \times \mathbb{R}^{m_D}$, where $u_C=(u_{C1},u_{C2}\IfTp{}{,\dots, u_{C,N}})$, \IfTp{$m_{C1}+m_{C2}=m_C$}{}, $u_D=(u_{D1},u_{D2}\IfTp{}{,\dots, u_{D,N}})$, and \IfTp{$m_{D1}+m_{D2}=m_D$}{}. {A solution to the hybrid system $\HS\IfIncd{_s}{} $ with input is defined as follows.}

\begin{definition}{Solution to $\HS\IfIncd{_s}{}$}
A hybrid signal $(\phi,u)$ {defines a solution pair to (\IfIncd{\ref{Hinc}}{\ref{Heq}})} 
  if $\phi \in \mathcal{X}$, $u=(u_C,u_D) \in \mathcal{U}$,   $\textup{dom}\phi=\textup{dom}u$, and 
  \begin{itemize}
    \item $(\phi(0,0),u_C(0,0)) \in \overline{C}$ or $(\phi(0,0),u_D(0,0)) \in D$,
    \item For each $j\in \mathbb{N}$ such that $I^j_{\phi}$ has a nonempty \NotAutomss{interior }$\textup{int} I^j_{\phi}$, we have, for all
    $t \in \textup{int} I^j_{\phi}$,
    {$(\phi(t,j),u_C(t,j)) \in C
    $}
    and, for almost all $t \in I^j_{\phi}$,
    \IfAutomss{\sj{$\frac{d}{dt}\phi(t,j)\IfIncd{\in}{=} F(\phi(t,j), u_C(t,j))$,}}{\begin{equation*}
        \frac{d}{dt}\phi(t,j)\IfIncd{\in}{=} F(\phi(t,j), u_C(t,j))
    \end{equation*}}
    \item For all $(t,j)\in \dom \phi$ such that $(t,j+1)\in \dom \phi$,  
    \IfAutomss{\sj{$(\phi(t,j),u_D(t,j)) {\in} D,$ $
	\phi(t,j+1) \IfIncd{\in}{=} G(\phi(t,j),u_D(t,j)) $.}}{\begin{eqnarray*}
	(\phi(t,j),u_D(t,j)) &\in& D \\
	\phi(t,j+1) &\IfIncd{\in}{=}& G(\phi(t,j),u_D(t,j)) 
\end{eqnarray*}}
\end{itemize}
A solution pair $(\phi, u)$ is a compact solution pair if $\phi$ is a compact hybrid arc\IfAutomss{\>\cite[Definition 2.5]{65}}{; see Definition \ref{htd}}.
\label{SolutiontocalH}
\end{definition}

Given a solution pair $(\phi,u)$, the component $\phi$ is referred to as the state trajectory.
In this article, the same symbols are used to denote input actions and their values. The context clarifies the meaning of $u$, as follows: ``the function $u$\pnn{,}'' ``the signal $u$\pnn{,}'' or ``the hybrid signal $u$'' that appears in ``the solution pair $(\phi,u)$'' refer to the input action, whereas  ``$u$'' refers to the input value as a point in $\mathbb{R}^{m_C} \times \mathbb{R}^{m_D}$ in any other case. \NotAutomss{The reader can replace ``the function $u$'' by ``$u_\phi$'', which is the input action yielding the state trajectory $\phi$.}

{A solution pair $(\phi,u)$ to $\HS\IfIncd{_s}{}$ from $\xi \in \reals^n$ is complete if $\dom (\phi,u)$ is unbounded. It is maximal if there is no solution $(\psi,w)$ from $\xi$ such that $\phi(t,j)=\psi(t,j)$ and $u(t,j)=w(t,j)$ for all $(t,j)\in \dom (\phi,u)$ and $\dom (\phi,u)$ is a proper subset of $\dom (\psi,w)$.} 
 We denote by $\hat{\mathcal{S}}_{\HS\IfIncd{_s}{}}(M)$ the set of solution pairs $(\phi,u)$ to $\HS\IfIncd{_s}{}$ as in (\IfIncd{\ref{Hinc}}{\ref{Heq}}) such that $\phi(0,0) \in M$. The set $\mathcal{S}_{\HS\IfIncd{_s}{}}(M) \subset \hat{\mathcal{S}}_{\HS\IfIncd{_s}{}}(M)$ denotes all maximal solution pairs from $M$. 
 {Given $\xi \in \mathbb{R}^n$, we define the set of input actions that yield maximal solutions to $\HS\IfIncd{_s}{}$ from $\xi$ 
 as 
	  $\mathcal{U}_{\HS\IfIncd{_s}{}}
    (\xi):= \{u :\exists(\phi,u) \in \hat{\mathcal{S}}_{\HS\IfIncd{_s}{}}(\xi)
  \}$.} 
  \IfPers{For a given $u \in \mathcal{U}$, we denote the set of maximal state trajectories to $\HS\IfIncd{_s}{} $ from $\xi$ for $u$ by $\mathcal{R}(\xi, u)=\{ \phi :(\phi,u) \in \mathcal{S}_\HS (\xi) \}$. We say $u$ renders a maximal trajectory $\phi$ to $\HS\IfIncd{_s}{} $ from $\xi$ if $\phi \in \mathcal{R}(\xi, u)$.}
For a given $u \in \mathcal{U}$, we denote the set of maximal state trajectories to $\HS\IfIncd{_s}{}$ from $\xi$ for $u$ by $\mathcal{R}(\xi, u)=\{ \phi :(\phi,u) \in \mathcal{S}_{\HS\IfIncd{_s}{}} (\xi) \}$. We say $u$ renders a maximal trajectory $\phi$ to $\HS\IfIncd{_s}{}$ from $\xi$ if $\phi \in \mathcal{R}(\xi, u)$.
 {A complete solution $(\phi,u)$ is discrete if $\dom (\phi,u) \subset \{0\} \times \nats$ and continuous if $\dom (\phi,u)\subset \realsgeq \times \{0\} $.}
  
We define the projections of $C \subseteq \reals^n \times \reals^{m_C}$ and $D \subseteq \reals^n \times \reals^{m_D}$ onto $\mathbb{R}^n$, respectively, as \pno{
$	\Pi(C):=\{ \xi \in  \mathbb{R}^n : \exists u_C \in \mathbb{R}^{m_C}  \text{ s.t. } 		(\xi,u_C) \in C \}
$ and 
$	\Pi(D):=\{ \xi \in  \mathbb{R}^n : \exists u_D \in \mathbb{R}^{m_D}   \text{ s.t. } 		(\xi,u_D) \in D \}. 
$}
We also define the set-valued maps that output the allowed input values \pno{at a given state $x$} as
\pno{
$ \Pi_u^C(x)=\{ u_C \in  \mathbb{R}^{m_C} : 
		(x,u_C) \in C \},
$ and 
$	\Pi_u^D(x)=\{ u_D \in  \mathbb{R}^{m_D} : 
		(x,u_D) \in D \}.
$} 
Moreover, {$\sup_t \dom \phi :=\sup\{ t \in \mathbb{R}_{\geq 0}: \exists j 
\>\text{s.t.} \> (t,j) \in \dom \phi\}$, 
$\sup_j \dom \phi :=\sup\{  j \in \mathbb{N}: \exists t 
\>\text{s.t.} \> (t,j) \in \dom \phi\}$}, and 
$\sup \dom \phi:= (\sup_t \dom \phi, \sup_j \dom \phi)$. Whenever $\dom \phi$ is compact, $\dom \phi \supset \max \dom \phi := \sup \dom \phi$.


\NotAutomss{
The following conditions guarantee uniqueness of solutions to $\HS\IfIncd{_s}{}$ as in (\IfIncd{\ref{Hinc}}{\ref{Heq}}) \cite[Proposition 2.11]{65}.
\begin{proposition}{Uniqueness of Solutions}
 Consider the hybrid system $\HS\IfIncd{_s}{} $ as in (\IfIncd{\ref{Hinc}}{\ref{Heq}}). For every $\xi \in \Pi(\overline{C}) \cup \Pi(D)$ and each $u \in \mathcal{U}$  there exists a unique maximal {solution $(\phi,u)$} with $\phi(0,0)=\xi$ provided that {the following holds:}
 \begin{enumerate}[label=\arabic*)]
\item   for every $\xi \in \Pi(\overline{C})\setminus \Pi(D)$ and  $T>0$, if two locally absolutely continuous {functions} $z_1,z_2:{I_z} \rightarrow \mathbb{R}^n$ and a {Lebesgue} measurable {function} $u_{z}:{I_z} \rightarrow \mathbb{R}^{m_C}$ {with $I_z$ of the form $I_z=[0,T)$ or $I_z = [0,T]$,} are such that, {for each $i \in \{1,2 \}$,} $\dot{z}_i(t) \IfIncd{\in}{=} F(z_i(t),u_{z}(t))$ for almost all $t \in {I_z}$, $(z_i(t), u_{z}(t)) \in C$ for all $t \in 
\textup{int}I_z$, and $z_i(0)=\xi$,  then $z_1(t)=z_2(t)$ for every $t \in {I_z}$;
%
  \item for every $(\xi,u_D) \in D$, $G(\xi, u_D)$ consists of one point.
  \end{enumerate}
  \label{UniquenessHu}
  \end{proposition}
}
  \IfIncd{}{
    \NotAutom{
  \begin{proof}
  Proceeding by contradiction, suppose there exist two maximal solutions to $\HS$, $(\phi_1,{u})$ and $(\phi_2,{u})$, with $u=(u_C,u_D)$, {$\dom\phi_1 = \dom\phi_2 = \dom u$}, and $\phi_1(0,0) = \phi_2(0,0)$ such that $\phi_1$ and $\phi_2$ are not identical over ${\dom u}
 $, {namely, there exists $(t^*,j^*) \in {\dom u}
 $ such that $\phi_1(t^*,j^*) \neq \phi_2(t^*,j^*)$.}  We have the following three cases:
  \begin{enumerate}[label=\alph*)]
	  \item If $(t^*,j^*) \in [t',t'') \times \{{j^*}\} \subset {\dom u}
	  $ for some  $t''>t'\geq 0$, then, the functions $z_1, z_2: I_z
	  \rightarrow \mathbb{R}^n$, and 
	  $u_{z}:I_z \rightarrow \mathbb{R}^{m_C}$, with {$T=t''-t'$} and $I_z$ of the form $I_z = [0,T)$, defined
	  {for each $i \in \{1,2 \}$ as} 
	${z}_i(t) = \phi_i ({t'+t},j)$ for all $t \in I_z$, and as $u_z(t) = u_C({t'+t},j)$ for all $t \in I_z$, 
	satisfy, {by Definition \ref{SolutiontocalH}}, $\dot{z}_i(t) \in F({z}_i(t),u_{z}(t))$ for almost all $t \in I_z$, and $(z_i(t), u_{z}(t)) \in C$ for all $t \in \textup{int}I_z$.
	However, at $t^* \in I_z$, $z_1(t^*)\neq z_2(t^*)$, 
    which contradicts item 1.
	  %
	  %
	    \item If $(t^*,j^*) \in [t',t''] \times \{{j^*}\} \subset {\dom u}
	  $ for some  $t''>t'\geq 0$, then, the functions $z_1, z_2: I_z
	  \rightarrow \mathbb{R}^n$, and 
	  $u_{z}:I_z \rightarrow \mathbb{R}^{m_C}$, with {$T=t''-t'$} and $I_z$ of the form $I_z = [0,T]$, defined
	  {for each $i \in \{1,2 \}$ as} 
	${z}_i(t) = \phi_i ({t'+t},j)$ for all $t \in I_z$, and as $u_z(t) = u_C({t'+t},j)$ for all $t \in I_z$, 
	satisfy, {by Definition \ref{SolutiontocalH}}, $\dot{z}_i(t) \in F({z}_i(t),u_{z}(t))$ for almost all $t \in [0,T)$, and $(z_i(t), u_{z}(t)) \in C$ for all $t \in \textup{int}I_z$.
	However, at $t^* \in I_z$, $z_1(t^*)\neq z_2(t^*)$, 
    which contradicts item 1.
	  %
	  %
	  %
	  %
	  \item If $(t^*,j^*)$ is such that 
	  $(t^*,j^*-1) \in \dom u
	  $,
	  then $(t^*,j^*-1)$ is a jump time of $\phi_1$ and $\phi_2$.
	  Then, by Definition \ref{SolutiontocalH},
for each $i \in \{1,2\},$ 
$\phi_i(t^*, j^*) = G(\phi_i(t^*,j^*-1),u_D(t^*, j^*-1))$.
Since, 
  $\phi_1(t^*,j^*) \neq \phi_2(t^*,j^*)$,
   and $\phi_1(t^*,  j^*-1)= \phi_2(t^*, j^*-1),
  G(\phi_1(t^*,j^*-1),u_D(t^*, j^*-1))$ takes more than one value, which contradicts item 2. 
	  %
  \end{enumerate}
Thus, provided items $1$ and $2$ hold, for every $\xi \in \Pi(\overline{C}) \cup \Pi(D)$ and each $u \in \mathcal{U}$  there exists a unique maximal {solution $(\phi,u)$} with $\phi(0,0)=\xi$. 
  \end{proof}
}
  }
\subsection{Hybrid Closed-loop Systems}
Given 
a hybrid system $\HS\IfIncd{_s}{}$ as in \eqref{Heq} and a function $\kappa:=(\kappa_C, \kappa_D)$ with $\kappa:\mathbb{R}^n \rightarrow \mathbb{R}^{m_C} \times \mathbb{R}^{m_D}$, 
 the autonomous hybrid system resulting from assigning $u=\kappa(x)$, namely, the hybrid closed-loop system, is given by
\begin{equation}
\HS_\kappa  \left\{ 
\begin{matrix}
\dot{x} &\IfIncd{\in}{=}& F(x,\kappa_C(x)) \quad  x \in C_\kappa \\ 
 x^+ &\IfIncd{\in}{=}& G(x,\kappa_D(x)) \quad x \in D_\kappa 
\end{matrix}
\right. 
\label{Hkeq}
\end{equation}
where $C_\kappa:=\{x \in \mathbb{R}^n: (x,\kappa_C(x)) \in C\}$ and $D_\kappa:=\{x \in \mathbb{R}^n: (x,\kappa_D(x)) \in D\}$.

A solution to the closed-loop hybrid system $\HS_\kappa$ is defined \IfAutomss{similar to Definition \ref{SolutiontocalH}.}{as follows.

\begin{definition}{Solution to $\HS_\kappa$} 
A hybrid arc $\phi$ defines a solution to the hybrid system $\HS_\kappa$ in (\ref{Hkeq}) 
  if 
\begin{itemize}
    \item $\phi(0,0) \in \overline{C_\kappa} \cup D_\kappa$,
    \item For each $j\in \mathbb{N}$ such that $I^J_{\phi}$ has a nonempty interior $ \textup{int} I^J_{\phi}$, we have, for all
    $t \in \textup{int} I^J_{\phi}$,
    \begin{equation*}
        \phi(t,j) \in C_\kappa
    \end{equation*}  
    and, for almost all $t \in I^J_{\phi}$,
    \begin{equation*}
        \frac{d}{dt}\phi(t,j)\IfIncd{\in}{=} F(\phi(t,j), \kappa_C(\phi(t,j)))
    \end{equation*}
    \item For all $(t,j) \in \dom \phi$ such that $(t,j+1)\in \dom \phi$,  
    \begin{eqnarray*}
	\phi(t,j) &\in& D_\kappa \\
	\phi(t,j+1) &\IfIncd{\in}{=}& G(\phi(t,j),\kappa_D({\phi}(t,j))) 
\end{eqnarray*}
\end{itemize}
A solution $\phi$ is a compact solution if $\phi$ is a compact hybrid arc.
\end{definition}
}
We denote by $\hat{\mathcal{S}}_{\HS_\kappa}(M)$ the set of solutions $\phi$ to $\HS_\kappa$ as in (\ref{Hkeq}) such that $\phi(0,0) \in M$. The set $\mathcal{S}_{\HS_\kappa}(M) \subset \hat{\mathcal{S}}_{\HS_\kappa}(M)$ denotes all maximal solutions from $M$. 
\IfIh{}{Given $\mathcal{T}\subset \realsgeq \times \nats$, let us denote by $\hat{\mathcal{S}}^\mathcal{T}_{\HS  }(M)\subset \hat{\mathcal{S}}_{\HS_\kappa}(M)$ the set of compact solutions with 
terminal time in $\mathcal{T}$, i.e., if $(\phi,u) \in \hat{\mathcal{S}}^\mathcal{T}_{\HS  }(M)$ and $\max \dom (\phi,u) =(T_\phi,J_{\phi^*})$, then $(T_\phi,J_{\phi^*})\in \mathcal{T}$.}


\section{ Two-player Zero-sum Hybrid Games}\label{sec: 2pzshg}
\subsection{Game Formulation}\label{Sec: Game form}
Following the formulation in \cite{basar1999dynamic}, for each $i\in\{1,2\}$, consider the $\mathit{i}$-th player $P_i$ with dynamics described by $\HS_i$ as in (\IfIncd{\ref{Hinc}}{\ref{Heq}}) with data $(C_i, F_i, D_i, G_i)$, state $x_i \in \mathbb{R}^{n_i}$, and input $u_{ i}=(u_{Ci}, u_{Di})\in \reals^{m_{Ci}} \times \reals^{m_{Di}}$, where {$C_i\subset \mathbb{R}^{n} \times \mathbb{R}^{m_{C}}$, $F_i: \mathbb{R}^{n} \times \mathbb{R}^{m_{C}} \IfIncd{\rightrightarrows}{\rightarrow}  \mathbb{R}^{n_i}$, $D_i \subset \mathbb{R}^{n}  \times \mathbb{R}^{m_{D}}$ and ${G_i: \mathbb{R}^{n}  \times \mathbb{R}^{m_{D}} \IfIncd{\rightrightarrows}{\rightarrow} \mathbb{R}^{n_i} }$}, with \IfTp{$n_1+n_2=n$}{}. 
We denote by $\mathcal{U}_i=\mathcal{U}_{Ci} \times \mathcal{U}_{Di}$ the set of hybrid inputs for $\HS_i$; see Definition 2.3.

Notice that each player's dynamics are described in terms {of} 
maps and sets defined in the entire state and input space rather than the individual spaces ($\reals^n$ and $\reals^m$ rather than $\reals^{n_i}$ and $\reals^{m_i}$, respectively). This allows to model the ability of each player's state to evolve according to the state variables and input of the other players.

\begin{definition}{Elements of a two-player zero-sum hybrid game} \IfTp{A two}{An $N$}-player zero-sum hybrid game is composed by 
\begin{enumerate}[label=\arabic*)]
\IfTp{}{\item A communication network described by a graph $\Gamma=(\mathcal{V,E,G})$ such that $\sum_{i= 1}^N n_{i} = n$, $\sum_{i= 1}^N m_{Ci} = m_C$ and $\sum_{i= 1}^N m_{Di} = m_D$, with $n_{i},m_{Ci} , m_{Di} \in \mathbb{N}$ for each $i \in \IfTp{\{1,2\}}{\mathcal{V}}$, where N is the number of players in the game.}
\item The state $x=(x_1, x_2\IfTp{}{, \dots, x_N}) \in \reals^n$, where, for each $i \in \IfTp{\{1,2\}}{\mathcal{V}}$, $x_i \in \reals^{n_i}$ is the state of player $P_i$.
\item The set of joint input actions $\mathcal{U}=\mathcal{U}_1 \times \mathcal{U}_2 \IfTp{}{\times \dots \times \mathcal{U}_N}$ with elements $u=(u_{ 1}, u_{ 2} \IfTp{}{,\dots, u_{ N}})$, where, for each $i \in \IfTp{\{1,2\}}{\mathcal{V}}$, $u_i=(u_{Ci},u_{Di})$ is a hybrid input. For each $i \in \IfTp{\{1,2\}}{\mathcal{V}}$, $P_i$ selects $u_i$ independently of $P_{\sj{-i}}$, who selects \IfTp{$u_{\sj{-i}}$}{any $u_j$, $j\neq i$}, namely, the joint input action $u$ has components $u_i$ that are independently chosen by each player\footnote{\sj{The subindex $-i$ refers to the player $P_{3-i}$.}}.
\item The dynamics of the game, described as in $(\IfIncd{\ref{Hinc}}{\ref{Heq}}
)$ and denoted by $\HS\IfIncd{_s}{}$, with data 
 \IfPers{Comment: Before we had
 \begin{equation*} \hspace{0.3cm}
 \begin{array}{r@{\>\>}c@{\>\>}l}
 C&:=&C_1 \times C_2 \IfTp{}{\times \dots \times   C_N}\\
 D&:=&\{(x,u_D) \in \mathbb{R}^n \times \mathbb{R}^{m_D}\IfConf{\hspace{-0.1cm}}{}:
 (x_i,u_{Di}) \in D_i,
  i \in \IfTp{\{1,2\}}{\mathcal{V}}  \}\\
 G(x,u_D)&:=&\{ \hat{G}_i(x,u_{D}):
 (x_i,u_{Di}) \in D_i, i \in \IfTp{\{1,2\}}{\mathcal{V}} \}
  \hspace{2.2cm}\IfConf{}{ \forall (x,u_D) \in D }
  \end{array}
 \end{equation*}
 We modify it, so the data of each player depends on the full state which implies that whether player $P_i$ flows or jumps (potentially) depends on the full state. Every example we present so far (in which player $P_2$ consists only of an input and no states) follows this structure. Now we have:
 }
\begin{equation*} \hspace{-0.2cm}\vspace{-0.2cm}
    \begin{array}{r@{\>\>}c@{\>\>}l}
    C&:=&{C_1 \cap C_2} \IfTp{}{\times \dots \times   C_N}\\
    F(x,u_C)&:=&(F_1(x, u_{C}),F_2(x, u_{C})\IfTp{}{, \dots, F_N(x, u_{C})}) \hspace{1cm} \IfConf{}{ \forall (x,u_C) \in C}\\
    D&:=&{ D_1 \cup D_2}
    \\
    G(x,u_D)&:=&\{ \hat{G}_i(x,u_{D}):
    {(x,u_{D})} \in D_i, i \in \IfTp{\{1,2\}}{\mathcal{V}} \}
     \hspace{2.2cm}\IfConf{}{ \forall (x,u_D) \in D }
     \end{array}
    \end{equation*}
where \IfTp{$\hat{G}_1(x,u_{D})=(G_1(x, u_{D}),I_{n_2})$, $\hat{G}_2(x,u_{D})=(I_{n_1},G_2(x, u_{D}))$}{\pnn{
$\hat{G}_i(x,u_{D})=(I_{n_1},I_{n_2}, \dots,G_i(x, u_{D}) ,\dots,I_{n_N})$}}, $u_C=(u_{C1},u_{C2})$, and $u_D=(u_{D1},u_{D2})$.
%
\item For each $i  \in \IfTp{\{1,2\}}{\mathcal{V}}$, a strategy space $\mathcal{K}_i$ of $P_i$ defined as a collection of mappings \IfTp{$\kappa_i:\mathbb{R}^{n\IfTp{}{_{in}}} \rightarrow \mathbb{R}^{m_{Ci}} \times \mathbb{R}^{m_{Di}}$}{$\gamma_i:\realsgeq \times \nats \times \mathbb{R}^{n_{in}} \rightarrow \mathbb{R}^{m_{Ci}} \times \mathbb{R}^{m_{Di}}$}\IfTp{}{, where $n_{in}=\sum_{j=1}^N n_j \mathcal{G}(j,i)$ denotes the size of the state information shared with $P_i$}. The strategy space of the game, namely {$\mathcal{K}=\mathcal{K}_1 \times \mathcal{K}_2$}, is the collection of mappings with elements \IfIh{$\kappa=(\kappa_1,\kappa_2\IfTp{}{,\dots,\kappa_N})$}{$\gamma=(\gamma_1,\gamma_2,\dots,\gamma_N)$}, where \IfIh{$\kappa_i\in \mathcal{K}_i$}{$\gamma_i\in \mathcal{K}_i$} for each $i \in \IfTp{\{1,2\}}{\mathcal{V}}$, such that every maximal solution $(\phi,u)$ to $\HS\IfIncd{_s}{}$ with input assigned as $\dom \phi \ni (t,j) \mapsto u_i(t,j)= \IfIh{\kappa}{\gamma}_i(\IfIh{}{t,j,}\phi(t,j))$ for each $i\in \IfTp{\{1,2\}}{\mathcal{V}}$ is complete. Each \IfIh{$\kappa_i \in \mathcal{K}_i$}{$\gamma_i \in \mathcal{K}_i$} is said to be a permissible pure\footnote{This, in contrast to when $\mathcal{K}_i$ is defined as a probability distribution, 
\sj{namely, when $ \mathcal{K}_i$ is 
the space of mixed strategies.}}\sj{feedback} strategy for $P_i$.
 \item  
%
A scalar-valued functional\IfIncd{\footnote{Given that we do not insist on having unique solutions, the cost $\mathcal{J}$ measures the largest cost of the solutions yielded to $\HS\IfIncd{_s}{} $ from $\xi$ by $u$. Thus, its arguments are \pnn{hybrid inputs as in Definition \ref{Def:HyInput}} and not solution pairs.}}{} \IfIh{$(\xi, u) \mapsto \mathcal{J}_i(\xi, u)$}{$\mathcal{J}_i:\mathbb{R}^n \times \mathcal{U} \rightarrow \mathbb{R}$} defined for each $i \in \IfTp{\{1,2\}}{\mathcal{V}}$, and called the cost associated to $P_i$. {For each $u\in \mathcal{U}$,} we refer to a single cost functional $\mathcal{J}{:=\mathcal{J}_1=-\mathcal{J}_2}$ as the cost associated to the \IfIncd{solutions}{unique solution} to $\HS\IfIncd{_s}{}$ from $\xi$ for $u$, 
and its structure is defined for each type of game. \IfIncd{It quantifies the largest cost over the solutions that $u$ yields to $\HS_\pnn{s}$ from $\xi$.}{}
\end{enumerate}
\label{elements}
\end{definition}
%
\IfPers{At times we present the joint input action from the perspective of the $i$-th player as $u=(u_1,\IfTp{u_2}{\dots, u_i,\dots , u_N})=(u_i,u_{-i})$, where $u_{-i}$ is the set of input actions of \IfTp{player $P_{-i}$}{players other than $P_i$}.} \IfPers{In the games studied in this work, player $P_i$ wants to minimize the cost $\mathcal{J}_i$. For the two-player game with $\mathcal{J}_1=-\mathcal{J}_2$, a gain by player $P_1$ implies a loss to player $P_2$, leading to the term ``zero-sum games." For this type of games, we refer to a single cost functional $\mathcal{J}=\mathcal{J}_1=-\mathcal{J}_2$. }
{
  \begin{remark}{Players' state} In scenarios where each player has its own dynamics, as in {pursue-evasion \cite{weintraub2020AnIT}, or target defense \cite{TDGame} }games, it is common to have a state associated to each player, namely $x_1$ for $P_1$ and $x_2$ for $P_2$, justifying the partition of the state $x$ in $x_1$ and $x_2$.  
  When the players do not have their own dynamics but can independently select an input, e.g., $P_1$ selects $u_1$ and $P_2$ selects $u_2$ to control a common state $x$, such state can be associated, without loss of generality, to either of the players, e.g., $x=x_1$ 
  with $n=n_1$ and $n_2=0$. \NotAutomss{This is illustrated in Example \ref{NumericalEx}.}
  \end{remark}
    }

{Notice that Definition \ref{elements} is general enough to cover games with a finite horizon, for which additional conditions 
specify the end of the game, e.g., a terminal set in the state space or fixed duration specifications \cite{leudo2022FH}.}

\NotAutomss{We say that a game formulation is in normal (or matrix) form when it describes only the correspondences between
strategies and costs. 
On the other hand, 
we refer to the mathematical description of a game to be in the Kuhn's extensive form if \pno{the formulation describes:
\begin{itemize}
\item the evolution of the game defined by its dynamics, 
\item the decision-making process defined by the strategies, 
\item the sharing of information between the players defined by the communication network, and
\item their outcome defined by the cost associated to each player.
\end{itemize}
}
If a game is formulated in a Kuhn's extensive form, then it admits a solution \cite{basar1999dynamic}. }
{From a given initial condition {$\xi$}, a given strategy $\kappa \in \mathcal{K}$ potentially leads to nonunique solutions\footnote{A given strategy $\kappa$ can lead to multiple input actions due to a nonempty $\Pi(C) \cap \Pi(D)$.} $(\phi^1,u^1), (\phi^2,u^2),\dots, (\phi^k,u^k)$ to $\HS\IfIncd{_s}{}$, where $u^l=\kappa(\phi^l)$ {and $\phi^l(0,0)=\xi$} for each $l \in \{1,2,\dots k\}$.
Thus, for the formulation in Definition \ref{elements} to be in Kuhn's extensive form \IfAutomss{\cite{basar1999dynamic}, so it admits a solution}{}, 
an appropriate cost definition is
required so 
each strategy $\kappa \in \mathcal{K}$ has a unique cost correspondence, namely, every solution $(\phi^l,u^l)$ 
with $u^l=\kappa(\phi^l)$, $l \in \{1,2,\dots, k\}$ is assigned the same cost.}
%
\IfPers{Given the elements of a \pn{zero-sum} game in the context of Definition \ref{elements}, the general solution of the game is given by the following definition.
\pn{\begin{definition}{Pareto-dominance solution}
An outcome of a game is Pareto dominated if some other outcome would make at least one player better off without hurting any other player. That is, some other outcome is weakly preferred by all players and strictly preferred by at least one player. If an outcome is not Pareto dominated by any other, then it is Pareto optimal, named after Vilfredo Pareto.
 \end{definition}
\begin{definition}{Pareto-optimal solution}For the case in which the player $i$ aims to minimize its own cost $\mathcal{J}_i: \mathbb{R}^n \times \mathcal{U} \rightarrow \mathbb{R}_{\geq 0}$
, $i \in \IfTp{\{1,2\}}{\mathcal{V}}$, (TO DEFINE IF LCi, LDi, Vi)  following structure (\ref{defJ}), we consider the scenario in which they cooperate, acting in full trust of each other, to attain a mutually beneficial solution. Here, the cost $\mathcal{J}_i$ depends on the full state $x$ and the input $u$ of the network $\HS  $.
For this scenario,  we say an action $u^* \in \mathcal{U}$ is a \textit{Pareto-optimal solution} if for each $\xi \in \Pi(\overline{C})\cup \Pi(D)$, there does not exist another action $u^p \in \mathcal{U}$ such that 
\begin{equation*}
\mathcal{J}_i(\xi,u^p) \leq \mathcal{J}_i(\xi,u^*)
\end{equation*}
for all $i \in \IfTp{\{1,2\}}{\mathcal{V}}$, with 
\begin{equation*}
\mathcal{J}_i(\xi,u^p) < \mathcal{J}_i(\xi,u^* )
\end{equation*}
for some $i \in \IfTp{\{1,2\}}{\mathcal{V}}$. 
In other words,  $u^* \in \mathcal{U}$ is a Pareto-optimal solution if for any action $u^p \in \mathcal{U}$, such that 
\begin{equation*}
\mathcal{J}_i(\xi,u^p) < \mathcal{J}_i(\xi,u^* )
\end{equation*}
for some $i \in \IfTp{\{1,2\}}{\mathcal{V}}$, there exists an $i \in \IfTp{\{1,2\}}{\mathcal{V}}$ such that
\begin{equation*}
\mathcal{J}_i(\xi,u^p) > \mathcal{J}_i(\xi,u^* )
\end{equation*}
Hence, if $u^*$ is a Pareto-optimal solution rendering costs $\mathcal{J}^*_i$, $i \in \IfTp{\{1,2\}}{\mathcal{V}}$, there does not exist another solution to the game for which the cost of some players is reduced without increasing the cost of others (\cite{Hespanha09}).
\end{definition}}

 \pn{\begin{definition}{Pure strategy Nash equilibrium} For a \IfTp{two}{N}-player noncooperative game, with \IfTp{}{communication network $\Gamma=(\mathcal{V,E,G})$ and }dynamics $\HS$ as in (\ref{Heq}),  
we say a strategy \IfIh{$\kappa^*=(\kappa_1^*, \IfTp{\kappa_2^*}{\dots, \kappa_N^*})\in \mathcal{K}$}{$\gamma^*=(\gamma_1^*, \dots, \gamma_N^*)\in \mathcal{K}$} 
 is \textit{a pure\footnote{This is in contrast to when it is defined as a probability distribution of selecting certain mappings, referred to as a mixed strategy.} strategy Nash equilibrium} if for any $\xi \in \Pi(\overline{C}) \cup \Pi(D)$,
 every $u^*=(u^*_1,\IfTp{u^*_2}{ \dots, u^*_N})$ rendering a maximal trajectory $\phi^*$ to $\HS  $ from $\xi$, 
with components defined as \IfIh{$\dom \phi^* \ni (t,j) \mapsto u^*_i(t,j)= \kappa_i^*(\IfIh{}{t,j,}\phi^*(t,j))$}{$\dom \phi^* \ni (t,j) \mapsto u^*_i(t,j)= \gamma_i^*(t,j,\phi^*(t,j))$}, for each $i\in \IfTp{\{1,2\}}{\mathcal{V}}$, 
satisfies
\IfConf{
\begin{eqnarray}
\mathcal{J}_1(\xi,u^*) 
&\leq&  \mathcal{J}_1(\xi,(u_{1}, u_{-1}^*)) \nonumber\\&& \hspace{1cm}  \forall u_1:\exists \phi : (\phi, (u_{1}, u_{-1}^*)) \in \mathcal{S}_{\HS  } (\xi),
\nonumber \\
\mathcal{J}_2(\xi,u^*) 
&\leq&  \mathcal{J}_2(\xi,(u_{2}, u_{-2}^*)) \nonumber\\ && \hspace{1cm}  \forall  u_2:\exists \phi :(\phi, (u_{2}, u_{-2}^*)) \in \mathcal{S}_{\HS  } (\xi), \nonumber \\
\IfTp{}{&\vdots&
\nonumber\\
\mathcal{J}(\xi,u^*) 
&\leq&  \mathcal{J}(\xi,(u_{N}, u_{-N}^*))\nonumber \\ && \hspace{1cm}  \forall  u_N:\exists \phi :(\phi, (u_{N}, u_{-N}^*)) \in \mathcal{S}_{\HS  } (\xi) \nonumber\\}
 \label{NashEqNCoopIneq}
\end{eqnarray}
\vspace{-0.9cm}
}
{
\begin{equation}
\begin{array}{r@{}l}
\mathcal{J}(\xi,u^*) 
&{}\leq  \mathcal{J}(\xi,(u_{1}, u_{2}^*\IfTp{}{, \dots, u_{N}^*})) \hspace{1cm}  \forall u_1:\exists \phi : (\phi, (u_{1}, u_{2}^*\IfTp{}{, \dots, u_{N}^*})) \in \mathcal{S}_{\HS  } (\xi),
 \\
\mathcal{J}(\xi,u^*) 
&{}\leq 
 \mathcal{J}(\xi,(u_{1}^*, u_{2}\IfTp{}{, \dots, u_{N}^*})) \hspace{1cm} \forall  u_2:\exists \phi :(\phi, (u_{1}^*, u_{2}\IfTp{}{, \dots, u_{N}^*})) \in \mathcal{S}_{\HS  } (\xi)\IfTp{}{, 
 \\
&{ }\hspace{0.2cm}\vdots
\\
\mathcal{J}(\xi,u^*) 
&{}\leq  \mathcal{J}(\xi,(u_{ 1}^*, u_{2}^*, \dots, u_{N}))\hspace{1cm} \forall  u_N:\exists \phi :(\phi, (u_{ 1}^*, u_{2}^*, \dots, u_{N})) \in \mathcal{S}_{\HS} (\xi)}.
\end{array}
 \label{NashEqNCoopIneq}
\end{equation}
 } 
 \label{NashEqCoop}
 \end{definition}}}
 %
 %
%
\subsection{Equilibrium Solution Concept}
Given the formulation of the elements of a 
hybrid game in Definition \ref{elements}, 
its solution is defined as follows.

\begin{definition}{\sj{Feedback} saddle-point equilibrium}\label{NashEqNonCoop}  Consider a two-player zero-sum game, with 
dynamics $\HS\IfIncd{_s}{}$ as in (\IfIncd{\ref{Hinc}}{\ref{Heq}}) with $\mathcal{J}_1=\mathcal{J}$, $\mathcal{J}_2=-\mathcal{J}$, for a given cost functional $\mathcal{J}:\mathbb{R}^n\times \mathcal{U} \rightarrow \mathbb{R}$. 
We say that a strategy $\IfIh{\kappa}{\gamma}=(\IfIh{\kappa}{\gamma}_1, \IfIh{\kappa}{\gamma}_2) \in \mathcal{K}$ 
is a \textit{\sj{feedback} saddle-point equilibrium} if for each $\xi \in \Pi(\overline{C}) \cup \Pi(D)$,
every 
 hybrid input $u^*=(u_1^*, u_2^*)$ 
{such that there exists $\phi^* \in \mathcal{R}(\xi,u^*)$,}
with components defined as $\dom \phi^* \ni (t,j) \mapsto u_i^*(t,j)= \IfIh{\kappa}{\gamma}_i(\IfIh{}{t,j,}\phi^*(t,j))$, 
$i\in \IfTp{\{1,2\}}{\mathcal{V}}$,
satisfies
\begin{eqnarray}
 \mathcal{J}(\xi,(u_{ 1}^*, u_{ 2})) \leq
\mathcal{J}(\xi, u^*) \leq
 \mathcal{J}(\xi,(u_{ 1}, u_{ 2}^*))
 \label{SaddlePointIneq}
\end{eqnarray}
%
{for all hybrid inputs $u_1$ and $u_2$ such that $\mathcal{R}(\xi, (u_{ 1}^*, u_{ 2}))$ and $\mathcal{R}(\xi, (u_{ 1}, u_{ 2}^*))$ are nonempty, \sj{ and every such hybrid input $u^*$ is strongly time consistent}.} 
 \end{definition}
Definition \ref{NashEqNonCoop} is a generalization of the classical pure strategy Nash equilibrium \cite[(6.3)]{basar1999dynamic} to the case where the players exhibit hybrid dynamics and opposite optimization goals.
 In words, we refer to the strategy {\IfIh{$\kappa=(\kappa_1,\kappa_2)$}{$\gamma^*$}} 
 as a \sj{feedback} saddle-point when 
a player $P_i$ cannot improve the cost $\J_i$ by playing any strategy different from \IfIh{$\kappa_i$}{$\gamma_i^*$} when \IfTp{the}{every } player $P_{\sj{-i}}$ is playing \NotAutomss{the strategy of the saddle-point,\>}\IfIh{$\kappa_{\sj{-i}}$}{$\gamma_{\sj{-i}}^*$}.     
%
{Condition \eqref{SaddlePointIneq} is verified over the set of inputs that define joint input actions $(u_1^*,u_2)$ and $(u_1,u_2^*)$, yielding at least one nontrivial solution to $\HS\IfIncd{_s}{}$ from $\xi$.}

\sj{
  \NotAutom{
  Let $u_{[(s,k),(r,l)]}$ denote the truncation of $u \in \mathcal{U}, \sup \dom u=(T,J)$, to the hybrid time interval $[(s,k),(r,l)]\subset \dom u$, and 
\begin{multline}
(\phi,u)^{\beta}_{[(s,k),(r,l)]} :=\\
 \big\{
(\phi,u)\in \mathcal S_\HS: 
u_{[(0,0),(s,k)]}=\beta(\phi_{[(0,0),(s,k)]}),\\
u_{[(r,l),(T,J)]}=\beta(\phi_{[(r,l),(T,J)]})
\big\}
\end{multline}
denote the set of solutions where the input in the intervals $[(0,0),(s,k)]$ and $[(r,l),(T,J)]$ is defined in terms of $\beta \in \mathcal{K}$. Then, we define time consistency as follows.

\begin{definition}{}
An input action $u^*$ that solves Problem ($\diamond$) is \emph{strongly time consistent} (STC) if its truncation to the interval $[(s,k),(T,J)]$,  $u^*_{[(s,k),(T,J)]}$, solves the truncated version of Problem ($\diamond$) over the set of solutions $(\phi,u)^{\beta}_{[(s,k),(T,J)]}$, for every $\beta \in \mathcal{K}$ and every $(s,k) \in ((0,0), (T,J)] \cap \dom (\phi,u)$.
\end{definition}
  }
{
\begin{remark}{Time consistency and subgame perfection}
The permissible strategies considered in this work have a feedback information structure, {in the sense} that they depend only on the current value of the state, and not on any past history of the
{\>values of the state or hybrid time}. 
Given $\xi \in \reals^n$, we say that an input action $u^*$ is \emph{strongly time consistent} if even when the past history of input values that led $\HS$ as in \eqref{Heq} to $\xi$ were not optimal, the action $u^*$ is still a solution for the remaining of the game (subgame), which is defined in the forthcoming Problem ($\diamond$), starting from $\xi$.
{When} this property holds for every state $\xi$ in $\Pi(C) \cup \Pi(D)$, we say that $u^*$ is \emph{subgame perfect}, see \cite{FERSHTMAN1989191}.
Then, under a strategy space that does not impose structural restrictions on the permissible strategies, (e.g., a linear dependence on the state) the saddle-point equilibrium strategy, when it exists, is said to be strongly time consistent if its components $\kappa_C$ and $\kappa_D$ lead to input actions that are strongly time consistent for each $\xi$ in {$\Pi(C) \cup \Pi(D)$}.
Notice that given the hybrid time horizon structure of the input actions considered in this work, the saddle-point equilibrium is time independent. This results in truncations of input actions not keeping {track} of previous hybrid time values, i.e., if there exists any past history of strategies that led to the current state, this is hidden for the evaluation of the saddle-point equilibrium at the current state, which results in preservation of optimality in the subgame, property known as \emph{permanent optimality} \cite[Section 5.6]{basar1999dynamic}.
\end{remark}
}
}
\IfPers{Notice that the saddle-point, as a solution to the zero-sum \IfTp{two}{$N$}-player game, is a strategy in $\mathcal{K}$,  though the concept of a solution to a hybrid system $\HS_\pnn{s}$ 
is a hybrid arc.

 \begin{remark}{Joint Action}
The pure strategy Nash equilibrium in Definition \ref{NashEqNonCoop} implies that for each $i \in \IfTp{\{1,2\}}{\mathcal{V}}$, the $i-$th component of each $u^*=(u^*_1, \IfTp{u^*_2}{\dots, u^*_N})$ 
satisfies
	\begin{equation}
u_i^* \in
	\argmin_{u_i:(\phi, (u_{i}, u_{\sj{-i}}^*)) \in \mathcal{S}_{\HS }  (\xi)}  \mathcal{J}_i(\xi,(u_{i}, u_{\sj{-i}}^*))
	\label{minprobNash}
	\end{equation}
%
for which the input action of \IfTp{}{any }player $P_{\sj{-i}}$ is fixed.
Thus, the joint input action is composed by solutions to the minimization problem (\ref{minprobNash}) for each $i\in \IfTp{\{1,2\}}{\mathcal{V}}$. 
\end{remark}}
%
\IfPers{Let $\mathcal{A} \subset \mathbb{R}^n$ be a closed set and consider
\begin{equation*}
	\mathcal{X}_\mathcal{A}:=\{ \phi \in \mathcal{X}: \displaystyle \lim \limits_{\substack{(t,j) \in \dom \phi\\
	(t,j)\rightarrow \sup \dom \phi}} \left | 
	\phi(t,j) \right|_\mathcal{A}=0 \}
\end{equation*}
\begin{equation*}
\mathcal{U}_\mathcal{A}(\xi) :=  \left\{ u\in \mathcal{U}: \exists  \phi \in \mathcal{R} (\xi,u) \cap \mathcal{X}_{\mathcal{A}} \right\}
\end{equation*}
To move to another section
}
%
\IfIh{}{ \section{One-player Infinite horizon Game}
Consider a one-player hybrid game with dynamics $\HS$ described by (\ref{Heq}) for given $(C,F,D,G)$. 
\begin{assumption}{}
The flow map $F$ is single valued and Lipschitz continuous in $\overline{C}$. The jump map $G$ is single valued.
 \label{AssLips}
\end{assumption}
 \pnn{Under Assumption \ref{AssLips}, conditions in Proposition \ref{UniquenessHu} are satisfied, rendering solutions to $\HS$ for a given $u$ from $\xi$ unique.}

Given $\xi \in \Pi(\overline{C}) \cup \Pi(D)$, a joint input action $u=(u_C, u_D)\in \mathcal{U}$ such that maximal solutions to $\HS  $ from $\xi$ for $u$ are complete, 
the stage cost for flows $L_C:\mathbb{R}^n  \times \mathbb{R}^{m_C} \rightarrow \mathbb{R}_{\geq 0}$, the stage cost for jumps $L_D:\mathbb{R}^n  \times \mathbb{R}^{m_D} \rightarrow \mathbb{R}_{\geq 0}$, and the terminal cost $q: \mathbb{R}^n \rightarrow \mathbb{R}$, we define the cost associated to 
the solution $(\phi,u)$ to $\HS$ from $\xi$, under Assumption \ref{AssLips}, as
\begin{multline}
\mathcal{J}(\xi, u) := 
 \sum_{j=0}^{\sup_j \dom \phi} \int_{t_{j}}^{t_{j+1}} L_C(\phi(t,j),u_{C}(t,j))dt + \\ 
 \sum_{j=0}^{\sup_j \dom \phi -1} \hspace{-0.4cm}L_D(\phi(t_{j+1},j),u_{D}(t_{j+1},j))
 + \underset{(t,j) \in \textup{dom}\phi}{\limsup_{t+j\rightarrow \infty} } q(\phi(t,j))\IfConf{\\}{}
\label{defJ}
\end{multline}
where $\{t_j\}_{j=0}^{\sup_j \dom \phi}$ is a nondecreasing sequence 
associated to the definition of the hybrid time domain of $\phi$\IfAutomss{\>\cite[Definition 2.3]{65}}{; see Definition \ref{htd}}.
%
%
%
%
%
%
%
%

The one-player infinite horizon version of the game is formalized as the following optimization problem.

\textit{Problem ($\star$):} 
Given $\xi \in \mathbb{R}^n$, under Assumption \ref{AssLips}, solve 
\begin{eqnarray}
\underset{ u \in \mathcal{U}_{\HS  }^\infty(\xi)
}{\text{minimize}}
&& \mathcal{J}(\xi, u)
\label{problemih}
\end{eqnarray}
where $\mathcal{U}_{\HS  }^\infty$ is the set of joint input actions yielding maximal complete solutions to $\HS$ as defined in Section 2.1.
\begin{assumption}{Existence of a solution to the game}
Given $\xi \in \Pi(\overline{C}) \cup \Pi(D)$, a solution to Problem ($\star$) exists, namely there exists an input action $u^*=(u_C^*,u_D^*) \in \mathcal{U}_{\HS  }^\infty (\xi)$ \pnn{such that $\mathcal{J}(\xi, u^*)<\infty$}, that attains the minimum in (\ref{problemih}), and as a consequence satisfies (\ref{NashEqNCoopIneq}) in Definition \ref{NashEqCoop}. In addition, there exists a strategy $\IfIh{\kappa}{\gamma} \in \mathcal{K}$ such that every complete solution to the closed-loop system \IfIh{$\HS_\kappa$}{$\HS_\gamma$} as in (\ref{Hkeq}), from $\xi$, has a cost \pnn{that is finite and} equal to the minimum in (\ref{problemih}). 
%
%
\label{SolExistence}
\end{assumption}
\begin{remark}{Optimal control solution}
The solution to Problem $(\star)$, when it exists, can be expressed in terms of a pure strategy Nash equilibrium $\IfIh{\kappa}{\gamma}^*$ as in Definition \ref{NashEqCoop}. Each $u^*$ rendering a trajectory $\phi^*$ such that $(\phi^*,u^*) \in \mathcal{S}_{\HS}^\infty (\xi)$ with $\dom \phi^* \ni (t,j) \mapsto u^*_i(t,j)= \IfIh{\kappa}{\gamma}_i^*(\IfIh{}{t,j,}\phi^*(t,j))$  for each $i\in \IfTp{\{1,2\}}{\mathcal{V}}$, satisfies
\begin{eqnarray*}
u^*\in  
\underset{ u \in \mathcal{U}_{\HS  }^\infty(\xi)
}{\arg\min}
&& \mathcal{J}(\xi, u)
\label{Nasheq1pih}
\end{eqnarray*}
and is it referred to as an optimal control. 
 The pure strategy Nash equilibrium $\IfIh{\kappa}{\gamma}^*$ coincides with the solution to the infinite horizon optimal control problem for hybrid systems, see \cite[Problem 1]{ferrante2019certifying}.
\end{remark}

\begin{definition}{Value function} Given $\xi \in \Pi(\overline{C}) \cup \Pi(D)$, \pnn{under Assumption \ref{AssLips},} the value function  at $\xi$ is given by 
\begin{equation}
\mathcal{J}^*(\xi):=\min_{
 u\in \mathcal{U}_{\HS  }^\infty(\xi)}
\mathcal{J}(\xi,u)
 \label{cost2goih}
\end{equation}
\label{ValueFunction}
\end{definition}


\IfPers{\textit{Hybrid Dynamic Programming}\\
Let us consider specific points of the solution to study the different behaviors in their evolution. For a given complete trajectory $\phi$, we denote $(T,J)= \sup \dom \phi$. 
Let us also denote by $(\phi^*,u^*)$ the solution that attains the infimum in (\ref{cost2goih}).\\
\textbf{Continuous Evolution}\\
Pick a $\xi \in \Pi(C)$ and a $h \in \mathbb{R}_{>0}$ such that every solution  from $\xi$ flows for at least $h$ seconds. 
Thus, the cost-to-go from $\xi$ is given by
\begin{multline*}
V(\xi):=\inf_{
(\phi, u)\in \mathcal{S}_{\HS  }(\xi)}
\int_{0}^{h} L_C(\phi(t,j),u_{C}(t,j))dt
+\int_{h}^{t_1} L_C(\phi(t,j),u_{C}(t,j))dt\\
+\sum_{j=1}^{J } \int_{t_{j}}^{t_{j+1}} L_C(\phi(t,j),u_{C}(t,j))dt +  
 \sum_{j=0}^{J-1} L_D(\phi(t_{j+1},j),u_{D}(t_{j+1},j))
+ \limsup_{(t,j)\rightarrow \sup \textup{dom} \phi} V(\phi(t,j))\\
= \inf_{(\phi, u)\in \hat{\mathcal{S}}^{(h,0)}_{\HS  }(\xi)} 
\left ( 
\int_{0}^{h} L_C(\phi(t,j),u_{C}(t,j))dt
\right.
+ 
\inf_{(\phi, u)\in \mathcal{S}_{\HS  }(\phi^*(h,0))}
 \\
\left( \left. \sum_{j=0}^{J} \int_{t_{j}}^{t_{j+1}} L_C(\phi(t,j),u_{C}(t,j))dt +  
 \sum_{j=0}^{J-1} L_D(\phi(t_{j+1},j),u_{D}(t_{j+1},j))
+ \limsup_{(t,j)\rightarrow \sup \textup{dom} \phi} V(\phi(t,j))\right) \right)
\end{multline*}
 The inner infimum in the previous expression is the cost-to-go from $\phi^*(h,0)$, so we can rewrite
\begin{equation*}
V(\xi)=
\inf_{(\phi, u)\in \hat{\mathcal{S}}^{(h,0)}_{\HS  }(\xi)} 
\int_{0}^{h} L_C(\phi(t,i),u_{C}(t,i))dt
+ V(\phi^*(h,0))
\end{equation*}
By dividing by $h$ and taking the limit, we obtain
\begin{equation*}
0= \lim_{h \rightarrow 0} 
\inf_{(\phi, u)\in \hat{\mathcal{S}}^{(h,0)}_{\HS  }(\xi)} 
\frac{1}{h}\int_{0}^{h} L_C(\phi(t,i),u_{C}(t,i))dt
+ \frac{V(\phi^*(h,0))-V(\phi^*(0,0))}{h}.
\end{equation*}
Thus, given that $V$ is time independent, one has the following point-wise condition to guarantee optimal flow 
%
\begin{equation}
0=  \inf_{u \in \Pi_u^C(x)}
\left( 
 L_C(x,u)
+\left\langle \nabla V(x),F(x,u_C) \right\rangle \right) 
\hspace{1cm} \forall x \in \Pi(C)  
\label{HJBih}
\end{equation}
\textbf{Discrete Evolution}\\
Pick now a $\xi \in \Pi(D)$. 
The cost-to-go from $\xi$ 
is given by
\begin{multline*}
V(\xi):=\inf_{
(\phi, u)\in \mathcal{S}_{\HS  }(\xi)}
L_D(\xi,u_{D}(0,0))\\ 
+\sum_{j=1}^{J} \int_{t_{j}}^{t_{j+1}} L_C(\phi(t,j),u_{C}(t,j))dt +  
 \sum_{j=1}^{J-1} L_D(\phi(t_{j+1},j),u_{D}(t_{j+1},j))
 + \limsup_{(t,j)\rightarrow \sup \textup{dom} \phi} V(\phi(t,j))
 \end{multline*}
 \begin{multline*}
= \inf_
{u_D(0, 0)}
\left( L_D(\xi,u_{D}(0,0))
+ 
\inf_{(\phi, u)\in 
\mathcal{S}_{\HS  }(\phi^*(0, 1))}
\right. \\
\left. \left( \sum_{j=0}^{J} \int_{t_{j}}^{t_{j+1}} L_C(\phi(t,j),u_{C}(t,j))dt +  
 \sum_{j=0}^{J-1} L_D(\phi(t_{j+1},j),u_{D}(t_{j+1},j))
 + \limsup_{(t,j)\rightarrow \sup \textup{dom} \phi} V(\phi(t,j))\right) \right)
\end{multline*}
The inner infimum is the minimum cost for a game starting at $\phi^*(0,1)=G(\phi^*(0,0),u^*_{D}(0,0))$, so we can rewrite
\begin{equation*}
V(\xi)
 = \inf_
 {u_D(0,0) \in \Pi_u^D(x)}
\left( L_D(\phi(0,0),u_{D}(0,0))
+ V(G(\phi(0,0),u_{D}(0,0)))
 \right)
\end{equation*}
Thus, one has the following point-wise condition to guarantee optimal jumps 
\begin{equation}
V(x)=  \inf_{u \in \Pi_u^D(x)}
\left( 
 L_D(x,u)
+ V(G(x,u)) \right) 
\hspace{1cm} \forall x \in \Pi(D)
\label{Bellmanih}
\end{equation}
}

\begin{remark}{Existence of optimal processes for hybrid infinite horizon optimal control problems} Sufficient conditions to guarantee the existence of a solution to a hybrid infinite horizon optimal control problem, as in Problem ($\star$), are provided in  \cite{GOEBEL2019153}. Given Problem ($\star$), a joint input action $u=(u_C, u_D)\in \mathcal{U}$ such that maximal solutions to $\HS  $ from $\xi \in \Pi(C) \cup \Pi(D)$ for $u$ are complete,  the stage cost for flows $L_C:\mathbb{R}^n  \times \mathbb{R}^{m_C} \rightarrow \mathbb{R}_{\geq 0}$, the stage cost for jumps $L_D:\mathbb{R}^n  \times \mathbb{R}^{m_D} \rightarrow \mathbb{R}_{\geq 0}$, and the terminal cost $q:\reals^n \rightarrow \reals$, one has that if
\begin{itemize}
\item $F(x,u_C)=f_0(x)+ \sum_{l=1}^{m_c} f_l(x)u_{C,l}$, where $f_l$ is continuous,  and there exists $M>0$ such that $|f_l(x)|\leq M(1 + |x|)$ for all $x \in \Pi(\overline{C})$, $l=0,1,\dots,m_C;$
\item $C=C_x \times U_C$ for a closed set $C_x \subset \reals^n$ and a closed and convex set $U_C \subset \reals^{m_C}$;
\item $L_C$ is lower semicontinuous, convex in $u_C$ for every $x \in C_x$, and either
$U_C$ is bounded or there exist $r>1,\eta>0, \delta \in \reals$ so that $L_C(x,u_C) \geq \eta |u_C|^r +\delta$ for every $(x,u_C) \in C$; 
\item $G$ is locally bounded;
\item $L_D$ is lower semicontinuous, and either $D=D_x \times U_D$ for a set $D_x \subset \reals^n$ and a bounded set $U_D \subset \reals^{m_D}$ or $L_D(x,u_D) \geq \beta(|u_D|)$ for every $(x,u_D) \in D$ for some $\beta:[0,\infty) \rightarrow [0,\infty)$ such that $\lim_{r \to \infty} \beta(r)=\infty$; 
\item \pn{The set of forward complete hybrid time domains starting at $(0, 0)$ is closed with respect to set convergence;}
\item The function $q$ is lower semicontinuous;
 \end{itemize}
 then $\mathcal{J}^*$ is lower semicontinuous, and if $\mathcal{J}^*(\xi)<\infty$, then there exists a solution to Problem $(\star)$. 
\end{remark}

\begin{theorem}{HJB for Problem ($\star$)}
Given a one-player hybrid game with dynamics $\HS$ as in (\ref{Heq}) with data $(C,F,D,G)$,
 and stage costs $L_C:\mathbb{R}^n  \times \mathbb{R}^{m_C} \rightarrow \mathbb{R}_{\geq 0} ,L_D:\mathbb{R}^n  \times \mathbb{R}^{m_D} \rightarrow \mathbb{R}_{\geq 0}$, and terminal cost $q: \mathbb{R}^n \rightarrow \mathbb{R}$\pnn{, suppose that Assumptions \ref{AssLips} and \ref{SolExistence} hold}. If there exists a function $V: \mathbb{R}^n \rightarrow \mathbb{R}$ that is continuously differentiable on a neighborhood of $\Pi(C)$, satisfies the \textit{Hamilton–Jacobi–Bellman-Isaacs hybrid equations}  
\eqs{
0=  \min_{u_C \in \Pi_u^C(x)}
\left\{ 
 L_C(x,u_C)
+\left\langle \nabla V(x),F(x,u_C) \right\rangle \right\} 
\hspace{1cm} }{ \forall x \in \Pi(\overline{C}),  
\label{HJBih}
}
\eqs{
V(x)=  \min_{u_D \in \Pi_u^D(x)}
\left\{ 
 L_D(x,u_D)
+ V(G(x,u_D)) \right\} 
\hspace{2cm} }{  \forall x \in \Pi(D) \NotConf{\hspace{0.7cm}}
\label{Bellmanih}
}
and, for each $\xi \in \Pi(\overline{C}) \cup \Pi(D)$,  each $(\phi, u)\in \mathcal{S}_{\HS  }^\infty(\xi)$ satisfies
\begin{equation}
\underset{(t,j) \in \textup{dom}\phi}{\limsup_{t+j\rightarrow \infty} } V(\phi(t,j))=\underset{(t,j) \in \textup{dom}\phi}{\limsup_{t+j\rightarrow \infty} } q(\phi(t,j))
\label{TerminalCond}
\end{equation}
then
\begin{equation}
\J^*(\xi)= V(\xi) \hspace{1cm} \forall \xi \in \Pi(\overline{C}) \cup \Pi(D)
\end{equation}
 \label{thHJBih}
 \end{theorem}
 \begin{corollary}{Stationary feedback law}
Given a one-player hybrid game with dynamics $\HS$ as in (\ref{Heq}) with data $(C,F,D,G)$,
 and stage costs $L_C:\mathbb{R}^n  \times \mathbb{R}^{m_C} \rightarrow \mathbb{R}_{\geq 0} ,L_D:\mathbb{R}^n  \times \mathbb{R}^{m_D} \rightarrow \mathbb{R}_{\geq 0}$, and terminal cost $q: \mathbb{R}^n \rightarrow \mathbb{R}$, \pnn{suppose that Assumptions \ref{AssLips} and \ref{SolExistence} hold}. If a continuously differentiable function $V$ satisfying (\ref{HJBih}), (\ref{Bellmanih}) exists, so a minimum of (\ref{HJBih}) and (\ref{Bellmanih}) exist, respectively in $\Pi_u^C(x)$ and $\Pi_u^D(x)$ for every $x \in \Pi(C) \cup \Pi(D)$, \IfIh{}{and given that $L_C,L_D, F, q,$ and $G$ are time-independent, }then the feedback law
$\kappa:=(\kappa_C,\kappa_D): \mathbb{R}^n \rightarrow \mathbb{R}^{m_C} \times \mathbb{R}^{m_D}$  with values
\eqs{
\kappa_C(x)\in\underset{u_C \in \Pi_u^C(x)}{\arg \min} \left\{ 
 L_C(x,u_C)
+\left\langle \nabla V(x),F(x,u_C) \right\rangle \right\} \hspace{1cm} }{  \forall x \in \Pi(\overline{C})
\label{kHJBeqih}
}
and
\eqs{
\kappa_D(x)\in\underset{u_D \in \Pi_u^D(x)}{\arg \min}  \left\{ 
 L_D(x,u_D)
+ V(G(x,u_D)) \right\} \hspace{2cm} }{  \forall x \in \Pi(D) 
\label{kBeqih}
}
defines an optimal control 
at every $x \in \Pi(\overline{C}) \cup \Pi(D)$. In addition, if for each $\xi \in \Pi(\overline{C}) \cup \Pi(D)$,  each $\phi \in \mathcal{S}_{\HS_\kappa  }^\infty(\xi)$ satisfies (\ref{TerminalCond}), then $\kappa$ is a pure strategy Nash equilibrium for the one-player infinite horizon game.
\label{FBLrem}
\end{corollary}
%
\begin{proof}
To show the claim we proceed as follows:
\begin{enumerate}
\item \pnn{Given the differentiable function $V$, select} a strategy $\IfIh{\kappa}{\gamma}$ attaining the minimum in (\ref{HJBih}) and (\ref{Bellmanih}) for each point in $\Pi(\overline{C}) \cup \Pi(D)$.
\item Pick an initial condition $\xi$ and evaluate the cost associated to any 
solution yield by $\IfIh{\kappa}{\gamma}$ from $\xi$. Show that this cost coincides with the value of the function $V$ at $\xi$. 
\item Lower bound the cost associated to any other solution from $\xi$ by the value of the function $V$ evaluated at $\xi$. 
\item By showing that the cost of any other solution from $\xi$ is not less than the cost of the solution yield by $\IfIh{\kappa}{\gamma}$ from $\xi$, we show optimality of $\IfIh{\kappa}{\gamma}$ in Problem ($\star$). 
\end{enumerate}
\NotConf{Proceeding as in item 1 above, 
given 
a function $V: \reals^n \to \reals$ that is continuously differentiable on $\Pi(C)$,
let $\IfIh{\kappa}{\gamma}:=(\IfIh{\kappa}{\gamma}_C,\IfIh{\kappa}{\gamma}_D) \in \mathcal{K}$ 
attain the minimum in (\ref{HJBih}) and (\ref{Bellmanih}) for every $x \in \Pi(\overline{C} \cap D)$, i.e.,
\begin{equation*}
\IfIh{\kappa}{\gamma}_C(x)\in\underset{u_C \in \Pi_u^C(x)}{\argmin} \left\{ 
L_C(x,u_C)
+\left\langle \nabla V(x),F(x,u_C) \right\rangle \right\}
\hspace{1cm} \forall x \in \Pi(\overline{C})
\label{kHJBeqihproofth}
\end{equation*}
and
\begin{equation*}
\IfIh{\kappa}{\gamma}_D(x)\in \underset{u_D \in \Pi_u^D(x)}{\argmin} \left\{ 
L_D(x,u_D) + V(G(x,u_D)) \right\}
\hspace{1cm} \forall x \in \Pi(D) 
\label{kBeqihproofth}
\end{equation*}

Continuing as in item 2 above, pick any $\xi \in \Pi(\overline{C}) \cup \Pi(D)$ and any $(\phi^*,u^*) \in \mathcal{S}_{\HS }^\infty (\xi)$ with $\dom \phi^* \ni (t,j) \mapsto u^*(t,j)= \IfIh{\kappa}{\gamma}(\IfIh{}{t,j,}\phi^*(t,j))$, 
and $(t_{J_{\phi^*}+1},J_{\phi^*}) \in \dom (\phi^*,u^*) $. 
\pnn{We show that there does not exist a complete solution to $\HS$ with cost smaller than the cost of $(\phi^*,u^*)$.}
Given that $V$ satisfies 
(\ref{HJBih}), and $\IfIh{\kappa}{\gamma}_C$ 
attains the minimum in (\ref{HJBih}), for each $j \in \nats$ such that $I_{\phi^*}^j=[t_j,t_{j+1}]$ has a nonempty interior int$I_{\phi^*}^j$, we have that for every $t \in \textup{int}I_{\phi^*}^j$, 
\begin{equation*}
0=L_C(\phi^*(t,j),\kappa_C(\phi^*(t,j)))+\left\langle \nabla V(\phi^*(t,j),\kappa_C(\phi^*(t,j))),F(\phi^*(t,j),\kappa_C(\phi^*(t,j))) \right\rangle
\end{equation*}
Given that  $\frac{d}{d t}V(\phi^*(t,j)):=\left\langle \nabla V(\phi^*(t,j),\kappa_C(\phi^*(t,j))),F(\phi^*(t,j),\kappa_C(\phi^*(t,j))) \right\rangle$ for every $(t,j):t \in \dom_t \phi^*$, 
by integrating over the interval $[t_j,t_{j+1}]$,
we obtain
\begin{eqnarray*}
0=  \int_{t_{j}}^{t_{j+1}} \left( L_C(\phi^*(t,j),\IfIh{\kappa}{\gamma}_{ C}(\phi^*(t,j)))
+ \frac{d}{d t}V(\phi^*(t,j)) \right)dt
\end{eqnarray*}
from where we have
\begin{eqnarray*}
0=  \int_{t_{j}}^{t_{j+1}} L_C(\phi^*(t,j),\IfIh{\kappa}{\gamma}_{ C}(\phi^*(t,j))) dt
+  V(\phi^*(t_{j+1},j))- V(\phi^*(t_j,j))
\end{eqnarray*}
%
%
%
Summing from $j=0$ to $j=J_{\phi^*}$, 
we obtain 
\begin{eqnarray*}
0
=\sum_{j=0}^{J_{\phi^*}}  \int_{t_{j}}^{t_{j+1}} L_C(\phi^*(t,j),\IfIh{\kappa}{\gamma}_{ C}(\phi^*(t,j))) dt
+\sum_{j=0}^{J_{\phi^*}} \left( V(\phi^*(t_{j+1},j))-V(\phi^*(t_j,j)) \right)
\end{eqnarray*}
Then, solving for $V$ at the initial condition $\phi^*(0,0)$, we obtain
%
%
\begin{eqnarray}
V(\phi^*(0,0))
=\sum_{j=0}^{J_{\phi^*}}  \int_{t_{j}}^{t_{j+1}} L_C(\phi^*(t,j),\IfIh{\kappa}{\gamma}_{ C}(\phi^*(t,j))) dt+ V(\phi^*(t_1,0))\nonumber \\
+\sum_{j=1}^{J_{\phi^*}} \left( V(\phi^*(t_{j+1},j))-V(\phi^*(t_j,j)) \right) \label{ContinCost2goih}
\end{eqnarray}
Given that $V$ satisfies 
(\ref{Bellmanih}) and $\IfIh{\kappa}{\gamma}_D$ 
attains the infimum in (\ref{Bellmanih}), for every $(t_{j+1},j) \in \dom \phi^*$ such that $(t_{j+1},j+1) \in \dom \phi^*$,  we have that
\begin{eqnarray*}
V(\phi^*(t_{j+1},j))&=&  \min_{
u_{D}(t_{j+1},j) \in \Pi_u(\phi^*(t_{j+1},j),D)}
\left\{ 
 L_D(\phi^*(t_{j+1},j),u_{ D}(t_{j+1},j))
+ V(G(\phi^*(t_{j+1},j), u_D(t_{j+1},j) \right\}\\ 
&=&  
 L_D(\phi^*(t_{j+1},j),\IfIh{\kappa}{\gamma}_{ D}(\phi^*(t_{j+1},j)))
+ V(G(\phi^*(t_{j+1},j), \IfIh{\kappa}{\gamma}_D(\phi^*(t_{j+1},j)))) \\
&=& L_D(\phi^*(t_{j+1},j),\IfIh{\kappa}{\gamma}_{ D}(\phi^*(t_{j+1},j)))
+ V(\phi^*(t_{j+1},j+1)) 
\end{eqnarray*}
where $\phi^*(t_{j+1},j)\in D_{\IfIh{\kappa}{\gamma}}$ is defined in (\ref{Hkeq}).
Summing both sides from $j=0$ to $j=J_{\phi^*}-1$, 
we obtain
\begin{eqnarray*}
\sum_{j=0}^{J_{\phi^*}-1} V(\phi^*(t_{j+1},j))
= \sum_{j=0}^{J_{\phi^*}-1} L_D(\phi^*(t_{j+1},j),\IfIh{\kappa}{\gamma}_{ D}(\phi^*(t_{j+1},j))) 
+\sum_{j=0}^{J_{\phi^*}-1} V( \phi^*(t_{j+1},j+1))
\end{eqnarray*}
Then, solving for $V$ at the first jump time $(t_1,0)$, we obtain
%
\begin{eqnarray}
V(\phi^*(t_{1},0))
= V( \phi^*(t_{1},1))
+ \sum_{j=0}^{J_{\phi^*}-1} L_D(\phi^*(t_{j+1},j),\IfIh{\kappa}{\gamma}_{ D}(\phi^*(t_{j+1},j))) \nonumber
\\
+\sum_{j=1}^{J_{\phi^*}-1}\left( V( \phi^*(t_{j+1},j+1))-V(\phi^*(t_{j+1},j)) \right)
\label{DiscreteCost2goih}
\end{eqnarray}
Given that $\phi^*(0,0)=\xi$, by substituting the right-hand side of (\ref{DiscreteCost2goih}) in  (\ref{ContinCost2goih}), we obtain
\begin{eqnarray}
V(\xi)
&=&\sum_{j=0}^{J_{\phi^*}}  \int_{t_{j}}^{t_{j+1}} L_C(\phi^*(t,j),\IfIh{\kappa}{\gamma}_{ C}(\phi^*(t,j))) dt+ V(\phi^*(t_1,0))\nonumber 
+\sum_{j=1}^{J_{\phi^*}} \left( V(\phi^*(t_{j+1},j))-V(\phi^*(t_j,j)) \right) \nonumber \\
&=&\sum_{j=0}^{J_{\phi^*}}  \int_{t_{j}}^{t_{j+1}} L_C(\phi^*(t,j),\IfIh{\kappa}{\gamma}_{ C}(\phi^*(t,j))) dt
+ \sum_{j=0}^{J_{\phi^*}-1}L_D(\phi^*(t_{j+1},j),\IfIh{\kappa}{\gamma}_{ D}(\phi^*(t_{j+1},j)))\nonumber \\ &&+ V(\phi^*(t_{1},1))+\sum_{j=1}^{J_{\phi^*}-1}\left( V(\phi^*(t_{j+1},j+1))-V(\phi^*(t_{j+1},j)) \right) \nonumber 
\\&&+\sum_{j=1}^{J_{\phi^*}} \left( V(\phi^*(t_{j+1},j))-V(\phi^*(t_j,j)) \right) \nonumber 
\end{eqnarray}
\pnn{Since 
\begin{eqnarray}
&& V(\phi^*(t_{1},1))+\sum_{j=1}^{J_{\phi^*}-1}\left( V(\phi^*(t_{j+1},j+1))-V(\phi^*(t_{j+1},j)) \right) \nonumber 
+\sum_{j=1}^{J_{\phi^*}} \left( V(\phi^*(t_{j+1},j))-V(\phi^*(t_j,j)) \right) \nonumber
\\ &=& V(\phi^*(t_{J_{\phi^*}+1},J_{\phi^*}))+V(\phi^*(t_{1},1))+\sum_{j=1}^{J_{\phi^*}-1}\left( V(\phi^*(t_{j+1},j+1)) \right) \nonumber 
-\sum_{j=1}^{J_{\phi^*}} \left(V(\phi^*(t_j,j)) \right) \nonumber
\\ &=&V(\phi^*(t_{J_{\phi^*}+1},J_{\phi^*})) \nonumber
\end{eqnarray}
then we have}
\begin{eqnarray}
V(\xi)
&=&\sum_{j=0}^{J_{\phi^*}}  \int_{t_{j}}^{t_{j+1}} L_C(\phi^*(t,j),\IfIh{\kappa}{\gamma}_{ C}(\phi^*(t,j))) dt
+ \sum_{j=0}^{J_{\phi^*}-1}L_D(\phi^*(t_{j+1},j),\IfIh{\kappa}{\gamma}_{ D}(\phi^*(t_{j+1},j))) \nonumber
\\ &&+V(\phi^*(t_{J_{\phi^*}+1},J_{\phi^*}))
\nonumber
\end{eqnarray}

By taking the limit when $(t_{J_{\phi^*}+1},J_{\phi^*}) \rightarrow \sup \dom \phi^*$, and given that (\ref{TerminalCond}) holds, we have 
\begin{eqnarray}
V(\xi)&=& \sum_{j=0}^{\sup_j \dom \phi^*}   \int_{t_{j}}^{t_{j+1}} L_C(\phi^*(t,j),\IfIh{\kappa}{\gamma}_{ C}(\phi^*(t,j))) dt
+ \sum_{j=0}^{\sup_j \dom \phi^* -1} L_D(\phi^*(t_{j+1},j),\IfIh{\kappa}{\gamma}_{ D}(\phi^*(t_{j+1},j))) \nonumber 
\\ &&
+\underset{(t,j) \in \textup{dom}\phi^*}{\limsup_{t+j\rightarrow \infty} } V(\phi^*(t,j))
\nonumber \\ 
&=& \mathcal{J}(\xi,u^*)
\label{OptimalCostih}
\end{eqnarray}
for any $(\phi^*,u^*) \in \mathcal{S}_{\HS  }^\infty (\xi)$ with $\dom \phi^* \ni (t,j) \mapsto u^*(t,j)= \IfIh{\kappa}{\gamma}(\IfIh{}{t,j,}\phi^*(t,j))$.

Continuing with item 3 as above, pick any $(\phi_s,u) \in \mathcal{S}_{\HS  }^\infty (\xi)$ with 
$\dom \phi_s \ni (t,j) \mapsto u(t,j)= \bar{\IfIh{\kappa}{\gamma}}(\IfIh{}{t,j,}\phi_s(t,j))$ for some $\bar{\IfIh{\kappa}{\gamma}} \in \mathcal{K}$ 
  and $(t_{J_{\phi_s}+1},J_{\phi_s}) \in \dom (\phi_s,u)$. Since $\bar{\IfIh{\kappa}{\gamma}}$ does not necessarily attain the minimum in (\ref{HJBih}), then,  for each $j \in \nats$ such that $I_{\phi_s}^j=[t_j,t_{j+1}]$ has a nonempty interior int$I_{\phi_s}^j$, we have that for every $t \in \textup{int}I_{\phi_s}^j$,
\begin{eqnarray*}
0
\leq   L_C(\phi_s(t,j),u_{ C}(t,j))
+\left\langle \nabla V(\phi_s(t,j),u_C(t,j)),F(\phi_s(t,j),u_C(t,j)) \right\rangle
\end{eqnarray*}
Given that  $\frac{d}{d t}V(\phi_s(t,j)):=\left\langle \nabla V(\phi_s(t,j),u_C(t,j)),F(\phi_s(t,j),u_C(t,j)) \right\rangle$ for every $(t,j):t \in \dom_t \phi_s$, 
by integrating over the interval $[t_j,t_{j+1}]$, 
we obtain
\begin{eqnarray*}
0 \leq \int_{t_{j}}^{t_{j+1}}  \left( L_C(\phi_s(t,j),u_{ C}(t,j))
+  \frac{d}{d t}V(\phi_s(t,j)) \right) dt
\end{eqnarray*}
from which we have
\begin{eqnarray*}
 V(\phi_s(t_j,j))
\leq \int_{t_{j}}^{t_{j+1}}   L_C(\phi_s(t,j),u_{ C}(t,j))dt
+V(\phi_s(t_{j+1},j))
\end{eqnarray*}
Summing both sides from  $j=0$ to $j=J_{\phi_s}$, 
we obtain 
\begin{eqnarray}
\sum_{j=0}^{J_{\phi_s}}V(\phi_s(t_j,j)) 
\leq \sum_{j=0}^{J_{\phi_s}}  \int_{t_{j}}^{t_{j+1}}  L_C(\phi_s(t,j),u_{C}(t,j))dt
+\sum_{j=0}^{J_{\phi_s}}V(\phi_s(t_{j+1},j))
\nonumber
\end{eqnarray}
Then, solving for $V$ at the initial time $(0,0)$, we obtain
%
%
\begin{eqnarray}
V(\phi_s(0,0))
\leq \sum_{j=0}^{J_{\phi_s}}  \int_{t_{j}}^{t_{j+1}}  L_C(\phi_s(t,j),u_{C}(t,j))dt+ V(\phi_s(t_1,0))\nonumber \\+\sum_{j=1}^{J_{\phi_s}} \left( V(\phi_s(t_{j+1},j))-V(\phi_s(t_j,j)) \right) 
\label{ContinCost2goNOih}
\end{eqnarray}
In addition, since $u$ defined by $\bar{\IfIh{\kappa}{\gamma}}$ does not necessarily attain the infimum in (\ref{Bellmanih}), then for every $(t_{j+1},j) \in \dom \phi_s $ such that   $(t_{j+1},j+1) \in \dom \phi_s$, we have
\begin{eqnarray*}
V(\phi_s(t_{j+1},j))
&\leq&   
L_D(\phi_s(t_{j+1},j),u_{ D}(t_{j+1},j))
+ V(G(\phi_s(t_{j+1},j),u_{D}(t_{j+1},j)))\\
&=&L_D(\phi_s(t_{j+1},j),u_{ D}(t_{j+1},j))
+ V( \phi_s(t_{j+1},j+1))
\end{eqnarray*}
Summing both sides from  $j=0$ to $j=J_{\phi_s}-1$, 
we obtain
\begin{eqnarray}
\sum_{j=0}^{J_{\phi_s}-1}V(\phi_s(t_{j+1},j)) 
\leq 
 \sum_{j=0}^{J_{\phi_s}-1}L_D(\phi_s(t_{j+1},j),u_{ D}(t_{j+1},j))
+\sum_{j=0}^{J_{\phi_s}-1} V(\phi_s(t_{j+1},j+1))
\nonumber
\end{eqnarray}
Then, solving for $V$ at the first jump time $(t_1,0)$, we obtain
%
\begin{eqnarray}
V(\phi_s(t_{1},0))
\leq V( \phi_s(t_{1},1))
+ \sum_{j=0}^{J_{\phi_s}-1}L_D(\phi_s(t_{j+1},j),u_{ D}(t_{j+1},j))\nonumber
\\
+\sum_{j=1}^{J_{\phi_s}-1}\left( V(\phi_s(t_{j+1},j+1))-V(\phi_s(t_{j+1},j)) \right)
\label{DiscreteCost2goNOih}
\end{eqnarray}
In addition, given that $\phi_s(0,0)=\xi$, upperbounding $V(\phi_s(t_{1},0))$ in (\ref{ContinCost2goNOih}) by the right-hand side of (\ref{DiscreteCost2goNOih}), we obtain 
\begin{eqnarray}
V(\xi)
&\leq& \sum_{j=0}^{J_{\phi_s}}  \int_{t_{j}}^{t_{j+1}}  L_C(\phi_s(t,j),u_{C}(t,j))dt+ V(\phi_s(t_1,0)) \nonumber 
\\&&+\sum_{j=1}^{J_{\phi_s}} \left( V(\phi_s(t_{j+1},j))-V(\phi_s(t_j,j)) \right) \nonumber \\
&\leq&
\sum_{j=0}^{J_{\phi_s}}  \int_{t_{j}}^{t_{j+1}}  L_C(\phi_s(t,j),u_{C}(t,j))dt
+ \sum_{j=0}^{J_{\phi_s}-1}L_D(\phi_s(t_{j+1},j),u_{D}(t_{j+1},j))\nonumber
\\&&+\sum_{j=1}^{J_{\phi_s}-1}\left( V(\phi_s(t_{j+1},j+1))-V(\phi_s(t_{j+1},j)) \right) \nonumber \\
&&+ V(\phi_s(t_{1},1))+\sum_{j=1}^{J_{\phi_s}} \left( V(\phi_s(t_{j+1},j))-V(\phi_s(t_j,j)) \right)\nonumber 
\nonumber
\end{eqnarray}
\pnn{Since 
\begin{eqnarray}
&& V(\phi_s(t_{1},1))+\sum_{j=1}^{J_{\phi_s}-1}\left( V(\phi_s(t_{j+1},j+1))-V(\phi_s(t_{j+1},j)) \right) \nonumber 
+\sum_{j=1}^{J_{\phi_s}} \left( V(\phi_s(t_{j+1},j))-V(\phi_s(t_j,j)) \right) \nonumber
\\ &=& V(t_{{J_{\phi_s}}+1},J_{\phi_s}) +V(\phi_s(t_{1},1))+\sum_{j=1}^{J_{\phi_s}-1}\left( V(\phi_s(t_{j+1},j+1)) \right) \nonumber 
-\sum_{j=1}^{J_{\phi_s}} \left(V(\phi_s(t_j,j)) \right) \nonumber
\\ &=&V(t_{{J_{\phi_s}}+1},J_{\phi_s}) \nonumber
\end{eqnarray}
then we have}
\begin{eqnarray}
V(\xi)
&=&\sum_{j=0}^{J_{\phi_s}}  \int_{t_{j}}^{t_{j+1}}  L_C(\phi_s(t,j),u_{C}(t,j))dt
+ \sum_{j=0}^{J_{\phi_s}-1}L_D(\phi_s(t_{j+1},j),u_{D}(t_{j+1},j))\nonumber
\\&&+V(\phi_s(t_{J_{\phi_s}+1},J_{\phi_s}))
\nonumber
\end{eqnarray}
By taking the limit when $(t_{{J_{\phi_s}}+1},J_{\phi_s}) \rightarrow \sup \dom \phi_s$, and given that (\ref{TerminalCond}) holds, we have
\begin{eqnarray}
V(\xi)&\leq&\sum_{j=0}^{\sup_j \dom \phi_s} \int_{t_{j}}^{t_{j+1}}  L_C(\phi_s(t,j),u_{C}(t,j))dt
+ \sum_{j=0}^{{\sup_j \dom \phi_s}-1}L_D(\phi_s(t_{j+1},j),u_{D}(t_{j+1},j))\nonumber
\\ &&
+\underset{(t,j) \in \textup{dom}\phi_s}{\limsup_{t+j\rightarrow \infty} } V(\phi_s(t,j))
\nonumber \\ 
&\leq& \mathcal{J}(\xi,u)
\label{BoundCostNO}
\end{eqnarray}
for $u$ defined by any $\bar{\IfIh{\kappa}{\gamma}} \in \mathcal{K}$. 

\pnn{Finally, proceeding as in item 4 above, under Assumption \ref{SolExistence},
by applying the $\min$ on each side of  (\ref{BoundCostNO}) over the set $\mathcal{U}_\HS^\infty (\xi)$, we obtain
\begin{equation*}
V(\xi) \leq  \underset{ u \in \mathcal{U}_{\HS  }^\infty(\xi)
}{\text{min}}\mathcal{J}(\xi, u).
\end{equation*}
Given that $V(\xi)=\mathcal{J}(\xi, u^
*)$ from (\ref{OptimalCostih}), we have that for any $\xi \in \Pi(\overline{C}) \cup \Pi(D)$, each 
$u^*$ defined by $\IfIh{\kappa}{\gamma}$ satisfies
\begin{equation*}
\mathcal{J}(\xi,u^*) \leq \underset{ u \in \mathcal{U}_{\HS  }^\infty(\xi)
}{\text{min}} \mathcal{J}(\xi, u)
\end{equation*}
which implies
\begin{equation}
\mathcal{J}(\xi,u^*) = \underset{ u \in \mathcal{U}_{\HS  }^\infty(\xi)
}{\text{min}} \mathcal{J}(\xi, u).
\label{LowerBoundd}
\end{equation}
Thus, from (\ref{OptimalCostih}) and (\ref{LowerBoundd}), $\kappa$ is the optimal strategy and $V(\xi)$ is the value function for $\HS$, as in Definition \ref{ValueFunction}.}}
\end{proof}
\begin{remark}{Equivalent conditions}
If there exists a continuously differentiable function $V: \mathbb{R}^n \rightarrow \mathbb{R}$ on $\Pi(C)$ that satisfies (\ref{HJBih}), (\ref{Bellmanih}) and a feedback law $\kappa:=(\kappa_C,\kappa_D): \mathbb{R}^n \rightarrow \mathbb{R}^{m_C} \times \mathbb{R}^{m_D}$ that satisfies (\ref{kHJBeqih}), (\ref{kBeqih}), then
\begin{eqnarray}
&& \IfAutom{\mathcal{L}_C(x,\kappa_C(x))}
  {L_C(x,\kappa_C(x))
+
\left\langle \nabla V(x),F(x,\kappa_C(x)) \right\rangle}=0  \hspace{1.2cm} \IfConf{\nonumber \\&& \hspace{6cm} }{} \forall x \in C_\kappa 
\label{kHJBeqiha}
\\
&&L_C(x,u_C)
+\left\langle \nabla V(x),F(x,u_C) \right\rangle \geq 0 \hspace{1.2cm}\IfConf{\nonumber \\&& \hspace{5.2cm} }{}\forall (x,u_C) \in C 
\label{kHJBeqihb}
\\
&& L_D(x,\kappa_D(x))
+ V(G(x,\kappa_D(x)))-V(x)=0  \hspace{1.2cm}\IfConf{\nonumber \\&& \hspace{6cm} }{}\forall x \in D_\kappa 
\label{kBeqihc}
\\
&& L_D(x,u_D)
+ V(G(x,u_D)) -V(x) \geq0 \hspace{1.2cm} \IfConf{\nonumber \\&& \hspace{5.2cm} }{}\forall (x,u_D) \in D
\label{kBeqihd}
\end{eqnarray}
In fact, from (\ref{kHJBeqih}), (\ref{kBeqih}) we have
\begin{multline}
\underset{u_C \in \Pi_u^C(x)}{\min} \left\{ 
 L_C(x,u_C)
+\left\langle \nabla V(x),F(x,u_C) \right\rangle \right\} \\
=L_C(x,\kappa_C(x))
+\left\langle \nabla V(x),F(x,\kappa_C(x))\right\rangle \hspace{0.9cm}\forall x \in \Pi(C) \IfConf{\\}{}
\label{kHJBeqihproof}
\end{multline}
and
\begin{multline}
\underset{u_D \in \Pi_u^D(x)}{\min}  \left\{ 
 L_D(x,u_D)
+ V(G(x,u_D)) \right\}\IfConf{\\}{}=L_D(x,\kappa_D(x))
+ V(G(x,\kappa_D(x))) \hspace{1cm} \forall x \in \Pi(D) 
\label{kBeqihproof}
\end{multline}
Thus, (\ref{kHJBeqiha}) and (\ref{kBeqihc}) are implications of (\ref{kHJBeqihproof}), (\ref{HJBih}) and (\ref{kBeqihproof}), (\ref{Bellmanih}) respectively. Similarly, (\ref{kHJBeqihb}) and (\ref{kBeqihd}) are derived straightforwardly from (\ref{HJBih}) and (\ref{Bellmanih}). 
\end{remark}
From \cite[Remark 3]{ferrante2019certifying}, if $C_\kappa=\Pi(C), D_\kappa=\Pi(D)$, then $V$ and $\kappa$ satisfying (\ref{kHJBeqiha}), (\ref{kHJBeqihb}), (\ref{kBeqihc}), and (\ref{kBeqihd}) satisfy also (\ref{HJBih}), (\ref{Bellmanih}), (\ref{kHJBeqih}), and (\ref{kBeqih}).

Next, we illustrate in some examples with linear dynamics and quadratic costs how Theorem \ref{thHJBih} provides conditions to characterize
the value function and the Nash equilibrium in a noncooperative infinite horizon one-player game covering both continuous, discrete and hybrid systems.
\begin{example}{Linear Quadratic Differential Games}
Consider a one-player hybrid game with state $x \in \reals^n$, input $u=(u_C, u_D) \in \reals^{m_C}\times \reals^{m_D}$, and dynamics $\HS$ as in (\ref{Heq}) with 
\begin{eqnarray}
C&:=&\mathbb{R}^n \times \mathbb{R}^{m_C} \nonumber
\\
F(x,u_c)&:=&A_C x+B_C u_C \nonumber
\label{LinearGame}
\end{eqnarray}
and jump set $D=\emptyset$, where $A_C \in \mathbb{R}^{n\times n}$ and $B_C \in \mathbb{R}^{n\times m_C}$. 
Consider the stage costs $L_C(x,u_C):=x^\top Q_C x+ u_C^\top R_C u_C$ and $L_D(x,u_D):=0$, the terminal cost $q(x):=x^\top P x$ and the function $V(x):=x^\top P x$, where $Q_C \in \mathbb{S}^n_+$, $R_C \in \mathbb{S}^{m_C}_+$ and  $P\in \mathbb{S}^n_+$. We refer to this game as an infinite horizon one-player linear quadratic differential game $(\mathcal{LQC})$. The function $V$ is continuously differentiable in $\reals^n$, and
\begin{eqnarray}
 \min_{u \in \Pi_u^C(x)}
\left\{L_C(x,u) +\left\langle \nabla V(x),F(x,u_C) \right\rangle \right\}
&=&
  \min_{u_C \in \mathbb{R}^{m_C}}
\left\{x^\top Q_C x+ u_C^\top R_C u_C +2x^\top P (A_C x+B_C u_C ) \right\}
\nonumber
\\&=&
  \min_{u_C \in \mathbb{R}^{m_C}}
\left\{ x^\top (Q_C+2 P A_C) x+ u_C^\top R_C u_C
+ 2 x^\top  P B_C u_C  \right\}
\nonumber
\\&=&0
\label{HJBLQGC}
\end{eqnarray}
is satisfied $ \forall x \in \reals^n$ (and attained by $\kappa_C(x)=-R_C^{-1}B_C^\top P x$)
if the continuous time algebraic Riccati equation (CARE)
\begin{equation}
0=Q_C +2PA_C- P B_C R_C^{-1}B_C ^\top P
\label{Riccati}
\end{equation}
is satisfied. 
Given that $V$ is continuously differentiable, and that (\ref{Bellmanih}) and (\ref{HJBih}) are satisfied thanks to  (\ref{HJBLQGC}) and $D=\emptyset$, from Theorem \ref{thHJBih}, the value function is
\begin{eqnarray}
\mathcal{J}^*(\xi) &:=& 
 \xi^\top P \xi.
\label{ContLCost}
\end{eqnarray}
for any $\xi \in \reals^n$.
%
\IfPers{Given that this is a purely continuous system, the calculation of the cost for when applying the feedback law $\kappa_C(x)=-R_C^{-1}B_C^\top P x$, which renders a solution $(\phi,u_C^*)$ with $\dom \phi \ni (t,j) \mapsto u_C^*(t,j)= \kappa_C(\phi(t,j))$, is as follows
\begin{equation*}
\mathcal{J}(\xi,u_C^*) =
 \int_{t_0}^{\infty} L_C(\phi(t,j),u_{C}^*(t,j))dt  + \underset{(t,j) \in \textup{dom}\phi}{\limsup_{t+j\rightarrow \infty} } q(\phi(t,j))
\end{equation*}
and, when (CARE) is satisfied, one has
\begin{eqnarray*}
 \int_{t_0}^{\infty} L_C(\phi(t,j),u_{C}^*(t,j))dt  &=& 
 \int_{0}^{\infty} \phi^\top(t,j) Q_C \phi(t,j) +{u_{C}^*}^\top(t,j) R_C u_{C}^*(t,j)dt  
 \\&=&  \int_{0}^{\infty} \phi^\top(t,j) (Q_C +P B_C R_C^{-1}B_C ^\top P)\phi(t,j) dt  
  \\&=&  \int_{0}^{\infty} \phi^\top(t,j) (2P B_C R_C^{-1}B_C ^\top P-2PA_C)\phi(t,j) dt  
  \\&=&-\int_{0}^{\infty} \frac{d V(\phi(t,j))}{d t} dt
    \\&=&V(\phi(0,0))-\limsup_{t\rightarrow \infty}V(\phi(t,0))
  \\&=& 
  \phi^\top(0,0) P \phi(0,0)- \limsup_{t\rightarrow \infty} \phi^\top(t,0) P \phi(t,0)
\end{eqnarray*}
which, given that $\phi(0,0)=\xi$ and $q(\phi(t,j))=\phi^\top(t,0) P \phi(t,0)$, implies
$	\mathcal{J}(\xi,u_C^*) =\xi^\top P \xi$, showing that the feedback law $\kappa_C$ attains the optimal cost.} 
From Corollary \ref{FBLrem}, the feedback law $\kappa_C$ is a pure strategy Nash equilibrium for the $(\mathcal{LQC})$ game.
\end{example}
\begin{example}{Linear Quadratic Difference Games}
Consider a one-player hybrid game with state $x \in \reals^n$, input $u=(u_C, u_D) \in \reals^{m_C}\times \reals^{m_D}$, and dynamics $\HS$ as in (\ref{Heq}) with
\begin{eqnarray}
D&=&\mathbb{R}^n \times \mathbb{R}^{m_D} \nonumber
\\
G(x,u_D)&=&A_D x+B_D u_D
\label{LinearDGame}
\end{eqnarray}
and flow set $C=\emptyset$, where $A_D \in \mathbb{R}^{n\times n}$ and $B_D \in \mathbb{R}^{n\times m_D}$. 
Consider the stage costs $L_C(x,u_C):=0$ and $L_D(x,u_D):=x^\top Q_D x+ u_D^\top R_D u_D$, the terminal cost $q(x):=x^\top P x$ and the function $V(x):=x^\top P x$, where $Q_D\in \mathbb{S}^n_+$, $R_D\in \mathbb{S}^{m_D}_+$ and  $P \in \mathbb{S}^n_+$. We refer to this game as an infinite horizon one player linear quadratic difference game $(\mathcal{LQD})$. The function $V$ is such that
\begin{eqnarray}
\min_{u \in \Pi_u^D(x)}
&&\left\{ 
 L_D(x,u)
+ V(G(x,u)) \right\} \nonumber
\\&=&
  \min_{u_D \in \mathbb{R}^{m_D}}
\left\{x^\top Q_D x+ u_D^\top R_D u_D +(A_D x+B_D u_D )^\top P (A_D x+B_D u_D ) \right\}
\nonumber
\\&=&
  \min_{u_D \in \mathbb{R}^{m_D}}
\left\{  x^\top (Q_D+A_D^\top P A_D) x+ u_D^\top (R_D +B_D^\top P B_D) u_D
+ 2 x^\top A_D^\top  P B_D u_D\right\}
\nonumber
\\&=&x^\top P x,
\label{HJBLQGD}
\end{eqnarray}
is satisfied $ \forall x \in \reals^n$ (and attained by $\kappa_D(x)=-(R_D+B_D^\top PB_D)^{-1}B_D^\top P A_D x$)
when the discrete-time algebraic Riccati equation (DARE)
\begin{equation}
P=Q_D+A_D^\top P A_D-A_D^\top P B_D (R_D+B_D^\top P B_D)^{-1}B_D^\top P A_D
\label{DARE}
\end{equation}
is satisfied. 
Given that $V$ is continuously differentiable, and that (\ref{Bellmanih}) and (\ref{HJBih}) are satisfied thanks to $C=\emptyset$ and (\ref{HJBLQGD}), from Theorem \ref{thHJBih}, the value function is

\begin{eqnarray}
\mathcal{J}^*(\xi) &:=& 
 \xi^\top P \xi
\label{ContLCost}
\end{eqnarray}
for any $\xi \in \reals^n$.
\IfPers{Given that this is a purely discrete system, the calculation of the cost for when applying the feedback law $\kappa_D(x)=-(R_D+B_D^\top PB_D)^{-1}B_D^\top P A_D x$, which renders a solution $(\phi,u_D^*)$ with $\dom \phi \ni (t,j) \mapsto u_D^*(t,j)= \kappa_D(\phi(t,j))$, is as follows
\begin{equation*}
\mathcal{J}(\xi,u_D^*) =
  \sum_{j=0}^{\infty} L_D(\phi(t_{j+1},j),u_{D}^*(t_{j+1},j))+ \underset{(t,j) \in \textup{dom}\phi}{\limsup_{t+j\rightarrow \infty} } q(\phi(t,j))
\end{equation*}
and, when (DARE) is satisfied, one has
\begin{eqnarray*}
 \sum_{j=0}^{\infty} L_D(\phi(t_{j+1},j),u_{D}^*(t_{j+1},j)) &=& 
 \sum_{j=0}^{\infty} \phi^\top(t,j)Q_D \phi(t,j) +{u_{D}^*}^\top(t,j) R_D u_{D}^*(t,j)  
 \\&=&  \sum_{j=0}^{\infty} \phi^\top(t,j) (Q_D +A_D^\top P B_D (R_D+B_D^\top P B_D)^{-1}B_D^\top P A_D)\phi(t,j)   
  \\&=&  \sum_{j=0}^{\infty} \phi^\top(t,j) (2A_D^\top P B_D (R_D+B_D^\top P B_D)^{-1}B_D^\top P A_D+P-A_D^\top PA_D)\phi(t,j)  
  \\&=&\sum_{j=0}^{\infty} V(\phi(t,j))-V(G(\phi(t,j),u_D^*(t,j))
  \\&=&V(\phi(0,0))-\limsup_{j\rightarrow \infty}V(\phi(0,j))
  \\&=& 
  \phi^\top(0,0) P \phi(0,0)- \limsup_{j\rightarrow \infty} \phi^\top(0,j) P \phi(0,j)
\end{eqnarray*}
which, given that $\phi(0,0)=\xi$ and $q(\phi(t,j))=\phi^\top(0,j) P \phi(0,j)$, implies
$	\mathcal{J}(\xi,u_D^*) =\xi^\top P \xi$, showing that the feedback law $\kappa_D$ attains the optimal cost.} 
From Corollary \ref{FBLrem}, the feedback law $\kappa_D$ is a pure strategy Nash equilibrium for the $(\mathcal{LQD})$ game.
\end{example}
\begin{example}{Numerical  Hybrid Game with Nonunique Solutions}
Consider a one-player hybrid game with state $x \in \reals$, input $u=(u_C, u_D) \in \reals \times \reals$, and dynamics $\HS$ as in (\ref{Heq}) described by
\begin{eqnarray}
C&:=&[0,2] \nonumber
\\
F(x,u_C)&:=&a x_1+b u_C \nonumber
\\
D&:=&\{1\} \nonumber
\\
G(x,u_D)&:=&0.5
\label{nuHGame}
\end{eqnarray}
where $a=-1,b=-1$. Consider the stage costs $L_C(x,u_C):=x^2Q_C+u_C^2 R_C $ and $L_D(x,u_D):=P(x^2-0.5^2)$, the terminal cost $q(x):=P x^2$ and the function $V(x):=P x^2$, where $Q_C=1$, $R_C=1.304$,and  $P =0.429$.  The function $V$ is such that
\begin{eqnarray}
 &&\min_{u_C \in \Pi_u^C(x)}
\left(L_C(x,u_C) +\left\langle \nabla V(x),F(x,u_C) \right\rangle \right) 
\nonumber\\&=&
\min_{u_C \in \mathbb{R}}
\left( x^\top (Q+ 2P a) x+ u_C^\top R_C u_C
+  2x^\top P b u_C \right)
 \nonumber \\ &=&
\min_{u_C \in \mathbb{R}}
\left( x^2 0.142+ 1.304u_C^2 
+  0.858x u_C \right)\nonumber \\
&=&0
\label{HJBLQGnusC}
\end{eqnarray}
holds $ \forall x \in \reals^n$ (and attained by $\kappa_C(x)=-R_C^{-1}b P x$),
leading $V(x)= Px^2$ as a solution to (\ref{HJBih}). In addition, the function $V$ is such that
\begin{eqnarray}
&&\min_{u \in \Pi_u^D(x)}
\left\{
 L_D(x,u)
+ V(G(x,u)) \right\} \nonumber
\\&=&
  \min_{u_D \in \mathbb{R}^{m_D}}
\left\{P(x^2-0.5^2)+P0.5^2 \right\} 
\nonumber
\\&=& P x^2,
\label{HJBLQGnusD}
\end{eqnarray}
 which leads $V(x)= Px^2$ as a solution to (\ref{Bellmanih}).
Thus, given that $V$ is continuously on a neighborhood of $\Pi(C)$, and that (\ref{HJBih}) and (\ref{Bellmanih}) hold thanks to (\ref{HJBLQGnusC}) and (\ref{HJBLQGnusD}), from Theorem \ref{thHJBih} the value function is
\begin{eqnarray}
\mathcal{J}^*(\xi) &:=& 
 P \xi^2
\label{HyLCost}
\end{eqnarray}
for any $\xi \in \Pi(\overline{C}) \cup \Pi(D)$.
Notice that, since $\Pi(\overline{C}) \cap \Pi(D) \neq \emptyset$, the set of maximal solutions from $x=2$ is given by $\mathcal{S}_\kappa(2)=\{\phi_\kappa,\phi_h\}$, where $\phi_\kappa$ has a domain $ \dom \phi_\kappa=\realsgeq \times \{0\}$, and is given by
\begin{equation*}
  \phi_\kappa(t,0)=  2 e^{-1.3289 t} 
\end{equation*}
In simple words, $\phi_\kappa $ flows from $x=2$ up to $x=0$. The maximal solution $\phi_h$ has a domain $ \dom \phi_h=[0,2.10747)\times \{0\} \cup [2.10747, \infty)\times \{1\}$, and is given by
\begin{equation*}
 \phi_h(t,0)= 2 e^{-1.3289 t} ,\hspace{0.5cm} \phi_h(t,1)=
0.5 e^{-1.3289 (t-2.10747)} 
\end{equation*}
In simple words, the trajectory $\phi_h$ flows from $x=2$ up to $x=1$, then it jumps to $x=0.5$ and flows up to $x=0$.\\
%
By denoting the corresponding input signals as $u_\kappa=\kappa(\phi_\kappa)$ and $u_h=\kappa(\phi_h)$, and given that the cost of flowing from $x=1$ to $x=0.5$ is equal to the cost of jumping, we have that the costs of the solutions $(\phi_\kappa,u_\kappa)$ and $(\phi_h,u_h)$ are equal to $P2^2$. 
%
%
Figure \ref{fig:my_label} illustrates this behavior.
\begin{figure}[H]
    \centering
    \includegraphics[width=15cm, center]{Figures/Counterexample_[211]Fx.png}
    \caption{Nonunique solutions attaining optimal cost. Continuous solution (green). Hybrid solution (blue-red).}
    \label{fig:my_label}
\end{figure}
This aligns to the Nash equilibrium Definition \ref{NashEqCoop}, and Assumption \ref{SolExistence}, since every maximal solution that $\kappa$ renders from $\xi=2$ attains the minimum value of the cost  $\mathcal{J}$ over the set of input actions that yield maximal solutions to $\HS$ with data (\ref{nuHGame}).
\label{NumericalEx}
\end{example}
\subsection{Linear Quadratic Hybrid Games with Periodic Jumps}
Given a time $\bar{T}\in \mathbb{R}$, consider a one-player hybrid game with state $x=(x_p,\tau)   \in \reals^n \times [0, \bar{T}]$, input $u=(u_C, u_D) \in \reals^{m_C}\times \reals^{m_D}$, and dynamics $\HS$ as in (\ref{Heq}) described by
\begin{eqnarray*}
C&:=&\mathbb{R}^n \times [0,\bar{T}] \times \mathbb{R}^{m_C}
\\
F&:=&(A_C x_p+B_C u_C,1)
\\
D&:=&\mathbb{R}^n \times \{\bar{T}\} \times \mathbb{R}^{m_D}
\\
G&:=&(A_D x_p+B_D u_D,0)
\end{eqnarray*}
where $A_C,A_D \in \mathbb{R}^{n\times n}$ and $B_C \in \mathbb{R}^{n\times m_C}$, $B_D \in \mathbb{R}^{n \times m_D}$.
Consider the stage costs $L_C(x,u_D):=x_p^\top Q_C x_p+ u_C^\top R_C u_C$ and $L_D(x,u_D):=x_p^\top Q_D x_p+ u_D^\top R_D u_D$, the terminal cost $q(x):=x_p^\top P(\tau) x_p$ and the function $V(x):=x_p^\top P(\tau) x_p$, where $Q_C,Q_D\in \mathbb{S}^n_+$, $R_C \in \mathbb{S}^{m_C}_+,R_D\in \mathbb{S}^{m_D}_+$ and  $P(\tau)\in \mathbb{S}^n_+ \forall \tau \in [0,\bar{T}]$. We refer to this game as an infinite horizon one player linear quadratic game with periodic jumps $(\mathcal{LQH})$.

\begin{corollary}{Hybrid Riccati equation} Given $\bar{T}\in \mathbb{R}$, $A_C,A_D \in \mathbb{R}^{n\times n}$, $B_C \in \mathbb{R}^{n\times m_C}$, $B_D \in \mathbb{R}^{n \times m_D}$, $Q_C,Q_D\in \mathbb{S}^n_+$, $R_C \in \mathbb{S}^{m_C}_+,R_D\in \mathbb{S}^{m_D}_+$,
suppose there exists a matrix $P:[0,\bar{T}]\rightarrow \mathbb{S}^n_+ $ continuously differentiable such that
\begin{equation}
-\frac{d }{d \tau}P(\tau) =Q_C+P(\tau) A_C+A_C^\top P(\tau)-P(\tau) B_C R_C^{-1}B_C^\top P(\tau)  \hspace{1cm} \forall \tau  \in (0,\bar{T})
\label{DiffRiccatiih} 
\end{equation}
and 
\begin{equation}
R_D+B_D^\top P B_D \in \mathbb{S}_{0+}^{m_D}
\label{lqeqinv}
\end{equation}
with 
\begin{equation}
P (\bar{T})=Q_D+A_D^\top P(0) A_D-A_D^\top P(0) B_D (R_D+B_D^\top P(0) B_D)^{-1}B_D^\top P(0) A_D
\label{Riccatiih} 
\end{equation}

Then, the feedback law $\kappa:=(\kappa_C,\kappa_D)$, with values
\begin{equation}
\kappa_C(x)=-R_C^{-1}B_C^\top P(\tau) x_p \hspace{1cm} \forall x \in \Pi(C)
\label{NashkCLQih}
\end{equation}
\begin{equation}
\kappa_D(x)=-(R_D+B_D^\top P(0) B_D)^{-1}B_D^\top P(0) A_D x_p \hspace{1cm} \forall x \in \Pi(D)
\label{NashkDLQih}
\end{equation}
is a Nash equilibrium for the $(\mathcal{LQH})$ game. In addition, given  $\xi=(\xi_p,\xi_\tau) \in \Pi(\overline{C}) \cup \Pi(D)$, the Nash equilibrium outcome is equal to $\xi_p^\top P(\xi_\tau) \xi_p$.
\end{corollary}
\begin{proof}
We show that when conditions (\ref{DiffRiccatiih})-(\ref{Riccatiih}) hold, by using the result in Theorem \ref{thHJBih},  the function $V$ is equal to the value function for this game. By using Corollary \ref{FBLrem}, we show that the feedback law as in (\ref{NashkCLQih}), (\ref{NashkDLQih}) attains such a cost. 
We can write (\ref{HJBih}) in Theorem \ref{thHJBih} as
\begin{eqnarray}
&&0=\min_{u \in \Pi_u^C(x)} \mathcal{L}_C(x,u_C),
 \nonumber
\\ 
&&\mathcal{L}_C(x,u_C)= 
 x_p^\top Q_C x_p+ u_C^\top R_C u_C
+2x_p^\top P(\tau) (A_C x_p+B_C u_C)  
+ x_p^\top \frac{{d} }{{d} \tau}P(\tau) x_p
\label{HJBihlq}
\end{eqnarray}
First, 
given that (\ref{DiffRiccatiih}) holds, and $x_p^\top(P(\tau) A_C+A_C^\top P(\tau))x_p=x_p^\top (2P(\tau)A)x_p$ for every $x\in C$, one has
\begin{equation}
\mathcal{L}_C(x,u_C)=
 x_p^\top (
 P(\tau) B_C R_C^{-1}B_C^\top P(\tau)) x_p+ u_C^\top R_C u_C
+ 2 x_p^\top  P(\tau) B_C u_C 
\end{equation}
The first order necessary condition for optimality
\begin{equation*}
\frac{\partial}{\partial u_C}\left.\left( 
 x_p^\top (
 P(\tau) B_C R_C^{-1}B_C^\top P(\tau)) x_p+ u_C^\top R_C u_C
+ 2 x_p^\top  P(\tau) B_C u_C  \right)\right|_{u_C^*} =0
\end{equation*}
is satisfied by the point
\begin{equation}
u_C^*=-R_C^{-1}B_C^\top P(\tau) x_p
\label{LQHPu*}
\end{equation}
Given that $R_C \in \mathbb{S}^{m_D}_+$, the second-order necessary condition for optimality 
\begin{equation*}
\frac{\partial^2}{\partial u_C^2}\left.\left( 
 x_p^\top (
  P(\tau) B_C R_C^{-1}B_C^\top P(\tau)) x_p+ u_C^\top R_C u_C
+ 2 x_p^\top  P(\tau) B_C u_C  \right)\right|_{u_C^*} \succeq 0,
\end{equation*}
holds, rendering $u_C^*$ in (\ref{LQHPu*}) as a minimizer of the optimization problem in (\ref{HJBihlq}).
In addition, 
it satisfies $\mathcal{L}_C(x,u_C^*)=0$, rendering $V(x)=x_p^\top P(\tau)x_p$ to be a solution of (\ref{HJBih}) in Theorem \ref{thHJBih}. 

On the other hand, we can write (\ref{Bellmanih}) in Theorem \ref{thHJBih} as
\begin{eqnarray}
&&x_p^\top P(\bar{T}) x_p =\min_{u \in \Pi_u^D(x)} \mathcal{L}_D(x,u_D),
 \nonumber
\\ 
&&\mathcal{L}_D(x,u_D)= 
x_p^\top Q_D x_p+ u_D^\top R_D u_D
+ (A_D x_p+B_D u_D)^\top P(0) (A_D x_p+B_D u_D) 
\label{Bellmanihlq}
\end{eqnarray}
which can be expanded as
\begin{equation}
\mathcal{L}_D(x,u_D)=  
 x_p^\top (Q_D+A_D^\top P(0) A_D) x_p+ u_D^\top (R_D +B_D^\top P(0) B_D) u_D
+ 2 x_p^\top A_D^\top  P(0) B_D u_D 
\end{equation}
The first order necessary condition for optimality
\begin{equation*}
\frac{\partial}{\partial u_D}\left.\left( 
 x_p^\top (Q_D+A_D^\top P(0) A_D) x_p+ u_D^\top (R_D +B_D^\top P(0) B_D) u_D
+ 2 x_p^\top A_D^\top  P(0) B_D u_D  \right)\right|_{u_D^*} =0
\end{equation*}
is satisfied by the point
\begin{equation}
u_D^*=-(R_D+B_D^\top P(0) B_D)^{-1}B_D^\top P(0) A_D x_p
\label{LQHPud*}
\end{equation}
Given that (\ref{lqeqinv}) holds, the second-order necessary condition for optimality 
\begin{equation*}
\frac{\partial^2}{\partial u_D^2}\left.\left( 
 x_p^\top (Q_D+A_D^\top P(0) A_D) x_p+ u_D^\top (R_D +B_D^\top P(0) B_D) u_D
+ 2 x_p^\top A_D^\top  P(0) B_D u_D  \right)\right|_{u_D^*} \succeq 0,
\end{equation*}
is satisfied, rendering $u_D^*$ as in (\ref{LQHPud*}) as a minimizer of the optimization problem in (\ref{Bellmanihlq}).
In addition,  $u_D^*$
satisfies $\mathcal{L}_D(x,u_D^*)=x_p^\top P(\bar{T})x_p$, with $P(\bar{T})$ as in (\ref{Riccatiih}),  rendering $V(x)=x_p^\top P(\tau)x_p$ to be a solution of (\ref{Bellmanih}) in Theorem \ref{thHJBih}. 

Then, from Theorem \ref{thHJBih} and Corollary (\ref{FBLrem}), for every $ x=(x_p,\tau) \in \Pi(\overline{C}) \cup \Pi(D)$ the value function is 
\begin{equation*}
\J^*(x)= x_p^\top P(\tau)x_p
\end{equation*}
and $\kappa=(\kappa_C,\kappa_D)$  as in (\ref{NashkCLQih}), (\ref{NashkDLQih}) is a pure strategy Nash equilibrium.
%
%
%
%
\end{proof}
%
%
\subsection{Lyapunov Functions for Hybrid Games}
Next, we introduce some definitions to present a result that connects optimality and stability properties for hybrid noncooperative games.
\begin{definition}{Class-$\mathcal{K}_\infty$ functions}
	A function $\rho: \mathbb{R}_{\geq 0} \rightarrow  \mathbb{R}_{\geq 0}$ is a class 
	$\mathcal{K}_\infty$ function, also written as $\rho \in \mathcal{K}_\infty$, if $\rho$ is
	 zero at zero, continuous, strictly increasing, and unbounded.
\end{definition}

\begin{definition}{Positive definite functions}
Consider a $\kappa:=(\kappa_C,\kappa_D): \mathbb{R}^n \rightarrow \mathbb{R}^{m_C} \times \mathbb{R}^{m_D}$. 
We say that a function $\rho: \reals_{\geq 0}\rightarrow \reals_{\geq 0}$ is positive definite, also written as $\rho \in \mathcal{PD}$, if $\rho(s)>0$ and $\rho(0)=0$. We say that a function $\rho: \reals^n \times \reals^{m_C} \rightarrow \reals_{\geq 0}$  is positive definite with respect to a set $\A \subset \reals^n$, in composition with $\kappa_C$, also written as $\rho \in \mathcal{PD}_{\kappa_C}(\A)$, if $\rho(s,\kappa_C(s)) > 0$ for all $s \notin \A$ and $\rho(\A, \kappa_C(\A)) = \{0\}$. Likewise, we say that a function $\rho: \reals^n \times \reals^{m_D} \rightarrow \reals_{\geq 0}$  is positive definite with respect to a set $\A \subset \reals^n$, in composition with $\kappa_D$, also written as $\rho \in \mathcal{PD}_{\kappa_D}(\A)$, if $\rho(s,\kappa_D(s)) > 0$ for all $s \notin \A$ and $\rho(\A, \kappa_D(\A)) = \{0\}$.
\end{definition}

\begin{corollary}{Nash equilibrium under the existence of a Lyapunov function} 
Consider a one-player hybrid game with closed-loop dynamics $\HS_{\kappa}$ as in (\ref{Hkeq}) with data $(C,F,D,G)$ and $\kappa:=(\kappa_C,\kappa_D): \mathbb{R}^n \rightarrow \mathbb{R}^{m_C} \times \mathbb{R}^{m_D}$ such that $C_\kappa = \Pi(C), D_\kappa = \Pi(D)$ and \pnn{suppose that Assumption \ref{AssLipsZ} holds.} Given a closed set $\mathcal{A} \subset \reals^n$, continuous stage costs $L_C:
\pnn{C} \rightarrow \mathbb{R}_{\geq 0}, L_D:
\pnn{D} \rightarrow \mathbb{R}_{\geq 0}$, such that
$L_C \in \mathcal{PD}_{\kappa_C}(\mathcal{A}\cap C_\kappa)$, $L_D \in \mathcal{PD}_{\kappa_D}(\mathcal{A}\cap D_\kappa)$, and a terminal cost $q:\reals^n \rightarrow \reals$, suppose there exist a function $V: \mathbb{R}^n \rightarrow \mathbb{R}$ that is continuously differentiable on an open set containing $\Pi(\overline{C})$,  satisfying (\ref{kHJBeqiha}), (\ref{kHJBeqihb}), (\ref{kBeqihc}) and (\ref{kBeqihd}), and such that for each $\xi \in \Pi(\overline{C}) \cup \Pi(D)$, for each $\phi \in \mathcal{S}_{\HS_{\kappa}  }^\infty(\xi)$, (\ref{TerminalCond}) holds. If there exist $\alpha_1,\alpha_2 \in \mathcal{K}_\infty$ such that for all $x \in \Pi(\overline{C}) \cup \Pi(D)$
\begin{equation}
\alpha_1(|x|_\A) \leq V(x) \leq \alpha_2(|x|_\A)
\label{alphabound}
\end{equation}
\pnn{and either
\begin{enumerate}
\item $\A \cap C_\kappa \neq \emptyset $ and $\A \cap D_\kappa \neq \emptyset $;
\item $\A \cap C_\kappa \neq \emptyset $, $\A \cap D_\kappa = \emptyset $, and there exists a continuous $\lambda \in \mathcal{PD}$ such that $L_D(x,\kappa_D(x)) \geq \lambda(|x|_\A)$ for all $x \in D_\kappa$;
\item $\A \cap C_\kappa = \emptyset $, $\A \cap D_\kappa \neq \emptyset $, and there exists a continuous $\lambda \in \mathcal{PD}$ such that $L_C(x,\kappa_D(x)) \geq \lambda(|x|_\A)$ for all $x \in C_\kappa$;
\end{enumerate}}
 then,   
\begin{equation}
\J^*(\xi)= V(\xi),
\end{equation}
for each $\xi \in \Pi(\overline{C}) \cup \Pi(D)$. The feedback law
$\kappa$ is the Nash equilibrium, and it renders $\A$ uniformly globally asymptotically stable for $\HS_{\kappa}$ in $\Pi(\overline{C}) \cup \Pi(D)$. 
\label{CorStability}
\end{corollary}
\begin{proof}
Given that $C_\kappa = \Pi(C), D_\kappa = \Pi(D)$, and $L_C,L_D,V,\kappa:=(\kappa_C,\kappa_D)$ are such that (\ref{kHJBeqiha}), (\ref{kHJBeqihb}), (\ref{kBeqihc}) and (\ref{kBeqihd}) hold, then, $V$ and $\kappa$ satisfy (\ref{HJBih}), (\ref{Bellmanih}), (\ref{kHJBeqih}), and (\ref{kBeqih}). 
Since (\ref{TerminalCond}) holds for each $\xi \in \Pi(\overline{C}) \cup \Pi(D)$, for each $\phi \in \mathcal{S}_{\HS_{\kappa}  }^\infty(\xi)$, we have from Theorem \ref{thHJBih} that $V$ is the value function as in (\ref{cost2goih}) for $\HS_{\kappa}$ at $\Pi(\overline{C}) \cup \Pi(D)$. 
From Corollary \ref{FBLrem} we have that the feedback law
 $\kappa:=(\kappa_C,\kappa_D)$ that satisfies (\ref{kHJBeqih}), (\ref{kBeqih}) is the pure strategy Nash equilibrium for this game. 
In addition, given that $V$ is a Lyapunov candidate for $\HS_{\kappa}$ \cite[Definition 3.16]{65}, 	the functions $\alpha_1, \alpha_2$ satisfy (\ref{alphabound}), and given that (\ref{kHJBeqihb}), (\ref{kBeqihd}) hold, then, from \cite[Theorem 3.18]{65} we have that $\A$ is uniformly globally asymptotically stable for $\HS_{\kappa}$. \pnn{Notice that the role of $\rho$ therein is played by the stage costs $L_C \in \mathcal{PD}_{\kappa_C}(\A \cap C_\kappa)$, $L_D \in \mathcal{PD}_{\kappa_D}(\A \cap D_\kappa)$ or the function $\lambda \in \mathcal{PD}$.}
 \end{proof}

\begin{example}{Linear Quadratic Differential Games}
Taking up the $(\mathcal{LQHC})$ game, recall that $\kappa:=(-R_C^{-1} B_C^\top P x,\cdot)$ and consider the case in which $(A_C-B_C R_C^{-1}B_C^\top P)$ is Hurwitz and $\A=\{0\}$. Given that $L_C \in \mathcal{PD}_{\kappa_C}(\A \cap C_\kappa)$, $V$ is continuously differentiable, and (\ref{kHJBeqiha}) and (\ref{kHJBeqihb}) hold, by setting, $\alpha_1(|x|_\A)=x^\top (P-I) x$ and $\alpha_2(|x|_\A)=x^\top (P+I) x$, from Corollary \ref{CorStability}, we have that $\kappa$ is the Nash equilibrium and renders $\A$ uniformly globally asymptotically stable for $\HS_{\kappa}$.
\end{example}
\begin{example}{Linear Quadratic Difference Games}
Taking up the $(\mathcal{LQHD})$ game, recall that $\kappa:=(\cdot,-(R_D+B_D^\top PB_D)^{-1}B_D^\top P A_D x)$ and consider the case in which $(A_D-B_D (R_D+B_D^\top PB_D)^{-1}B_D^\top P A_D )$ is Schur and $\A=\{0\}$. Given that $L_D \in \mathcal{PD}_{\kappa_D}(\A \cap D_\kappa)$, and  (\ref{kBeqihc}), and (\ref{kBeqihd}) hold, by setting $\alpha_1(|x|_\A)=x^\top (P-I) x$ and $\alpha_2(|x|_\A)=x^\top (P+I) x$, from Corollary \ref{CorStability}, we have that $\kappa$ is the Nash equilibrium and renders $\A$ uniformly globally asymptotically stable for $\HS_{\kappa}$.
\end{example}
\begin{example}{Numerical Hybrid Game with Nonunique Solutions}
Taking up Example \ref{NumericalEx}, recall that $a<0$ and $\kappa:=(\kappa_C,\kappa_D)=(-R_C^{-1}bPx,\cdot)$ and let $\A=\{0\}$. Given that $L_C \in \mathcal{PD}_{\kappa_C }(\A \cap C_\kappa)$, \pnn{$L_D \in \mathcal{PD}_{\kappa_D}(\A \cap D_\kappa)$ vacuously}, and (\ref{kHJBeqiha}),  (\ref{kHJBeqihb}), (\ref{kBeqihc}), and (\ref{kBeqihd}) hold, by setting $\alpha_1(|x|_\A)=x^\top (P-I) x$ and $\alpha_2(|x|_\A)=x^\top (P+I) x$, 
and \pnn{for $|x|_\A \mapsto \lambda(|x|_\A)=P \frac{x^2}{2}$, $L_D(x,\kappa_D(x)) \geq \lambda(|x|_\A)$ for all $x \in D_\kappa$,}  
from Corollary \ref{CorStability}, we have that $\kappa$ is the Nash equilibrium and renders $\A$ uniformly globally asymptotically stable for $\HS$ with data as in (\ref{nuHGame}). 
\end{example}}
%
%
\IfPers{Let $\mathcal{A} \subset \mathbb{R}^n$ be a closed set and consider
\begin{equation*}
	\mathcal{X}_\mathcal{A}:=\{ \phi \in \mathcal{X}: \displaystyle \lim \limits_{\substack{(t,j) \in \dom \phi\\
	(t,j)\rightarrow \sup \dom \phi}} \left | 
	\phi(t,j) \right|_\mathcal{A}=0 \}
\end{equation*}
\begin{equation*}
\mathcal{U}_\mathcal{A}(\xi) :=  \left\{ u\in \mathcal{U}: \exists  \phi \in \mathcal{R} (\xi,u) \cap \mathcal{X}_{\mathcal{A}} \right\}
\end{equation*}
To move to another section
}
%
\subsection{Problem Statement}\label{subsec: probstat}
We formulate an 
optimization problem to solve a two-player zero-sum hybrid game with \sj{variable} terminal time and a terminal set, \pno{and provide sufficient} conditions to characterize the solution. 

Following the formulation in Definition \ref{elements}, {consider a two-player zero-sum hybrid game with dynamics $\HS\IfIncd{_s}{}$ described by (\IfIncd{\ref{Hinc}}{\ref{Heq}}) \IfTp{}{with $N=2$ }
with data $(C,F,D,G)$. }
Let the closed set $X \subset \Pi(C) \cup \Pi(D)$ be the terminal constraint set. We say that a solution $(\phi,u)$ to $\HS$ is \emph{feasible} if there exists $(T,J) \in \dom (\phi,u)$ such that $\phi(T,J) \in X$. In addition, \pn{we make $(T,J)$ to be both the terminal time of $(\phi,u)$ and the first time at which $\phi$ reaches} $X$, i.e., there does not exist $(t,j) \in \dom \phi$ with $t+j < T+J$ 
such that $\phi(t,j) \in X$ \pn{and $(T,J) = \max \dom (\phi, u)$; hence $\dom \phi$ is compact.}\footnote{\pn{When $X = \emptyset$, 
the requirement that $\phi$ belongs to $X$ is not enforced, hence, there is no terminal constraint and the two-player zero-sum hybrid game {evolves} over {an} infinite (hybrid) horizon when $\dom \phi$ is unbounded\sj{, i.e., $\phi$ is complete.}}}
\IfIncd{}{Uniqueness of solutions for a given input implies a unique correspondence {from cost} to control {action}, which allows this type of games to be \textit{well-defined}, so 
{that an} equilibrium solution is defined \cite[Remark 5.3]{basar1999dynamic}.
This justifies the following assumption.
\begin{assumption}{}
  The flow map $F$ and the flow set $C$ are such that solutions to {$\dot{x}=F(x,u_C)$ $(x,u_C)\in C$} are unique for each input $u_C$. The jump map $G$ is single valued.
 \label{AssLipsZ}
\end{assumption}
Sufficient conditions to guarantee that Assumption \ref{AssLipsZ} holds include Lipschitz continuity of the flow map $F$, provided it is a single-valued function.
Under Assumption \ref{AssLipsZ}, \NotAutomss{the conditions in Proposition \ref{UniquenessHu} are satisfied, so }solutions to $\HS$ are unique\footnote{{Under Assumption \ref{AssLipsZ}, the domain of the input $u$ specifies whether from points \pn{in} $\Pi(C)\cap \Pi(D)$ a jump or flow occur.}} for each $u\in \mathcal{U}$.}

Given $\xi \in \Pi(C) \cup \Pi(D)$, a joint input action $u=(u_C, u_D)\in \mathcal{U}
$, 
the stage cost for flows $L_C:\mathbb{R}^n  \times \mathbb{R}^{m_C} \rightarrow \mathbb{R}_{\geq 0}$, the stage cost for jumps $L_D:\mathbb{R}^n  \times \mathbb{R}^{m_D} \rightarrow \mathbb{R}_{\geq 0}$, and the terminal cost $q: \mathbb{R}^n \rightarrow \mathbb{R}$, we define the cost associated to the solution\IfIncd{s}{ $(\phi,u)$} to $\HS\IfIncd{_s}{}$ 
from $\xi$, 
\IfIncd{}{under Assumption \ref{AssLipsZ},} as 
\IfIncd{\begin{equation}
\label{defJTNCincinc}
\mathcal{J}(\xi,u) := 
\sup_{\phi \in \R_s (\xi,u)}
\widetilde{\mathcal{J}}(\phi,u)
\end{equation}
where\footnote{\pn{Notice that $\mathcal{J}$ depends on the initial condition $\xi$ and input $u$, while $\widetilde{\mathcal{J}}$ depends on the solution pair $(\phi,u)$ with $\phi(0,0)=\xi$.}}}{}
%
\begin{equation}\label{defJTNCinc}
  \begin{split}
  &{\mathcal{J}(\xi,u)}:=
  \sum_{j=0}^
  {\mathclap{\substack{{\sup_j \dom \phi}} }}\quad
   \int_{t_{j}}^{t_{j+1}} L_C(\phi(t,j),u_{C}(t,j))dt  
  \\& \hspace{1.5cm}
  + \sum_{j=0}^
  {\mathclap{\substack{{\sup_j \dom \phi -1}} }} \quad
  L_D(\phi(t_{j+1},j),u_{D}(t_{j+1},j))
  \\& \hspace{1.5cm}   +\underset{(t,j) \in \textup{dom}\phi}{\limsup_{t+j\rightarrow \sup_t \dom \phi + \sup_j \dom \phi} } q(\phi(t,j))
\end{split}
\end{equation}
where {$t_{\sup_j \dom \phi+1}=\sup_t \dom \phi$ defines the upper limit of the last integral, }
and $\{t_j\}_{j=0}^{\sup_j \dom \phi}$ is a nondecreasing sequence
associated to the definition of the hybrid time domain of {$(\phi,u)$}\IfAutomss{\>\cite[Definition 2.3]{65}}{; see Definition \ref{htd}}. 
\IfIncd{The cost is defined as the largest cost of solutions from $\xi$.}{}

\pno{
\pn{When $X$ is nonempty,} the set $\mathcal{S}_{\HS}^X(\xi) \subset {\mathcal{S}}_{\HS}(\xi)$ denotes all maximal solutions from $\xi$ that reach $X$ at their terminal time. When $X$ is empty, $\mathcal{S}_{\HS}^X(\xi)$ is the set of complete solutions from $\xi$. We define the set of input actions that yield maximal solutions to $\HS\IfIncd{_s}{}$ from $\xi$ entering $X$
 as 
	  $\mathcal{U}_{\HS\IfIncd{_s}{}}^X
    (\xi):= \{u :
    \exists (\phi,u) \in 
    \mathcal{S}_{\HS\IfIncd{_s}{}}^X(\xi)
  \}$.}
{The feasible set $\mathcal{M}\subset \Pi(C) \cup \Pi(D)$} is the set of states $\xi$ such that there exists 
$(\phi, u) \in \hat{\mathcal{S}}_{\mathcal{H}  }^X(\xi) $ with $\phi(T,J)\in X $, where $(T,J)$ is the terminal time of $\dom (\phi,u)$\NotAutomss{, namely, $(T,J) = {\max \dom \phi}$}.

We \IfAutomss{now}{are ready to } formulate the two-player zero-sum game.

\color{black}
  \textit{Problem ($\diamond$):} 
Given \pn{the terminal set $X$, the feasible set $\mathcal{M}\NotAutomss{\subset \Pi(\overline{C}) \cup \Pi(D)}$,} and $\xi \in \mathcal{M}$, under Assumption \ref{AssLipsZ}, solve
\begin{eqnarray}\label{problemzsih}
 \underset{u=(u_1,u_2)  \in  \pn{\mathcal{U}_\HS^X (\xi)} }
 { \underset{u_{1}}{\textup{minimize}}\>\> \underset{u_{2}}{\textup{maximize}}}
 && \mathcal{J}(\xi, u)
 \end{eqnarray}
\sj{over the space of feedback strategies.}
\begin{remark}{Infinite horizon games}
 { \pn{When the terminal set $X$ is empty and maximal solutions are complete,}
  Problem $(\diamond)$ reduces to an infinite horizon hybrid game as in \cite{leudohygames}, 
  as stated in footnote $4$.
  %
  In this case, the feasible set satisfies $\mathcal{M} = \Pi(C) \cup \Pi(D)$ and, for each $\xi \in \mathcal{M}$, the set of complete solutions \pn{$\mathcal{S}_{\HS}^X(\xi)$} is nonempty. 
  For \IfAutomss{such}{infinite horizon\>} games, 
the set 
$\mathcal{U}_{\HS\IfIncd{_s}{}}^X$ in Problem $(\diamond)$ denotes \IfAutomss{the inputs that yield}{all joint input actions yielding}\> maximal complete solutions to $\HS\IfIncd{_s}{}$.}
\end{remark}



{
\begin{remark}{\IfIh{}{Pure strategy }\sj{Feedback} saddle-point equilibrium and min-max  
\sj{input action}}
A solution to Problem ($\diamond$), when it exists, can be expressed in terms of the \sj{feedback} 
saddle-point equilibrium $\IfIh{\kappa=(\kappa_1,\kappa_2)}{\gamma}$ for the two-player zero-sum 
game. {Each $u^*=(u^*_1, u^*_2)$ that renders a state trajectory $\phi^*\in \mathcal{R}(\xi,u^*)$}, 
with components defined as $\dom \phi^* \ni (t,j) \mapsto u^*_i(t,j)= \IfIh{\kappa}{\gamma}_i(\IfIh{}{t,j,}\phi^*(t,j))$ for each $i\in \IfTp{\{1,2\}}{\mathcal{V}}$, satisfies
\NotAutom{\begin{eqnarray*}
u^*=  
\underset{u=(u_1,u_2) \in \mathcal{U}_{\HS\IfIncd{_s}{}}^\infty(\xi)}
{\arg \min_{u_{1}} \max_{u_{2}}}
 \mathcal{J}(\xi, u)
=
 \underset{u=(u_1,u_2) \in \mathcal{U}_{\HS\IfIncd{_s}{}}^\infty(\xi)}
{\arg  \max_{u_{2}}\min_{u_{1}}}
 \mathcal{J}(\xi, u)
\end{eqnarray*}}
\begin{eqnarray*}
  u^*=
  \underset{\mathclap
  {\substack{u=(u_1,u_2) \in \pn{\mathcal{U}_\HS^X (\xi)}
  }}}
  {\arg \min_{u_{1}} \max_{u_{2}}}
  \>\> \mathcal{J}(\xi, u)
  =
  \underset{\mathclap
  {\substack{u=(u_1,u_2) \in \pn{\mathcal{U}_\HS^X (\xi)}
  }}}{
  {\arg  \max_{u_{2}}\min_{u_{1}}}}
  \>\> \mathcal{J}(\xi, u)
  \end{eqnarray*}
and it is referred to as a min-max \sj{input action} at $\xi$.
\end{remark}

\begin{definition}{Value function}
\NotAutom{  
  Given $\xi \in \Pi(\overline{C}) \cup \Pi(D)$, 
the value function at $\xi$, {when it exists,} is given by 
\begin{equation}
  \begin{split}
\mathcal{J}^*(\xi):=
\underset{u=(u_1,u_2) \in \mathcal{U}_{\HS\IfIncd{_s}{}}^\infty(\xi)}
{\min_{u_{1}} \max_{u_{2}}} \>\>
\mathcal{J}(\xi,u)
\\
=\underset{u=(u_1,u_2) \in \mathcal{U}_{\HS\IfIncd{_s}{}}^\infty(\xi)}
{\max_{u_{2}} \min_{u_{1}} } \>\>
\mathcal{J}(\xi,u)
  \end{split}
 \label{cost2gozsih}
\end{equation}}
Given \pn{the terminal set $X$, the feasible set $\mathcal{M}\subset \Pi(\overline{C}) \cup \Pi(D)$,} and $\xi \in \mathcal{M}$, under Assumption \ref{AssLipsZ}, the value function at $\xi$ is given by 
\begin{equation}
\mathcal{J}^*(\xi):=
\underset{\mathclap
{\substack{u=(u_1,u_2) \in \pn{\mathcal{U}_\HS^X (\xi)}
}}}
{\min_{u_{1}} \max_{u_{2}}}\>\>
\mathcal{J}(\xi, u)
=
\underset{\mathclap
{\substack{u=(u_1,u_2) \in \pn{\mathcal{U}_\HS^X (\xi)}
}}}
{\max_{u_{2}}\min_{u_{1}}} \>\>
\mathcal{J}(\xi, u)
 \label{cost2gozsih}
\end{equation}
\label{ValueFunctionz}
\end{definition}}
%
%
%
%
%
\vspace{-0.4cm}
\IfIncd{
\subsection{Cost Upper Bounds}
In general, the cost evaluation tools employed in approaches based on dynamic programming fall short to characterize strategies to attain a saddle point equilibrium solution for a two-player zero-sum game with dynamics given by hybrid inclusions. 
\pnn{The classical conditions involved therein do not guarantee the existence of a
lower value function when player $P_2$ plays optimally.}
Thus,{ to avoid setting an ill-defined game, in this section, \pnn{we 
provide sufficient conditions to find an upper value function, namely, to solve the following problem.

\textit{Problem ($\star$):} 
Given $\xi \in \mathbb{R}^n$ and $u_1^* \in \mathcal{U}_1$, 
solve
\begin{eqnarray}
\underset{(u_1^*,u_2) \in \mathcal{U}_{\HS\IfIncd{_s}{} }^\infty(\xi)}
{ \underset{u_{2}}{\textup{maximize}}}
&& \mathcal{J}(\xi, (u_1^*,u_2))
\label{problemzsih}
\end{eqnarray}
where $u_1^*$ is the optimal action selected by player $P_1$
and $\mathcal{U}_{\HS\IfIncd{_s}{} }^\infty$ is the set of joint input actions yielding maximal solutions to $\HS\IfIncd{_s}{} $ complete, as defined in Section 2.1.} }
\IfPers{In this section, we provide sufficient conditions to solve Problem $(\diamond)$ via finding a control strategy that minimizes the largest-possible (due to the set-valued dynamics) cost under the maximizing adversarial action and evaluating it.} 
First, we provide pointwise conditions that allow to upper bound the cost for a given control action.

\begin{proposition}{Upper Bound for a given Control Action}
  Given a system with dynamics $\HS\IfIncd{_s}{} $ as in (\IfIncd{\ref{Hinc}}{\ref{Heq}}) \IfTp{}{with $N=2$, }with data $(C,F,D,G)$,
  stage costs $L_C:\mathbb{R}^n  \times \mathbb{R}^{m_C} \rightarrow \mathbb{R}_{\geq 0}$ {and} $L_D:\mathbb{R}^n  \times \mathbb{R}^{m_D} \rightarrow \mathbb{R}_{\geq 0}$, and terminal cost $q: \mathbb{R}^n \rightarrow \mathbb{R}$, suppose there exists a function $V: \mathbb{R}^n \rightarrow \mathbb{R}$ that is continuously differentiable on a neighborhood of $\Pi(C)$ such that
  \begin{equation}
    \begin{medsize}
      \begin{split}
    0\geq  
   L_C(x,u_C) 
  +{\underset{f \in F(x,u_C)}{\sup}} \left\langle \nabla V(x),f \right\rangle 
  \hspace{1cm} \forall (x,u_{C}) \in C, 
\end{split}
\end{medsize}
\label{ubcolinc}
  \end{equation}
  \begin{equation}
    \begin{medsize}
      \begin{split}
        V(x)\geq 
         L_D(x,u_D)
        + {\underset{g \in G(x,u_D)}{\sup}}
        V(g)
        \hspace{1cm} \forall (x,u_{D}) \in D .
      \end{split}
    \end{medsize}
    \label{ubdolind}
  \end{equation}
  Let $(\phi,u)$ be any solution to $\HS\IfIncd{_s}{} $ from $\xi \in \Pi(\overline{C}) \cup \Pi(D)$. Then, 
  \begin{equation}
  \begin{medsize}
    {\widetilde{\mathcal{J}}}(\phi,u) \leq V(\xi).
    \end{medsize}
  \end{equation} 
  \label{UpperBoundOpenLoop}
\end{proposition}
\begin{proof}
From \eqref{ubcolinc}, \pnn{and given a solution $(\phi,u)$ to $\HS\IfIncd{_s}{} $ from $\xi \in \Pi(\overline{C}) \cup \Pi(D)$}, for each $j \in \nats$ such that $I_{\phi}^j=[t_j,t_{j+1}]$ has a nonempty interior int$I_{\phi}^j$, we have, for all $t \in \textup{int}I_{\phi}^j$,
\begin{equation}
\begin{medsize}
\begin{split}
  L_C({\phi(t,j),u_C(t,j)}) 
 +\frac{d}{d t}V(\phi(t,j)) \hspace{4cm}\\\leq
    L_C(\phi(t,j),u_C(t,j)) 
 +{\underset{f \in F({\phi(t,j),u_C(t,j)})}{\sup}}\left\langle \nabla V(\phi(t,j)),f \right\rangle \leq 0
 \end{split}
 \end{medsize}
 \label{eq:PrHamiltonC1}
\end{equation}
\NotAutom{by integrating over the interval $[t_j,t_{j+1}]$, 
we obtain
\begin{eqnarray*}
0\geq  \int_{t_{j}}^{t_{j+1}} \left( L_C(\phi(t,j),u_{ C}(t,j))
+ \frac{d}{d t}V(\phi(t,j)) \right )dt
\end{eqnarray*}
from where we have
\begin{multline*}
0\geq  \int_{t_{j}}^{t_{j+1}} L_C(\phi(t,j),u_{ C}(t,j)) dt
\\+  V(\phi^*(t_{j+1},j))- V(\phi(t_j,j))
\end{multline*}
Summing from $j=0$ to $j=J_{\phi}$, 
we obtain 
\begin{multline*}
0
\geq\sum_{j=0}^{J_{\phi}}  \int_{t_{j}}^{t_{j+1}} L_C(\phi(t,j),u_{ C}(t,j)) dt
\\+\sum_{j=0}^{J_{\phi}} \left( V(\phi(t_{j+1},j))-V(\phi(t_j,j)) \right)
\end{multline*}
Then, solving for $V$ at the initial condition $\phi(0,0)$, we obtain
\begin{equation} \label{ContinCost2goihzol}
\begin{medsize}
\begin{split} 
V(\phi(0,0))
\geq\sum_{j=0}^{J_{\phi}}  \int_{t_{j}}^{t_{j+1}} L_C(\phi(t,j),u_{ C}(t,j)) dt
\\  + V(\phi(t_1,0))
+\sum_{j=1}^{J_{\phi}} \left( V(\phi(t_{j+1},j))-V(\phi(t_j,j)) \right) 
\end{split}
\end{medsize}
\end{equation}}
In addition, from \eqref{ubdolind}, for every {$(t_{j+1},j) \in \dom \phi$ such that $(t_{j+1},j+1) \in \dom \phi$}, we have 
\begin{equation}
\begin{medsize}
\begin{split}
  L_D(\phi(t_{j+1},&j),u_D(t_{j+1},j))
        + V(\phi(t_{j+1},j+1)) 
        - V(\phi(t_{j+1},j))
        \\
        \leq& L_D(\phi(t_{j+1},j),u_D(t_{j+1},j))
        \\
        &+ {\underset{g \in G(\phi(t_{j+1},j),u_D(t_{j+1},j))}{\sup}}
        V(g) - V(\phi(t_{j+1},j))
        \leq 0
        \label{eq:PrHamiltonD1}
\end{split}
\end{medsize}
\end{equation}
\NotAutom{Summing from $j=0$ to $j=J_{\phi}-1$,
we obtain
\begin{equation}
\begin{medsize}
\begin{split} 
\sum_{j=0}^{J_{\phi}-1}V(\phi(t_{j+1},j)) 
\geq 
\sum_{j=0}^{J_{\phi}-1}L_D(\phi(t_{j+1},j),u_{ D}(t_{j+1},j))\nonumber
\\
+\sum_{j=0}^{J_{\phi}-1} V(\phi(t_{j+1},j+1))
\nonumber
\end{split}
\end{medsize}
\end{equation}
Then, solving for $V$ at the first jump time, we obtain
\begin{eqnarray} \nonumber
&&V(\phi(t_{1},0))
\geq V( \phi(t_{1},1))
\\&&
+ \sum_{j=0}^{J_{\phi}-1}L_D(\phi(t_{j+1},j),u_{ D}(t_{j+1},j)) \label{DiscreteCost2goNOihzwincol}
\\&&
+\sum_{j=1}^{J_{\phi}-1}\left( V(\phi(t_{j+1},j+1))-V(\phi(t_{j+1},j)) \right) \nonumber
\end{eqnarray}
In addition, given that $\phi(0,0)=\xi$, lower bounding $V(\phi(t_{1},0))$ in (\ref{ContinCost2goihzol}) by the right-hand side of (\ref{DiscreteCost2goNOihzwincol}), we obtain 
\begin{equation*}
\begin{medsize}
\begin{split} 
V(\xi)
\geq& \sum_{j=0}^{J_{\phi}} \int_{t_{j}}^{t_{j+1}} L_C(\phi(t,j),u_{C}(t,j))dt+ V(\phi(t_1,0)) 
\\&
+\sum_{j=1}^{J_{\phi}} \left( V(\phi(t_{j+1},j))-V(\phi(t_j,j)) \right)  \\
\geq&
\sum_{j=0}^{J_{\phi}} \int_{t_{j}}^{t_{j+1}} L_C(\phi(t,j),u_{C}(t,j))dt
\\&
+ \sum_{j=0}^{J_{\phi}-1}L_D(\phi(t_{j+1},j),u_{D}(t_{j+1},j))
 \\
&+\sum_{j=1}^{J_{\phi}-1}\left( V(\phi(t_{j+1},j+1))-V(\phi(t_{j+1},j)) \right) 
\\&
+ V(\phi(t_{1},1))+\sum_{j=1}^{J_{\phi}} \left( V(\phi(t_{j+1},j))-V(\phi(t_j,j)) \right) 
\end{split}
\end{medsize}
\end{equation*}
Since 
\begin{eqnarray*} \nonumber
&& V(\phi(t_{1},1))
\\&&
+\sum_{j=1}^{J_{\phi}-1}\left( V(\phi(t_{j+1},j+1))-V(\phi(t_{j+1},j)) \right) \nonumber 
\\&&
+\sum_{j=1}^{J_{\phi}} \left( V(\phi(t_{j+1},j))-V(\phi(t_j,j)) \right) \nonumber
\\ &=& V(\phi(t_{{J_{\phi}}+1},J_{\phi})) +V(\phi(t_{1},1))\nonumber
\\&&
+\sum_{j=1}^{J_{\phi}-1}\left( V(\phi(t_{j+1},j+1)) \right) \nonumber 
-\sum_{j=1}^{J_{\phi}} \left(V(\phi(t_j,j)) \right) \nonumber
\\ &=&V(\phi(t_{{J_{\phi}}+1},J_{\phi})) 
\end{eqnarray*}
then we have
\begin{eqnarray*}\nonumber
V(\xi)
&\geq&\sum_{j=0}^{J_{\phi}} \int_{t_{j}}^{t_{j+1}} L_C(\phi(t,j),u_{C}(t,j))dt
\\&&
+ \sum_{j=0}^{J_{\phi}-1}L_D(\phi(t_{j+1},j),u_{D}(t_{j+1},j))
\\&&+V(\phi(t_{J_{\phi}+1},J_{\phi}))
\nonumber
\end{eqnarray*}
By taking the limit when $(t_{J_{\phi}+1},J_{\phi}) \rightarrow \sup \dom \phi$, given the definition of the cost (\ref{defJTNCinc}), we have
\begin{eqnarray} \nonumber
V(\xi)&\geq&\sum_{j=0}^{\sup_j \dom \phi} \int_{t_{j}}^{t_{j+1}} L_C(\phi(t,j),u_{C}(t,j))dt
\\&&
+ \sum_{j=0}^{{\sup_j \dom \phi}-1}L_D(\phi(t_{j+1},j),u_{D}(t_{j+1},j))
\nonumber \\
&& 
+\underset{(t,j) \in \textup{dom}\phi}{\limsup_{t+j\rightarrow \infty} } V(\phi(t,j))
\nonumber \\ 
&=&J(\phi,u)
\label{BoundCostNOzwinc}
\end{eqnarray}}
\IfAutom{Then, thanks to \eqref{eq:PrHamiltonC1} and \eqref{eq:PrHamiltonD1}, from Proposition \ref{Pp:ProofDerivation}, we have that 
$V(\xi) \geq {\widetilde{\mathcal{J}}}(\phi, u)$, with ${\widetilde{\mathcal{J}}}$ defined as in \eqref{defJTNCinc}.
}{}
\end{proof}
Based on Proposition \ref{UpperBoundOpenLoop}, we introduce the main result of the section with sufficient conditions to 
{characterize an optimal control feedback strategy to}
\pno{attain and} evaluate the upper bound cost without computing solutions. 
\begin{ttheorem}{Cost Upper Bound and Evaluation under Optimal Control}
  Given a 
  system with dynamics $\HS\IfIncd{_s}{} $ as in (\IfIncd{\ref{Hinc}}{\ref{Heq}}) \IfTp{}{with $N=2$, }with data $(C,F,D,G)$,
   stage costs $L_C:\mathbb{R}^n  \times \mathbb{R}^{m_C} \rightarrow \mathbb{R}_{\geq 0}$ {and} $L_D:\mathbb{R}^n  \times \mathbb{R}^{m_D} \rightarrow \mathbb{R}_{\geq 0}$, and terminal cost $q: \mathbb{R}^n \rightarrow \mathbb{R}$, {suppose the following hold:} 
   \begin{enumerate}
     \item There exists a function $V: \mathbb{R}^n \rightarrow \mathbb{R}$ that is continuously differentiable on a neighborhood of $\Pi(C)$ and a feedback law
   $\kappa:=(\kappa_C,\kappa_D){=((\kappa_{C1},\kappa_{C2}),(\kappa_{D1},\kappa_{D2}))}: \mathbb{R}^n \rightarrow \mathbb{R}^{m_C} \times \mathbb{R}^{m_D}$ {such that $F(x,\kappa_C(x))$ and 
   $G(x,\kappa_D(x))$ are compact for every $x$ such that $(x,\kappa_C(x))\in C$ and $(x,\kappa_D(x))\in D$, respectively,}  
   and 
     that satisfy 
  \begin{equation}
  \begin{medsize}
  \begin{split}
    0= 
   L_C(x,\kappa_C(x)) 
  +
  {\underset{f \in F(x,\kappa_C(x))}{\max}}
  \left\langle \nabla V(x),f \right\rangle 
  \IfConf{\\}{}
  \hspace{1cm} \forall x: (x,\kappa_C(x)) \in C, 
  \end{split}
  \end{medsize}
  \label{HJBzsihinc}
  \end{equation}
  \begin{equation}
    \begin{medsize}
      \begin{split}
    0\geq  
   L_C(x,\kappa_{C1}(x),u_{C2}) 
  +
  {\underset{f \in F(x,\kappa_{C1}(x),u_{C2})}{\sup}}
  \left\langle \nabla V(x),f \right\rangle 
  \IfConf{\\}{}
  \hspace{1cm} \forall (x,u_{C2}):(x,\kappa_{C1}(x),u_{C2}) \in C, 
\end{split}
\end{medsize}
\label{HJBzsihincb}
  \end{equation}
  \begin{equation}
  \begin{medsize}
  \begin{split}
    V(x)= 
    L_D(x,\kappa_D(x))
  +{\underset{g \in G(x,\kappa_D(x))}{\max}}
  V(g)
   \IfConf{\\}{}
   \hspace{1cm} \forall x: (x,\kappa_D(x)) \in D,
  \end{split}
  \end{medsize}
  \label{Bellmanzsihinc}
  \end{equation}
  \begin{equation}
    \begin{medsize}
      \begin{split}
        V(x)\geq 
         L_D(x,\kappa_{D1}(x),u_{D2})
        +
        {\underset{g \in G(x,\kappa_{D1}(x),u_{D2})}{\sup}}
        V(g) 
        \hspace{0.2cm}\\ 
        \hspace{0.1cm} \forall (x,u_{D2}):(x,\kappa_{D1}(x),u_{D2}) \in D ,
      \end{split}
    \end{medsize}
    \label{Bellmanzsihincb}
  \end{equation}
  \item For each $\xi \in \Pi(\overline{C}) \cup \Pi(D)$, each $(\phi, u)\in \mathcal{S}_{\HS\IfIncd{_s}{}   }^\infty(\xi)$ satisfies
  \begin{equation}
  \underset{(t,j) \in \textup{dom}\phi}{\limsup_{t+j\rightarrow \infty} } V(\phi(t,j))=\underset{(t,j) \in \textup{dom}\phi}{\limsup_{t+j\rightarrow \infty} } q(\phi(t,j))\pno{,}
  \label{TerminalCondCGinc}
  \end{equation}
  \end{enumerate}
  Then  \pn{ 
    \begin{equation}
    \begin{medsize}
      \begin{split}
    \mathcal{J}(\xi,u^*) 
    =\underset{(u_1^*,u_2) \in \mathcal{U}_{\HS\IfIncd{_s}{}   }^\infty(\xi)}
    { 
    \max_{u_{2}}}
    \mathcal{J}(\xi,(u_1^*,u_2))
    =V(\xi) 
      \end{split}
    \end{medsize}
  \label{ResultValueinc}
  \end{equation}
  for all $\xi \in \Pi(\overline{C}) \cup \Pi(D)$, where $u^*=(u_1^*,u_2^*)$ is defined by $\kappa$, namely, 
  there exists $(\phi^*,u^*)\in \mathcal{S}_{\HS\IfIncd{_s}{}   }^\infty(\xi)$ such that $\dom \phi^* \ni (t,j) \mapsto u^*(t,j)= \IfIh{\kappa}{\gamma}(\IfIh{}{t,j,}\phi^*(t,j))$.}
  \label{thHJBszihinc} 
  \end{ttheorem}
  \begin{proof}
  To show the claim we proceed as follows:
  \begin{enumerate}[label=\alph*)]
  \item Pick an initial condition $\xi$ and evaluate the cost associated to the 
  solutions yield by $\IfIh{\kappa=(\kappa_C,\kappa_D)}{\gamma}$ 
  from $\xi$. Find an upper bound for this cost. 
  \item Upper bound the cost associated to the solutions from $\xi$ when $P_{1}$ plays $\IfIh{\kappa}{\gamma}_1:=(\IfIh{\kappa}{\gamma}_{C1},\IfIh{\kappa}{\gamma}_{D1})$ by the value of the function $V$ evaluated at $\xi$. 
  \item By showing that the cost of the solutions from $\xi$ when $P_{1}$ plays $\IfIh{\kappa}{\gamma}_1$ are not greater than the cost of any solution yield by $\IfIh{\kappa=(\kappa_1,\kappa_2)}{\gamma}$ from $\xi$,
  we {establish \eqref{ResultValueinc}}.
  \end{enumerate}
  Proceeding as in item a above, {following \cite{ferrante2019certifying,ferrante2018optimality},} pick any $\xi \in \Pi(\overline{C}) \cup \Pi(D)$ and any solution $(\phi^*,u^*)$ to 
  \begin{equation}
  \begin{medsize}
{\HS_{\textup{max}}}\left\{ 
\begin{matrix}
\dot{x} &\in& 
\pnn{\underset{f \in F(x,u_{C})}{\textup{argmax}}}\left\langle \nabla V(x),f \right\rangle \quad  (x,u_C) \in C \\ 
x^+ &\in& 
\pnn{\underset{g \in G(x,u_D)}{\textup{argmax}}}V(g)\hspace{1.2cm}   (x,u_D) \in D
\end{matrix} \right.
\end{medsize}
\label{OptimalSystProof}
  \end{equation}
  {that thanks to Lemma \ref{Pp:sol2Hmaxsol2Hk} is also a solution to $\HS_{s}$,} 
  with $\dom \phi^* \ni (t,j) \mapsto u^*(t,j)= \IfIh{\kappa}{\gamma}(\IfIh{}{t,j,}\phi^*(t,j))$. 
  %
  Given that $V$ and $\kappa$ satisfy 
  (\ref{HJBzsihinc}), 
  for each $j \in \nats$ such that $I_{\phi^*}^j=[t_j,t_{j+1}]$ has a nonempty interior int$I_{\phi^*}^j$, 
  we have:
  \begin{enumerate}[i)]
    \item 
   for all $t \in \textup{int}I_{\phi^*}^j$, 
  \begin{equation}
  \begin{medsize}
  \begin{split}
   0
  =&  L_C(\phi^*(t,j),\IfIh{\kappa}{\gamma}_{ C}(\phi^*(t,j)))
  \\&\hspace{1cm}+{\underset{f \in F(\phi^*(t,j),\IfIh{\kappa}{\gamma}_C(\phi^*(t,j)))}{\textup{max}}}\left\langle \nabla V(\phi^*(t,j)), f \right\rangle 
  \end{split}
  \end{medsize}
  \label{eq:PrHamiltonC2}
  \end{equation}
  and $\phi^*(t,j) \in C_{\IfIh{\kappa}{\gamma}}$, as in (\ref{Hkeq}). 
\end{enumerate}
  Given that {$(\phi^*,u^*)$ is a solution to $\HS_{\textup{max}}$, and $V$ is continuously differentiable on a neighborhood of $\Pi(C)$, we can express its total derivative as}
  $$\frac{d}{d t}V(\phi^*(t,j))=\pnn{\underset{f \in F(\phi^*(t,j),\IfIh{\kappa}{\gamma}_C(\phi^*(t,j)))}{\textup{argmax}}}\left\langle \nabla V(\phi^*(t,j)), f \right\rangle$$
   for every $(t,j)\in \textup{int}(I_{\phi^*}^j)\times \{j\}$ with $\textup{int}(I_{\phi^*}^j)$ nonempty. 
  \NotConf{By integrating over the interval $[t_j,t_{j+1}]$, 
  we obtain
  \begin{eqnarray*}
  0=  \int_{t_{j}}^{t_{j+1}} \left( L_C(\phi^*(t,j),\IfIh{\kappa}{\gamma}_{ C}(\phi^*(t,j)))
  + \frac{d}{d t}V(\phi^*(t,j)) \right )dt
  \end{eqnarray*}
  from where we have
  \begin{multline*}
  0=  \int_{t_{j}}^{t_{j+1}} L_C(\phi^*(t,j),\IfIh{\kappa}{\gamma}_{ C}(\phi^*(t,j))) dt
  \\+  V(\phi^*(t_{j+1},j))- V(\phi^*(t_j,j))
  \end{multline*}
  %
  %
  %
  Summing from $j=0$ to $j=J_{\phi^*}$, 
  we obtain 
  \begin{multline*}
  0
  =\sum_{j=0}^{J_{\phi^*}}  \int_{t_{j}}^{t_{j+1}} L_C(\phi^*(t,j),\IfIh{\kappa}{\gamma}_{ C}(\phi^*(t,j))) dt
  \\+\sum_{j=0}^{J_{\phi^*}} \left( V(\phi^*(t_{j+1},j))-V(\phi^*(t_j,j)) \right)
  \end{multline*}
  Then, solving for $V$ at the initial condition $\phi^*(0,0)$, we obtain
  %
  %
  \begin{equation} \label{ContinCost2goihzinc}
  \begin{medsize}
  \begin{split} 
  V(\phi^*(0,0))
  =\sum_{j=0}^{J_{\phi^*}}  \int_{t_{j}}^{t_{j+1}} L_C(\phi^*(t,j),\IfIh{\kappa}{\gamma}_{ C}(\phi^*(t,j))) dt
  \\  + V(\phi^*(t_1,0))
  +\sum_{j=1}^{J_{\phi^*}} \left( V(\phi^*(t_{j+1},j))-V(\phi^*(t_j,j)) \right) 
  \end{split}
  \end{medsize}
  \end{equation}
  }
 Given that $V$ and $\IfIh{\kappa}{\gamma}_D$ satisfy 
  (\ref{Bellmanzsihinc}), we have: 
  \begin{enumerate}[i)]\setcounter{enumi}{1}
    \item 
  for every $(t_{j+1},j) \in \dom \phi^*$ such that $(t_{j+1},j+1) \in \dom \phi^*$, 
  \begin{equation}
  \begin{medsize}
  \begin{split} 
  V(\phi^*(t_{j+1},j))
  = 
  L_D(\phi^*(t_{j+1},j),\IfIh{\kappa}{\gamma}_{ D}(\phi^*(t_{j+1},j)))
  \hspace{1.2cm} 
  \\   \hspace{0.5cm} 
  +{\underset{g \in G(\phi^*(t_{j+1},j), \IfIh{\kappa}{\gamma}_D(\phi^*(t_{j+1},j)))}{\textup{max}}} V(g) \\
  {=} L_D(\phi^*(t_{j+1},j),\IfIh{\kappa}{\gamma}_{ D}(\phi^*(t_{j+1},j)))
  + V(\phi^*(t_{j+1},j+1)) 
  \end{split}
  \end{medsize}
  \label{eq:PrHamiltonD2}
  \end{equation}
  where $\phi^*(t,j)\in D_{\IfIh{\kappa}{\gamma}}$ as in (\ref{Hkeq}).
\end{enumerate}
\NotAutom{  Summing both sides from $j=0$ to $j=J_{\phi^*}-1$, 
  we obtain
  \begin{equation*}
  \begin{medsize}
  \begin{split} 
  \sum_{j=0}^{J_{\phi^*}-1} V(\phi^*(t,j))
  =& \sum_{j=0}^{J_{\phi^*}-1} L_D(\phi^*(t,j),\IfIh{\kappa}{\gamma}_{ D}(\phi^*(t,j))) 
  \\
  &+\sum_{j=0}^{J_{\phi^*}-1} V( \phi^*(t,j+1))
  \end{split}
  \end{medsize}
  \end{equation*}
  Then, solving for $V$ at the first jump time, we obtain
  %
  \begin{eqnarray} \nonumber
  &&V(\phi^*(t_{1},0))
  = 
  \\ \nonumber
  &&V( \phi^*(t_{1},1))
  + \sum_{j=0}^{J_{\phi^*}-1} L_D(\phi^*(t,j),\IfIh{\kappa}{\gamma}_{ D}(\phi^*(t,j))) 
  \\ \label{DiscreteCost2goihzinc}
  &&+\sum_{j=1}^{J_{\phi^*}-1}\left( V( \phi^*(t,j+1))-V(\phi^*(t,j)) \right)
  \end{eqnarray}
  Given that $\phi^*(0,0)=\xi$, by substituting the right-hand side of (\ref{DiscreteCost2goihzinc}) in  (\ref{ContinCost2goihzinc}), we obtain
  \begin{equation} 
  \begin{medsize}
  \begin{split} 
  V(\xi)
  =&\sum_{j=0}^{J_{\phi^*}}  \int_{t_{j}}^{t} L_C(\phi^*(t,j),\IfIh{\kappa}{\gamma}_{ C}(\phi^*(t,j))) dt
  \nonumber \\&
  + V(\phi^*(t_1,0))\nonumber 
  +\sum_{j=1}^{J_{\phi^*}} \left( V(\phi^*(t,j))-V(\phi^*(t_j,j)) \right) \nonumber \\
  =&\sum_{j=0}^{J_{\phi^*}}  \int_{t_{j}}^{t} L_C(\phi^*(t,j),\IfIh{\kappa}{\gamma}_{ C}(\phi^*(t,j))) dt
  \nonumber \\&
  + \sum_{j=0}^{J_{\phi^*}-1}L_D(\phi^*(t,j),\IfIh{\kappa}{\gamma}_{ D}(\phi^*(t,j)))
  + V(\phi^*(t_{1},1))
  \nonumber \\&
  +\sum_{j=1}^{J_{\phi^*}-1}\left( V(\phi^*(t,j+1))-V(\phi^*(t,j)) \right) 
  \nonumber \\&
  +\sum_{j=1}^{J_{\phi^*}} \left( V(\phi^*(t,j))-V(\phi^*(t_j,j)) \right) \nonumber 
  \end{split}
  \end{medsize}
  \end{equation}
  Since 
  \begin{eqnarray}
  && V(\phi^*(t_{1},1))
  \nonumber \\&&
  +\sum_{j=1}^{J_{\phi^*}-1}\left( V(\phi^*(t,j+1))-V(\phi^*(t,j)) \right) 
  \nonumber \\&&
  +\sum_{j=1}^{J_{\phi^*}} \left( V(\phi^*(t,j))-V(\phi^*(t_j,j)) \right) \nonumber
  \\ &=& V(\phi^*(t_{J_{\phi^*}+1},J_{\phi^*}))+V(\phi^*(t_{1},1))
  \nonumber \\&&
  +\sum_{j=1}^{J_{\phi^*}-1}\left( V(\phi^*(t,j+1)) \right) \nonumber 
  -\sum_{j=1}^{J_{\phi^*}} \left(V(\phi^*(t_j,j)) \right) \nonumber
  \\ &=&V(\phi^*(t_{J_{\phi^*}+1},J_{\phi^*})) 
  \end{eqnarray}
  then we have
  \begin{eqnarray}
  V(\xi)
  &=&\sum_{j=0}^{J_{\phi^*}}  \int_{t_{j}}^{t} L_C(\phi^*(t,j),\IfIh{\kappa}{\gamma}_{ C}(\phi^*(t,j))) dt
  \nonumber \\&&
  + \sum_{j=0}^{J_{\phi^*}-1}L_D(\phi^*(t,j),\IfIh{\kappa}{\gamma}_{ D}(\phi^*(t,j))) 
  \\&&+V(\phi^*(t_{J_{\phi^*}+1},J_{\phi^*})) \nonumber
  \end{eqnarray}
  By taking the limit when $(t_{J_{\phi^*}+1},J_{\phi^*}) \rightarrow \sup \dom \phi^*$, and given that (\ref{TerminalCondCGinc}) holds, we have 
  \begin{eqnarray}
  V(\xi)&=& \sum_{j=0}^{\sup_j \dom \phi^*}   \int_{t_{j}}^{t} L_C(\phi^*(t,j),\IfIh{\kappa}{\gamma}_{ C}(\phi^*(t,j))) dt
  \nonumber \\&&
  + \sum_{j=0}^{\sup_j \dom \phi^* -1} L_D(\phi^*(t,j),\IfIh{\kappa}{\gamma}_{ D}(\phi^*(t,j))) 
  \nonumber
  \\ && 
  +\underset{(t,j) \in \textup{dom}\phi^*}{\limsup_{t+j\rightarrow \infty} } V(\phi^*(t,j))
  \nonumber \\ 
  &=& J(\phi^*,u^*)
  \label{WCaseCostihzinc}
  \end{eqnarray}}
  \IfConf{Now, when $(\phi^*,u^*)$ is complete with $\dom \phi^* \ni (t,j) \mapsto u^*(t,j)= \IfIh{\kappa}{\gamma}(\IfIh{}{t,j,}\phi^*(t,j))$, thanks to \eqref{eq:PrHamiltonC2} and \eqref{eq:PrHamiltonD2}, from Proposition \ref{Pp:ProofDerivation} and Corollary \ref{Cor:ProofDerivation}, we have that 
  \begin{equation}
  V(\xi) = \widetilde{\mathcal{J}}(\phi^*, u^*).
    \label{WCaseCostihzinc}
  \end{equation}
  }{}
  \pnn{Using \eqref{HJBzsihinc}, for all $x$ such that $(x,\kappa_C(x)) \in C$, it holds
  \begin{equation}
    0\geq
    L_C(x,\kappa_C(x)) 
   +\left\langle \nabla V(x),f \right\rangle 
   \quad \forall f \in F(x,\kappa_C(x))
  \end{equation}
  and using \eqref{Bellmanzsihinc}, for all $x$ such that $(x,\kappa_D(x)) \in D$, we have
  \begin{equation}
    V(x)\geq 
    L_D(x,\kappa_D(x))
  + V(g)
  \quad \forall g \in G(x,\kappa_D(x)).
  \end{equation}
  Thus, for any arbitrary $(\phi, u^*)\in \mathcal{S}_{\HS_{s}}^\infty(\xi)$, we have from Proposition \ref{UpperBoundOpenLoop} and \eqref{WCaseCostihzinc} that
  \begin{equation}
    \pnn{\widetilde{\mathcal{J}}}(\phi,u^*) \leq \pnn{\widetilde{\mathcal{J}}}(\phi^*,u^*)=V(\xi)  
  \end{equation}
  which also implies that when the control action $u^*$ is defined by the feedback law $\kappa$, the largest cost of solutions from $\xi$ satisfies 
  \begin{equation}
    \J(\xi,u^*) = V(\xi).  
	\label{OptimalCostihzinc}
  \end{equation}
  }
  Proceeding with item b as above, pick any $(\phi_w,u^w) \in \mathcal{S}_{\HS\IfIncd{_s}{}  }^w (\xi)$ where $\mathcal{S}_{\HS\IfIncd{_s}{}  }^w (\xi)(\subset \mathcal{S}_{\HS\IfIncd{_s}{}  } (\xi))$ is the set of solutions $(\phi,u)$ with $u=(u_{1},u_{2})$, 
  $\dom \phi \ni (t,j) \mapsto u_1(t,j)= \IfIh{\kappa}{\gamma}_1(\IfIh{}{t,j,}\phi(t,j))$ for ${\IfIh{\kappa}{\gamma}}_1:=(\IfIh{\kappa}{\gamma}_{C1}, \IfIh{\kappa}{\gamma}_{D1}) $
  , and $\dom \phi \ni (t,j) \mapsto u_2(t,j)= \bar{\IfIh{\kappa}{\gamma}}_2(\IfIh{}{t,j,}\phi(t,j))$ for some $\bar{\IfIh{\kappa}{\gamma}}_2 \in \mathcal{K}_2$. 
  Since $\bar{\IfIh{\kappa}{\gamma}}_2$ does not necessarily attain the upper bound in 
  \eqref{HJBzsihincb}, then, 
  for each $j \in \nats$ such that $I_{\phi_w}^j=[t_j, t_{j+1}]$ has a nonempty interior int$I_{\phi_w}^j$, we have that for every $t \in \textup{int}I_{\phi_w}^j$,
  \begin{equation*}
\begin{medsize}
\begin{split}
  0
  \geq L_C(\phi_w(t,j),u_{ C}^w(t,j))+\left\langle \nabla V(\phi_w(t,j)),f \right\rangle \hspace{2cm}
  \\ 
   \forall f \in F(\phi_w(t,j),u_C^w(t,j))
\end{split}
\end{medsize}
  \end{equation*}
  which implies
  \begin{equation} 
  0 \geq L_C(\phi_w(t,j),u_{ C}^w(t,j))+\frac{d}{d t}V(\phi_w(t,j))\IfAutom{.}{,}
  \label{eq:PrHamiltonC3}
  \end{equation}
  \NotConf{ and 
  by integrating over the interval $[t_j,t_{j+1}]$, 
  we obtain
  \begin{eqnarray*}
  0 \geq \int_{t_{j}}^{t_{j+1}} \left( L_C(\phi_w(t,j),u_{ C}^w(t,j))
  + \frac{d}{d t}V(\phi_w(t,j)) \right) dt
  \end{eqnarray*}
  from which we have
  \begin{eqnarray*}
  V(\phi_w(t_j,j))
  \geq \int_{t_{j}}^{t_{j+1}} L_C(\phi_w(t,j),u_{ C}^w(t,j))dt
  \\
  +V(\phi_w(t_{j+1},j))
  \end{eqnarray*}
  Summing both sides from $j=0$ to $j=J_{\phi_w}$, 
  we obtain 
  \begin{eqnarray}
  \sum_{j=0}^{J_{\phi_w}}V(\phi_w(t_j,j)) 
  \geq \sum_{j=0}^{J_{\phi_w}} \int_{t_{j}}^{t_{j+1}} L_C(\phi_w(t,j),u_{C}^w(t,j))dt
  \nonumber \\
  +\sum_{j=0}^{J_{\phi_w}}V(\phi_w(t_{j+1},j))
  \nonumber
  \end{eqnarray}
  Then, solving for $V$ at the initial condition $\phi_w(0,0)$, we obtain
  \begin{equation}\label{ContinCost2goNOihzwinc}
  \begin{medsize}
  \begin{split} 
  &V(\phi_w(0,0))
  \geq \sum_{j=0}^{J_{\phi_w}} \int_{t_{j}}^{t_{j+1}} L_C(\phi_w(t,j),u_{C}^w(t,j))dt
  \\&
  + V(\phi_w(t_1,0))+\sum_{j=1}^{J_{\phi_w}} \left( V(\phi_w(t_{j+1},j))-V(\phi_w(t_j,j)) \right) 
  \end{split}
  \end{medsize}
  \end{equation}
  }
In addition, since $\bar{\IfIh{\kappa}{\gamma}}_2$ does not necessarily attain the upper bound in \eqref{Bellmanzsihincb}, then, 
for every $(t_{j+1},j) \in \dom \phi_w $ such that $(t_{j+1},j+1) \in \dom \phi_w$, we have
   \begin{equation*}
\begin{medsize}
\begin{split}
  V(\phi_w(t_{j+1},j))
  \geq 
  L_D(\phi_w(t_{j+1},j),u_{ D}^w(t_{j+1},j))
  + V(g)
  \hspace{1cm}
  \\
   \quad \forall g\in G(\phi_w(t_{j+1},j),u_{D}^w(t_{j+1},j))
    \end{split}
  \end{medsize}
  \end{equation*}
  and 
     \begin{equation}
\begin{medsize}
\begin{split}
    V(\phi_w(t_{j+1},j))
    \geq&
     L_D(\phi_w(t_{j+1},j),u_{ D}^w(t_{j+1},j))
    \\&
    + V( \phi_w(t_{j+1},j+1))
      \end{split}
  \end{medsize}
   \label{eq:PrHamiltonD3}
  \end{equation}
  \NotConf{Summing both sides from $j=0$ to $j=J_{\phi_w}-1$, 
  we obtain
  \begin{equation}
  \begin{medsize}
  \begin{split} 
  \sum_{j=0}^{J_{\phi_w}-1}V(\phi_w(t,j)) 
  \geq 
  \sum_{j=0}^{J_{\phi_w}-1}L_D(\phi_w(t,j),u_{ D}^w(t,j))\nonumber
  \\
  +\sum_{j=0}^{J_{\phi_w}-1} V(\phi_w(t,j+1))
  \nonumber
  \end{split}
  \end{medsize}
  \end{equation}
  Then, solving for $V$ at the first jump time, we obtain
  %
  \begin{eqnarray} \nonumber
  &&V(\phi_w(t_{1},0))
  \geq V( \phi_w(t_{1},1))
  \\&&
  + \sum_{j=0}^{J_{\phi_w}-1}L_D(\phi_w(t,j),u_{ D}^w(t,j)) \label{DiscreteCost2goNOihzwinc}
  \\&&
  +\sum_{j=1}^{J_{\phi_w}-1}\left( V(\phi_w(t,j+1))-V(\phi_w(t,j)) \right) \nonumber
  \end{eqnarray}
  In addition, given that $\phi_w(0,0)=\xi$, lower bounding $V(\phi_w(t_{1},0))$ in (\ref{ContinCost2goNOihzwinc}) by the right-hand side of (\ref{DiscreteCost2goNOihzwinc}), we obtain 
  \begin{equation*}
  \begin{medsize}
  \begin{split} 
  V(\xi)
  \geq& \sum_{j=0}^{J_{\phi_w}} \int_{t_{j}}^{t} L_C(\phi_w(t,j),u_{C}^w(t,j))dt+ V(\phi_w(t_1,0)) 
  \\&
  +\sum_{j=1}^{J_{\phi_w}} \left( V(\phi_w(t,j))-V(\phi_w(t_j,j)) \right)  \\
  \geq&
  \sum_{j=0}^{J_{\phi_w}} \int_{t_{j}}^{t} L_C(\phi_w(t,j),u_{C}^w(t,j))dt
  \\&
  + \sum_{j=0}^{J_{\phi_w}-1}L_D(\phi_w(t,j),u_{D}^w(t,j))
   \\
  &+\sum_{j=1}^{J_{\phi_w}-1}\left( V(\phi_w(t,j+1))-V(\phi_w(t,j)) \right) 
  \\&
  + V(\phi_w(t_{1},1))+\sum_{j=1}^{J_{\phi_w}} \left( V(\phi_w(t,j))-V(\phi_w(t_j,j)) \right) 
  \end{split}
  \end{medsize}
  \end{equation*}
  Since 
  \begin{eqnarray*} \nonumber
  && V(\phi_w(t_{1},1))
  \\&&
  +\sum_{j=1}^{J_{\phi_w}-1}\left( V(\phi_w(t,j+1))-V(\phi_w(t,j)) \right) \nonumber 
  \\&&
  +\sum_{j=1}^{J_{\phi_w}} \left( V(\phi_w(t,j))-V(\phi_w(t_j,j)) \right) \nonumber
  \\ &=& V(\phi_w(t_{{J_{\phi_w}}+1},J_{\phi_w})) +V(\phi_w(t_{1},1))\nonumber
  \\&&
  +\sum_{j=1}^{J_{\phi_w}-1}\left( V(\phi_w(t,j+1)) \right) \nonumber 
  -\sum_{j=1}^{J_{\phi_w}} \left(V(\phi_w(t_j,j)) \right) \nonumber
  \\ &=&V(\phi_w(t_{{J_{\phi_w}}+1},J_{\phi_w})) 
  \end{eqnarray*}
  then we have
  \begin{eqnarray*}\nonumber
  V(\xi)
  &\geq&\sum_{j=0}^{J_{\phi_w}} \int_{t_{j}}^{t} L_C(\phi_w(t,j),u_{C}^w(t,j))dt
  \\&&
  + \sum_{j=0}^{J_{\phi_w}-1}L_D(\phi_w(t,j),u_{D}^w(t,j))
  \\&&+V(\phi_w(t_{J_{\phi_w}+1},J_{\phi_w}))
  \nonumber
  \end{eqnarray*}
  By taking the limit when $(t_{J_{\phi_w}+1},J_{\phi_w}) \rightarrow \sup \dom \phi_w$, and given that (\ref{TerminalCondCGinc}) holds, we have
  \begin{eqnarray} \nonumber
  V(\xi)&\geq&\sum_{j=0}^{\sup_j \dom \phi_w} \int_{t_{j}}^{t_{j+1}} L_C(\phi_w(t,j),u_{C}^w(t,j))dt
  \\&&
  + \sum_{j=0}^{{\sup_j \dom \phi_w}-1}L_D(\phi_w(t_{j+1},j),u_{D}^w(t_{j+1},j))
  \nonumber \\
  && 
  +\underset{(t,j) \in \textup{dom}\phi_w}{\limsup_{t+j\rightarrow \infty} } V(\phi_w(t,j))
  \nonumber \\ 
  &=&J(\phi_w,u^w)
  \label{BoundCostNOzwinc}
  \end{eqnarray}
  }
  \IfConf{Now, when $(\phi_w,u^w)$ is complete, with $u^w=(u_1^w,u_2^w)$, $u_1^w$ defined by $\IfIh{\kappa}{\gamma}_1$, 
  and $u_2^w$ defined by any $\bar{\IfIh{\kappa}{\gamma}}_2 \in \mathcal{K}_2$, thanks to \eqref{eq:PrHamiltonC3} and \eqref{eq:PrHamiltonD3}, from Proposition \ref{Pp:ProofDerivation} and Corollary \ref{Cor:ProofDerivation}, we have that 
  \begin{equation}
  V(\xi) \geq  {\widetilde{\mathcal{J}}}(\phi_w,u^w).
    \label{BoundCostNOzwinc}
  \end{equation}
  }{}
  \\\\
  Finally, 
  by proceeding as in item c above,
  by applying the supremum on each side of (\ref{BoundCostNOzwinc})
  over the set $\mathcal{S}_{\HS\IfIncd{_s}{} }^\infty (\xi)$, we obtain
  \begin{equation*}
  V(\xi) \geq \underset{\pn{u_2:}(\phi_w, ({\kappa_1(\phi_w)},u_2)) \in \mathcal{S}_{\HS\IfIncd{_s}{}   }^\infty(\xi)
  }{\text{sup}}
  {\widetilde{\mathcal{J}}}(\phi_w, (\kappa_1(\phi_w),u_2))=:\underline{V}(\xi).
  \end{equation*}
  Given that $V(\xi)=\mathcal{J}(\xi, u^*)$ 
  from (\ref{OptimalCostihzinc}), we have that for any $\xi \in \Pi(\overline{C}) \cup \Pi(D)$, \pnn{each 
  $(\phi^*,u^*)\in \mathcal{S}_{\HS\IfIncd{_s}{}  }^\infty(\xi)$ with $u^*=(\kappa_1(\phi^*),\kappa_2(\phi^*))$ satisfies}
  \begin{equation}
    \underline{V}(\xi)
  \leq
  \mathcal{J}(\xi,u^*) 
   \label{Boundsszinc}
  \end{equation}
  Since  $(\phi^*,u^*) \in \mathcal{S}_{\HS_s  }^w(\xi) \cap \mathcal{S}_{\HS_s  }^\infty(\xi)$
  , we have
  \begin{equation}
    \underline{V}(\xi)
  =
  \underset{(\phi^*, (\kappa_1(\phi^*),\kappa_2(\phi^*))) \in \mathcal{S}_{\HS\IfIncd{_s}{}  }^\infty(\xi)
  }{\text{sup}} 
  \mathcal{J}(\xi, (\kappa_1(\phi^*),\kappa_2(\phi^*))).
  \label{supJk}
  \end{equation}
Given that the supremum in \eqref{supJk} is attained by $\underline{V}(\xi)$, \eqref{Boundsszinc} leads to
  \pn{
  \begin{equation}
  \mathcal{J}(\xi,u^*) 
  =\underset{(u_1^*,u_2) \in \mathcal{U}_{\HS\IfIncd{_s}{}   }^\infty(\xi)}
  { 
  \max_{u_{2}}}
  \mathcal{J}(\xi,(u_1^*,u_2))
  \label{LowerBounddzinc}
  \end{equation}}
  Thus, from (\ref{WCaseCostihzinc}) and (\ref{LowerBounddzinc}), $V(\xi)$ is the upper bound for the cost of solutions to $\HS\IfIncd{_s}{} $ 
  and the strategy $\kappa=(\kappa_1,\kappa_2)$ attains it.
\end{proof}
}{}
\NotAutom{
\pn{\section{Existence of Optimal Control Laws}}

Based on the formulation of the previous section, we establish the following results to guarantee the existence of a minmax control for the finite-time horizon and infinite-time horizon problems.
Such results
are based on assumptions that include the existence for the continuous-time case, the existence for the discrete-time case, and mild conditions on the endpoint penalties. 
In the finite-horizon case, boundedness and closedness with respect to set convergence of the set of admissible hybrid time domains is assumed.

\pn{Personal Note: The set of admissible hybrid time domains is the domain set of any admissible control action.
The proof uses a minmax sequence to construct an optimal solution pair.
If a limiting hybrid time domain exists for hybrid time domains of a minmax sequence, then standard non-hybrid results imply that an optimal solution pair
on that limiting hybrid time domain can be deduced.
}

\pn{To move here definition of a Minmax sequence. 
To define: limiting hybrid time domain.}

{
\begin{conjecture}{Existence of optimal open-loop control for hybrid finite-horizon minmax problems} 
Given Problem ($\diamond_f$), a joint input action $u=(u_C, u_D)=((u_{C1},u_{C2}),(u_{D1},u_{D2}))\in \U_\HS$,  the stage cost for flows $L_C:\mathbb{R}^n  \times \mathbb{R}^{m_C} \rightarrow \mathbb{R}_{\geq 0}$, the stage cost for jumps $L_D:\mathbb{R}^n  \times \mathbb{R}^{m_D} \rightarrow \mathbb{R}_{\geq 0}$, the terminal cost $q:\reals^n \rightarrow \reals$, if

\begin{itemize}
\pn{\item[(EC1)] $F(x,u_C)=f_0(x)+ \sum_{l=1}^{m_c} f_l(x)u_{C_l}$, where each $f_l$ is continuous and has linear growth, i.e., there exists an $M>0$ such that $|f_l(x)|\leq M(1 + |x|)$ for all $x \in \Pi(\overline{C})$, $l\in \{0,1,\dots,m_C\};$
\item[(EC2)] $C=C_x \times U_C$ for a closed set $C_x \subset \reals^n$ and a closed and convex set $U_C \subset \reals^{m_C}$;
\item[(EC3)] $L_C$ is lower semicontinuous, convex in $u_C$ for every $x \in C_x$, and either
$U_C$ is bounded or there exist $r>1,\eta>0, \delta \in \reals$ so that $L_C(x,u_C) \geq \eta |u_C|^r +\delta$ for every $(x,u_C) \in C$; 
\item[(DB1)] $G$ is locally bounded;
\item[(DB2)] $L_D$ is lower semicontinuous, and either $D=D_x \times U_D$ for a set $D_x \subset \reals^n$ and a bounded set $U_D \subset \reals^{m_D}$ or $L_D(x,u_D) \geq \beta(|u_D|)$ for every $(x,u_D) \in D$ for some $\beta:[0,\infty) \rightarrow [0,\infty)$ such that $\lim_{r \to \infty} \beta(r)=\infty$; 
\item[A3.81] The set of 
compact hybrid time domains starting at $(0, 0)$ is bounded and is closed with respect to set convergence;
\item[A3.82] The function $q$ is lower semicontinuous.}
 \end{itemize}
  then $\mathcal{J}^*$ is lower semicontinuous, and if $\mathcal{J}^*(\xi)<\infty$, then there exists a solution to Problem $(\diamond_f)$. 
\end{conjecture}
\begin{proof} To complete.

\begin{center}\pnn{
--------------------One-player finite-horizon single-valued maps-----------------}\end{center}
\IfPers{
DEFINITIONS with Goebel'19 notation 

$L_C:\reals^n\times \reals^{m_C} \rightarrow \realsgeq$ and $L_D:\reals^n\times \reals^{m_D} \rightarrow \realsgeq$.
Hybrid running cost $\Gamma$ defined for compact solution pairs.
Terminal cost $l:\reals^n\times \reals^n \rightarrow\realsgeq$.

Problem ($\star$) (Free initial and terminal conditions)
\begin{equation*}
\textup{minimize  } l(\phi(t_a,j_a),\phi(t_b,j_b))+\Gamma(\phi,u)
\label{Problem1Ih}
\end{equation*}
over all solution pairs with compact hybrid time domain.

A solution pair to $\HS$ is feasible to Problem ($\star$) if
\begin{enumerate}
	\item its domain is a compact hybrid time domain, and 
	\item its cost is finite.
\end{enumerate}


\textbf{Remark 1.} Iff $\min \left\{ l(\phi(t_a,j_a),\phi(t_b,j_b))+\Gamma(\phi,u)\right\}$ is finite, 
\begin{enumerate}
\item The set of solution pairs to $\HS$ that are feasible to Problem($\star$) is not empty, and
	\item there  exists a minimizing sequence, namely, a sequence of solution pairs $(\phi_i,u_i)$ that are feasible to Problem($\star$) such that	$l(\phi_i(t_b, j_b))+\Gamma (\phi_i,u_i)$ 	converge to  $\min \left\{ l(\phi(t_a,j_a),\phi(t_b,j_b))+\Gamma(\phi,u)\right\}$ . 
	\end{enumerate}

ASSUMPTIONS 

For the existence of solutions to Problem($\star$)

\begin{enumerate}
\item Minimizers of integral costs exist for the continuous part of (1) and the minima depend in a (lower semi)continuous way on the endpoints ($\xi, \eta$) of optimal solutions.

	Consider the function $S:[0,\infty)\times \reals^n \times \reals^n \rightarrow [0,\infty]$, where
	\begin{equation}
	S(\tau,	\xi, \eta)=\inf \left\{\int_0^\tau L_C(\phi(t),u(t)) dt : \phi(0)=\xi,\phi(\tau)=\eta  \right\}
	\label{Si}
	\end{equation}
	where the infimum is taken over all integrable $u : [0, \tau ] \rightarrow \reals^{m_C}$ and absolutely continuous $\phi : 	[0, \tau] \rightarrow \reals^n$ satisfying the continuous part of the control dynamics (1):
	$(\phi(t),u(t)) \in C$ for all $t\in (0,\tau),$     
	$\dot{\phi}(t) = F(\phi(t),u(t))$ for almost all $t\in (0,\tau)$.

	\textbf{Assumption 3.1.} The function $S$ in (\ref{Si}) is (lower semi)continuous and, if $S(\tau , \xi , \eta) < \infty$, then the infimum defining $S(\tau , \xi , \eta)$ is attained.

\item There exist minimizers for the discrete part of $\HS$ as in (1):

		\textbf{Assumption 3.2.} The set $D$ is closed. The mapping $G$ is (outer semi)continuous relative to $D$. The 	function $L_D$ is (lower semi)continuous.

\item Conditions on the terminal cost and the set of compact hybrid time domains. 
\pn{For hybrid time domains, set convergence reduces to convergence of endpoints of intervals $I_j$, for all relevant $j$.} Example 5.3 HS book.

		\textbf{Assumption 3.3} The set of compact hybrid time domains is bounded and is closed with respect to set convergence. The function $l$ is (lower semi)continuous.

\textbf{Definition. (Boundedness of the set of compact hybrid time domains)} The set of compact hybrid time domains is bounded if the union of all the compact hybrid time domains is a bounded subset of $\reals^2$.

\textbf{Definition. (Closedness of the set of compact hybrid time domains)} The set of hybrid time domains is closed with respect to set convergence if every sequence of compact hybrid time domains converges to a compact hybrid time domain.

\textbf{Example:} Given a compact $K\subset \reals^2$, the set of all compact hybrid time domains is contained in $K$.

\item 	\textbf{Property 3.4}  The sequence $(\phi_i, u_i)$ of solution pairs to $\HS$ is such that the sequences $\phi_i(t^i_j , j), \phi_i(t^i_{j+1}, j)$, and $u_i(t^i_{j+1}, j)$ are uniformly bounded.

In other words, the property requires that there exists $r > 0$ so that, for every $i$, for every $\inf_j \dom \phi_i \leq j \leq \sup_j \dom \phi_i$,

$$||\phi_i(t^i_j , j)|| < r, ||\phi_i(t^i_{j+1}, j)|| < r, ||u_i(t^i_{j+1}, j)|| < r,$$
where $\{t_j^i\}_{j=0}^{\sup_j \dom \phi}$ is a nondecreasing sequence associated to the definition of the hybrid time domain of $(\phi^i,u^i)$\IfAutomss{\>\cite[Definition 2.3]{65}}{; see Definition \ref{htd}}.

\end{enumerate}

	\textbf{Theorem 3.5} Suppose that Assumptions 3.1–3.3 hold. If there exists a minimizing sequence for Problem($\star$) that has Property 3.4, then
there exists a solution pair to $\HS$ that solves Problem($\star$), namely, an optimal solution pair.

\textbf{Definition. \pn{[cite Clarke]}} If $C=C_x \times U_C$, the pair $(F, U_C)$ is said to be finitely generated if 
$F$ has the form  
$$F(x,u_C)=f_0(x)+f(x)u_C=f_0(x)+ \sum_{l=1}^{m_c} f_l(x)u_{C_l},$$
where $f$ is a function whose values are $n\times m_C$ matrices (whose columns are the vectors $f_1,f_2,\dots,f_{m_C}$), and $u_C=(u_{C_1},u_{C_2},\dots,u_{C_{m_C}})\in \reals^{m_C}$.

Thus, a finitely generated system corresponds to (certain) linear combinations of a finite family of vector fields $\{f_l: 1 \leq i \leq m_C \}$, added to the \textit{drift} term $f_0$. When $F$
has this form, we also say that $F$ is  \textit{affine in the control variable}.

	\textbf{Proposition 5.2.} If $C=C_x \times U_C$ and
	\begin{itemize}
	\item[(EC1)] $(F,U_C)$ is \textit{finitely generated}, where each $f_l$ is continuous, and has linear growth: there exists $M>0$ such that $|f_l(x)|\leq M(1 + |x|)$ for all $x \in \Pi(\overline{C})$, $l=0,1,\dots,m_C;$
\item[(EC2)] the set $C_x \subset \reals^n$ is closed and the set $U_C \subset \reals^{m_C}$ is closed and convex;
\item[(EC3)] $L_C$ is lower semicontinuous, convex in $u_C$ for every $x \in C_x$, and either
\begin{itemize}
\item $U_C$ is bounded, or 
\item there exist $r>1,\eta>0, \delta \in \reals$ such that
\end{itemize}
 $$L_C(x,u_C) \geq \eta |u_C|^r +\delta\textup{ for every }(x,u_C) \in C;$$ 
\end{itemize}
then Assumption 3.1 holds.

\textbf{Remark.} Condition (DB) implies Assumption 3.2.

	\textbf{Proposition 5.3}  If $C=C_x \times U_C$, conditions (EC1), (EC2), (EC3) hold, 
	\begin{itemize}
	\item[(DB1)] $G$ is locally bounded;
\item[(DB2)] $L_D$ is lower semicontinuous, and either $D=D_x \times U_D$ for a set $D_x \subset \reals^n$ and a bounded set $U_D \subset \reals^{m_D}$ or $L_D(x,u_D) \geq \beta(|u_D|)$ for every $(x,u_D) \in D$ for some $\beta:[0,\infty) \rightarrow [0,\infty)$ such that $\lim_{r \to \infty} \beta(r)=\infty$; 
	\end{itemize}
	and a sequence $(\phi_i,u_i)$ of compact solution pairs is such that each of the sequences $\phi_i( \inf \dom \phi_i), \dom(\phi_i,u_i)$ and $\Gamma(\phi_i,u_i)$ is uniformly bounded, then the sequence $(\phi_i,u_i)$ has Property 3.4.

Given $(T,J) \in \reals_{\geq 0} \times \nats_{\geq 0}$, define the value function at $\xi \in \reals^n$ as
$$\mathcal{J}^* (\xi)= \inf \{\Gamma (\phi, u) + l(\phi(T,J)): \phi(0,0)=\xi \} $$
where the infimum is taken over all solution pairs to $\HS$  
with compact hybrid time domain and such that $\sup \dom \phi = (T,J)$.

\pnn{
	\textbf{Corollary 3.6 (Finite-horizon value function for fixed initial condition and terminal cost):} Suppose that  (EC1), (EC2), (EC3), (DB), and the set of compact hybrid time domains satisfies Assumption 3.3.
Then, the value function is (lower semi)continuous and if $\mathcal{J}^*(\xi)<\infty$, then a solution to Problem($\star$) exists, i.e. there exists an optimal solution pair attaining the infimum in Problem($\star$).}

This holds only for those $\xi$ such that the set of solution pairs to $\HS$ from $\xi$ that are feasible to Problem($\star$) is not empty.

\pn{Which of these conditions can be relaxed given that there is a fixed initial state and fixed terminal time?}

	}
DEFINITIONS 

$L_C:\reals^n\times \reals^{m_C} \rightarrow \realsgeq$ and $L_D:\reals^n\times \reals^{m_D} \rightarrow \realsgeq$.
Terminal cost $q: \reals^n \rightarrow\realsgeq$.
Given a compact solution pair $(\phi, u)$ to $\HS$ from $\xi$ with terminal time $(T,J$), its cost associated is defined as

\begin{equation}
\label{defJex}
\hspace{-0.25cm}
\begin{aligned}
&\mathcal{J}(\xi,u) := 
 \sum_{j=0}^{J} \int_{t_{j}}^{t_{j+1}} L_C(\phi(t,j),u_{C}(t,j))dt  \\&+ 
 \sum_{j=0}^{J -1} 
 L_D(\phi(t_{j+1},j),u_{D}(t_{j+1},j))
 + 
 q(\phi(T,J))
\end{aligned}
\end{equation}
where $\{t_j\}_{j=0}^{\sup_j \dom \phi}$ is a nondecreasing sequence associated to the definition of the hybrid time domain of $\phi$\IfAutomss{\>\cite[Definition 2.3]{65}}{; see Definition \ref{htd}}.

Problem ($\star_f$) (Fixed initial state and terminal time)
\begin{equation*}
\textup{minimize  } \mathcal{J}(\xi,u)
\label{Problem1Ih}
\end{equation*}
over all input actions with compact hybrid time domain and terminal time $(T,J)$.

An input action is feasible to Problem ($\star$) if
\begin{enumerate}
	\item its domain is a compact hybrid time domain,
	\item its terminal time is $(T,J)$, and 
	\item its cost is finite.
\end{enumerate}


Given $(T,J) \in \reals_{\geq 0} \times \nats_{\geq 0}$, define the value function at $\xi \in \reals^n$ as
$$\mathcal{J}^* (\xi)= \inf \mathcal{J}(\xi,u) $$
where the infimum is taken over all input actions 
with compact hybrid time domain and terminal time $(T,J)$.

\textbf{Remark 2.} Iff $ \mathcal{J}^*(\xi)$ is finite, 
\begin{enumerate}
\item The set of input actions that are feasible to Problem($\star_f$) is not empty, and
	\item there  exists a minimizing sequence, namely, a sequence of solutions $(\phi_i,u_i)$ to $\HS$ from $\xi$, such that $u_i$ are feasible to Problem($\star$) and 	$\mathcal{J} (\xi_i,u_i)$ 	converge to  $\mathcal{J}^*(\xi)$. 
	\end{enumerate}

ASSUMPTIONS 

For the existence of solutions to Problem($\star_f$)

\begin{enumerate}
\item Minimizers of integral costs exist for the continuous part of (1) and the minima depend in a (lower semi)continuous way on the endpoints ($\xi, \eta$) of optimal solutions.

	Consider the function $S:[0,\infty)\times \reals^n \times \reals^n \rightarrow [0,\infty]$, where
	\begin{equation}
	S(\tau,	\xi, \eta)=\inf \left\{\int_0^\tau L_C(\phi(t),u(t)) dt : \phi(0)=\xi,\phi(\tau)=\eta  \right\}
	\label{Si}
	\end{equation}
	where the infimum is taken over all integrable $u : [0, \tau ] \rightarrow \reals^{m_C}$ and absolutely continuous $\phi : 	[0, \tau] \rightarrow \reals^n$ satisfying the continuous part of the control dynamics (1):
	
	$(\phi(t),u(t)) \in C$ for all $t\in (0,\tau),$     
	
	$\dot{\phi}(t) = F(\phi(t),u(t))$ for almost all $t\in (0,\tau)$.

	\textbf{Assumption 3.1.} The function $S$ in (\ref{Si}) is (lower semi)continuous and, if $S(\tau , \xi , \eta) < \infty$, then 	the infimum defining $S(\tau , \xi , \eta)$ is attained.

\item There exist minimizers for the discrete part of $\HS$ as in (1):

		\textbf{Assumption 3.2.} The set $D$ is closed. The mapping $G$ is (outer semi)continuous relative to $D$. The 	function $L_D$ is (lower semi)continuous.

\item Conditions on the terminal cost and the set of compact hybrid time domains. 
\pn{For hybrid time domains, set convergence reduces to convergence of endpoints of intervals $I_j$, for all relevant $j$.} Example 5.3 HS book.

		\textbf{Assumption 3.3}  The function $q$ is 	(lower semi)continuous.

\textbf{Remark.} The set of compact hybrid time domains starting at $(0,0)$ and ending at $(T,J)$ is bounded and is closed with respect to set convergence. Given a compact $K\subset \reals^2$, the set of all compact hybrid time domains contained in $K$ is bounded and is closed with respect to set convergence. 

\textbf{Definition. (Boundedness of the set of compact hybrid time domains):}  The set of compact hybrid time domains starting at $(0,0)$ and ending at $\mathcal{T}$ is bounded if the union of all the compact hybrid time domains starting at $(0,0)$ and ending at $\mathcal{T}$ is a bounded subset of $\reals^2$.

\textbf{Closedness of the set of compact hybrid time domains:} The set of compact hybrid time domains starting at $(0,0)$ and ending at $\mathcal{T}$ is closed with respect to set convergence if every sequence of compact hybrid time domains starting at $(0,0)$ and ending at $\mathcal{T}$ converges to a compact hybrid time domain that starts at $(0,0)$ and ends at $\mathcal{T}$.

\item 	\textbf{Property 3.4}  The sequence $(\phi_i, u_i)$ of solutions to $\HS$ from $\xi$ is such that the sequences $\phi_i(t^i_j , j), \phi_i(t^i_{j+1}, j)$, and $u_i(t^i_{j+1}, j)$ are uniformly bounded.

In other words, the property requires that there exists $r > 0$ so that, for every $i$, for every 
$0  \leq j \leq \sup_j \dom \phi_i$,

$$||\phi_i(t^i_j , j)|| < r, ||\phi_i(t^i_{j+1}, j)|| < r, ||u_i(t^i_{j+1}, j)|| < r,$$
where $\{t_j^i\}_{j=0}^{\sup_j \dom \phi}$ is a nondecreasing sequence associated to the definition of the hybrid time domain of $(\phi^i,u^i)$\IfAutomss{\>\cite[Definition 2.3]{65}}{; see Definition \ref{htd}}.

\end{enumerate}

\pnn{\textbf{Theorem 3.5} Suppose that Assumptions 3.1–3.3 hold. Given $\xi \in \reals^n$, if there exists a minimizing sequence for Problem($\star_f$) that has Property 3.4, then
there exists a solution pair to $\HS$
from $\xi$ that solves Problem($\star_f$), namely, an optimal solution pair.}

\textbf{Definition. \pn{[cite Clarke]}} If $C=C_x \times U_C$, the pair $(F, U_C)$ is said to be finitely generated if 
$F$ has the form  
$$F(x,u_C)=f_0(x)+f(x)u_C=f_0(x)+ \sum_{l=1}^{m_c} f_l(x)u_{C_l},$$
where $f$ is a function whose values are $n\times m_C$ matrices (whose columns are the vectors $f_1,f_2,\dots,f_{m_C}$), and $u_C=(u_{C_1},u_{C_2},\dots,u_{C_{m_C}})\in \reals^{m_C}$.

Thus, a finitely generated system corresponds to (certain) linear combinations of a finite family of vector fields $\{f_l: 1 \leq i \leq m_C \}$, added to the \textit{drift} term $f_0$ . When $F$
has this form, we also say that $F$ is  \textit{affine in the control variable}.

	\textbf{Proposition 5.2.} If $C=C_x \times U_C$ and
	\begin{itemize}
	\item[(EC1)] $(F,U_C)$ is \textit{finitely generated}, where each $f_l$ is continuous,  and has linear growth: there exists $M>0$ such that $|f_l(x)|\leq M(1 + |x|)$ for all $x \in \Pi(\overline{C})$, $l=0,1,\dots,m_C;$
\item[(EC2)] the set $C_x \subset \reals^n$ is closed and the set $U_C \subset \reals^{m_C}$ is closed and convex;
\item[(EC3)] $L_C$ is lower semicontinuous, convex in $u_C$ for every $x \in C_x$, and either
\begin{itemize}
\item $U_C$ is bounded, or 
\item there exist $r>1,\eta>0, \delta \in \reals$ such that
\end{itemize}
 $$L_C(x,u_C) \geq \eta |u_C|^r +\delta\textup{ for every }(x,u_C) \in C;$$ 
\end{itemize}
then Assumption 3.1 holds.

\textbf{Remark.} Condition (DB) implies Assumption 3.2.

	\textbf{Proposition 5.3}  If $C=C_x \times U_C$, conditions (EC1), (EC2), (EC3) hold, 
	\begin{itemize}
	\item[(DB1)] $G$ is locally bounded;
\item[(DB2)] $L_D$ is lower semicontinuous, and either $D=D_x \times U_D$ for a set $D_x \subset \reals^n$ and a bounded set $U_D \subset \reals^{m_D}$ or $L_D(x,u_D) \geq \beta(|u_D|)$ for every $(x,u_D) \in D$ for some $\beta:[0,\infty) \rightarrow [0,\infty)$ such that $\lim_{r \to \infty} \beta(r)=\infty$; 
	\end{itemize}
	and a sequence $(\phi_i,u_i)$ of compact solution pairs to $\HS$ from $\xi$ with terminal time in $\mathcal{T}$ is such that \pn{the sequence $\mathcal{J}(\xi,u_i)$ is uniformly bounded,} then the sequence $(\phi_i,u_i)$ has Property 3.4.

\pn{Note: The term uniformly makes reference to every element of the sequence being bounded by the same constant. 
}

\pnn{
	\textbf{Corollary 3.6:} Suppose that  (EC1), (EC2), (EC3), (DB), and Assumption 3.3 holds.
Then, the value function $\mathcal{J}^*$ is (lower semi)continuous and given $\xi \in \reals^n$, \pn{if $\mathcal{J}^*(\xi)<\infty$}, then a solution to Problem($\star_f$) exists, i.e. there exists an optimal input action attaining the infimum in Problem($\star_f$).}

\pn{Note: Can the cost be infinite for the finite horizon case?}

This holds only for those $\xi$ such that the set of solution pairs to $\HS$ from $\xi$ that are feasible to Problem($\star$) is not empty.

{
\begin{center}
\pn{
-----------------------------Two-player Finite-horizon-----------------------------
}\end{center}

\textbf{DEFINITIONS}

$L_C$ and $L_D$ map to $\realsgeq$.
Cost $\Gamma_\infty$ defined over limit to infinity of $\Gamma$.

\textbf{Problem ($\diamond$)}
\begin{equation}
\underset{(\phi,(u_1,u_2)) \in \mathcal{S}_{\mathcal{H}  }^\infty}
{ \min_{u_{1}} \max_{u_{2}}} \Gamma_\infty(\phi,(u_1,u_2))
\label{Problem1Ih}
\end{equation}

Optimization over all solution pairs with forward complete hybrid time domain.

A solution pair to $\HS$ feasible to Problem ($\diamond$) is such that 
\begin{enumerate}
	\item its domain is a forward complete hybrid time domain, and 
	\item $\Gamma_\infty(\phi,u)< \infty$.
\end{enumerate}

If the value function is finite, 
\begin{enumerate}
	\item The set of solution pairs to $\HS$ that are feasible is not empty,
	\item an optimizer sequence is a sequence of solution pairs $(\phi_i,u_i)$ that are feasible to Problem($\diamond$) such that $\Gamma_\infty (\phi_i,u_i)$ converges to the value function. 
\end{enumerate}

\textbf{ASSUMPTIONS }

For the existence of solutions to Problem($\diamond$)
\begin{enumerate}
\item Optimizers  in the min-max sense of integral terms of the cost exist for the continuous part of (1) and the optima depend in a (lower semi)continuous way on the initial point ($\xi$) of optimal solutions.

	Consider the function $S_\infty:\reals^n \rightarrow [0,\infty]$, defined as
	\begin{equation}
	S_\infty(\xi)=\underset{
	}
{ \inf_{u_{1}} \sup_{u_{2}}} \left\{\int_0^\infty L_C(\phi(t),u(t)) dt | \phi(0)=\xi \right\}
	\label{Sinfz}
	\end{equation}
	where the infimum and supremum are taken over all locally integrable $u_1 : [0, \infty] \rightarrow \reals^{m_1}, u_2 : [0, \infty] \rightarrow \reals^{m_2}$, respectively, and locally absolutely continuous $\phi :[0, \infty] \rightarrow \reals^n$ satisfying the continuous part of the control dynamics (1),

$(\phi (t),u(t)) \in C$ for all $t \in (0, \infty), \hspace{1cm} \dot{\phi}(t) \in F(\phi(t),u(t))$ for almost all $t \in (0 , \infty)$  
	
\pn{	\textbf{Assumption 3.1z} The functions $S_i,S_s$ are (lower semi)continuous and, if $S_i(\tau , \xi , \eta) < \infty, S_s(\tau , \xi , \eta) > -\infty$, then 	the infimum defining $S_i(\tau , \xi , \eta)$ and the supremum defining $S_s(\tau , \xi , \eta)$ attained.
}

\item There exist optimizers of the discrete part of the $\HS$ as in (1):

		\textbf{Assumption 3.2} The set D is closed. The mapping $G$ is (outer semi)continuous relative to $D$. The 	function $L_D$ is (lower semi)continuous.
		
\textbf{	Assumption 3.7z} The function $S_\infty$ in (\ref{Sinf}) is (lower semi)continuous and, if $S_\infty(\xi) < \infty$, then the inf sup defining $S_\infty(\xi)$ is attained, and 
for \pn{
	\begin{equation*}
	S_i(\tau,	\xi, \eta)=\inf_{u_1} \left\{\int_0^\tau L_C(\phi(t),u(t)) dt | \phi(0)=\xi,\phi(\tau)=\eta  \right\},
	\label{S}
	\end{equation*}
		\begin{equation*}
	S_s(\tau,	\xi, \eta)=\sup_{u_2} \left\{\int_0^\tau L_C(\phi(t),u(t)) dt | \phi(0)=\xi,\phi(\tau)=\eta  \right\}
	\label{S}
	\end{equation*}
 solved over all integrable $u_1, u_2 : [0, \tau ] \rightarrow \reals^m$ and absolutely continuous $\phi : 	[0, \tau] \rightarrow \reals^n$ satisfying the continuous part of the control dynamics (1), we have
 
\begin{equation}
\lim_{\tau \nearrow  \infty} \sup_{\xi' \rightarrow \xi} \sup_\eta  S_s(\tau,	\xi', \eta) \leq S_\infty(\xi) \leq \lim_{\tau \nearrow  \infty} \inf_{\xi' \rightarrow \xi} \inf_\eta  S_i(\tau,	\xi', \eta)
\end{equation}}
\item Regarding the set of forward complete hybrid time domains. 

		\textbf{Assumption 3.8} The set of forward complete hybrid time domains starting at $(0,0)$ is closed with respect to set convergence. The function q is 	(lower semi)continuous.

\pn{I can say that $\mathcal{S}_{\mathcal{H}  }^\infty$ is such that the set of domains of solutions therein is closed with respect to set convergence. What do you think?}

\item 	\textbf{Property 3.9}  The sequence $(\phi_i, u_i)$ of solution pairs to $\HS$ is such that the sequences $\phi_i(t^i_j , j), \phi_i(t^i_{j+1}, j)$, and $u_i(t^i_{j+1}, j)$ are locally uniformly bounded.


\end{enumerate}

\pnn{
	\textbf{Theorem 3.10z} Suppose that Assumptions 3.1z, 3.2, 3.7z and 3.8 hold. If there exists an optimizing sequence for Problem($\diamond$) that has Property 3.9, then
there exists a solution pair to $\HS$ that solves for Problem($\diamond$), namely, an optimal pair.}

\pn{	\textbf{Proposition 5.5z} Sufficient conditions (EC1), (EC2) and (EC3) for Assumption 3.7z.}

Condition (DB) imply Assumption 3.2.

\pn{	\textbf{Proposition 5.6}  Sufficient conditions (EC1), (EC2), (EC3) and (DB) ensure that optimizing sequences have Property 3.9.}

\pnn{
	\textbf{Corollary 3.11 (Infinite horizon value function: fixed initial condition and no terminal cost):}  (EC1), (EC2), (EC3), (DB), and Assumption 3.8.
Then, the value function is (lower semi)continuous and if $-\infty<V(\xi)<\infty$, then a solution to Problem($\diamond$) exists, i.e. there exists an optimal solution pair to $\HS$ from $\xi$ attaining the infimum in Problem($\diamond$).}

The following results  are used to prove Theorem 3.10z.

\pn{\textbf{Lemma 5.4z} Suppose that Assumptions 3.1z, 3.2 and 3.7z hold. If $(\phi_i,u_i)$ are forward complete solutions to $\HS$ such that
\begin{itemize}
\item  $\liminf_{i\rightarrow \infty} \Gamma_\infty(\phi_i,u_i) < \infty$:
\item $\inf \dom(\phi_i,u_i) = (0, 0) \in \dom(\phi_i,u_i)$ for all $i \in \nats$ and $\lim_{i\rightarrow \infty} \phi_i(0,0)=\xi$;
\item $(\phi_i,u_i)$ satisfy Property 3.9;

then there exists a forward complete solution pair $(\phi,u)$ to $\HS$ so that $\phi(0,0)=\xi$ and
\end{itemize}
\begin{equation}
\Gamma_\infty(\phi,u) \leq \liminf_{i\rightarrow \infty} \Gamma_\infty(\phi_i,u_i)
\end{equation}}
}
\end{proof}}
\IfPers{
\begin{theorem}{Existence of optimal open-loop control for hybrid finite and infinite horizon minmax problems} 
Given Problem ($\diamond$), a joint input action $u=(u_C, u_D)=((u_{C1},u_{C2}),(u_{D1},u_{D2}))\in \U$ such that maximal solutions to $\mathcal{H}  $ from $\xi \in \Pi(C) \cup \Pi(D)$ for $u$ are complete,  the stage cost for flows $L_C:\mathbb{R}^n  \times \mathbb{R}^{m_C} \rightarrow \mathbb{R}_{\geq 0}$, the stage cost for jumps $L_D:\mathbb{R}^n  \times \mathbb{R}^{m_D} \rightarrow \mathbb{R}_{\geq 0}$, the terminal cost $q:\reals^n \rightarrow \reals$, if

\begin{itemize}
\pn{\item[(EC1)] $F(x,u_C)=f_0(x)+ \sum_{l=1}^{m_c} f_l(x)u_{C_l}$, where $f_l$ is continuous,  and there exists $M>0$ such that $|f_l(x)|\leq M(1 + |x|)$ for all $x \in \Pi(\overline{C})$, $l=0,1,\dots,m_C;$
\item[(EC2)] $C=C_x \times U_C$ for a closed set $C_x \subset \reals^n$ and a closed and convex set $U_C \subset \reals^{m_C}$;
\item[(EC3)] $L_C$ is lower semicontinuous, convex in $u_C$ for every $x \in C_x$, and either
$U_C$ is bounded or there exist $r>1,\eta>0, \delta \in \reals$ so that $L_C(x,u_C) \geq \eta |u_C|^r +\delta$ for every $(x,u_C) \in C$; 
\item[(DB1)] $G$ is locally bounded;
\item[(DB2)] $L_D$ is lower semicontinuous, and either $D=D_x \times U_D$ for a set $D_x \subset \reals^n$ and a bounded set $U_D \subset \reals^{m_D}$ or $L_D(x,u_D) \geq \beta(|u_D|)$ for every $(x,u_D) \in D$ for some $\beta:[0,\infty) \rightarrow [0,\infty)$ such that $\lim_{r \to \infty} \beta(r)=\infty$; 
\item[A3.81] The set of forward complete hybrid time domains starting at $(0, 0)$ is closed with respect to set convergence;
\item[A3.82] The function $q$ is lower semicontinuous.}
 \end{itemize}
  then $\mathcal{J}^*$ is lower semicontinuous, and if $\mathcal{J}^*(\xi)<\infty$, then there exists a solution to Problem $(\diamond)$. 
\end{theorem}
\begin{proof}

\begin{center}
\pnn{
-----------------------------One-player Infinite horizon-----------------------------
}\end{center}

\textbf{DEFINITIONS}

$L_C$ and $L_D$ map to $\realsgeq$.
Cost $\Gamma_\infty$ defined over limit to infinity of $\Gamma$.
No terminal cost.

\textbf{Problem ($\circ$):}
\begin{equation}
\textup{minimize  } l(\phi(0,0))+\Gamma_\infty(\phi,u)
\label{Problem2Ih}
\end{equation}

Minimization over all solution pairs with forward complete hybrid time domain.

A solution pair to $\HS$ feasible to Problem ($\circ$) is such that 
\begin{enumerate}
	\item its domain is a forward complete hybrid time domain, and 
	\item its cost is finite.
\end{enumerate}

If the value function is finite, 
\begin{enumerate}
	\item The set of solution pairs to $\HS$ that are feasible is not empty,
	\item a minimizing sequence is a sequence of solution pairs $(\phi_i,u_i)$ that are feasible to Problem($\circ$) such that $\Gamma_\infty (\phi_i,u_i)$ converge to the value function. 
\end{enumerate}

\textbf{ASSUMPTIONS}

For the existence of solutions to Problem($\circ$)
\begin{enumerate}
\item Minimizers of integral terms of the cost exist for the continuous part of (1) and the minima depend in a (lower semi)continuous way on the initial point ($\xi$) of optimal solutions.

	Consider the function $S_\infty:\reals^n \rightarrow [0,\infty]$, where
	\begin{equation}
	S_\infty(\xi)=\inf \left\{\int_0^\infty L_C(\phi(t),u(t)) dt | \phi(0)=\xi \right\}
	\label{Sinf}
	\end{equation}
	where the infimum is taken over all LOCALLY integrable $u : [0, \infty] \rightarrow \reals^m$ and LOCALLY absolutely continuous $\phi :[0, \infty] \rightarrow \reals^n$ satisfying the continuous part of the control dynamics (1),

$(\phi (t),u(t)) \in C$ for all $t \in (0, \infty), \hspace{1cm} \dot{\phi}(t) \in F(\phi(t),u(t))$ for almost all $t \in (0 , \infty)$  

\textbf{Assumption 3.7.} The function $S_\infty$ in (\ref{Sinfz}) is (lower semi)continuous and, if $S_\infty(\xi) < \infty$, then the infimum defining $S_\infty(\xi)$ is attained, and 

\begin{equation}
S_\infty(\xi) \leq \lim_{\tau \nearrow  \infty} \inf_{\xi' \rightarrow \xi} \inf_\eta  S(\tau,	\xi', \eta)
\end{equation}

\item Regarding the set of forward complete hybrid time domains. 

		\textbf{Assumption 3.8} The set of forward complete hybrid time domains starting at $(0,0)$ is closed with respect to set convergence. The function q is 	(lower semi)continuous.

	Every sequence of forward complete hybrid time domains converges to a forward complete hybrid time domain when the set of such is closed with respect to set convergence.
	
	\pn{WHY IS THIS AN ASSUMPTION? IT IS NOT THE CASE? HOW TO PROVE IT? I can impose it in the set over which we optimize.}

\item 	\textbf{Property 3.9}  The sequence $(\phi_i, u_i)$ of solution pairs to $\HS$ is such that the sequences $\phi_i(t^i_j , j), \phi_i(t^i_{j+1}, j)$, and $u_i(t^i_{j+1}, j)$ are locally uniformly bounded.

In other words, the property requires that for all $\tau>0$, there exists $r > 0$ so that, for every $i$, for every $\inf_j \dom \phi_i \leq j \leq \sup_j \dom \phi_i$, and times $t^i_j,t^i_{j+1}$ satisfying $t+j<\tau$,

$$||\phi_i(t^i_j , j)|| < r, ||\phi_i(t^i_{j+1}, j)|| < r, ||u_i(t^i_{j+1}, j)|| < r,$$
where $\{t_j^i\}_{j=0}^{\sup_j \dom \phi}$ is a nondecreasing sequence associated to the definition of the hybrid time domain of $(\phi^i,u^i)$;\IfAutomss{\>\cite[Definition 2.3]{65}}{; see Definition \ref{htd}}.
\end{enumerate}

\pnn{
	\textbf{Theorem 3.10.} Suppose that Assumptions 3.1, 3.2, 3.7 and 3.8 hold. If there exists a minimizing sequence for Problem($\circ$) that has Property 3.9, then
there exists a solution pair to $\HS$ that solves for Problem($\circ$), namely, an optimal pair.}

	\textbf{Proposition 5.5.} Sufficient conditions (EC1), (EC2) and (EC3) for Assumption 3.7.

Condition (DB) imply Assumption 3.2.

	\textbf{Proposition 5.6}  Sufficient conditions (EC1), (EC2), (EC3) and (DB) ensure that minimizing sequences have Property 3.9.

\pnn{
	\textbf{Corollary 3.11 (Infinite horizon value function: fixed initial condition and no terminal cost):}  (EC1), (EC2), (EC3), (DB), and Assumption 3.8.
Then, the value function is (lower semi)continuous and if $V(\xi)<\infty$, then a solution to Problem($\circ$) exists, i.e. there exists an optimal solution pair to $\HS$ from $\xi$ attaining the infimum in Problem($\circ$).}

\begin{center}
\pnn{
-----------------------------Two-player Infinite horizon-----------------------------
}\end{center}

\textbf{DEFINITIONS}

$L_C$ and $L_D$ map to $\realsgeq$.
Cost $\Gamma_\infty$ defined over limit to infinity of $\Gamma$.

\textbf{Problem ($\diamond$)}
\begin{equation}
\underset{(\phi,(u_1,u_2)) \in \mathcal{S}_{\mathcal{H}  }^\infty}
{ \min_{u_{1}} \max_{u_{2}}} \Gamma_\infty(\phi,(u_1,u_2))
\label{Problem1Ih}
\end{equation}

Optimization over all solution pairs with forward complete hybrid time domain.

A solution pair to $\HS$ feasible to Problem ($\diamond$) is such that 
\begin{enumerate}
	\item its domain is a forward complete hybrid time domain, and 
	\item $\Gamma_\infty(\phi,u)< \infty$.
\end{enumerate}

If the value function is finite, 
\begin{enumerate}
	\item The set of solution pairs to $\HS$ that are feasible is not empty,
	\item an optimizer sequence is a sequence of solution pairs $(\phi_i,u_i)$ that are feasible to Problem($\diamond$) such that $\Gamma_\infty (\phi_i,u_i)$ converges to the value function. 
\end{enumerate}

\textbf{ASSUMPTIONS }

For the existence of solutions to Problem($\diamond$)
\begin{enumerate}
\item Optimizers  in the min-max sense of integral terms of the cost exist for the continuous part of (1) and the optima depend in a (lower semi)continuous way on the initial point ($\xi$) of optimal solutions.

	Consider the function $S_\infty:\reals^n \rightarrow [0,\infty]$, defined as
	\begin{equation}
	S_\infty(\xi)=\underset{
	}
{ \inf_{u_{1}} \sup_{u_{2}}} \left\{\int_0^\infty L_C(\phi(t),u(t)) dt | \phi(0)=\xi \right\}
	\label{Sinfz}
	\end{equation}
	where the infimum and supremum are taken over all locally integrable $u_1 : [0, \infty] \rightarrow \reals^{m_1}, u_2 : [0, \infty] \rightarrow \reals^{m_2}$, respectively, and locally absolutely continuous $\phi :[0, \infty] \rightarrow \reals^n$ satisfying the continuous part of the control dynamics (1),

$(\phi (t),u(t)) \in C$ for all $t \in (0, \infty), \hspace{1cm} \dot{\phi}(t) \in F(\phi(t),u(t))$ for almost all $t \in (0 , \infty)$  
	
\pn{	\textbf{Assumption 3.1z} The functions $S_i,S_s$ are (lower semi)continuous and, if $S_i(\tau , \xi , \eta) < \infty, S_s(\tau , \xi , \eta) > -\infty$, then 	the infimum defining $S_i(\tau , \xi , \eta)$ and the supremum defining $S_s(\tau , \xi , \eta)$ attained.
}

\item There exist optimizers of the discrete part of the $\HS$ as in (1):

		\textbf{Assumption 3.2} The set D is closed. The mapping $G$ is (outer semi)continuous relative to $D$. The 	function $L_D$ is (lower semi)continuous.
		
\textbf{	Assumption 3.7z} The function $S_\infty$ in (\ref{Sinf}) is (lower semi)continuous and, if $S_\infty(\xi) < \infty$, then the inf sup defining $S_\infty(\xi)$ is attained, and 
for \pn{
	\begin{equation*}
	S_i(\tau,	\xi, \eta)=\inf_{u_1} \left\{\int_0^\tau L_C(\phi(t),u(t)) dt | \phi(0)=\xi,\phi(\tau)=\eta  \right\},
	\label{S}
	\end{equation*}
		\begin{equation*}
	S_s(\tau,	\xi, \eta)=\sup_{u_2} \left\{\int_0^\tau L_C(\phi(t),u(t)) dt | \phi(0)=\xi,\phi(\tau)=\eta  \right\}
	\label{S}
	\end{equation*}
 solved over all integrable $u_1, u_2 : [0, \tau ] \rightarrow \reals^m$ and absolutely continuous $\phi : 	[0, \tau] \rightarrow \reals^n$ satisfying the continuous part of the control dynamics (1), we have
 
\begin{equation}
\lim_{\tau \nearrow  \infty} \sup_{\xi' \rightarrow \xi} \sup_\eta  S_s(\tau,	\xi', \eta) \leq S_\infty(\xi) \leq \lim_{\tau \nearrow  \infty} \inf_{\xi' \rightarrow \xi} \inf_\eta  S_i(\tau,	\xi', \eta)
\end{equation}}
\item Regarding the set of forward complete hybrid time domains. 

		\textbf{Assumption 3.8} The set of forward complete hybrid time domains starting at $(0,0)$ is closed with respect to set convergence. The function q is 	(lower semi)continuous.

\pn{I can say that $\mathcal{S}_{\mathcal{H}  }^\infty$ is such that the set of domains of solutions therein is closed with respect to set convergence. What do you think?}

\item 	\textbf{Property 3.9}  The sequence $(\phi_i, u_i)$ of solution pairs to $\HS$ is such that the sequences $\phi_i(t^i_j , j), \phi_i(t^i_{j+1}, j)$, and $u_i(t^i_{j+1}, j)$ are locally uniformly bounded.


\end{enumerate}

\pnn{
	\textbf{Theorem 3.10z} Suppose that Assumptions 3.1z, 3.2, 3.7z and 3.8 hold. If there exists an optimizing sequence for Problem($\diamond$) that has Property 3.9, then
there exists a solution pair to $\HS$ that solves for Problem($\diamond$), namely, an optimal pair.}

\pn{	\textbf{Proposition 5.5z} Sufficient conditions (EC1), (EC2) and (EC3) for Assumption 3.7z.}

Condition (DB) imply Assumption 3.2.

\pn{	\textbf{Proposition 5.6}  Sufficient conditions (EC1), (EC2), (EC3) and (DB) ensure that optimizing sequences have Property 3.9.}

\pnn{
	\textbf{Corollary 3.11 (Infinite horizon value function: fixed initial condition and no terminal cost):}  (EC1), (EC2), (EC3), (DB), and Assumption 3.8.
Then, the value function is (lower semi)continuous and if $-\infty<V(\xi)<\infty$, then a solution to Problem($\diamond$) exists, i.e. there exists an optimal solution pair to $\HS$ from $\xi$ attaining the infimum in Problem($\diamond$).}

The following results  are used to prove Theorem 3.10z.

\pn{\textbf{Lemma 5.4z} Suppose that Assumptions 3.1z, 3.2 and 3.7z hold. If $(\phi_i,u_i)$ are forward complete solutions to $\HS$ such that
\begin{itemize}
\item  $\liminf_{i\rightarrow \infty} \Gamma_\infty(\phi_i,u_i) < \infty$:
\item $\inf \dom(\phi_i,u_i) = (0, 0) \in \dom(\phi_i,u_i)$ for all $i \in \nats$ and $\lim_{i\rightarrow \infty} \phi_i(0,0)=\xi$;
\item $(\phi_i,u_i)$ satisfy Property 3.9;

then there exists a forward complete solution pair $(\phi,u)$ to $\HS$ so that $\phi(0,0)=\xi$ and
\end{itemize}
\begin{equation}
\Gamma_\infty(\phi,u) \leq \liminf_{i\rightarrow \infty} \Gamma_\infty(\phi_i,u_i)
\end{equation}}
\end{proof}
}
}
  \section{
Hamilton–Jacobi–Bellman-Isaacs Equations for Two-player Zero-sum Hybrid Games}

The following result provides sufficient conditions to characterize the value function, and the feedback law that attains it. It addresses the solution to Problem $(\diamond)$ 
showing that the optimizer is the saddle-point equilibrium.
It involves the feasible set $\mathcal{M}$ to reduce the set over which the sufficient conditions need to be checked. When $\mathcal{M}$ is not known, it could {just} be replaced by $\reals^n$.
\begin{ttheorem}{Hamilton–Jacobi–Bellman-Isaacs (HJBI) for Problem ($\diamond$)}
\NotAutom{
  Given a two-player zero-sum hybrid game with dynamics $\HS$ as in (\ref{Heq}) \IfTp{}{with $N=2$, }with data $(C,F,D,G)$, {satisfying Assumption \ref{AssLipsZ}},
 stage costs $L_C:\mathbb{R}^n  \times \mathbb{R}^{m_C} \rightarrow \mathbb{R}_{\geq 0}$ {and} $L_D:\mathbb{R}^n  \times \mathbb{R}^{m_D} \rightarrow \mathbb{R}_{\geq 0}$, and terminal cost $q: \mathbb{R}^n \rightarrow \mathbb{R}$, {suppose the following hold:} 
 \begin{enumerate}[label=\arabic*)]
   \item There exists a function $V: \mathbb{R}^n \rightarrow \mathbb{R}$ that is continuously differentiable on a neighborhood of\> $\Pi(C)$ that satisfies the \textit{Hamilton–Jacobi–Bellman-Isaacs hybrid equations given as} 
{
\begin{equation}
\begin{split}
0&= \hspace{-0.6cm}  \underset{u_C=(u_{C1},u_{C2}) \in \Pi_u^C(x)}
{\min_{u_{C1}} \max_{u_{C2}}}  
{\mathcal{L}_C(x,u_C)}
\\&= \hspace{-0.6cm} \underset{u_C=(u_{C1},u_{C2}) \in \Pi_u^C(x)}
{ \max_{u_{C2}}\min_{u_{C1}}}
{\mathcal{L}_C(x,u_C)}
\hspace{0.4cm}\forall x \in \Pi(C),  
\end{split}
\label{HJBzsih}
\end{equation}}
where ${\mathcal{L}_C(x,u_C):=  L_C(x,u_C)+\left\langle \nabla V(x),F(x,u_C) \right\rangle }$,

\begin{equation}
\begin{split}
V(x)&=\hspace{-0.6cm}\underset{u_D=(u_{D1},u_{D2}) \in \Pi_u^D(x)}
{\min_{u_{D1}} \max_{u_{D2}}}  
{\mathcal{L}_D(x,u_D)}
\\&=\hspace{-0.6cm}\underset{u_D=(u_{D1},u_{D2}) \in \Pi_u^D(x)}
{ \max_{u_{D2}}\min_{u_{D1}}}  
{\mathcal{L}_D(x,u_D)}
\hspace{0.5cm} \forall x \in \Pi(D),
\end{split}
\label{Bellmanzsih}
\end{equation}
where ${\mathcal{L}_D(x,u_D):=  L_D(x,u_D) + V(G(x,u_D)) }$.
\item For each $\xi \in \Pi(\overline{C}) \cup \Pi(D)$, each $(\phi, u)\in \mathcal{S}_{\HS  }^\infty(\xi)$ satisfies
\begin{equation}
\underset{(t,j) \in \textup{dom}\phi}{\limsup_{t+j\rightarrow \infty} } V(\phi(t,j))=\underset{(t,j) \in \textup{dom}\phi}{\limsup_{t+j\rightarrow \infty} } q(\phi(t,j)){.}
\label{TerminalCondCG}
\end{equation}
\end{enumerate}
Then   
\begin{equation}
\J^*(\xi)= V(\xi) \hspace{2cm} \forall \xi \in \Pi(\overline{C}) \cup \Pi(D),
\label{ResultValue}
\end{equation}
and any 
feedback law
$\kappa:=(\kappa_C,\kappa_D)\NotConf{=((\kappa_{C1},\kappa_{C2}),(\kappa_{D1},\kappa_{D2}))}: \mathbb{R}^n \rightarrow \mathbb{R}^{m_C} \times \mathbb{R}^{m_D}$  with values 
{
\begin{equation}
  \hspace{-0.2cm}
\kappa_{C}(x)\in 
{\arg \min_{u_{C1}} \max_{u_{C2}}}_
{\mathclap{\substack{\\ \\ {u_C=(u_{C1},u_{C2}) \in \Pi_u^C(x)}} }}\
{\mathcal{L}_C(x,u_C)}
\hspace{0.7cm}
\forall x \in \Pi(C)
\label{kHJBeqzsihc}
\end{equation}}
and
\begin{equation}
\kappa_{D}(x)\in 
{\arg \min_{u_{D1}} \max_{u_{D2}}}_
{\mathclap{\substack{\\ \\ {u_D=(u_{D1},u_{D2}) \in \Pi_u^D(x)}} }}\
{\mathcal{L}_D(x,u_D)}
\hspace{0.7cm}
\forall x \in \Pi(D) 
\label{kBeqzsihc}
\end{equation}
is a pure strategy saddle-point equilibrium for the two-player zero-sum hybrid game with {infinite horizon and} $\mathcal{J}_1=\mathcal{J}$, $\mathcal{J}_2=-\mathcal{J}$.

}
%
%
\pno{Given a two-player zero-sum hybrid game with dynamics $\mathcal{H}$ as in (\ref{Heq}) \IfTp{}{with $N=2$, }with data $(C,F,D,G)$ {satisfying Assumption \ref{AssLipsZ}},
 stage costs $L_C:\mathbb{R}^n  \times \mathbb{R}^{m_C} \rightarrow \mathbb{R}_{\geq 0} ,L_D:\mathbb{R}^n  \times \mathbb{R}^{m_D} \rightarrow \mathbb{R}_{\geq 0}$, terminal cost $q: \mathbb{R}^n \rightarrow \mathbb{R}$,
 \pn{(potentially empty) terminal set $X$, and feasible set $\mathcal{M}$}, suppose the following hold:
  \begin{enumerate}[label=\arabic*)]
    \item There exists a function $V: \mathbb{R}^n \rightarrow \mathbb{R}$ that is continuously differentiable on a neighborhood of\> $\Pi(C)$ {and} that satisfies the \textit{Hamilton–Jacobi–Bellman-Isaacs (HJBI) hybrid equations given as} 
    \begin{equation}
      \begin{split}
      0&= 
      {\min_{u_{C1}} \max_{u_{C2}}}_
      {\mathclap{\substack{\\ \\  u_C=(u_{C1},u_{C2}) \in \Pi_u^C(x) \qquad}}}
      \quad
      {\mathcal{L}_C(x,u_C)}
      \\&= 
      { \max_{u_{C2}}\min_{u_{C1}}}_
      {\mathclap{\substack{\\ \\  u_C=(u_{C1},u_{C2}) \in \Pi_u^C(x) \qquad}}}
        \quad
      {\mathcal{L}_C(x,u_C)}
      \hspace{0.4cm}\forall x \in \Pi(C) \pnn{\cap \mathcal{M}},  
      \end{split}
      \label{HJBzsih}
      \end{equation}
      where ${\mathcal{L}_C(x,u_C):=  L_C(x,u_C)+\left\langle \nabla V(x),F(x,u_C) \right\rangle }$,
      
      \begin{equation}
      \begin{split}\hspace*{-0.6cm}
      V(x)&=
      {\min_{u_{D1}} \max_{u_{D2}}}_
      {\mathclap{\substack{\\ \\  u_D=(u_{D1},u_{D2}) \in \Pi_u^D(x) \qquad}}} \quad 
      {\mathcal{L}_D(x,u_D)}
      \\&=
      { \max_{u_{D2}}\min_{u_{D1}}}_
      {\mathclap{\substack{\\ \\  u_D=(u_{D1},u_{D2}) \in \Pi_u^D(x) \qquad}}} \quad  
      {\mathcal{L}_D(x,u_D)}
      \hspace{0.4cm} \forall x \in \Pi(D) \pnn{\cap \mathcal{M}},
      \end{split}
      \label{Bellmanzsih}
      \end{equation}
      where ${\mathcal{L}_D(x,u_D):=  L_D(x,u_D) + V(G(x,u_D)) }$.
    \item For each $\xi \in \mathcal{M} $, each $(\phi, u)\in \pn{{\mathcal{S}}_{\HS}^X(\xi)}$ satisfies\footnote{The boundary condition \eqref{TerminalCondCG} matches the value of $V$ to the terminal cost $q$ at the final value of $\phi$.}
    {
\begin{equation}
  \hspace{0.1cm}
\limsup_{\mathclap{\substack{\\{t+j\rightarrow \sup_t \dom \phi + \sup_j \dom \phi}
\\ (t,j) \in \textup{dom}\phi \hspace{0cm} }}} \quad
V(\phi(t,j))
\quad
=
\quad
\limsup_{\mathclap{\substack{\\{t+j\rightarrow \sup_t \dom \phi + \sup_j \dom \phi}
\\ (t,j) \in \textup{dom}\phi \hspace{0cm} }}} \>
q(\phi(t,j))
\label{TerminalCondCG}
\end{equation}}
\end{enumerate}
Then   
\begin{equation}
\J^*(\xi)= V(\xi) \hspace{2cm} \forall \xi \in \Pi(\overline{C}) \cup \Pi(D),
\label{ResultValue}
\end{equation}
and any 
feedback law
$\kappa:=(\kappa_C,\kappa_D)\NotConf{=((\kappa_{C1},\kappa_{C2}),(\kappa_{D1},\kappa_{D2}))}: \mathbb{R}^n \rightarrow \mathbb{R}^{m_C} \times \mathbb{R}^{m_D}$  with values 
\begin{equation}
  \hspace{-0.2cm}
\kappa_{C}(x)\in 
{\arg \min_{u_{C1}} \max_{u_{C2}}}_
{\mathclap{\substack{\\ \\ {u_C=(u_{C1},u_{C2}) \in \Pi_u^C(x)}\hspace*{1cm}} }}\
{\mathcal{L}_C(x,u_C)}
\hspace{0.7cm}
\forall x \in \Pi(C)\pnn{\cap \mathcal{M}}
\label{kHJBeqzsihc}
\end{equation}
\NotAutomss{and}
\begin{equation}
\kappa_{D}(x)\in 
{\arg \min_{u_{D1}} \max_{u_{D2}}}_
{\mathclap{\substack{\\ \\ {u_D=(u_{D1},u_{D2}) \in \Pi_u^D(x)}\hspace*{1cm}} }}\
{\mathcal{L}_D(x,u_D)}
\hspace{0.7cm}
\forall x \in \Pi(D) \pnn{\cap \mathcal{M}}
\label{kBeqzsihc}
\end{equation}
is a pure strategy saddle-point equilibrium for {Problem~$(\diamond)$} with $\mathcal{J}_1=\mathcal{J}$, $\mathcal{J}_2=-\mathcal{J}$, {where $\mathcal{J}$ is as in \eqref{defJTNCinc}}.
}
\label{thHJBszih} 
\end{ttheorem}

\NotConf{\begin{proofs}
To show the claim we \pno{apply cost evaluation tools built upon dynamic programming approaches} and proceed as follows:
\begin{enumerate}[label=\arabic*)]
\item Pick an initial condition $\xi$ and evaluate the cost associated to any 
solution yielded by $\IfIh{\kappa=(\kappa_C,\kappa_D)}{\gamma}$, with values as in (\ref{kHJBeqzsihc}) and (\ref{kBeqzsihc}), from $\xi$. Show that this cost coincides with the value of the function $V$ at $\xi$. 
\item Lower bound the cost associated to any solution from $\xi$ when $P_{2}$ plays $\IfIh{\kappa}{\gamma}_2:=(\IfIh{\kappa}{\gamma}_{C2},\IfIh{\kappa}{\gamma}_{D2})$ by the value of the function $V$ evaluated at $\xi$. 
\item Upper bound the cost associated to any solution from $\xi$ when $P_{1}$ plays $\IfIh{\kappa}{\gamma}_1:=(\IfIh{\kappa}{\gamma}_{C1},\IfIh{\kappa}{\gamma}_{D1})$ by the value of the function $V$ evaluated at $\xi$. 
\item By showing that the cost of any solution from $\xi$ when $P_{1}$ plays $\IfIh{\kappa}{\gamma}_1$ is not less than the cost of any solution yielded by $\IfIh{\kappa}{\gamma}$ from $\xi$,
and by showing that the cost of any solution from $\xi$ when $P_{2}$ plays $\IfIh{\kappa}{\gamma}_2$ is not larger than the cost of any solution yielded by $\IfIh{\kappa}{\gamma}$ from $\xi$, we show optimality of $\IfIh{\kappa}{\gamma}$ in Problem ($\diamond$) in the min-max sense. 
\end{enumerate}
\end{proofs}}
\IfIncd{The proof is presented in the Appendix.}{}
\IfAutomss{W}{Notice that w}hen the players select the optimal strategy, the value function equals the function $V$. 
\sj{The result does not require computing solutions to $\HS$}, at the price of {finding} the function $V$ satisfying the conditions therein.

\pn{The terminal set $X$ determines the size of the compact hybrid time domain of the solutions considered in Theorem \ref{thHJBszih}. Based on reachability tools, given \NotAutomss{a terminal set\>}$X$, the feasible set $\mathcal{M}$ can be computed for certain class of systems.}
\pno{  
When \NotAutomss{the feasible set\>}$\mathcal{M}$ is known a priori, the set of states for which equations (\ref{HJBzsih}) and (\ref{Bellmanzsih}) need to be enforced could be smaller than \NotAutomss{the sets of states studied }in the infinite horizon counterpart. 
}
%
\subsection{Proof of Theorem \ref{thHJBszih}}
  Before we present the proof of Theorem \ref{thHJBszih}, we present the following results providing sufficient conditions to bound and exactly evaluate the cost of the game. These results are instrumental on guaranteeing that the saddle-point equilibrium is attained and in evaluating the value function of the game.
\begin{proposition}{Time-dependent conditions for upper bound}
  Consider 
  $(\phi,u) \in \pn{\mathcal{S}_\HS^X} 
  (\xi)$
   with $u=(u_C, u_D)$, such that
%
 \begin{enumerate}[1)]
\item for each $j \in \nats$ such that $I_{\phi}^j
$  has a nonempty interior\footnote{\pno{When $j={\sup_j \dom \phi}\in \nats$ and $\sup_t \dom \phi = \infty$, we define $t_{j+1}:=\infty$.}} int$I_{\phi}^j$, 
\begin{equation}
    L_C(\phi(t,j),u_C(t,j)) 
   +\frac{d}{d t}V(\phi(t,j))  \leq
 0 \quad \>\>
  \forall t \in \textup{int}I_{\phi}^j
\label{eq:upperboundhamiltonian}
\end{equation}
 \NotAutomss{and} 
\item  for every $(t_{j+1},j) \in \dom \phi$ such that $(t_{j+1},j+1) \in \dom \phi$, 
\begin{equation}
\begin{split}
  L_D(\phi(t_{j+1},j),u_D(t_{j+1},j))
  \hspace{2cm}
  \\
        + V(\phi(t_{j+1},j+1))-V(\phi(t_{j+1},j)) 
         \leq 0.
\end{split}
         \label{upperbounddiscrete}
\end{equation}
\end{enumerate}
\vspace{-0.6cm}
Then
\begin{equation}
\begin{split}
   \sum_{j=0}^{\sup_j \dom \phi} \int_{t_{j}}^{t_{j+1}} L_C(\phi(t,j),u_{C}(t,j))dt
   \hspace{0.3cm}
    &
 \\ 
 + \sum_{j=0}^{{\sup_j \dom \phi}-1}L_D(\phi(t_{j+1},j),u_{D}(t_{j+1},j))&
  \\
 +\underset{(t,j) \in \textup{dom}\phi}{\limsup_{t+j\rightarrow \sup_t \dom \phi + \sup_j \dom \phi}} V(\phi(t,j)) 
 &\leq V(\xi).
 \end{split}
 \label{eq:costboundVxi}
 \end{equation}
 \label{Pp:ProofDerivation}
\end{proposition}
\pno{\vspace*{-0.4cm}
\begin{proof}
Given a 
$(\phi,u)\in \mathcal{S}_\HS
 (\xi)$, where $\{t_j\}_{j=0}^{\sup_j \dom \phi}$ is a nondecreasing sequence associated with the hybrid time domain of $(\phi,u)$ as in \IfAutomss{\cite[Definition 2.3]{65}}{Definition \ref{htd}},
{
for each $j \in \nats$ such that $I_{\phi}^j
$ has a nonempty interior int$I_{\phi}^j$, 
}
by integrating \eqref{eq:upperboundhamiltonian} over 
$I_{\phi}^j$, 
we obtain

\begin{equation*}
\begin{split}
0\geq  \int_{t_{j}}^{t_{j+1}} \left( L_C(\phi(t,j),u_{ C}(t,j))
+ \frac{d}{d t}V(\phi(t,j)) \right )dt
\end{split}
\end{equation*}
from where we have
\begin{equation*}
\begin{split}
0\geq  \int_{t_{j}}^{t_{j+1}} L_C(\phi(t,j),u_{ C}(t,j)) dt
\hspace{3cm}
\\
+  V(\phi(t_{j+1},j))- V(\phi(t_j,j))
\end{split}
\end{equation*}
{Pick $(t^*,j^*) \in \dom (\phi,u)$.} Summing from $j=0$ to $j={j^*}$ 
we obtain 
\begin{equation*}
\begin{split}
0
\geq\sum_{j=0}^{{j^*}}  \int_{t_{j}}^{t_{j+1}} L_C(\phi(t,j),u_{ C}(t,j)) dt
\hspace{3cm}
\\+\sum_{j=0}^{{j^*}} \left( V(\phi(t_{j+1},j))-V(\phi(t_j,j)) \right)
\end{split}
\end{equation*}
Then, solving for $V$ at the initial condition $\phi(0,0)$, we obtain

\begin{equation} \label{ContinCost2goihzol}
\begin{split} 
&  \hspace{-0.5cm}V(\phi(0,0))
\geq\sum_{j=0}^{{j^*}}  \int_{t_{j}}^{t_{j+1}} L_C(\phi(t,j),u_{ C}(t,j)) dt
\\  +& V(\phi(t_1,0))
+\sum_{j=1}^{{j^*}} \left( V(\phi(t_{j+1},j))-V(\phi(t_j,j)) \right) 
\end{split}
\end{equation}
%
In addition, {if $j^*>0$,} 
adding \eqref{upperbounddiscrete} from $j=0$ to $j={j^*}-1$, 
we obtain
\begin{equation}
\begin{split} 
\sum_{j=0}^{{j^*}-1}V(\phi(t_{j+1},j)) 
\geq 
\sum_{j=0}^{{j^*}-1}L_D(\phi(t_{j+1},j),u_{ D}(t_{j+1},j))\nonumber
\\
+\sum_{j=0}^{{j^*}-1} V(\phi(t_{j+1},j+1))
\nonumber
\end{split}
\end{equation}
Then, solving for $V$ at the first jump time, we obtain
\begin{eqnarray}
  & V(\phi(t_{1},0)) \hspace{-0.1cm}
\geq \hspace{-0.1cm} V( \phi(t_{1},1)) \nonumber
+\hspace{-0.2cm}
 \sum_{j=0}^{{j^*}-1} \hspace{-0.15cm}
 L_D(\phi(t_{j+1},j),u_{ D}(t_{j+1},j)) 
\\
&\hspace{0.5cm}+\sum_{j=1}^{{j^*}-1}\left( V(\phi(t_{j+1},j+1))-V(\phi(t_{j+1},j)) \right) 
\label{DiscreteCost2goNOihzwincol}
\end{eqnarray}
In addition, given that $\phi(0,0)=\xi$, lower bounding $V(\phi(t_{1},0))$ in (\ref{ContinCost2goihzol}) by the right-hand side of (\ref{DiscreteCost2goNOihzwincol}), we obtain 
\begin{align*} 
V(\xi)
\geq& 
  \sum_{j=0}^{{j^*}} \int_{t_{j}}^{t_{j+1}} L_C(\phi(t,j),u_{C}(t,j))dt+ V(\phi(t_1,0)) 
\\&
+\sum_{j=1}^{{j^*}} \left( V(\phi(t_{j+1},j))-V(\phi(t_j,j)) \right)  
\NotAutomss{
\\
\geq&
\sum_{j=0}^{{j^*}} \int_{t_{j}}^{t_{j+1}} L_C(\phi(t,j),u_{C}(t,j))dt
\\&
+ \sum_{j=0}^{{j^*}-1}L_D(\phi(t_{j+1},j),u_{D}(t_{j+1},j))
 \\
&+\sum_{j=1}^{{j^*}-1}\left( V(\phi(t_{j+1},j+1))-V(\phi(t_{j+1},j)) \right) 
\\&
+ V(\phi(t_{1},1))+\sum_{j=1}^{{j^*}} \left( V(\phi(t_{j+1},j))-V(\phi(t_j,j)) \right) 
}
\end{align*}
Since 
\begin{eqnarray*} \nonumber
&& \hspace{-0.2cm}V(\phi(t_{1},1))
+\sum_{j=1}^{{j^*}-1}\left( V(\phi(t_{j+1},j+1))-V(\phi(t_{j+1},j)) \right) \nonumber 
\\&&
+\sum_{j=1}^{{j^*}} \left( V(\phi(t_{j+1},j))-V(\phi(t_j,j)) \right) \nonumber
\NotAutomss{
  \\
 &=& V(\phi(t_{{{j^*}}+1},{j^*})) +V(\phi(t_{1},1))\nonumber
\\&&
+\sum_{j=1}^{{j^*}-1}\left( V(\phi(t_{j+1},j+1)) \right) \nonumber 
-\sum_{j=1}^{{j^*}} \left(V(\phi(t_j,j)) \right) \nonumber
}\IfAutomss{=}{
\\ &=&}
V(\phi(t_{{{j^*}}+1},{j^*})) 
\end{eqnarray*}
then we have
{
\begin{eqnarray*}\nonumber
V(\xi)
&\geq&V(\phi(t_{{j^*}+1},{j^*}))+\sum_{j=0}^{{j^*}} \int_{t_{j}}^{t_{j+1}} \hspace{-0.4cm}L_C(\phi(t,j),u_{C}(t,j))dt
\\&&
+ \sum_{j=0}^{{j^*}-1}L_D(\phi(t_{j+1},j),u_{D}(t_{j+1},j))
\nonumber
\end{eqnarray*}}
{By taking the limit when {$t_{{j^*}+1}+{j^*} \rightarrow \sup_t \dom \phi + \sup_j \dom \phi$, 
we establish \eqref{eq:costboundVxi}.}}
Notice that if $j^*=0$, the solution $(\phi,u)$ is continuous and \eqref{eq:costboundVxi} reduces to 
\begin{equation*}\nonumber
V(\xi)
\geq 
{\limsup}_
{\mathclap{\substack{\\ \\  
  t^* \rightarrow \sup_t \dom \phi \quad \quad \quad
} }}
\quad
\int_{t_{0}}^{t^*} 
\hspace*{-0.2cm}
L_C(\phi(t,0),u_{C}(t,0))dt
+V(\phi(t^*,0)).
\nonumber
\end{equation*}
On the other hand, if $t_{j^*+1}=0$ for all $j^*$, {the solution $(\phi,u)$ is discrete} and 
\eqref{eq:costboundVxi} reduces to 
\begin{equation*}\nonumber
V(\xi)
\geq
{\limsup}_
{\mathclap{\substack{\\ \\  
  j^* \rightarrow \sup_j \dom \phi \quad \quad \quad
} }}
\quad
 \sum_{j=0}^{{j^*}-1}
 \hspace{-0.15cm}
 L_D(\phi(0,j),u_{D}(0,j))
+V(\phi(0,{j^*})).
\nonumber
\end{equation*}
%
\vspace{-0.5cm}
\end{proof}
}{}
The following corollary is immediate from the proof of Proposition \ref{Pp:ProofDerivation}.
\begin{corollary}{Change of Signs}
  If the inequalities in the conditions in Proposition \ref{Pp:ProofDerivation} are inverted, {namely, if ``$\leq$'' in \eqref{eq:upperboundhamiltonian} and \eqref{upperbounddiscrete} is replaced with ``$\geq$'',} then \eqref{eq:costboundVxi} holds with the inequality inverted. Likewise, 
  if the conditions in Proposition \ref{Pp:ProofDerivation} hold {with} equalities, then \eqref{eq:costboundVxi} holds {with} equality.
  \label{Cor:ProofDerivation}
\end{corollary}
We are ready to prove Theorem \ref{thHJBszih}. 
To show it, we proceed as follows:
\begin{enumerate}[label=\alph*)]
\item Pick an initial condition $\xi \pn{\in \mathcal{M}}$ and evaluate the cost associated to any 
solution from $\xi$ yielded by $\IfIh{\kappa=(\kappa_C,\kappa_D)}{\gamma}$, with values as in (\ref{kHJBeqzsihc}) and (\ref{kBeqzsihc}). Show that this cost coincides with the value of the function $V$ at $\xi$. 
\item Lower bound the cost associated to any 
solution from $\xi$ when $P_{2}$ plays $\IfIh{\kappa}{\gamma}_2:=(\IfIh{\kappa}{\gamma}_{C2},\IfIh{\kappa}{\gamma}_{D2})$ by the value of the function $V$ evaluated at $\xi$. 
\item Upper bound the cost associated to any 
solution from $\xi$ when $P_{1}$ plays $\IfIh{\kappa}{\gamma}_1:=(\IfIh{\kappa}{\gamma}_{C1},\IfIh{\kappa}{\gamma}_{D1})$ by the value of the function $V$ evaluated at $\xi$. 
\item By showing that the cost of any 
solution from $\xi$ when $P_{1}$ plays $\IfIh{\kappa}{\gamma}_1$ is not {greater} than the cost of any solution yield by $\IfIh{\kappa}{\gamma}$ from $\xi$,
and by showing that the cost of any solution from $\xi$ when $P_{2}$ plays $\IfIh{\kappa}{\gamma}_2$ is not {less} than the cost of any solution yield by $\IfIh{\kappa}{\gamma}$ from $\xi$, we show optimality of $\IfIh{\kappa}{\gamma}$ in Problem ($\diamond$) in the min-max sense. 
\end{enumerate}
Proceeding as in item {a} above, pick any $\xi \in \pn{\mathcal{M}}$ 
and any {$(\phi^*,u^*) \in \pn{\mathcal{S}_{\HS  }^X}
 (\xi)$} with $\dom \phi^* \ni (t,j) \mapsto u^*(t,j)= \IfIh{\kappa}{\gamma}(\IfIh{}{t,j,}\phi^*(t,j))$. 
\IfPers{and $(t_{J_{\phi^*}+1},J_{\phi^*}) \in \dom (\phi^*,u^*) $} 
We show that $(\phi^*,u^*)$ is optimal in the min-max sense.
%
Given that $V$ satisfies 
(\ref{HJBzsih}), and $\IfIh{\kappa}{\gamma}_C$ 
is as in (\ref{kHJBeqzsihc}), for each $j \in \nats$ such that $I_{\phi^*}^j=[t_j,t_{j+1}]$ has a nonempty interior int$I_{\phi^*}^j$, we have, for all $t \in \textup{int}I_{\phi^*}^j$, 
{
\begin{equation*}
\begin{split}
0
\IfPers{=&  \underset{u_{ C}(t,j)=(u_{C1},u_{C2}) \in \Pi_u(\phi^*(t,j),C)}
{\min_{u_{C1}} \max_{u_{C2}}}  
\left\{ 
L_C(\phi^*(t,j),u_{ C}(t,j))
\right.
\\& \hspace{2.8cm}
\left.
+\left\langle \nabla V(\phi^*(t,j)),F(\phi^*(t,j),u_C(t,j)) \right\rangle \right\}\\ 
}
=&  L_C(\phi^*(t,j),\IfIh{\kappa}{\gamma}_{ C}(\phi^*(t,j)))
\\&\hspace{1.3cm}+\left\langle \nabla V(\phi^*(t,j)), F(\phi^*(t,j),\IfIh{\kappa}{\gamma}_C(\phi^*(t,j))) \right\rangle
\end{split}
\end{equation*}}
and $\phi^*(t,j) \in C_{\IfIh{\kappa}{\gamma}}$, as in (\ref{Hkeq}). Given that $V$ is continuously differentiable on a neighborhood of $\Pi(C)$, we can express its total derivative along $\phi^*$ as
\begin{equation}
\frac{d}{d t}V(\phi^*(t,j))=\left\langle \nabla V(\phi^*(t,j)), F(\phi^*(t,j),\IfIh{\kappa}{\gamma}_C(\phi^*(t,j))) \right\rangle
\label{eq:PrHamiltonC4}
\end{equation}
for every $(t,j)$ $\in \textup{int}(I^j_{\phi^*}) \times \{j\}$ with $\textup{int}(I^j_{\phi^*})$ 
nonempty\IfAutom{.}{, and 
by integrating over the interval $[t_j,t_{j+1}]$, 
we obtain
\begin{eqnarray*}
0=  \int_{t_{j}}^{t_{j+1}} \left( L_C(\phi^*(t,j),\IfIh{\kappa}{\gamma}_{ C}(\phi^*(t,j)))
+ \frac{d}{d t}V(\phi^*(t,j)) \right )dt
\end{eqnarray*}
from where we have
\begin{multline*}
0=  \int_{t_{j}}^{t_{j+1}} L_C(\phi^*(t,j),\IfIh{\kappa}{\gamma}_{ C}(\phi^*(t,j))) dt
\\+  V(\phi^*(t_{j+1},j))- V(\phi^*(t_j,j))
\end{multline*}
%
%
%
Summing from $j=0$ to $j=J_{\phi^*}$, 
we obtain 
\begin{multline*}
0
=\sum_{j=0}^{J_{\phi^*}}  \int_{t_{j}}^{t_{j+1}} L_C(\phi^*(t,j),\IfIh{\kappa}{\gamma}_{ C}(\phi^*(t,j))) dt
\\+\sum_{j=0}^{J_{\phi^*}} \left( V(\phi^*(t_{j+1},j))-V(\phi^*(t_j,j)) \right)
\end{multline*}
Then, solving for $V$ at the initial condition $\phi^*(0,0)$, we obtain
%
%
\begin{equation} \label{ContinCost2goihz}
\begin{medsize}
\begin{split} 
V(\phi^*(0,0))
=\sum_{j=0}^{J_{\phi^*}}  \int_{t_{j}}^{t_{j+1}} L_C(\phi^*(t,j),\IfIh{\kappa}{\gamma}_{ C}(\phi^*(t,j))) dt
\\  + V(\phi^*(t_1,0))
+\sum_{j=1}^{J_{\phi^*}} \left( V(\phi^*(t_{j+1},j))-V(\phi^*(t_j,j)) \right) 
\end{split}
\end{medsize}
\end{equation}
}
Given that $V$ satisfies 
(\ref{Bellmanzsih}) and $\IfIh{\kappa}{\gamma}_D$ 
is as in (\ref{kBeqzsihc}), for every $(t_{j+1},j) \in \dom \phi^*$ such that $(t_{j+1},j+1) \in \dom \phi^*$, we have that
{
\begin{equation}
\begin{split} 
V(\phi^*(t_{j+1},j))
\IfPers{
=&\underset{u_{D}(t_{j+1},j)=(u_{D1},u_{D2}) \in \Pi_u(\phi^*(t_{j+1},j),D)}
{\min_{u_{D1}} \max_{u_{D2}}}  
\\&
\left\{ 
L_D(\phi^*(t_{j+1},j),u_{ D}(t_{j+1},j))
\right. \\ & \left. \hspace{2cm} 
+V(G(\phi^*(t_{j+1},j), u_D(t_{j+1},j) \right\}\\ 
}
=& \> 
L_D(\phi^*(t_{j+1},j),\IfIh{\kappa}{\gamma}_{ D}(\phi^*(t_{j+1},j)))
\\ &  \hspace{-0.7cm} 
+ V(G(\phi^*(t_{j+1},j), \IfIh{\kappa}{\gamma}_D(\phi^*(t_{j+1},j)))) \\
=&\> L_D(\phi^*(t_{j+1},j),\IfIh{\kappa}{\gamma}_{ D}(\phi^*(t_{j+1},j)))
\\ &  \hspace{1.7cm} 
+ V(\phi^*(t_{j+1},j+1)) 
\end{split}
\label{eq:PrHamiltonD4}
\end{equation}}
where $\phi^*(t_{j+1},j)\in D_{\IfIh{\kappa}{\gamma}}$ is defined in (\ref{Hkeq}).
\NotAutom{Summing both sides from $j=0$ to $j=J_{\phi^*}-1$, 
we obtain
\begin{equation*}
\begin{medsize}
\begin{split} 
\sum_{j=0}^{J_{\phi^*}-1} V(\phi^*(t_{j+1},j))
=& \sum_{j=0}^{J_{\phi^*}-1} L_D(\phi^*(t_{j+1},j),\IfIh{\kappa}{\gamma}_{ D}(\phi^*(t_{j+1},j))) 
\\
&+\sum_{j=0}^{J_{\phi^*}-1} V( \phi^*(t_{j+1},j+1))
\end{split}
\end{medsize}
\end{equation*}
Then, solving for $V$ at the first jump time, we obtain
%
\begin{eqnarray} \nonumber
&&V(\phi^*(t_{1},0))
= 
\\ \nonumber
&&V( \phi^*(t_{1},1))
+ \sum_{j=0}^{J_{\phi^*}-1} L_D(\phi^*(t_{j+1},j),\IfIh{\kappa}{\gamma}_{ D}(\phi^*(t_{j+1},j))) 
\\ \label{DiscreteCost2goihz}
&&+\sum_{j=1}^{J_{\phi^*}-1}\left( V( \phi^*(t_{j+1},j+1))-V(\phi^*(t_{j+1},j)) \right)
\end{eqnarray}
Given that $\phi^*(0,0)=\xi$, by substituting the right-hand side of (\ref{DiscreteCost2goihz}) in  (\ref{ContinCost2goihz}), we obtain
\begin{equation} 
\begin{medsize}
\begin{split} 
V(\xi)
=&\sum_{j=0}^{J_{\phi^*}}  \int_{t_{j}}^{t_{j+1}} L_C(\phi^*(t,j),\IfIh{\kappa}{\gamma}_{ C}(\phi^*(t,j))) dt
\nonumber \\&
+ V(\phi^*(t_1,0))\nonumber 
+\sum_{j=1}^{J_{\phi^*}} \left( V(\phi^*(t_{j+1},j))-V(\phi^*(t_j,j)) \right) \nonumber \\
=&\sum_{j=0}^{J_{\phi^*}}  \int_{t_{j}}^{t_{j+1}} L_C(\phi^*(t,j),\IfIh{\kappa}{\gamma}_{ C}(\phi^*(t,j))) dt
\nonumber \\&
+ \sum_{j=0}^{J_{\phi^*}-1}L_D(\phi^*(t_{j+1},j),\IfIh{\kappa}{\gamma}_{ D}(\phi^*(t_{j+1},j)))
+ V(\phi^*(t_{1},1))
\nonumber \\&
+\sum_{j=1}^{J_{\phi^*}-1}\left( V(\phi^*(t_{j+1},j+1))-V(\phi^*(t_{j+1},j)) \right) 
\nonumber \\&
+\sum_{j=1}^{J_{\phi^*}} \left( V(\phi^*(t_{j+1},j))-V(\phi^*(t_j,j)) \right) \nonumber 
\end{split}
\end{medsize}
\end{equation}
Since 
\begin{eqnarray}
&& V(\phi^*(t_{1},1))
\nonumber \\&&
+\sum_{j=1}^{J_{\phi^*}-1}\left( V(\phi^*(t_{j+1},j+1))-V(\phi^*(t_{j+1},j)) \right) 
\nonumber \\&&
+\sum_{j=1}^{J_{\phi^*}} \left( V(\phi^*(t_{j+1},j))-V(\phi^*(t_j,j)) \right) \nonumber
\\ &=& V(\phi^*(t_{J_{\phi^*}+1},J_{\phi^*}))+V(\phi^*(t_{1},1))
\nonumber \\&&
+\sum_{j=1}^{J_{\phi^*}-1}\left( V(\phi^*(t_{j+1},j+1)) \right) \nonumber 
-\sum_{j=1}^{J_{\phi^*}} \left(V(\phi^*(t_j,j)) \right) \nonumber
\\ &=&V(\phi^*(t_{J_{\phi^*}+1},J_{\phi^*})) 
\end{eqnarray}
then we have
\begin{eqnarray}
V(\xi)
&=&\sum_{j=0}^{J_{\phi^*}}  \int_{t_{j}}^{t_{j+1}} L_C(\phi^*(t,j),\IfIh{\kappa}{\gamma}_{ C}(\phi^*(t,j))) dt
\nonumber \\&&
+ \sum_{j=0}^{J_{\phi^*}-1}L_D(\phi^*(t_{j+1},j),\IfIh{\kappa}{\gamma}_{ D}(\phi^*(t_{j+1},j))) 
\\&&+V(\phi^*(t_{J_{\phi^*}+1},J_{\phi^*})) \nonumber
\end{eqnarray}
By taking the limit when $(t_{J_{\phi^*}+1},J_{\phi^*}) \rightarrow \sup \dom \phi^*$, and given that (\ref{TerminalCondCG}) holds, we have 
\begin{eqnarray}
V(\xi)&=& \sum_{j=0}^{\sup_j \dom \phi^*}   \int_{t_{j}}^{t_{j+1}} L_C(\phi^*(t,j),\IfIh{\kappa}{\gamma}_{ C}(\phi^*(t,j))) dt
\nonumber \\&&
+ \sum_{j=0}^{\sup_j \dom \phi^* -1} L_D(\phi^*(t_{j+1},j),\IfIh{\kappa}{\gamma}_{ D}(\phi^*(t_{j+1},j))) 
\nonumber
\\ && 
+\underset{(t,j) \in \textup{dom}\phi^*}{\limsup_{t+j\rightarrow \infty} } V(\phi^*(t,j))
\nonumber \\ 
&=& \mathcal{J}(\xi,u^*)
\label{OptimalCostihz}
\end{eqnarray}
}
 \IfAutom{Now, 
 given that $(\phi^*,u^*)$ is 
  maximal
  with $\dom \phi^* \ni (t,j) \mapsto u^*(t,j)= \IfIh{\kappa}{\gamma}(\IfIh{}{t,j,}\phi^*(t,j))$,
 thanks to \eqref{eq:PrHamiltonC4} and \eqref{eq:PrHamiltonD4},
 from Corollary \ref{Cor:ProofDerivation} {and \eqref{TerminalCondCG}}, we have that 
  \begin{equation}
  V(\xi) = \mathcal{J}(\xi,u^*).
    \label{OptimalCostihz}
  \end{equation}
  }{}
%
Continuing with item {b} as above, pick any 
%
$(\phi_s,u^s) \in \mathcal{S}_{\HS  }^s(\xi)$, where
$\mathcal{S}_{\HS  }^s (\xi) (\subset \pn{\mathcal{S}_{\HS}^X(\xi)})$ is the set of solutions $(\phi,u)$ with $u=(u_1,u_2)$, $\dom \phi \ni (t,j) \mapsto u_1(t,j)= \bar{\IfIh{\kappa}{\gamma}}_1(\IfIh{}{t,j,}\phi(t,j))$ for some $\bar{\IfIh{\kappa}{\gamma}}_1 \in \mathcal{K}_1$, $\dom \phi \ni (t,j) \mapsto u_2(t,j)= \IfIh{\kappa}{\gamma}_2(t,j,\phi(t,j))$ for ${\IfIh{\kappa}{\gamma}}_2:=(\IfIh{\kappa}{\gamma}_{C2}, \IfIh{\kappa}{\gamma}_{D2}) $ as in (\ref{kHJBeqzsihc}) and (\ref{kBeqzsihc}). 
Since $\bar{\IfIh{\kappa}{\gamma}}_1$ does not necessarily attain the minimum in (\ref{HJBzsih}), then,  for each $j \in \nats$ such that $I_{\phi_s}^j=[t_j, t_{j+1}]$ has a nonempty interior int$I_{\phi_s}^j$, we have for every $t \in \textup{int}I_{\phi_s}^j$,
{
\begin{equation*}
\begin{split}
0&
\leq  L_C(\phi_s(t,j),u_{ C}^s(t,j))
\\& \hspace{2cm}
+\left\langle \nabla V(\phi_s(t,j)),F(\phi_s(t,j),u_C^s(t,j)) \right\rangle.
\end{split}
\end{equation*}}
%
Similarly to \eqref{eq:PrHamiltonC4}, we have

\begin{equation} 
\begin{split}
\frac{d}{d t}V(\phi_s(t,j)):=\left\langle \nabla V(\phi_s(t,j)),F(\phi_s(t,j),u_C^s(t,j)) \right\rangle
\end{split}
\label{eq:PrHamiltonC5}
\end{equation}
 for every $(t,j) \in \textup{int}(I^j_{\phi_s}) \times \{j\}$ with $\textup{int}(I^j_{\phi_s})$ 
 nonempty{\IfAutom{.}{, and
by integrating over the interval $[t_j,t_{j+1}]$,  
we obtain
\begin{eqnarray*}
0 \leq \int_{t_{j}}^{t_{j+1}}  \left( L_C(\phi_s(t,j),u_{ C}^s(t,j))
+  \frac{d}{d t}V(\phi_s(t,j)) \right) dt
\end{eqnarray*}
from which we have
\begin{equation*}
\begin{medsize}
V(\phi_s(t_j,j))
\leq \int_{t_{j}}^{t_{j+1}}   L_C(\phi_s(t,j),u_{ C}^s(t,j))dt
+V(\phi_s(t_{j+1},j))
\end{medsize}
\end{equation*}
Summing both sides from  $j=0$ to $j=J_{\phi_s}$, 
we obtain 
\begin{eqnarray}
\sum_{j=0}^{J_{\phi_s}}V(\phi_s(t_j,j)) 
\leq \sum_{j=0}^{J_{\phi_s}}  \int_{t_{j}}^{t_{j+1}}  L_C(\phi_s(t,j),u_{C}^s(t,j))dt
\nonumber \\
+\sum_{j=0}^{J_{\phi_s}}V(\phi_s(t_{j+1},j))
\nonumber
\end{eqnarray}
Then, solving for $V$ at the initial condition $\phi_s(0,0)$, we obtain
%
%
\begin{equation}\label{ContinCost2goNOihz}
\begin{medsize}
\begin{split} 
&V(\phi_s(0,0))
\leq \sum_{j=0}^{J_{\phi_s}}  \int_{t_{j}}^{t_{j+1}}  L_C(\phi_s(t,j),u_{C}^s(t,j))dt
\\
&+ V(\phi_s(t_1,0))+\sum_{j=1}^{J_{\phi_s}} \left( V(\phi_s(t_{j+1},j))-V(\phi_s(t_j,j)) \right) 
\end{split}
\end{medsize}
\end{equation}
}
In addition, since 
{$\bar{\IfIh{\kappa}{\gamma}}_1$} does not necessarily attain the minimum in (\ref{Bellmanzsih}), then for every $(t_{j+1},j) \in \dom \phi_s $ such that   $(t_{j+1},j+1) \in \dom \phi_s$, we have
{
\begin{equation}
\begin{split} 
V(\phi_s(t_{j+1},j))
\leq&  \NotAutom{
  \> 
L_D(\phi_s(t_{j+1},j),u_{ D}^s(t_{j+1},j))
\\&\hspace{0.7cm}
+ V(G(\phi_s(t_{j+1},j),u_{D}^s(t_{j+1},j)))\\
{\leq}&
}
\>L_D(\phi_s(t_{j+1},j),u_{ D}^s(t_{j+1},j))
\\&\hspace{1.5cm}
+ V( \phi_s(t_{j+1},j+1)).
\end{split}
\label{eq:PrHamiltonD5}
\end{equation}}
\NotAutom{
Summing both sides from  $j=0$ to $j=J_{\phi_s}-1$, 
we obtain
\begin{eqnarray}
\sum_{j=0}^{J_{\phi_s}-1}V(\phi_s(t_{j+1},j)) 
\leq 
\sum_{j=0}^{J_{\phi_s}-1}L_D(\phi_s(t_{j+1},j),u_{ D}^s(t_{j+1},j))
\nonumber\\
+\sum_{j=0}^{J_{\phi_s}-1} V(\phi_s(t_{j+1},j+1))
\nonumber
\end{eqnarray}
Then, solving for $V$ at the first jump time, we obtain
%
\begin{equation}\label{DiscreteCost2goNOihz}
\begin{medsize}
\begin{split} 
V(\phi_s(t_{1},0))
\leq V( \phi_s(t_{1},1))
+ \sum_{j=0}^{J_{\phi_s}-1}L_D(\phi_s(t_{j+1},j),u_{ D}^s(t_{j+1},j))
\\
+\sum_{j=1}^{J_{\phi_s}-1}\left( V(\phi_s(t_{j+1},j+1))-V(\phi_s(t_{j+1},j)) \right)
\end{split}
\end{medsize}
\end{equation}
In addition, given that $\phi_s(0,0)=\xi$, upper bounding $V(\phi_s(t_{1},0))$ in (\ref{ContinCost2goNOihz}) by the right-hand side of (\ref{DiscreteCost2goNOihz}), we obtain 
\begin{equation*}
\begin{medsize}
\begin{split}
V(\xi)  
\leq& \sum_{j=0}^{J_{\phi_s}}  \int_{t_{j}}^{t_{j+1}}  L_C(\phi_s(t,j),u_{C}^s(t,j))dt+ V(\phi_s(t_1,0))
\\&
+\sum_{j=1}^{J_{\phi_s}} \left( V(\phi_s(t_{j+1},j))-V(\phi_s(t_j,j)) \right)  \\
\leq&
\sum_{j=0}^{J_{\phi_s}}  \int_{t_{j}}^{t_{j+1}}  L_C(\phi_s(t,j),u_{C}^s(t,j))dt
\\&+ \sum_{j=0}^{J_{\phi_s}-1}L_D(\phi_s(t_{j+1},j),u_{D}^s(t_{j+1},j))
 \\
&+\sum_{j=1}^{J_{\phi_s}-1}\left( V(\phi_s(t_{j+1},j+1))-V(\phi_s(t_{j+1},j)) \right)  
\\&+ V(\phi_s(t_{1},1))+\sum_{j=1}^{J_{\phi_s}} \left( V(\phi_s(t_{j+1},j))-V(\phi_s(t_j,j)) \right) 
\end{split}
\end{medsize}
\end{equation*}
Since 
\begin{eqnarray}
&& V(\phi_s(t_{1},1))\nonumber
\\&&+\sum_{j=1}^{J_{\phi_s}-1}\left( V(\phi_s(t_{j+1},j+1))-V(\phi_s(t_{j+1},j)) \right) \nonumber 
\\&&
+\sum_{j=1}^{J_{\phi_s}} \left( V(\phi_s(t_{j+1},j))-V(\phi_s(t_j,j)) \right) \nonumber
\\ &=& V(\phi_s(t_{{J_{\phi_s}}+1},J_{\phi_s})) +V(\phi_s(t_{1},1))\nonumber
\\&&
+\sum_{j=1}^{J_{\phi_s}-1}\left( V(\phi_s(t_{j+1},j+1)) \right)
-\sum_{j=1}^{J_{\phi_s}} \left(V(\phi_s(t_j,j)) \right) \nonumber
\\ &=&V(\phi_s(t_{{J_{\phi_s}}+1},J_{\phi_s}))
\end{eqnarray}
then we have
\begin{eqnarray}
V(\xi)
&\leq&\sum_{j=0}^{J_{\phi_s}}  \int_{t_{j}}^{t_{j+1}}  L_C(\phi_s(t,j),u_{C}^s(t,j))dt\nonumber
\\&&
+ \sum_{j=0}^{J_{\phi_s}-1}L_D(\phi_s(t_{j+1},j),u_{D}^s(t_{j+1},j))
\\&&+V(\phi_s(t_{J_{\phi_s}+1},J_{\phi_s}))
\nonumber
\end{eqnarray}
By taking the limit when $(t_{{J_{\phi_s}}+1},J_{\phi_s}) \rightarrow \sup \dom \phi_s$, and given that (\ref{TerminalCondCG}) holds, we have
\begin{eqnarray}
V(\xi)&\leq&\sum_{j=0}^{\sup_j \dom \phi_s} \int_{t_{j}}^{t_{j+1}}  L_C(\phi_s(t,j),u_{C}^s(t,j))dt\nonumber
\\&&\nonumber
+ \sum_{j=0}^{{\sup_j \dom \phi_s}-1}L_D(\phi_s(t_{j+1},j),u_{D}^s(t_{j+1},j))
\\&&
+\underset{(t,j) \in \textup{dom}\phi_s}{\limsup_{t+j\rightarrow \infty} } V(\phi_s(t,j))
\nonumber \\ 
&\leq& \mathcal{J}(\xi,u^s)
\label{BoundCostNOz}
\end{eqnarray}
}
\IfAutom{Now, 
given that $(\phi_s, \sj{u^s})$ is 
maximal, with $u^s=(u_1^s,u_2^s)$, $u_1^s$ defined by any $\bar{\IfIh{\kappa}{\gamma}}_1 \in \mathcal{K}_1$, and $u_2^s$ defined by $\IfIh{\kappa}{\gamma}_2$ as in (\ref{kHJBeqzsihc}) and (\ref{kBeqzsihc}), thanks to \eqref{eq:PrHamiltonC5} and \eqref{eq:PrHamiltonD5}, 
 from Proposition~\ref{Pp:ProofDerivation} {and \eqref{TerminalCondCG}}, we have 
  \begin{equation}
  V(\xi) \leq \mathcal{J}(\xi,u^s).
    \label{BoundCostNOz}
  \end{equation}
  }{}
%
Proceeding with item {c} as above, pick any 
$(\phi_w,u^w) \in {\mathcal{S}_{\HS}^w(\xi)}$, where  
$\mathcal{S}_{\HS}^w(\xi) (\subset \pn{\mathcal{S}_{\HS}^X(\xi)})$ is the set of solutions $(\phi,u)$ with $u=(u_{1},u_{2})$, $\dom \phi \ni (t,j) \mapsto u_1(t,j)= \IfIh{\kappa}{\gamma}_1(\IfIh{}{t,j,}\phi(t,j))$ for ${\IfIh{\kappa}{\gamma}}_1:=(\IfIh{\kappa}{\gamma}_{C1}, \IfIh{\kappa}{\gamma}_{D1}) $ as in (\ref{kHJBeqzsihc}) and (\ref{kBeqzsihc}), $\dom \phi \ni (t,j) \mapsto u_2(t,j)= \bar{\IfIh{\kappa}{\gamma}}_2(\IfIh{}{t,j,}\phi(t,j))$ for some $\bar{\IfIh{\kappa}{\gamma}}_2 \in \mathcal{K}_2$.
Since $\bar{\IfIh{\kappa}{\gamma}}_2$ does not necessarily attain the maximum in (\ref{HJBzsih}), then, for each $j \in \nats$ such that $I_{\phi_w}^j=[t_j, t_{j+1}]$ has a nonempty interior int$I_{\phi_w}^j$, we have for every $t \in \textup{int}I_{\phi_w}^j$,
{
\begin{equation*}
\begin{split}
0&
\geq L_C(\phi_w(t,j),u_{ C}^w(t,j))
\\& \hspace{2cm}
+\left\langle \nabla V(\phi_w(t,j)),F(\phi_w(t,j),u_C^w(t,j)) \right\rangle.
\end{split}
\end{equation*}
}
%
{Similarly to \eqref{eq:PrHamiltonC4}, we have} 
\begin{equation}
 \frac{d}{d t}V(\phi_w(t,j)):=\left\langle \nabla V(\phi_w(t,j)),F(\phi_w(t,j),u_C^w(t,j)) \right\rangle
 \label{eq:PrHamiltonC6}
 \end{equation}
  for every $(t,j) \in \textup{int}(I^j_{\phi_w}) \times \{j\}$ with $ \textup{int}(I^j_{\phi_w}) $ 
  nonempty{\IfAutom{.}{, and
by integrating over the interval $[t_j,t_{j+1}]$, 
we obtain
\begin{eqnarray*}
0 \geq \int_{t_{j}}^{t_{j+1}} \left( L_C(\phi_w(t,j),u_{ C}^w(t,j))
+ \frac{d}{d t}V(\phi_w(t,j)) \right) dt
\end{eqnarray*}
from which we have
\begin{eqnarray*}
V(\phi_w(t_j,j))
\geq \int_{t_{j}}^{t_{j+1}} L_C(\phi_w(t,j),u_{ C}^w(t,j))dt
\\
+V(\phi_w(t_{j+1},j))
\end{eqnarray*}
Summing both sides from $j=0$ to $j=J_{\phi_w}$, 
we obtain 
\begin{eqnarray}
\sum_{j=0}^{J_{\phi_w}}V(\phi_w(t_j,j)) 
\geq \sum_{j=0}^{J_{\phi_w}} \int_{t_{j}}^{t_{j+1}} L_C(\phi_w(t,j),u_{C}^w(t,j))dt
\nonumber \\
+\sum_{j=0}^{J_{\phi_w}}V(\phi_w(t_{j+1},j))
\nonumber
\end{eqnarray}
Then, solving for $V$ at the initial condition $\phi_w(0,0)$, we obtain
%
%
\begin{equation}\label{ContinCost2goNOihzw}
\begin{medsize}
\begin{split} 
&V(\phi_w(0,0))
\geq \sum_{j=0}^{J_{\phi_w}} \int_{t_{j}}^{t_{j+1}} L_C(\phi_w(t,j),u_{C}^w(t,j))dt
\\&
+ V(\phi_w(t_1,0))+\sum_{j=1}^{J_{\phi_w}} \left( V(\phi_w(t_{j+1},j))-V(\phi_w(t_j,j)) \right) 
\end{split}
\end{medsize}
\end{equation}
}
In addition, since 
$\bar{\IfIh{\kappa}{\gamma}}_2$ does not necessarily attain the maximum in (\ref{Bellmanzsih}), then for every $(t_{j+1},j) \in \dom \phi_w $ such that $(t_{j+1},j+1) \in \dom \phi_w$, we have
{
\begin{equation}
\begin{split} 
V(\phi_w(t_{j+1},j))
\geq& \NotAutom{
  \>
L_D(\phi_w(t_{j+1},j),u_{ D}^w(t_{j+1},j))
\\& \hspace{0.7cm}
+ V(G(\phi_w(t_{j+1},j),u_{D}^w(t_{j+1},j)))
\\
=&
}
\>L_D(\phi_w(t_{j+1},j),u_{ D}^w(t_{j+1},j))
\\&\hspace{1.5cm}
+ V( \phi_w(t_{j+1},j+1)).
\end{split}
\label{eq:PrHamiltonD6}
\end{equation}}
\NotAutom{
Summing both sides from $j=0$ to $j=J_{\phi_w}-1$, 
we obtain
\begin{equation}
\begin{medsize}
\begin{split} 
\sum_{j=0}^{J_{\phi_w}-1}V(\phi_w(t_{j+1},j)) 
\geq 
\sum_{j=0}^{J_{\phi_w}-1}L_D(\phi_w(t_{j+1},j),u_{ D}^w(t_{j+1},j))\nonumber
\\
+\sum_{j=0}^{J_{\phi_w}-1} V(\phi_w(t_{j+1},j+1))
\nonumber
\end{split}
\end{medsize}
\end{equation}
Then, solving for $V$ at the first jump time, we obtain
%
\begin{eqnarray} \nonumber
&&V(\phi_w(t_{1},0))
\geq V( \phi_w(t_{1},1))
\\&&
+ \sum_{j=0}^{J_{\phi_w}-1}L_D(\phi_w(t_{j+1},j),u_{ D}^w(t_{j+1},j)) \label{DiscreteCost2goNOihzw}
\\&&
+\sum_{j=1}^{J_{\phi_w}-1}\left( V(\phi_w(t_{j+1},j+1))-V(\phi_w(t_{j+1},j)) \right) \nonumber
\end{eqnarray}
In addition, given that $\phi_w(0,0)=\xi$, lower bounding $V(\phi_w(t_{1},0))$ in (\ref{ContinCost2goNOihzw}) by the right-hand side of (\ref{DiscreteCost2goNOihzw}), we obtain 
\begin{equation}
\begin{medsize}
\begin{split} 
V(\xi)
\geq& \sum_{j=0}^{J_{\phi_w}} \int_{t_{j}}^{t_{j+1}} L_C(\phi_w(t,j),u_{C}^w(t,j))dt+ V(\phi_w(t_1,0)) 
\\&
+\sum_{j=1}^{J_{\phi_w}} \left( V(\phi_w(t_{j+1},j))-V(\phi_w(t_j,j)) \right)  \\
\geq&
\sum_{j=0}^{J_{\phi_w}} \int_{t_{j}}^{t_{j+1}} L_C(\phi_w(t,j),u_{C}^w(t,j))dt
\\&
+ \sum_{j=0}^{J_{\phi_w}-1}L_D(\phi_w(t_{j+1},j),u_{D}^w(t_{j+1},j))
 \\
&+\sum_{j=1}^{J_{\phi_w}-1}\left( V(\phi_w(t_{j+1},j+1))-V(\phi_w(t_{j+1},j)) \right) 
\\&
+ V(\phi_w(t_{1},1))+\sum_{j=1}^{J_{\phi_w}} \left( V(\phi_w(t_{j+1},j))-V(\phi_w(t_j,j)) \right) 
\end{split}
\end{medsize}
\end{equation}
Since 
\begin{eqnarray} \nonumber
&& V(\phi_w(t_{1},1))
\\&&
+\sum_{j=1}^{J_{\phi_w}-1}\left( V(\phi_w(t_{j+1},j+1))-V(\phi_w(t_{j+1},j)) \right) \nonumber 
\\&&
+\sum_{j=1}^{J_{\phi_w}} \left( V(\phi_w(t_{j+1},j))-V(\phi_w(t_j,j)) \right) \nonumber
\\ &=& V(\phi_w(t_{{J_{\phi_w}}+1},J_{\phi_w})) +V(\phi_w(t_{1},1))\nonumber
\\&&
+\sum_{j=1}^{J_{\phi_w}-1}\left( V(\phi_w(t_{j+1},j+1)) \right) \nonumber 
-\sum_{j=1}^{J_{\phi_w}} \left(V(\phi_w(t_j,j)) \right) \nonumber
\\ &=&V(\phi_w(t_{{J_{\phi_w}}+1},J_{\phi_w})) 
\end{eqnarray}
then we have
\begin{eqnarray}\nonumber
V(\xi)
&\geq&\sum_{j=0}^{J_{\phi_w}} \int_{t_{j}}^{t_{j+1}} L_C(\phi_w(t,j),u_{C}^w(t,j))dt
\\&&
+ \sum_{j=0}^{J_{\phi_w}-1}L_D(\phi_w(t_{j+1},j),u_{D}^w(t_{j+1},j))
\\&&+V(\phi_w(t_{J_{\phi_w}+1},J_{\phi_w}))
\nonumber
\end{eqnarray}
By taking the limit when $(t_{J_{\phi_w}+1},J_{\phi_w}) \rightarrow \sup \dom \phi_w$, and given that (\ref{TerminalCondCG}) holds, we have
\begin{eqnarray} \nonumber
V(\xi)&\geq&\sum_{j=0}^{\sup_j \dom \phi_w} \int_{t_{j}}^{t_{j+1}} L_C(\phi_w(t,j),u_{C}^w(t,j))dt
\\&&
+ \sum_{j=0}^{{\sup_j \dom \phi_w}-1}L_D(\phi_w(t_{j+1},j),u_{D}^w(t_{j+1},j))
\nonumber \\
&& 
+\underset{(t,j) \in \textup{dom}\phi_w}{\limsup_{t+j\rightarrow \infty} } V(\phi_w(t,j))
\nonumber \\ 
&\geq& \mathcal{J}(\xi,u^w)
\label{BoundCostNOzw}
\end{eqnarray}
}
 \IfAutom{Now, 
 given that $(\phi_w,u^w)$ is 
 maximal, 
 %
 with  $u^w=(u_1^w,u_2^w)$, $u_1^w$ defined by $\IfIh{\kappa}{\gamma}_1$ as in (\ref{kHJBeqzsihc}) and (\ref{kBeqzsihc}), and $u_2^w$ defined by any $\bar{\IfIh{\kappa}{\gamma}}_2 \in \mathcal{K}_2$,}{}
 thanks to \eqref{eq:PrHamiltonC6} and \eqref{eq:PrHamiltonD6},
  from Corollary \ref{Cor:ProofDerivation} {and \eqref{TerminalCondCG}}, we have 
  \begin{equation}
  V(\xi) \geq \mathcal{J}(\xi,u^w).
    \label{BoundCostNOzw}
  \end{equation}
  }{}
%
{Finally, 
by proceeding as in item {d} above,
by applying the infimum on \NotAutomss{each side of }(\ref{BoundCostNOz}) 
over the set $\sj{\mathcal{S}_{\HS}^s}
(\xi)$, we obtain
\begin{equation}\label{a}
V(\xi) \leq  
\sj{    \underset{(\phi_s, \sj{u^s}) \in {\mathcal{S}_{\mathcal{H}}^s}(\xi)
}{\text{inf}}\mathcal{J}(\xi, \sj{u^s}) }
{=: \overline V(\xi)}
\end{equation}
\sj{Notice that the infimum in \eqref{a} is attained in ${\mathcal{S}_{\mathcal{H}}^s}(\xi)$ if there exists $(\phi_s, \sj{u^s}) \in {\mathcal{S}_{\mathcal{H}}^s}(\xi)$ such that $\mathcal{J}(\xi, \sj{u^s}) = V(\xi)$.}
By applying the supremum on \NotAutomss{each side of }(\ref{BoundCostNOzw})
over the set $\sj{\mathcal{S}_{\HS}^w}
 (\xi)$, \NotAutomss{we obtain}
\begin{equation}\label{b}
V(\xi) \geq 
\sj{\underset{(\phi_w, \sj{u^w}) \in {\mathcal{S}_{\mathcal{H}}^w}(\xi)}{\text{sup}}\mathcal{J}(\xi, \sj{u^w}) }
{=: \underline{ V}(\xi)}.
\end{equation}
\sj{Notice that the supremum in \eqref{b} is attained in ${\mathcal{S}_{\mathcal{H}}^w}(\xi)$ if there exists $(\phi_w, \sj{u^w}) \in {\mathcal{S}_{\mathcal{H}}^w}(\xi)$ such that $\mathcal{J}(\xi, \sj{u^w}) = V(\xi)$, at which $\underline{ V}(\xi)=V(\xi)$.}
Given that $V(\xi)=\mathcal{J}(\xi, u^*)$ from (\ref{OptimalCostihz}), we have that for any \pn{$\xi \in \mathcal{M}$}
, {each 
$(\phi^*,u^*)\in \pn{\mathcal{S}_{\HS }^X}
(\xi)$ with $u^*=(\kappa_1(\phi^*),\kappa_2(\phi^*))$ satisfies}
\begin{equation}
{\underline{V}(\xi)}
\leq
\mathcal{J}(\xi,u^*) 
\leq 
{\overline{V}(\xi)}
\label{Boundssz}
\end{equation}
Thanks to 
 $(\phi^*,u^*) \in \mathcal{S}_{\HS}^s(\xi) \cap \mathcal{S}_{\HS}^w(\xi)
 \pn{(\subset \mathcal{S}_{\HS }^X (\xi)) }
 $,} we have
 
\begin{equation}
\sj{
\overline{V}(\xi)=
\underset{(\phi_s, \sj{u^s}) \in {\mathcal{S}_{\mathcal{H}}^s}
(\xi)
}{\text{inf}}\mathcal{J}(\xi, \sj{u^s}) =\mathcal{J}(\xi, u^*)= V(\xi)
}
\label{supJk}
\end{equation}
\NotAutomss{and}
\begin{equation}
\sj{    \underline{V}(\xi) = \underset{(\phi_w, \sj{u^w}) \in {\mathcal{S}_{\mathcal{H}}^w}
    (\xi)
    }{\text{sup}}\mathcal{J}(\xi, \sj{u^w})=\mathcal{J}(\xi, u^*) = V(\xi)
    }
\label{infJk}
\end{equation}}
Since the infimum and supremum are attained in 
\eqref{supJk} and \eqref{infJk}
, respectively, \sj{by $(\phi^*,u^*)$,} \eqref{Boundssz}
leads to

\begin{equation}
\mathcal{J}(\xi,u^*) 
=\underset{(u_1,u_2) \in \pn{\mathcal{U}_{\HS}^X}
(\xi)}
{ \min_{u_{1}} \max_{u_{2}}}
\mathcal{J}(\xi,(u_1,u_2))
\label{LowerBounddz}
\end{equation}
Thus, from (\ref{OptimalCostihz}) and (\ref{LowerBounddz}), $V(\xi)$ is the value function for $\HS$, as in Definition \ref{ValueFunctionz} and from (\ref{Boundssz}), $\kappa$ is the saddle-point equilibrium as in Definition \ref{NashEqNonCoop}.
\eop\smallskip\vskip 3 pt
%
\pno{
  \begin{remark}{Connections between Theorem \ref{thHJBszih} and Problem $(\diamond)$}
  Given $\xi \in (\Pi(\overline{C}) \cup \Pi(D)) \cap \mathcal{M}$,
  if there exist a function $V$ satisfying the conditions in Theorem~\ref{thHJBszih},
  then
  a solution to Problem ($\diamond$) exists, namely there is an optimizer input action $u^*=(u_C^*,u_D^*)=((u_{C1}^*,u_{C2}^*),(u_{D1}^*,u_{D2}^*)) \in \pn{\mathcal{U}_{\HS  }^X} (\xi)$ that satisfies \eqref{SaddlePointIneq}, 
  and $V$ is {the} value function as in Definition \ref{ValueFunctionz}. 
  In addition, notice that
  the strategy $\kappa=(\kappa_C, \kappa_D)\in \mathcal{K}$ with elements as in (\ref{kHJBeqzsihc}) and (\ref{kBeqzsihc}) is such that every maximal solution to the closed-loop system $\HS_{\IfIh{\kappa}{\gamma}} $ from $\xi$ has a cost that is equal to the min-max in (\ref{problemzsih}), which is equal to the max-min.
  \label{SolExistenceZ}
  \end{remark}}
  
  \pno{
  \begin{remark}{Existence of a value function}
  Theorem~\ref{thHJBszih} does not explicitly rely on regularity conditions over the stage costs, {flow and jump} maps, convexity of $\J$, or compactness of the set of inputs \pn{$\mathcal{U}_\HS^X$}. 
  Sufficient conditions to guarantee the existence of a solution to Problem~$(\diamond)$ 
  are not {currently} available in the literature. One could expect that, as in any converse results, 
  guaranteeing the existence of a value function satisfying \eqref{HJBzsih} and \eqref{Bellmanzsih} would require the data of the system and the game to satisfy {certain} regularity properties. 
  \IfAutomss{Such is the case in optimal control problems \cite{GOEBEL2019153}.}
  {In the context of optimal control such regularity is required to guarantee existence \cite{GOEBEL2019153}.}
  %
  \end{remark}}
  
  \begin{remark}{Computation of the function $V$}
  \sj{In Mayer-type games with dynamics defined by hybrid finite-state automata as in \cite{tomlin2000game,ding2011toward,LevelSetMIanMitchell2000}, reachability-based approaches allow to synthesize safety controllers and compute the value function through the satisfaction of HJB conditions.
  }
  {In \sj{several other} cases, computing the saddle-point equilibrium strategy and the function $V$ satisfying the HJBI hybrid equations is difficult. 
  %
  This is a challenge already present in the certification of asymptotic stability. 
  However,} the complexity associated to the computation of a Lyapunov function does not diminish the contribution that the sufficient conditions for stability have had in the field. 
  In the same spirit, a contribution of Theorem \ref{thHJBszih}, as {an important step in games with dynamics defined as in \cite{65}, is in providing sufficient conditions that characterize value functions and saddle-point equilibria for such systems, similar to the results for continuous-time and discrete-time systems already available in the literature; \NotAutomss{see, } e.g., \cite{basar1999dynamic}.}
  \end{remark}
%
%
%
%
\sj{
  \begin{remark}{Time invariance of saddle-point equilibrium\NotAutomss{ and function $V$}}
  In general, games with fixed horizon define time-varying strategies \cite{leudo2022FH}, with some exceptions in the case in which the system is autonomous and the stage and terminal costs are stationary - see \cite[Remark 5.5]{basar1999dynamic}.
Problems with variable terminal time are in general non-stationary, due to the end time being optimized. This leads to optimal strategies and value function depending explicitly on time. 
{However}, when a 
terminal set 
is set, though there is a variable ({potentially} unbounded) terminal time, the optimality conditions are stationary, ({when the dynamics, the stage costs, and the terminal cost are time invariant}). Notice that for such a case, the terminal time becomes implicit rather than explicit because the game ends when the state enters the terminal set. {Hence,} the time at which the game ends is not fixed a priori but rather determined by the evolution of the state. In other words, the time dependency is replaced by a state-dependent stopping time. 
Thus, the value function and optimal strategies depend only on the state. 
  \end{remark}
  }
\IfIncd{}{\NotAutom{
\section{Linear Quadratic Hybrid Games
}
In this section, we consider a special case of our result that emerges {in practical scenarios with} hybrid systems with linear flow and jump maps and periodic jumps, {as in noise attenuation of cyber-physical systems,  see, e.g., \cite{ferrante2019certifying,possieri2020linear,possieri2017mathcal}}. We introduce a state variable $\tau$ that plays the role of a timer. Once $\tau$ reaches a fixed threshold $\bar{T}$, it triggers a jump in the state and resets $\tau$ to {zero. More precisely,
given $\bar{T}\in \mathbb{R}$, we consider a hybrid}
system with state $x=(x_p,\tau){=((x_{p1},x_{p2}),\tau)} \in \reals^n \times [0, \bar{T}]$, input $u=(u_C, u_D)=((u_{C1},u_{C2}),$
$(u_{D1},u_{D2})) \in \reals^{m_C}\times \reals^{m_D}$, and dynamics $\HS$ as in (\ref{Heq})\IfTp{,}{ with $N=2$ and} {defined} by
\begin{equation}\label{eq: CDFG priodic}
\begin{array}{rll}
C&:=&\mathbb{R}^n \times [0,\bar{T}] \times \mathbb{R}^{m_C}
\\
F(x,u_C)&:=&(A_C x_p+B_C u_C,1) 
\quad \forall(x,u_C) \in C
\\
D&:=&\mathbb{R}^n \times \{\bar{T}\} \times \mathbb{R}^{m_D}
\\
G(x,u_D)&:=&(A_D x_p+B_D u_D,0)
\quad \forall(x,u_D) \in D
\end{array}
\end{equation}
with $A_C=\left[\begin{smallmatrix}A_{C1} & 0 \\ 0 & A_{C2}\end{smallmatrix}\right], 
B_C=\left[ \begin{smallmatrix}B_{C1} \>B_{C2}\end{smallmatrix}\right],
A_D = \left[\begin{smallmatrix}A_{D1} & 0 \\ 0 & A_{D2}\end{smallmatrix}\right],$ and 
$ B_D = \left[ \begin{smallmatrix}B_{D1} \>B_{D2}\end{smallmatrix}\right] $.
The input $u_1=(u_{C1},u_{D1})$ is assigned by {player} $P_1$ and the input $u_2=(u_{C2},u_{D2})$ is assigned by {player} $P_2$. The problem of finding conditions for $u_1$ to minimize a cost functional $\mathcal{J}$ in the presence of the action $u_2$ that seeks to maximize it, is formulated as {the} two-player zero-sum game {as in Section \ref{Sec: Game form}.} Thus, \textcolor{black}{by solving Problem $(\diamond)$ for any $\xi \in \Pi(C)\cup \Pi(D)$}, 
the control objective is achieved. 
{Using Theorem \ref{thHJBszih},} the following result presents a tool for the solution of the optimal control problem for hybrid systems with linear maps and periodic jumps under an {adversarial} action.

\begin{corollary}{Hybrid Riccati equation for periodic jumps} Given {a hybrid system $\HS$ as in \eqref{Heq} defined by $(C,F,D,G)$ as in \eqref{eq: CDFG priodic}, let} $\bar{T}\in \mathbb{R}$, 
and,
{with the aim of pursuing minimum energy and distance to the origin, consider the cost functions $L_C(x,u_C):=x_p^\top Q_C x_p+ u_{C1}^\top R_{C1} u_{C1}+ u_{C2}^\top R_{C2} u_{C2}$, $L_D(x,u_D):=x_p^\top Q_D x_p+  u_{D1}^\top R_{D1} u_{D1}+ u_{D2}^\top R_{D2} u_{D2} $, and terminal cost $q(x):=x_p^\top P(\tau) x_p$ \textcolor{black}{defining $\J$ as in (\ref{defJTNCinc}), } {with} $Q_C,Q_D\in \mathbb{S}^n_+$, $R_{C1} \in \mathbb{S}^{m_{C_1}}_+$, $-R_{C2} \in \mathbb{S}^{m_{C_2}}_+$, $R_{D1} \in \mathbb{S}^{m_{D_1}}_+$, and $-R_{D2} \in \mathbb{S}^{m_{D_2}}_+$.} Suppose there exists a matrix {function} $P:[0,\bar{T}]\rightarrow \mathbb{S}^n_+ $ {that is} continuously differentiable {and} such that
{
\begin{equation}
\begin{split}
-\frac{d}{d\tau}P(\tau) =&-P(\tau)(  B_{C2} R_{C2}^{-1}B_{C2}^\top+B_{C1} R_{C1}^{-1}B_{C1}^\top) P(\tau)
\\&\hspace{0.1cm} +Q_C+P(\tau) A_C+A_C^\top P(\tau) \hspace{0.5cm} \forall \tau  \in (0,\bar{T}),
\end{split}
\label{DiffRiccatizsih}
\end{equation}}
\begin{equation}
\begin{array}{rrl}
-R_{D2} -B_{D2}^\top P(0) B_{D2},\>&
R_{D1} +B_{D1}^\top P(0) B_{D1}&
\in \mathbb{S}_{0+}^{m_D},
\end{array}
\label{zlqeqinv}
\end{equation}
the matrix $R_v=\left[\begin{smallmatrix}
R_{D1}+B_{D1}^\top P(0) B_{D1}
&
B_{D1}^\top P(0) B_{D2}
 \\
B_{D2}^\top P(0) B_{D1}
&
R_{D2}+B_{D2}^\top P(0) B_{D2}
\end{smallmatrix}\right] $ is invertible, and 
\begin{multline}
P (\bar{T})=Q_D+A_D^\top P(0) A_D \IfConf{\\}{}
-
\begin{bmatrix}
A_D^\top P(0) B_{D1}
 &
A_D^\top P(0) B_{D2}
\end{bmatrix}
R_v
^{-1}
\begin{bmatrix}
B_{D1}^\top P(0) A_D
\\
B_{D2}^\top P(0) A_D
\end{bmatrix}
\label{Riccatizsih} 
\end{multline}
{ where $A_C, B_{C1}, B_{C2}, A_D, B_{D1}$, and $B_{D2}$ are defined below \eqref{eq: CDFG priodic}.}
Then, the feedback law $\kappa:=(\kappa_C,\kappa_D)$, with values
\begin{eqnarray}
\kappa_C(x)=(-R_{C1}^{-1}B_{C1}^\top P(\tau) x_p,-R_{C2}^{-1}B_{C2}^\top P(\tau) x_p) \nonumber
\\ \hspace{0.4cm} \forall x \in \Pi(C),
\label{NashkCLQzsih}
\end{eqnarray}
\begin{equation}
\kappa_D(x)=-R_v
^{-1}
\begin{bmatrix}
B_{D1}^\top P(0) A_D
\\
B_{D2}^\top P(0) A_D
\end{bmatrix} x_p \hspace{0.5cm} \forall x \in \Pi(D)
\label{NashkDLQzsih}
\end{equation}
is the pure strategy saddle-point equilibrium for the two-player zero-sum hybrid game with periodic jumps. In addition, for each $x=(x_p,\tau) \in \Pi(\overline{C}) \cup \Pi(D)$, the value function is equal to $V(x):=x_p^\top P(\tau) x_p$.
\label{HREqPJ}
\end{corollary}
\begin{proof}
We show that when conditions (\ref{DiffRiccatizsih})-(\ref{Riccatizsih}) hold, by using Theorem \ref{thHJBszih} with \pn{$X=\emptyset$}, the value function is equal to the function $V$ and with 
 the feedback law $\kappa:=(\kappa_C,\kappa_D)$ with values as in (\ref{NashkCLQzsih}) and (\ref{NashkDLQzsih}), such a cost is attained. 
We can write (\ref{HJBzsih}) in Theorem \ref{thHJBszih} as
\begin{eqnarray}
&&0=
\underset{u_C=(u_{C1},u_{C2}) \in \Pi_u^C(x)}
{\min_{u_{C1}} \max_{u_{C2}}}  \mathcal{L}_C(x,u_C),
 \nonumber
\\ 
&&\mathcal{L}_C(x,u_C)= 
 x_p^\top Q_C x_p+ u_{C1}^\top R_{C1} u_{C1}
 + u_{C2}^\top R_{C2} u_{C2}
\nonumber \\
&&\hspace{0.9cm} {
+2x_p^\top P(\tau) (A_C x_p+B_C u_C)  
+ x_p^\top \frac{d}{d \tau}P(\tau) x_p}
\label{HJBzsihlq}
\end{eqnarray}
First, given that (\ref{DiffRiccatizsih}) holds, and $x_p^\top(P(\tau) A_C+A_C^\top P(\tau))x_p={2x_p^\top P(\tau)A_Cx_p}$ for every $x\in \Pi(C)$, one has
\begin{equation}
\begin{split}
\mathcal{L}_C(x,u_C)= \hspace{6.5cm} \nonumber\\
 x_p^\top P(\tau)(  B_{C2} R_{C2}^{-1}B_{C2}^\top+B_{C1} R_{C1}^{-1}B_{C1}^\top)P(\tau)  x_p
 \nonumber \\
 + u_{C1}^\top R_{C1} u_{C1}+ u_{C2}^\top R_{C2} u_{C2}
+ 2 x_p^\top  P(\tau) B_C u_C 
\end{split}
\end{equation}
The {first-order} necessary conditions for optimality
{
\IfAutom{
\begin{equation*}
\frac{\partial}{\partial u_{C1}}
\mathcal{L}_C(x,u_C)  \Big|_{u_{C}^*} =0,
\quad
\frac{\partial}{\partial u_{C2}}
\mathcal{L}_C(x,u_C) \Big|_{u_{C}^*} =0
\end{equation*}
{for all $(x,u_C)\in C$}}
{\begin{multline*}
\frac{\partial}{\partial u_{C1}} \left( 
 x_p^\top P(\tau)(  B_{C2} R_{C2}^{-1}B_{C2}^\top+B_{C1} R_{C1}^{-1}B_{C1}^\top)P(\tau)  x_p\right. \\ + u_{C1}^\top R_{C1} u_{C1}+ u_{C2}^\top R_{C2} u_{C2}
\\\left.\left.+ 2 x_p^\top  P(\tau)  (B_{C1} u_{C1} + B_{C2} u_{C2}) \right)\right|_{u_{C}^*} =0
\end{multline*}
\begin{multline*}
\frac{\partial}{\partial u_{C2} }\left( 
 x_p^\top P(\tau)(  B_{C2} R_{C2}^{-1}B_{C2}^\top+B_{C1} R_{C1}^{-1}B_{C1}^\top)P(\tau)  x_p\right.+\\ u_{C1}^\top R_{C1} u_{C1}+ u_{C2}^\top R_{C2} u_{C2}
 \\ \left.\left.+ 2 x_p^\top  P(\tau) (B_{C1} u_{C1} + B_{C2} u_{C2}) \right)\right|_{u_{C}^*} =0
\end{multline*}}
}
are satisfied by the point $u_C^*=(u_{C1}^*,u_{C2}^*)$, with values
\begin{equation}
u_{C1}^*=-R_{C1}^{-1}B_{C1}^\top P(\tau) x_p, 
\quad
u_{C2}^*=-R_{C2}^{-1}B_{C2}^\top P(\tau) x_p
\label{zLQHPu2*}
\end{equation}
{for each $x=(x_p,\tau)\in \Pi(C)$.}
Given that $R_{C1},-R_{C2} \in \mathbb{S}^{m_D}_+$, the second-order sufficient conditions for optimality 
%
{
\IfAutom{
\begin{equation*}
\frac{\partial^2}{\partial u_{C1}^2} 
 \mathcal{L}_C(x,u_C) \Big|_{u_C^*} \succeq 0,
 \quad
\frac{\partial^2}{\partial u_{C2}^2}
 \mathcal{L}_C(x,u_C)\Big|_{u_C^*} \preceq 0,
\end{equation*}
{for all $(x,u_C)\in C$}
}
{
\begin{multline*}
\frac{\partial^2}{\partial u_{C1}^2} \left( 
 x_p^\top P(\tau)(  B_{C2} R_{C2}^{-1}B_{C2}^\top+B_{C1} R_{C1}^{-1}B_{C1}^\top)P(\tau)  x_p\right. \\ + u_{C1}^\top R_{C1} u_{C1}+ u_{C2}^\top R_{C2} u_{C2}
 \\\left.\left.
+ 2 x_p^\top  P(\tau)  (B_{C1} u_{C1} + B_{C2} u_{C2}) \right)\right|_{u_{C1}^*} \succeq 0
\end{multline*}
\begin{multline*}
\frac{\partial^2}{\partial u_{C2}^2 }\left( 
 x_p^\top P(\tau)(  B_{C2} R_{C2}^{-1}B_{C2}^\top+B_{C1} R_{C1}^{-1}B_{C1}^\top)P(\tau)  x_p\right. \\ + u_{C1}^\top R_{C1} u_{C1}+ u_{C2}^\top R_{C2} u_{C2}
 \\\left.\left.+ 2 x_p^\top  P(\tau) (B_{C1} u_{C1} + B_{C2} u_{C2}) \right)\right|_{u_{C2}^*} \preceq 0
\end{multline*}}
}
hold, rendering $u_{C}^*$ as in 
\eqref{zLQHPu2*} as an optimizer of the min-max problem in (\ref{HJBzsihlq}).
In addition, it satisfies $\mathcal{L}_C(x,u_C^*)=0$, making $V(x)=x_p^\top P(\tau)x_p$ a solution {to} (\ref{HJBzsih}) in Theorem \ref{thHJBszih}. 

On the other hand, we can write (\ref{Bellmanzsih}) in Theorem \ref{thHJBszih} as
\begin{eqnarray}
&&x_p^\top P(\bar{T}) x_p =\underset{u_D=(u_{D1},u_{D2}) \in \Pi_u^D(x)}
{\min_{u_{D1}} \max_{u_{D2}}}  \mathcal{L}_D(x,u_D),
 \nonumber
\\ 
&&\mathcal{L}_D(x,u_D)= 
x_p^\top Q_D x_p+ u_{D1}^\top R_{D1} u_{D1}+ u_{D2}^\top R_{D2} u_{D2} \hspace{0.4cm}
\nonumber\\
&&\hspace{1.3cm} {
+ (A_D x_p+B_D u_D)^\top P(0) (A_D x_p+B_D u_D) }
\label{Bellmanzsihlq}
\end{eqnarray}
\NotAutom{which can be expanded as
{\small
\begin{multline}
\mathcal{L}_D(x,u_D)=  
 x_p^\top (Q_D+A_D^\top P(0) A_D) x_p + 2 x_p^\top A_D^\top  P(0) B_D u_D \\
 + u_{D1}^\top (R_{D1} +B_{D1}^\top P(0) B_{D1}) u_{D1}
 + u_{D2}^\top (R_{D2} +B_{D2}^\top P(0) B_{D2}) u_{D2}\\
 + u_{D1}^\top (B_{D1}^\top P(0) B_{D2}) u_{D2}
 + u_{D2}^\top (B_{D2}^\top P(0) B_{D1}) u_{D1}
\end{multline}}
}
{Similar to the case along flows, the first-order} necessary conditions for optimality
\NotAutom{
\begin{equation*}
\frac{\partial}{\partial u_{D1}}\left.
\mathcal{L}_D(x,u_D)  \right|_{u_{D}^*} =0,
\quad
\frac{\partial}{\partial u_{D2}}\left.
\mathcal{L}_D(x,u_D) \right|_{u_{D}^*} =0
\end{equation*}
}
are satisfied by the point $u_D^*=(u_{D1}^*,u_{D2}^*)$, such that, for each $x_p \in \Pi(D)$,
%
{
\begin{equation}
{\tiny
\begin{split}
u_{D}^*= \hspace{7.5cm}
\\-\begin{bmatrix}
R_{D1}+B_{D1}^\top P(0) B_{D1}
&
B_{D1}^\top P(0) B_{D2}
 \\
B_{D2}^\top P(0) B_{D1}
&
R_{D2}+B_{D2}^\top P(0) B_{D2}
\end{bmatrix} 
^{-1}
\begin{bmatrix}
B_{D1}^\top P(0) A_D
\\
B_{D2}^\top P(0) A_D
\end{bmatrix} x_p
\end{split}}
\label{zLQHPud*}
\end{equation}}
Given that (\ref{zlqeqinv}) holds, the second-order sufficient conditions for optimality 
\NotAutom{
\begin{equation*}
\frac{\partial^2}{\partial u_{D1}^2}\left. 
 \mathcal{L}_D(x,u_D) \right|_{u_D^*} \succeq 0,
 \quad
\frac{\partial^2}{\partial u_{D2}^2}\left.
 \mathcal{L}_D(x,u_D)\right|_{u_D^*} \preceq 0,
\end{equation*}
}
are satisfied, rendering $u_D^*$ as in (\ref{zLQHPud*}) as an optimizer of the \NotAutomss{min-max }problem in (\ref{Bellmanzsihlq}).
In addition, $u_D^*$ satisfies $\mathcal{L}_D(x,u_D^*)=x_p^\top P(\bar{T})x_p$, with $P(\bar{T})$ as in (\ref{Riccatizsih}), making $V(x)=x_p^\top P(\tau)x_p$ a solution of (\ref{Bellmanzsih})\NotAutomss{\> in Theorem \ref{thHJBszih}}. 

Then, given that $V$ is continuously differentiable on a neighborhood of $\Pi(C)$ and {that} Assumption \ref{AssLipsZ} holds, by applying Theorem \ref{thHJBszih}, in particular from (\ref{ResultValue}), for every $ \xi=(\xi_p,\xi_\tau) \in \Pi(\overline{C}) \cup \Pi(D)$ the value function is $\J^*(\xi)=\J(\xi,((u_{C1}^*,u_{D1}^*),(u_{C2}^*,u_{D2}^*))= \xi_p^\top P(\xi_\tau)\xi_p$. From (\ref{kHJBeqzsihc}) and (\ref{kBeqzsihc}), the feedback law $\kappa=(\kappa_C,\kappa_D)$ with values as in (\ref{NashkCLQzsih}) and (\ref{NashkDLQzsih}) is a pure strategy saddle-point equilibrium.
\end{proof}
}

\NotAutom{By following the same modeling approach and imposing conditions of the hybrid time domains, games for switching systems can be covered by Corollary \ref{HREqPJ}. By selecting appropriate stage costs, optimality is encoded in the satisfaction of the infinitesimal conditions instead of in the knowledge of specific solutions/trajectories. Note that for switching systems, the function $V$ might be independent of the timer state if the stage costs are independent of it as well.}

}
\section{Asymptotic Stability for Hybrid Games}
{We} present a result that connects optimality and asymptotic stability for \NotAutomss{two-player zero-sum }hybrid games. First, we introduce \IfAutomss{{positive definite functions.}}{definitions of some classes of functions.
\begin{definition}{Class-$\mathcal{K}_\infty$ functions}
	A function $\alpha: \mathbb{R}_{\geq 0} \rightarrow  \mathbb{R}_{\geq 0}$ is a class-$\mathcal{K}_\infty$ function, also written as $\alpha \in \mathcal{K}_\infty$, if $\alpha$ is
	 zero at zero, continuous, strictly increasing, and unbounded.
\end{definition}}

\begin{definition}{Positive definite functions}
{A} function $\rho: \reals_{\geq 0}\rightarrow \reals_{\geq 0}$ is positive definite, also written as $\rho \in \mathcal{PD}$, if $\rho(s)>0$ for all $s>0$ and $\rho(0)=0$.
{A} function $\rho: \reals^n \times \reals^{m} \rightarrow \reals_{\geq 0}$  is positive definite with respect to a set $\A \subset \reals^n$, in composition with $\kappa:\mathbb{R}^n \rightarrow \mathbb{R}^{m}$, also written as $\rho \in \mathcal{PD}_{\kappa}(\A)$, if $\rho(x,\kappa(x)) > 0$ for all $x \in \reals^n\setminus \A$ and $\rho(\A, {\kappa}(\A)) = \{0\}$. 
\label{PDsets}
\end{definition}

\NotAutomss{
\begin{definition}{
  Pre-asymptotic stability}
  A {closed} set $\mathcal{A} \subset \mathbb{R}^n$ is locally pre-asymptotically stable for a {hybrid closed-loop system} $\HS_\kappa$ as in (\ref{Hkeq}) if {it is} 
  \begin{itemize}
  \item  stable for $\HS_\kappa$, i.e., 
  if for every $\varepsilon>0$ there exists $\delta>0$ such that every solution $\phi$ to $\HS_\kappa$ with $|\phi(0,0)|_\mathcal{A} \leq \delta$ satisfies $|\phi(0,0)|_\mathcal{A} \leq \varepsilon$ for all $(t,j) \in \dom \phi$; {and}
  \item \pn{locally pre-attractive for $\HS_\kappa$, i.e., 
  there exists $\mu>0$ such that every solution $\phi$ to $\HS_\kappa$ with $|\phi(0,0)|_\mathcal{A} \leq \mu$ is bounded and, if $\phi$ is complete, then also $\lim_{t+j \rightarrow \infty} |\phi(t,j)|_\mathcal{A} =0$.} 
  \end{itemize}
  \end{definition}  
}
In the next result, we provide alternative conditions to those in Theorem \ref{thHJBszih} for the solution to Problem $(\diamond)$. 
\begin{lemma}{Equivalent conditions}\label{Lemma:EquivCond}
 Given $\HS_{\kappa}$ as in (\ref{Hkeq})  
{with data} $(C,F,D,G)$, \pn{the terminal set $X$, the feasible set $\mathcal{M}\subset \Pi(\overline{C}) \cup \Pi(D)$,} and {feedback} $\kappa:=(\kappa_C,\kappa_D)=((\kappa_{C1},\kappa_{C2}),(\kappa_{D1},\kappa_{D2})): \mathbb{R}^n \rightarrow \mathbb{R}^{m_C} \times \mathbb{R}^{m_D}$ {that satisfies \eqref{kHJBeqzsihc} and \eqref{kBeqzsihc}}, if there exists a function $V: \mathbb{R}^n \rightarrow \mathbb{R}$ that is continuously differentiable on a neighborhood of $\Pi(C)$ 
 such that\footnote{Notice that  $C_\kappa = \Pi(C)$ {and} $ D_\kappa = \Pi(D)$ when $\kappa_C(x)\in \Pi_u^C(x)$ for all $x\in \Pi(C)$ and $\kappa_D(x)\in \Pi_u^D(x)$ for all $x\in \Pi(D)$. \NotAutomss{In words, the feedback law $\kappa$ defining the {hybrid} closed-loop system $\HS_\kappa$ does not render input actions outside $C$ or $D$.}} $C_\kappa = \Pi(C)$ {and} $D_\kappa = \Pi(D)$, 
then (\ref{HJBzsih}), (\ref{Bellmanzsih}), (\ref{kHJBeqzsihc}), and (\ref{kBeqzsihc}) are satisfied 
if and only if
\begin{equation}
  \IfAutom{\mathcal{L}_C(x,\kappa_C(x))}{L_C(x,\kappa_C(x))
+\left\langle \nabla V(x),F(x,\kappa_C(x)) \right\rangle}
=0 \hspace{0.5cm} \forall x \in C_\kappa \pn{\cap \mathcal M} ,
\label{kHJBeqihza}
\end{equation}
{
\begin{equation}
\begin{split}
  \IfAutom{&\mathcal{L}_C(x,(u_{C1},\kappa_{C2}(x)))}{
 L_C(x,(u_{C1},\kappa_{C2}(x)))
+\left\langle \nabla V(x),F(x,(u_{C1},\kappa_{C2}(x))) \right\rangle
  }
\geq 0 \hspace{0.5cm}\\ \IfAutom{&\hspace{1.5cm}}{}\forall (x,u_{C1}): \> (x,(u_{C1},\kappa_{C2}(x)))\in C\pn{\cap \mathcal M},
\end{split}
\label{kHJBeqihzb}
\end{equation}
\begin{equation}
    \begin{split}
  \IfAutom{&\mathcal{L}_C(x,(\kappa_{C1}(x),u_{C2}))}{
 L_C(x,(\kappa_{C1}(x),u_{C2}))
+
\left\langle \nabla V(x),F(x,(\kappa_{C1}(x),u_{C2})) \right\rangle
  }
\leq 0 \hspace{0.5cm} \\ \IfAutom{&\hspace{1.5cm}}{} \forall (x,u_{C2}): \> (x,(\kappa_{C1}(x),u_{C2}))\in C\pn{\cap \mathcal M},
\end{split}
\label{kHJBeqihzc}
\end{equation}}
\begin{equation}
  \IfAutom{\mathcal{L}_D(x,\kappa_D(x))}{
 L_D(x,\kappa_D(x))
+ V(G(x,\kappa_D(x)))}=V(x) \hspace{0.5cm} \forall x \in D_\kappa\pn{\cap \mathcal M}, 
\label{kBeqihzd}
\end{equation}
{
\begin{equation}
\begin{array}{r@{}l}
  \IfAutom{&\mathcal{L}_D(x,(u_{D1},\kappa_{D2}(x)))}{ L_D(x,(u_{D1},\kappa_{D2}(x)))
+ V(G(x,(u_{D1},\kappa_{D2}(x))))}  \geq V(x) \\  \IfAutom{&\hspace{1cm}}{} \hspace{0.2cm} \forall (x,u_{D1}):\>(x,(u_{D1},\kappa_{D2}(x))) \in D\pn{\cap \mathcal M},
\end{array}
\label{kBeqihze}
\end{equation}
\begin{equation}
\begin{array}{r@{}l}
  \IfAutom{&\mathcal{L}_D(x,(\kappa_{D1}(x),u_{D2}))}{
 L_D(x,(\kappa_{D1}(x),u_{D2}))
+ V(G(x,(\kappa_{D1}(x),u_{D2})))
  }
\leq V(x) \\  \IfAutom{&\hspace{1cm}}{} \hspace{0.2cm} \forall (x,u_{D2}):\>(x,(\kappa_{D1}(x),u_{D2})) \in D\pn{\cap \mathcal M}.
\end{array}
\label{kBeqihzf}
\end{equation}}
\end{lemma}
\NotAutomss{{The proof is presented in the Appendix.}}
\NotConf{
  \par\noindent\textbf{Proof of Lemma \ref{Lemma:EquivCond}.}
  $(\rightarrow)$
From (\ref{kHJBeqzsihc}) and (\ref{kBeqzsihc}) we have
\begin{equation}
\underset{u_C=(u_{C1},u_{C2}) \in \Pi_u^C(x)}
{ \underset{u_{C1}}{\min}\> \underset{u_{C2}}{\max}} \left\{ 
L_C(x,u_C) \right.
+  \left.\left\langle \nabla V(x),F(x,u_C) \right\rangle \right\} =L_C(x,\kappa_C(x))
+\left\langle \nabla V(x) ,F(x,\kappa_C(x)) \right\rangle  \hspace{0.49cm} \forall x \in \Pi(C)
\label{kHJBeqihzproof}
\end{equation}
and

\begin{equation}
\underset{u_D=(u_{D1},u_{D2}) \in \Pi_u^D(x)}
{ \underset{u_{D1}}{\min} \underset{u_{D2}}{\max}} \left\{ 
L_D(x,u_D)
\right. + \left. V(G(x,u_D)) \right\}=L_D(x,\kappa_D(x))
+ V(G(x,\kappa_D(x))) \hspace{1.2cm} \forall x \in \Pi(D)
\label{kBeqihzproof}
\end{equation}
Thus, (\ref{HJBzsih}) and (\ref{kHJBeqihzproof}) imply 
\begin{equation}
L_C(x,\kappa_C(x))
+\left\langle \nabla V(x) ,F(x,\kappa_C(x)) \right\rangle=0  \hspace{0.49cm} \forall x \in \Pi(C), 
\label{kHJBeqihzproofa}
\end{equation}
while, (\ref{Bellmanzsih})  and (\ref{kBeqihzproof}) imply  
\begin{equation}
L_D(x,\kappa_D(x))
+ V(G(x,\kappa_D(x))) =V(x) \hspace{0.4cm} \forall x \in \Pi(D).
\label{kBeqihzproofa}
\end{equation}
From (\ref{kHJBeqihzproofa}) and (\ref{HJBzsih}), we have 
\begin{equation}
\underset{u_{C1}:(u_{C1},\kappa_{C2}(x)) \in \Pi_u^C(x)}
{ \underset{u_{C1}}{\min} 
} \left\{ 
L_C(x,(u_{C1},\kappa_{C2}(x))) \right.
+  \left.\left\langle \nabla V(x),F(x,(u_{C1},\kappa_{C2}(x))) \right\rangle \right\}  \geq 0 \hspace{0.4cm}\forall x \in \Pi(C)
\label{kHJBeqihzproofb}
\end{equation}
and
\begin{equation}
\underset{u_{C2}:(\kappa_{C1}(x), u_{C2}) \in \Pi_u^C(x)}
{ 
\underset{u_{C2}}{\max}
} \left\{ 
L_C(x,(\kappa_{C1}(x), u_{C2})) \right.
+  \left.\left\langle \nabla V(x),F(x,(\kappa_{C1}(x), u_{C2})) \right\rangle \right\}  \leq 0 \hspace{0.4cm} \forall x \in \Pi(C)
\label{kHJBeqihzproofc}
\end{equation}
which imply (\ref{kHJBeqihzb}) and (\ref{kHJBeqihzc}), respectively. Likewise, 
From (\ref{kBeqihzproofa}) and (\ref{Bellmanzsih}), we have 
\begin{equation}
\underset{u_{D1}:(u_{D1},\kappa_{D2}(x)) \in \Pi_u^D(x)}
{ \underset{u_{D1}}{\min} 
} \left\{ 
L_D(x,(u_{D1},\kappa_{D2}(x))) \right.
+  \left.    V(G(x,(u_{D1},\kappa_{D2}(x))))   \right\} \geq V(x) \hspace{0.4cm} \forall x \in \Pi(D)
\label{kBeqihzproofb}
\end{equation}
and
\begin{equation}
\underset{u_{D2}:(\kappa_{D1}(x), u_{D2}) \in \Pi_u^D(x)}
{ 
\underset{u_{D2}}{\max}
} \left\{ 
L_D(x,(\kappa_{D1}(x), u_{D2})) \right.
+ \left.  V(G(x,(\kappa_{D1}(x), u_{D2}))) \right\} \leq V(x) \hspace{0.4cm} \forall x \in \Pi(D)
\label{kBeqihzproofc}
\end{equation}
which imply (\ref{kBeqihze}) and (\ref{kBeqihzf}), respectively.\\
$(\leftarrow)$ Given $V$ and $\kappa:=(\kappa_C,\kappa_D)=((\kappa_{C1},\kappa_{C2}),(\kappa_{D1},\kappa_{D2}))$ such that (\ref{kHJBeqihza})-(\ref{kBeqihzf}) are satisfied, and such that $C_\kappa = \Pi(C), D_\kappa = \Pi(D)$, let us prove that $V$ and $\kappa$ satisfy (\ref{HJBzsih}), (\ref{Bellmanzsih}), (\ref{kHJBeqzsihc}), and (\ref{kBeqzsihc}). From (\ref{kHJBeqihza}) and (\ref{kHJBeqihzb}) we have 
\begin{equation}
\begin{split}
\underset{u_{C1}:(u_{C1},\kappa_{C2}(x)) \in \Pi_u^C(x)}
{ \underset{u_{C1}}{\min} 
} \left\{ 
L_C(x,(u_{C1},\kappa_{C2}(x))) \right.
+  \left.\left\langle \nabla V(x),F(x,(u_{C1},\kappa_{C2}(x))) \right\rangle \right\}  
\\ =L_C(x,\kappa_C(x))
+\left\langle \nabla V(x) ,F(x,\kappa_C(x)) \right\rangle =0 \hspace{0.4cm} \forall x \in \Pi(C)
\end{split}
\label{kHJBeqihzproof1left}
\end{equation}
and from (\ref{kHJBeqihza}) and (\ref{kHJBeqihzc}) we have 
\begin{equation}
\begin{split}
\underset{u_{C2}:(\kappa_{C1}(x), u_{C2}) \in \Pi_u^C(x)}
{ 
\underset{u_{C2}}{\max}
} \left\{ 
L_C(x,(\kappa_{C1}(x), u_{C2})) \right.
+  \left.\left\langle \nabla V(x),F(x,(\kappa_{C1}(x), u_{C2})) \right\rangle \right\}  
\\=L_C(x,\kappa_C(x))
+\left\langle \nabla V(x) ,F(x,\kappa_C(x)) \right\rangle =0 \hspace{0.4cm} \forall x \in \Pi(C)
\end{split}
\label{kHJBeqihzproof2left}
\end{equation}
Thus, (\ref{kHJBeqihzproof1left}) and (\ref{kHJBeqihzproof2left}) imply (\ref{HJBzsih}) and  (\ref{kHJBeqzsihc}).
Similarly, from (\ref{kBeqihzd}) and (\ref{kBeqihze}) we have 
\begin{equation}
\begin{split}
\underset{u_{D1}:(u_{D1},\kappa_{D2}(x)) \in \Pi_u^D(x)}
{ \underset{u_{D1}}{\min} 
} \left\{ 
L_D(x,(u_{D1},\kappa_{D2}(x))) \right.
+  \left.    V(G(x,(u_{D1},\kappa_{D2}(x))))   \right\} 
\\=L_D(x,\kappa_D(x))
+ V(G(x,\kappa_D(x))) =V(x) \hspace{0.4cm} \forall x \in \Pi(D)
\end{split}
\label{kBeqihzproof1left}
\end{equation}
and from (\ref{kBeqihzd}) and (\ref{kBeqihzf}) we have 
\begin{equation}
\begin{split}
\underset{u_{D2}:(\kappa_{D1}(x), u_{D2}) \in \Pi_u^D(x)}
{ 
\underset{u_{D2}}{\max}
} \left\{ 
L_D(x,(\kappa_{D1}(x), u_{D2})) \right.
+ \left.  V(G(x,(\kappa_{D1}(x), u_{D2}))) \right\}
\\ =L_D(x,\kappa_D(x))
+ V(G(x,\kappa_D(x))) =V(x) \hspace{0.4cm} \forall x \in \Pi(D)
\end{split}
\label{kBeqihzproof2left}
\end{equation}
Thus, (\ref{kBeqihzproof1left}) and (\ref{kBeqihzproof2left}) imply (\ref{Bellmanzsih}) and  (\ref{kBeqzsihc}).
}
\IfAutomss{{Under additional properties, the conditions above establish pre-asymptotic stability (pAS) of a set $\mathcal{A}$.
We define pAS for hybrid systems as in \cite[Definition  3.1
]{220}
}
}{}

\begin{ttheorem}{{Saddle-point equilibrium under the existence of a Lyapunov function}}\label{zCorStability}
Consider a two-player zero-sum hybrid game with
\pno{dynamics $\HS$ as in \eqref{Heq}
with data $(C,F,D,G)$ {satisfying Assumption \ref{AssLipsZ}},
and $\kappa:=(\kappa_C,\kappa_D): \mathbb{R}^n \rightarrow \mathbb{R}^{m_C} \times \mathbb{R}^{m_D}$
defining the closed-loop dynamics $\HS_{\kappa}$ as in (\ref{Hkeq})}  
such that $C_\kappa = \Pi(C)$ {and} $D_\kappa = \Pi(D)$. 
Given \pn{the terminal set $X$, the feasible set $\mathcal{M}\subset \Pi(\overline{C}) \cup \Pi(D)$, and} a closed set $\mathcal{A} \subset \Pi(C) \cup \Pi(D)$, continuous functions $L_C:
C \rightarrow \mathbb{R}_{\geq 0}$ and $L_D:
D \rightarrow \mathbb{R}_{\geq 0}$ 
defining the stage costs for flows and jumps, respectively, {and} $q:\reals^n \rightarrow \reals$ defining the terminal cost, suppose there exists a function $V: \mathbb{R}^n \rightarrow \mathbb{R}$ that is continuously differentiable on an open set containing $\overline{C_\kappa}$,  satisfying (\ref{kHJBeqihza})-(\ref{kBeqihzf}), and such that for each $\xi \in \pn{(\overline{C_\kappa} \cup D_\kappa)\cap \mathcal M}$, each $\phi \in \pn{\mathcal{S}_{\HS_{\kappa}  }^X}(\xi)$ satisfies (\ref{TerminalCondCG}). If 
one of the following conditions\NotAutomss{\footnote{The subindex in the set of positive definite functions $\mathcal{PD}_*$ denotes the feedback law that the functions in the set are composed with to satisfy the properties in Definition \ref{PDsets}.}} 
holds\IfAutomss{\footnote{{Class-$\mathcal{K}_\infty$ functions are defined as usual, see, e.g., \cite{65}.}}}

\begin{enumerate}[label=\arabic*)]
\item $L_C \in \mathcal{PD}_{\kappa_C}(\mathcal{A})$ and $L_D \in \mathcal{PD}_{\kappa_D}(\mathcal{A} )$;
\item $L_D \in \mathcal{PD}_{\kappa_D}(\mathcal{A})$ and there exists a continuous {function} $\eta \in \mathcal{PD}$ such that $L_C(x,\kappa_D(x)) \geq \eta(|x|_\A)$ for all $x \in \pn{C_\kappa\cap \mathcal{M}}$;
\item $L_C \in \mathcal{PD}_{\kappa_C}(\mathcal{A})$ and there exists a continuous {function} $\eta \in \mathcal{PD}$ such that $L_D(x,\kappa_D(x)) \geq \eta(|x|_\A)$ for all $x \in \pn{D_\kappa\cap \mathcal{M}}$;
\item {$L_C \equiv 0, L_D \in \mathcal{PD}_{\kappa_D}(\mathcal{A})$, and for each $r > 0$, there exist $\gamma_r\in \mathcal K_\infty$ and $N_r \geq 0$ such that for every solution $\phi \in \pn{\mathcal{S}_{\HS_\kappa}^X (\xi)}$, 
$|\phi(0, 0)|_\mathcal{A} \in (0, r]$, $(t, j) \in \dom \phi$, $t+j \geq T$ imply $j \geq \gamma_r (T)-N_r$
;}
\item {$L_C \in \mathcal{PD}_{\kappa_C}(\mathcal{A}), L_D \equiv 0$, and for each $r > 0$, there exist $\gamma_r\in \mathcal K_\infty$ and $N_r \geq 0$ such that for every solution $\phi \in \pn{\mathcal{S}_{\HS_\kappa}^X (\xi)}$, 
 $|\phi(0, 0)|_\mathcal{A} \in (0, r]$, $(t, j) \in \dom \phi$, $t+j \geq T$ imply $t \geq \gamma_r (T)-N_r$
;}
\item {$L_C(x,\kappa_C(x)) \geq -\lambda_C V(x)$ for all $x \in C_\kappa$, $L_D (x,\kappa_D(x)) \geq (1-e^{\lambda_D}) V(x)$ for all $x \in D_\kappa$, and there exist $\gamma > 0$ and $M>0$ such that, for each solution $\phi \in \pn{\mathcal{S}_{\HS_\kappa}^X (\xi)}$, 
 $(t, j) \in \dom \phi$ implies $\lambda_C t+ \lambda_D j \leq M-\gamma(t+j)$;}
\end{enumerate}
 then   
$\J^*(\xi)= V(\xi)$ for all $\xi \in \pn{(\overline{C_\kappa} \cup D_\kappa)\cap \mathcal{M}}$.
Furthermore, the feedback law
$\kappa$ is the saddle-point equilibrium (see Definition \ref{NashEqNonCoop}) and it renders $\A$ 
{
pre-asymptotically stable} 
for $\HS_{\kappa}$ \pn{with basin of attraction containing the largest sublevel set of $V$ contained in $\mathcal{M}$.}
\end{ttheorem}
%
%
%
\textbf{Proof.}
 Since{,} by assumption{,} we have that $C_\kappa = \Pi(C), D_\kappa = \Pi(D)$, and $V,\kappa:=(\kappa_C,\kappa_D)=((\kappa_{C1},\kappa_{C2}),$ $(\kappa_{D1},\kappa_{D2}))$ are such that (\ref{kHJBeqihza})-(\ref{kBeqihzf}) hold, then, thanks to Lemma \ref{Lemma:EquivCond}, $V$ and $\kappa$ satisfy (\ref{HJBzsih}), (\ref{Bellmanzsih}), (\ref{kHJBeqzsihc}), and (\ref{kBeqzsihc}). 
  Since\NotAutomss{\> in addition,} for each $\xi \in \pn{(\overline{C_\kappa} \cup D_\kappa) \cap \mathcal M}$, each $\phi \in \pn{\mathcal{S}_{\HS_{\kappa}  }^X}(\xi)$ satisfies (\ref{TerminalCondCG}), we have from Theorem \ref{thHJBszih} that $V$ is the value function as in (\ref{cost2gozsih}) for $\HS_{\kappa}$ at $\pn{(\overline{C_\kappa} \cup D_\kappa) \cap \mathcal M}$ 
  and the feedback law
  $\kappa$ with values (\ref{kHJBeqzsihc}), (\ref{kBeqzsihc}) is the saddle-point equilibrium for this game. 
{
Then, $V$ is a Lyapunov candidate for $\HS_{\kappa}$ \cite[Def. 3.16]{65} since $\overline{C_\kappa} \cup D_\kappa \subset \dom V=\reals^n$ and $V$ is continuously differentiable on an open set containing $\overline{C_\kappa}$}. From (\ref{kHJBeqihza}) {and} (\ref{kBeqihzd}), we have
  {
    \begin{equation}\label{kHJBeqihzabound}
    \begin{split}
  &\left\langle \nabla V(x),F(x,\kappa_C(x)) \right\rangle\leq- L_C(x,\kappa_C(x))\\
  &\hspace{5cm}
  \forall x \in C_\kappa \pn{\cap \mathcal M},
  \end{split}
  \end{equation}}
  \vspace*{-0.3cm}
  {
  \begin{equation}\label{kBeqihzbound}
  \begin{split}
  &V(G(x,\kappa_D(x)))-V(x)\leq-L_D(x,\kappa_D(x))  \\
  &\hspace{5cm}\forall x \in D_\kappa \pn{\cap \mathcal M}.
    \end{split} 
  \end{equation}}
  Moreover, {if}
  \begin{enumerate}[label=\alph*)]
  \item 
  Item 1, item 4, or item 5 above hold, {define}
  \begin{equation*}
  \begin{split}
  &\rho(x,\kappa(x)):=\\
  &\quad \quad\begin{medsize}\begin{cases}
  L_C(x,\kappa_C(x)) \quad &\textup{if}\> x\in C_\kappa \setminus D_\kappa
  \\
  \min\{L_C(x,\kappa_C(x)), L_D(x,\kappa_D(x))\}\quad &\textup{if}\> x\in C_\kappa\cap D_\kappa  
  \\
    L_D(x,\kappa_D(x)) \quad &\textup{if}\> x\in D_\kappa \setminus C_\kappa 
  \end{cases} \end{medsize}
  \end{split}
  \end{equation*}
  
  \item 
  Item 2 above holds, define
  \begin{equation*}
  \hspace{-0.4cm}
  \begin{split}
  \rho(x,\kappa(x)):=
  \begin{medsize}
  \begin{cases}
  \eta(|x|_\mathcal{A}) \quad &\textup{if}\> x\in C_\kappa \setminus D_\kappa
  \\
  \min\{\eta(|x|_\mathcal{A}), L_D(x,\kappa_D(x))\}\quad &\textup{if}\> x\in C_\kappa\cap D_\kappa  
  \\
    L_D(x,\kappa_D(x)) \quad &\textup{if}\> x\in D_\kappa \setminus C_\kappa
  \end{cases}  \end{medsize}
  \end{split}
  \end{equation*}
  
  \item 
 Item 3 above holds, {define}

  \begin{equation*}
  \hspace{-0.4cm}
  \begin{split}
  \rho(x,\kappa(x)):=
  \begin{medsize}
  \begin{cases}
  L_C(x,\kappa_C(x)) \quad &\textup{if}\> x\in C_\kappa \setminus D_\kappa
  \\
  \min\{L_C(x,\kappa_C(x)), \eta(|x|_\mathcal{A})\}\quad &\textup{if}\> x\in C_\kappa\cap D_\kappa  
  \\
    \eta(|x|_\mathcal{A}) \quad &\textup{if}\> x\in D_\kappa \setminus C_\kappa
  \end{cases}  \end{medsize}
  \end{split}
  \end{equation*}
  
\NotAutomss{
\item  Item 6 above holds, define
\begin{equation*}
  \begin{medsize}
    \begin{split}
    &\rho(x,\kappa(x)):=
    \begin{cases}
    \lambda_C V(x)
    \quad &\textup{if}\> x\in C_\kappa \setminus D_\kappa
    \\
    \min\{\lambda_C V(x),
    e^{\lambda_D} V(x)
    \}\quad &\textup{if}\> x\in C_\kappa\cap D_\kappa  
    \\
    e^{\lambda_D} V(x) 
    \quad &\textup{if}\> x\in D_\kappa \setminus C_\kappa
    \end{cases}
    \end{split}
    \end{medsize}
    \end{equation*}}
  \end{enumerate}
  \noindent Thus, 
  given that from (\ref{kHJBeqihzabound}) and (\ref{kBeqihzbound}), for each case above the function $\rho$ satisfies 
  \begin{equation}\label{kHJBeqihzaboundrho}
  \left\langle \nabla V(x),F(x,\kappa_C(x)) \right\rangle\leq- \rho(x,\kappa(x))\quad\forall x \in C_\kappa\pn{\cap \mathcal M},
  \end{equation}
  \begin{equation}\label{kBeqihzboundrho}
  V(G(x,\kappa_D(x)))-V(x)\leq-\rho(x,\kappa(x))\quad\forall x \in D_\kappa\pn{\cap \mathcal M}.
  \end{equation}
  \noindent {Thanks to 
  \cite[Theorem 3.19]{220}{, the set} $\A$ is 
  pAS for $\HS_{\kappa}$.} 
  \IfPers{
  On the other hand, if 
  \begin{enumerate}[label=\alph*)]
    \setcounter{enumi}{3}
  \item  Item 4 above holds, define
  \begin{equation*}
    \begin{medsize}
      \begin{split}
      &\rho(x,\kappa(x)):=
      \begin{cases}
      0
      \quad &\textup{if}\> x\in C_\kappa \setminus D_\kappa
      \\
      \min\{0,
       L_D(x,\kappa_D(x))\}\quad &\textup{if}\> x\in C_\kappa\cap D_\kappa  
      \\
        L_D(x,\kappa_D(x)) \quad &\textup{if}\> x\in D_\kappa \setminus C_\kappa
      \end{cases}
      \end{split}
      \end{medsize}
      \end{equation*}
    \item  Item 5 above holds, define
  \begin{equation*}
    \begin{medsize}
      \begin{split}
      &\rho(x,\kappa(x)):=
      \begin{cases}
      L_C(x,\kappa_C(x)) 
      \quad &\textup{if}\> x\in C_\kappa \setminus D_\kappa
      \\
      \min\{L_C(x,\kappa_C(x)),
      0
      \}\quad &\textup{if}\> x\in C_\kappa\cap D_\kappa  
      \\
      0
      \quad &\textup{if}\> x\in D_\kappa \setminus C_\kappa
      \end{cases}
      \end{split}
      \end{medsize}
      \end{equation*}
\item  Item 6 above holds, define
\begin{equation*}
  \begin{medsize}
    \begin{split}
    &\rho(x,\kappa(x)):=
    \begin{cases}
    \lambda_C V(x)
    \quad &\textup{if}\> x\in C_\kappa \setminus D_\kappa
    \\
    \min\{\lambda_C V(x),
    e^{\lambda_D} V(x)
    \}\quad &\textup{if}\> x\in C_\kappa\cap D_\kappa  
    \\
    e^{\lambda_D} V(x) 
    \quad &\textup{if}\> x\in D_\kappa \setminus C_\kappa
    \end{cases}
    \end{split}
    \end{medsize}
    \end{equation*}
    \end{enumerate}
  \noindent Thus, given the functions $\alpha_1, \alpha_2$ satisfying (\ref{zalphabound}),
  and given that from (\ref{kHJBeqihzabound}) and (\ref{kBeqihzbound}), the function $\rho$ satisfies \eqref{kHJBeqihzaboundrho} and \eqref{kBeqihzboundrho}, 
  thanks to \cite[Proposition 3.24]{65} (respectively,  \cite[Proposition 3.27]{65}, and \cite[Proposition 3.29]{65}), the set $\A$ is uniformly globally asymptotically stable for $\HS_{\kappa}$.}
  \IfAutomss{The case when item 6 holds follows directly from item e in \cite[Theorem 3.19]{220}.}{}
  \eop\smallskip\vskip 3 pt
\NotAutom{
\section{Saddle-Point Equilibria for Variable-terminal Time Hybrid Games}

In this section, we formulate a finite-horizon optimization problem to solve a two-player zero-sum hybrid game with free terminal time and a fixed terminal set, and provide the sufficient conditions to attain the solution. 

Following the formulation in Definition \ref{elements}, consider a two-player zero-sum game with dynamics $\mathcal{H} $ described by (\ref{Heq}) with data $(C,F,D,G)$. 
Let the closed set $X \subset \Pi(C) \cup \Pi(D)$ be the terminal constraint set, and we say a solution $(\phi,u)$ to $\HS$ is \emph{feasible} if there exists a finite $(T,J) \in \dom (\phi,u)$ such that $\phi(T,J) \in X$. In addition, let such $(T,J)$ be the terminal time of $(\phi,u)$ and the first time at which $\phi$ enters $X$, i.e., there does not exist any $(t,j) \in \dom \phi$ with $t+j < T+J$ 
 such that $\phi(t,j) \in X$.
%
 Since a well-defined game requires a unique correspondence from cost to input action, we require uniqueness of solutions to $\HS$ for a given input.
 Under Assumption \ref{AssLipsZ}, the conditions in Proposition \ref{UniquenessHu} are satisfied, rendering solutions to $\mathcal{H}$ for a given $u\in \mathcal{U}$ from $\xi$ unique.
 
 Given $\xi \in \Pi(C) \cup \Pi(D)$, a joint input action $u=(u_C, u_D)\in \U$, the stage cost for flows $L_C:\mathbb{R}^n  \times \mathbb{R}^{m_C} \rightarrow \mathbb{R}_{\geq 0}$, the stage cost for jumps $L_D:\mathbb{R}^n  \times \mathbb{R}^{m_D} \rightarrow \mathbb{R}_{\geq 0}$, and the terminal cost $q: \mathbb{R}^n \rightarrow \mathbb{R}$, 
we define the cost associated to the solution $(\phi,u)$ to $\mathcal{H}$ from $\xi$ that reaches $X$ for the first time at finite $(T,J) \in \realsgeq \times \nats$, under Assumption \ref{AssLipsZ}, as
\begin{equation}\label{defJvt}
  \begin{split}
  &{\mathcal{J}(\xi,u)}:=
  \sum_{j=0}^{J} \int_{t_{j}}^{t_{j+1}} L_C(\phi(t,j),u_{C}(t,j))dt  
  \\& \hspace{1.5cm}
  + 
  \sum_{j=0}^{J -1} 
  L_D(\phi(t_{j+1},j),u_{D}(t_{j+1},j))
  \\& \hspace{1.5cm}   +q(\phi(T,J))
\end{split}
\end{equation}
where $(T,J)= \max \dom \phi  $, $t_{J+1}=T$, and $\{t_j\}_{j=0}^{\sup_j \dom \phi}$ is a nondecreasing sequence associated to the definition of the hybrid time domain of $(\phi,u)$\IfAutomss{\>\cite[Definition 2.3]{65}}{; see Definition \ref{htd}}.

The feasible set $\mathcal{M}\subset \Pi(C) \cup \Pi(D)$ is the set of states $\xi$ such that there exists 
$(\phi, u) \in \hat{\mathcal{S}}_{\mathcal{H}  }(\xi) $ with $\phi(t,j)|_{\max \dom \phi}\in X$.

The solution to the two-player zero-sum game consists of solving the following problem.

\textit{Problem ($\star_\mathcal{M}$):} 
Given $\mathcal{M}\subset \Pi(\overline{C}) \cup \Pi(D)$ and $\xi \in \mathcal{M}$, under Assumption \ref{AssLipsZ}, solve
\begin{eqnarray}
\underset{u=(u_1,u_2)  \in  \mathcal{U} }
{ \underset{u_{1}}{\textup{minimize}}\>\> \underset{u_{2}}{\textup{maximize}}}
&& \mathcal{J}(\xi, u)
\\ \text{subject to}&&  
\phi(t,j)|_{\max \dom \phi} \in X \nonumber
\label{problemvt}
\end{eqnarray}
where $\phi$ is the 
maximal
state trajectory rendered by $u$ to $\HS$ from $\xi$.
\begin{remark}{\pnn{Terminal set} saddle-point equilibrium and min-max control}
A solution to Problem ($\star_\mathcal{M}$), when it exists, can be expressed in terms of the pure strategy saddle-point equilibrium \pnn{$\kappa=(\kappa_1,\kappa_2)$} for the two-player zero-sum variable-terminal-time game. 
Each $u^*=(u^*_1, u^*_2)$ that renders a state trajectory $\phi^*\in {\mathcal{R}} (\xi,u^*)$, with components defined as $\dom \phi^* \ni (t,j) \mapsto u^*_i(t,j)= \kappa_i(\phi^*(t,j))$  for each $i\in \IfTp{\{1,2\}}{\mathcal{V}}$, satisfies
\begin{eqnarray*}
u^*=
\underset{\mathclap
{\substack{u=(u_1,u_2) \in \mathcal{U}
\\ \phi(t,j)|_{\max \dom \phi}\in X
\\ \phi \in \mathcal{R}(\xi,u)
}}}
{\arg \min_{u_{1}} \max_{u_{2}}}
\>\> \mathcal{J}(\xi, u)
=
\underset{\mathclap
{\substack{u=(u_1,u_2) \in \mathcal{U}
\\ \phi(t,j)|_{\max \dom \phi}\in X
\\ \phi \in \mathcal{R}(\xi,u)
}}}{
{\arg  \max_{u_{2}}\min_{u_{1}}}}
\>\> \mathcal{J}(\xi, u)
\end{eqnarray*}
and it is referred to as a min-max control at $\xi$.
\end{remark}

\begin{definition}{Value function} 
Given $\mathcal{M}\subset \Pi(\overline{C}) \cup \Pi(D)$ and $\xi \in \mathcal{M}$, under Assumption \ref{AssLipsZ}, the value function at $\xi$ is given by 
\begin{equation}
\mathcal{J}^*(\xi):=
\underset{\mathclap
{\substack{u=(u_1,u_2) \in \mathcal{U}
\\ \phi(t,j)|_{\max \dom \phi}\in X
\\ \phi \in \mathcal{R}(\xi,u)
}}}
{\min_{u_{1}} \max_{u_{2}}}
\mathcal{J}(\xi, u)
=
\underset{\mathclap
{\substack{u=(u_1,u_2) \in \mathcal{U}
\\ \phi(t,j)|_{\max \dom \phi}\in X
\\ \phi \in \mathcal{R}(\xi,u)
}}}
{\max_{u_{2}}\min_{u_{1}}} 
\mathcal{J}(\xi, u)
 \label{cost2govt}
\end{equation}
\end{definition}
The following result provides sufficient conditions to characterize the value function, and the feedback law that attains it. It addresses the solution to Problem $(\star_\mathcal{M})$ for each $\xi \in (\Pi(\overline{C}) \cup \Pi(D)) \cap \mathcal{M}$ showing that the optimizer is the saddle-point equilibrium.
\begin{ttheorem}{Hamilton-Jacobi-Bellman-Isaacs (HJBI) for Problem ($\star_\mathcal{M}$)}
Given a two-player zero-sum hybrid game with dynamics $\mathcal{H}$ as in (\ref{Heq}) \IfTp{}{with $N=2$, }with data $(C,F,D,G)$, {satisfying Assumption \ref{AssLipsZ}},
 stage costs $L_C:\mathbb{R}^n  \times \mathbb{R}^{m_C} \rightarrow \mathbb{R}_{\geq 0} ,L_D:\mathbb{R}^n  \times \mathbb{R}^{m_D} \rightarrow \mathbb{R}_{\geq 0}$, terminal cost $q: \mathbb{R}^n \rightarrow \mathbb{R}$,
  and the feasible set $\mathcal{M}  \subset \Pi(C) \cup \Pi(D)$, suppose the following hold:
  \begin{enumerate}[label=\arabic*)]
    \item There exists a function $V: \mathbb{R}^n \rightarrow \mathbb{R}$ that is continuously differentiable on a neighborhood of\> $\Pi(C)$ that satisfies the \textit{Hamilton–Jacobi–Bellman-Isaacs hybrid equations given as} 
    \begin{equation}
      \begin{split}
      0&= \hspace{-0.6cm}  \underset{u_C=(u_{C1},u_{C2}) \in \Pi_u^C(x)}
      {\min_{u_{C1}} \max_{u_{C2}}}  
      {\mathcal{L}_C(x,u_C)}
      \\&= \hspace{-0.6cm} \underset{u_C=(u_{C1},u_{C2}) \in \Pi_u^C(x)}
      { \max_{u_{C2}}\min_{u_{C1}}}
      {\mathcal{L}_C(x,u_C)}
      \hspace{0.4cm}\forall x \in \Pi(C) \cap \mathcal{M},  
      \end{split}
      \label{HJBvt}
      \end{equation}
      where ${\mathcal{L}_C(x,u_C):=  L_C(x,u_C)+\left\langle \nabla V(x),F(x,u_C) \right\rangle }$,
      \begin{equation}
      \begin{split}
      V(x)&=\hspace{-0.6cm}\underset{u_D=(u_{D1},u_{D2}) \in \Pi_u^D(x)}
      {\min_{u_{D1}} \max_{u_{D2}}}  
      {\mathcal{L}_D(x,u_D)}
      \\&=\hspace{-0.6cm}\underset{u_D=(u_{D1},u_{D2}) \in \Pi_u^D(x)}
      { \max_{u_{D2}}\min_{u_{D1}}}  
      {\mathcal{L}_D(x,u_D)}
      \hspace{0.5cm} \forall x \in \Pi(D) \cap \mathcal{M},
      \end{split}
      \label{Bellmanvt}
      \end{equation}
      where ${\mathcal{L}_D(x,u_D):=  L_D(x,u_D) + V(G(x,u_D)) }$.
    \item For each $\xi \in \mathcal{M} $, each $(\phi, u)\in \hat{\mathcal{S}}_{\HS}(\xi)$ satisfies
\begin{equation}
V(\phi(t,j))=
q(\phi(t,j))
\hspace{1cm} \pnn{\forall (t,j) \in \textup{dom}\phi\>\> \textup{ s.t. }\>\> \phi(t,j) \in X}
\label{TerminalCondCGvt}
\end{equation}
\end{enumerate}
Then   
\begin{equation}
\J^*(\xi)= V(\xi) \hspace{2cm} \forall \xi \in \Pi(\overline{C}) \cup \Pi(D),
\label{ResultValue}
\end{equation}
and any 
feedback law
$\kappa:=(\kappa_C,\kappa_D)\NotConf{=((\kappa_{C1},\kappa_{C2}),(\kappa_{D1},\kappa_{D2}))}: \mathbb{R}^n \rightarrow \mathbb{R}^{m_C} \times \mathbb{R}^{m_D}$  with values 
{
\begin{equation}
  \hspace{-0.2cm}
\kappa_{C}(x)\in 
{\arg \min_{u_{C1}} \max_{u_{C2}}}_
{\mathclap{\substack{\\ \\ {u_C=(u_{C1},u_{C2}) \in \Pi_u^C(x)}} }}\
{\mathcal{L}_C(x,u_C)}
\hspace{0.7cm}
\forall x \in \Pi(C)\cap \mathcal{M}
\label{kHJBeqvt}
\end{equation}}
and
{
\begin{equation}
\kappa_{D}(x)\in 
{\arg \min_{u_{D1}} \max_{u_{D2}}}_
{\mathclap{\substack{\\ \\ {u_D=(u_{D1},u_{D2}) \in \Pi_u^D(x)}} }}\
{\mathcal{L}_D(x,u_D)}
\hspace{0.7cm}
\forall x \in \Pi(D) \cap \mathcal{M}
\label{kBeqvt}
\end{equation}
}
is a pure strategy saddle-point equilibrium for the two-player zero-sum variable-terminal time hybrid game with $\mathcal{J}_1=\mathcal{J}$, $\mathcal{J}_2=-\mathcal{J}$.
\label{thHJBszvt} 
\end{ttheorem}
\pnn{As in proof of Theorem (\ref{thHJBszih}), this conjecture relies on cost evaluation tools. When the feasible set $\mathcal{M}$ is known a priori, the set of states for which equations (\ref{HJBvt}) and (\ref{Bellmanvt}) need to be enforced could potentially be smaller than the sets of states studied in the infinite horizon counterpart. 
}
}
\section{Applications}
We illustrate {Theorem \ref{thHJBszih} in a disturbance rejection and a security problem by recasting them }as zero-sum hybrid games.
\subsection{
Robust hybrid LQR with aperiodic jumps}
\IfAutom{
In this section, \sj{we study a special case that emerges in applications featuring sample-and-hold control implementations and intermittent information scenarios --  see, e.g., \cite{163,274,ferrante2019certifying,possieri2020linear,possieri2017mathcal}.}
\sj{We introduce a state variable $\tau$ that plays the role of a timer. Once $\tau$ reaches an element in a 
threshold {set $\{T_1,T_2\}$ with $0 \leq T_1\leq T_2$,} 
 it potentially\footnote{When $T_1<T_2$, solutions can either evolve via flow or jump when $\tau = T_1$. A sequence $\{T_i\}_{i=1}^N$ can be handled similarly.} triggers a jump in the state and resets $\tau$ to zero. More precisely,
given $\bar{T}\in \mathbb{R}$, we consider a hybrid
system with state $x=(x_p,\tau){=((x_{p1},x_{p2}),\tau)} \in \reals^n \times [0, {T_2}
]$,}
input $u=(u_C, u_D)=((u_{C1},u_{C2}),$
$(u_{D1},u_{D2})) \in \reals^{m_C}\times \reals^{m_D}$, and dynamics $\HS$ as in (\ref{Heq})\IfTp{,}{ with $N=2$ and} {defined} by
\begin{equation}\label{eq: CDFG priodic}
\begin{array}{rll}
C&:=&\mathbb{R}^n \times [0,{T_2}
] \times \mathbb{R}^{m_C}
\\
F(x,u_C)&:=&(A_C x_p+B_C u_C,1) 
\quad \forall(x,u_C) \in C
\\
D&:=&\mathbb{R}^n \times \{{T_1,T_2}
\} \times \mathbb{R}^{m_D}
\\
G(x,u_D)&:=&(A_D x_p+B_D u_D,0)
\quad \forall(x,u_D) \in D
\end{array}
\end{equation}
with $A_C=\left[\begin{smallmatrix}A_{C1} & 0 \\ 0 & A_{C2}\end{smallmatrix}\right], 
B_C=\left[ \begin{smallmatrix}B_{C1} \>B_{C2}\end{smallmatrix}\right],
A_D = \left[\begin{smallmatrix}A_{D1} & 0 \\ 0 & A_{D2}\end{smallmatrix}\right],$ and 
$ B_D = \left[ \begin{smallmatrix}B_{D1} \>B_{D2}\end{smallmatrix}\right] $.
{Here,} the input $u_1:=(u_{C1},u_{D1})$ plays the role of the control  and is assigned by {player} $P_1$, and $u_2:=(u_{C2},u_{D2})$ is the disturbance input, which is assigned by {player} $P_2$. 
The problem of 
 {minimizing} the effect of the {worst-case} disturbance $u_2$ in the cost of complete solutions to $\HS$ is formulated as a two-player zero-sum game as in Section \ref{Sec: Game form}. Thus, by solving Problem $(\diamond)$ for every $\xi \in \Pi(C)\cup \Pi(D)$, 
the control objective is achieved. 
\IfPers{Note that the control objective does not include attenuation of $u_2$.} 

The following result presents a tool for the solution of optimal control problems for hybrid systems with linear maps and aperiodic jumps under an {adversarial} action.

\begin{proposition}{Hybrid Riccati equation for disturbance rejection with aperiodic jumps} Given a hybrid system $\HS$ as in \eqref{Heq} defined by $(C,F,D,G)$ as in \eqref{eq: CDFG priodic}
  \sj{with state $x=(x_p,\tau)
  \in \reals^n \times [0, {T_2}
  ]$, let $0\leq T_1\leq \bar{T}\leq T_2$,} 
  and,
  {with the aim of pursuing minimum energy and distance to the origin, consider the cost functions $L_C(x,u_C):=x_p^\top Q_C x_p+ u_{C1}^\top R_{C1} u_{C1}+ u_{C2}^\top R_{C2} u_{C2}$, $L_D(x,u_D):=x_p^\top Q_D x_p+  u_{D1}^\top R_{D1} u_{D1}+ u_{D2}^\top R_{D2} u_{D2} $, and terminal cost $q(x):=x_p^\top P(\tau) x_p$ \textcolor{black}{defining $\J$ as in (\ref{defJTNCinc}),} {with} $Q_C,Q_D\in \mathbb{S}^n_+$, $R_{C1} \in \mathbb{S}^{m_{C_1}}_+$, $-R_{C2} \in \mathbb{S}^{m_{C_2}}_+$, $R_{D1} \in \mathbb{S}^{m_{D_1}}_+$, and $-R_{D2} \in \mathbb{S}^{m_{D_2}}_+$.} Suppose there exists a matrix {function} $P:[0,{T_2}
  ]\rightarrow \mathbb{S}^n_+ $ {that is} continuously differentiable {and} such that
  {
  \begin{equation}
  \begin{split}
  -\frac{d}{d \tau}P(\tau) =&-P(\tau)(  B_{C2} R_{C2}^{-1}B_{C2}^\top+B_{C1} R_{C1}^{-1}B_{C1}^\top) P(\tau)
  \\&\hspace{-0.1cm} +Q_C+P(\tau) A_C+A_C^\top P(\tau) \hspace{0.5cm} \forall \tau  \in (0,{T_2}
  ),
  \end{split}
  \label{DiffRiccatizsih}
  \end{equation}}
  \begin{equation}
  \begin{array}{rrl}
  -R_{D2} -B_{D2}^\top P(0) B_{D2},\>&
  R_{D1} +B_{D1}^\top P(0) B_{D1}&
  \in \mathbb{S}_{0+}^{m_D},
  \end{array}
  \label{zlqeqinv}
  \end{equation}
  the matrix \IfAutomss{{$R_v:=R_D+[B_{D1}\>B_{D2}]^\top P(0) [B_{D1}\>B_{D2}]$, with $R_D:=\textup{diag}(R_{D1},R_{D2})$,}}{$R_v=\left[\begin{smallmatrix}
  R_{D1}+B_{D1}^\top P(0) B_{D1}
  &
  B_{D1}^\top P(0) B_{D2}
   \\
  B_{D2}^\top P(0) B_{D1}
  &
  R_{D2}+B_{D2}^\top P(0) B_{D2}
  \end{smallmatrix}\right]$} is invertible, and 
  \begin{equation}
    \hspace{-0.7cm}
    \begin{split}
  &P (\bar{T})=Q_D+A_D^\top P(0) A_D \IfConf{\\& \hspace{0.4cm}}{}
  - \IfAutomss{{A_D^\top P(0) [B_{D1}\>B_{D2}]}}{\left[
  \begin{smallmatrix}
  A_D^\top P(0) B_{D1}
  &
  A_D^\top P(0) B_{D2}
  \end{smallmatrix}
  \right]}
  R_v
  ^{-1} \IfAutomss{{[B_{D1}\>B_{D2}]^\top P(0) A_D}}{\left[
  \begin{smallmatrix}
  B_{D1}^\top P(0) A_D
  \\
  B_{D2}^\top P(0) A_D
  \end{smallmatrix}
  \right]}
\end{split}
  \label{Riccatizsih} 
  \end{equation}
  {at each $\bar T \in \{T_1,T_2\}$,} where $A_C, B_{C1}, B_{C2}, A_D, B_{D1}$, and $B_{D2}$ are defined below \eqref{eq: CDFG priodic}.
  Then, the feedback law $\kappa:=(\kappa_C,\kappa_D)$, with values
  \begin{eqnarray}
  \kappa_C(x)=(-R_{C1}^{-1}B_{C1}^\top P(\tau) x_p,-R_{C2}^{-1}B_{C2}^\top P(\tau) x_p) \nonumber
  \\ \hspace{0.4cm} \forall x \in \Pi(C),
  \label{NashkCLQzsih}
  \end{eqnarray}
  \begin{equation}
  \kappa_D(x)=-R_v
  ^{-1}
\IfAutomss{{[B_{D1}\>B_{D2}]^\top P(0) A_D}}{  \begin{bmatrix}
  B_{D1}^\top P(0) A_D
  \\
  B_{D2}^\top P(0) A_D
  \end{bmatrix}} x_p \hspace{0.5cm} \forall x \in \Pi(D)
  \label{NashkDLQzsih}
  \end{equation}
  is a pure strategy saddle-point equilibrium for the two-player zero-sum hybrid game with periodic jumps. In addition, for each $x=(x_p,\tau) \in \Pi(\overline{C}) \cup \Pi(D)$, the value function is equal to $V(x):=x_p^\top P(\tau) x_p$.
  \label{HREqPJ}
  \end{proposition}
  \begin{proof}
    We show that when \NotAutomss{conditions\>}(\ref{DiffRiccatizsih})-(\ref{Riccatizsih}) hold, by using Theorem \ref{thHJBszih}, the value function is equal to \NotAutomss{the function\>}$V$ and 
     the feedback law $\kappa:=(\kappa_C,\kappa_D)$ with values \NotAutomss{as\>}in (\ref{NashkCLQzsih}) and (\ref{NashkDLQzsih}), such a cost is attained. 
    We can write (\ref{HJBzsih}) in Theorem \ref{thHJBszih} as
    \begin{eqnarray}
    &&0=
    \underset{u_C=(u_{C1},u_{C2}) \in \Pi_u^C(x)}
    {\min_{u_{C1}} \max_{u_{C2}}}  \mathcal{L}_C(x,u_C),
     \nonumber
    \\ 
    &&\mathcal{L}_C(x,u_C)= 
     x_p^\top Q_C x_p+ u_{C1}^\top R_{C1} u_{C1}
     + u_{C2}^\top R_{C2} u_{C2}
    \nonumber \\
    &&\hspace{0.8cm} {
    +2x_p^\top P(\tau) (A_C x_p+B_C u_C)  
    + x_p^\top \frac{{d}}{{d} \tau}P(\tau)  x_p}
    \label{HJBzsihlq}
    \end{eqnarray}
    First, thanks to (\ref{DiffRiccatizsih}) and $x_p^\top(P(\tau) A_C+A_C^\top P(\tau))x_p={2x_p^\top P(\tau)A_Cx_p}$, for every $x\in \Pi(C)$, one has 
%
    \begin{equation}
    \begin{split}
    \mathcal{L}_C(x,u_C)= \hspace{6.5cm} \nonumber\\
     x_p^\top P(\tau)(  B_{C2} R_{C2}^{-1}B_{C2}^\top+B_{C1} R_{C1}^{-1}B_{C1}^\top)P(\tau)  x_p
     \nonumber \\
     + u_{C1}^\top R_{C1} u_{C1}+ u_{C2}^\top R_{C2} u_{C2}
    + 2 x_p^\top  P(\tau) B_C u_C 
    \end{split}
    \end{equation}

    The {first-order} necessary conditions for optimality
    {
    \IfAutom{
      \IfAutomss{{$\frac{\partial}{\partial u_{C1}}
      \mathcal{L}_C(x,u_C)  \big|_{u_{C}^*} =0,
      $ $\>\frac{\partial}{\partial u_{C2}}
      \mathcal{L}_C(x,u_C) \big|_{u_{C}^*} =0$}}{
    \begin{equation*}
    \frac{\partial}{\partial u_{C1}}
    \mathcal{L}_C(x,u_C)  \Big|_{u_{C}^*} =0,
    \quad
    \frac{\partial}{\partial u_{C2}}
    \mathcal{L}_C(x,u_C) \Big|_{u_{C}^*} =0
    \end{equation*}
      }
    {for all $(x,u_C)\in C$}}
    {\begin{multline*}
    \frac{\partial}{\partial u_{C1}} \left( 
     x_p^\top P(\tau)(  B_{C2} R_{C2}^{-1}B_{C2}^\top+B_{C1} R_{C1}^{-1}B_{C1}^\top)P(\tau)  x_p\right. \\ + u_{C1}^\top R_{C1} u_{C1}+ u_{C2}^\top R_{C2} u_{C2}
    \\\left.\left.+ 2 x_p^\top  P(\tau)  (B_{C1} u_{C1} + B_{C2} u_{C2}) \right)\right|_{u_{C}^*} =0
    \end{multline*}
    \begin{multline*}
    \frac{\partial}{\partial u_{C2} }\left( 
     x_p^\top P(\tau)(  B_{C2} R_{C2}^{-1}B_{C2}^\top+B_{C1} R_{C1}^{-1}B_{C1}^\top)P(\tau)  x_p\right.+\\ u_{C1}^\top R_{C1} u_{C1}+ u_{C2}^\top R_{C2} u_{C2}
     \\ \left.\left.+ 2 x_p^\top  P(\tau) (B_{C1} u_{C1} + B_{C2} u_{C2}) \right)\right|_{u_{C}^*} =0
    \end{multline*}}
    }
    are satisfied by the point $u_C^*=(u_{C1}^*,u_{C2}^*)$, with values
    \begin{equation}
    u_{C1}^*=-R_{C1}^{-1}B_{C1}^\top P(\tau) x_p, 
    \quad
    u_{C2}^*=-R_{C2}^{-1}B_{C2}^\top P(\tau) x_p
    \label{zLQHPu2*}
    \end{equation}
    {for each $x=(x_p,\tau)\in \Pi(C)$.}
    Since $R_{C1},-R_{C2} \in \mathbb{S}^{m_D}_+$, the second-order sufficient conditions for optimality 
    %
    \IfAutom{
      \IfAutomss{{$\frac{\partial^2}{\partial u_{C1}^2} 
      \mathcal{L}_C(x,u_C) \big|_{u_C^*} \succeq 0,
      $ $\>\frac{\partial^2}{\partial u_{C2}^2}
      \mathcal{L}_C(x,u_C)\big|_{u_C^*} \preceq 0,$}}{
    \begin{equation*}
    \frac{\partial^2}{\partial u_{C1}^2} 
     \mathcal{L}_C(x,u_C) \Big|_{u_C^*} \succeq 0,
     \quad
    \frac{\partial^2}{\partial u_{C2}^2}
     \mathcal{L}_C(x,u_C)\Big|_{u_C^*} \preceq 0,
    \end{equation*}}
    {hold for all $(x,u_C)\in C$,}
    }
    {
    \begin{multline*}
    \frac{\partial^2}{\partial u_{C1}^2} \left( 
     x_p^\top P(\tau)(  B_{C2} R_{C2}^{-1}B_{C2}^\top+B_{C1} R_{C1}^{-1}B_{C1}^\top)P(\tau)  x_p\right. \\ + u_{C1}^\top R_{C1} u_{C1}+ u_{C2}^\top R_{C2} u_{C2}
     \\\left.\left.
    + 2 x_p^\top  P(\tau)  (B_{C1} u_{C1} + B_{C2} u_{C2}) \right)\right|_{u_{C1}^*} \succeq 0
    \end{multline*}
    \begin{multline*}
    \frac{\partial^2}{\partial u_{C2}^2 }\left( 
     x_p^\top P(\tau)(  B_{C2} R_{C2}^{-1}B_{C2}^\top+B_{C1} R_{C1}^{-1}B_{C1}^\top)P(\tau)  x_p\right. \\ + u_{C1}^\top R_{C1} u_{C1}+ u_{C2}^\top R_{C2} u_{C2}
     \\\left.\left.+ 2 x_p^\top  P(\tau) (B_{C1} u_{C1} + B_{C2} u_{C2}) \right)\right|_{u_{C2}^*} \preceq 0
    \end{multline*}
    hold}
    rendering $u_{C}^*$ as in 
    \eqref{zLQHPu2*} as an optimizer of the min-max problem in (\ref{HJBzsihlq}).
    In addition, it satisfies $\mathcal{L}_C(x,u_C^*)=0$, making $V(x)=x_p^\top P(\tau)x_p$ a solution {to} (\ref{HJBzsih}) in Theorem \ref{thHJBszih}. 
    
    On the other hand, we can write (\ref{Bellmanzsih}) \NotAutomss{in Theorem \ref{thHJBszih}\>} as
    \begin{eqnarray}
    &&x_p^\top P(\bar{T}) x_p =\underset{u_D=(u_{D1},u_{D2}) \in \Pi_u^D(x)}
    {\min_{u_{D1}} \max_{u_{D2}}}  \mathcal{L}_D(x,u_D),
     \nonumber
    \\ 
    &&\mathcal{L}_D(x,u_D)= 
    x_p^\top Q_D x_p+ u_{D1}^\top R_{D1} u_{D1}+ u_{D2}^\top R_{D2} u_{D2} \hspace{0.4cm}
    \nonumber\\
    &&\hspace{1.3cm} {
    + (A_D x_p+B_D u_D)^\top P(0) (A_D x_p+B_D u_D) }
    \label{Bellmanzsihlq}
    \end{eqnarray}
    \NotAutom{which can be expanded as
    {\small
    \begin{multline}
    \mathcal{L}_D(x,u_D)=  
     x_p^\top (Q_D+A_D^\top P(0) A_D) x_p + 2 x_p^\top A_D^\top  P(0) B_D u_D \\
     + u_{D1}^\top (R_{D1} +B_{D1}^\top P(0) B_{D1}) u_{D1}
     + u_{D2}^\top (R_{D2} +B_{D2}^\top P(0) B_{D2}) u_{D2}\\
     + u_{D1}^\top (B_{D1}^\top P(0) B_{D2}) u_{D2}
     + u_{D2}^\top (B_{D2}^\top P(0) B_{D1}) u_{D1}
    \end{multline}}
    }
    {Similar to the case along flows, the first-order} necessary conditions for optimality
    \NotAutom{
    \begin{equation*}
    \frac{\partial}{\partial u_{D1}}\left.
    \mathcal{L}_D(x,u_D)  \right|_{u_{D}^*} =0,
    \quad
    \frac{\partial}{\partial u_{D2}}\left.
    \mathcal{L}_D(x,u_D) \right|_{u_{D}^*} =0
    \end{equation*}
    }
    are satisfied by the point $u_D^*=(u_{D1}^*,u_{D2}^*)$, such that, for each $x_p \in \Pi(D)$,
    %
    %
    \begin{equation}
      \IfAutomss{
      \begin{split}
      u_{D}^*= - R_v ^{-1}
	{[B_{D1}\>B_{D2}]^\top P(0) A_D} 
    	 x_p,
      \end{split}}{
      \begin{split}
      u_{D}^*= \hspace{7.5cm}
      \\{\tiny-\begin{bmatrix}
      R_{D1}+B_{D1}^\top P(0) B_{D1}
      &
      B_{D1}^\top P(0) B_{D2}
       \\
      B_{D2}^\top P(0) B_{D1}
      &
      R_{D2}+B_{D2}^\top P(0) B_{D2}
      \end{bmatrix} 
      ^{-1}
      \begin{bmatrix}
      B_{D1}^\top P(0) A_D
      \\
      B_{D2}^\top P(0) A_D
      \end{bmatrix} x_p}
      \end{split}}
    \label{zLQHPud*}
    \end{equation}
    Thanks to (\ref{zlqeqinv}), the second-order sufficient conditions for optimality 
    \NotAutom{
    \begin{equation*}
    \frac{\partial^2}{\partial u_{D1}^2}\left. 
     \mathcal{L}_D(x,u_D) \right|_{u_D^*} \succeq 0,
     \quad
    \frac{\partial^2}{\partial u_{D2}^2}\left.
     \mathcal{L}_D(x,u_D)\right|_{u_D^*} \preceq 0,
    \end{equation*}
    }
    are satisfied, rendering $u_D^*$ as in (\ref{zLQHPud*}) as an optimizer of the \NotAutomss{min-max }problem in (\ref{Bellmanzsihlq}).
    In addition, $u_D^*$ satisfies $\mathcal{L}_D(x,u_D^*)=x_p^\top P(\bar{T})x_p$ with {$\bar T \in \{T_1, T_2\}$ and} $P(\bar{T})$ as in (\ref{Riccatizsih}), making $V(x)=x_p^\top P(\tau)x_p$ a solution of (\ref{Bellmanzsih})\NotAutomss{\>in Theorem \ref{thHJBszih}}. 
    
    Then, given that $V$ is continuously differentiable on a neighborhood of $\Pi(C)$ and {that} Assumption \ref{AssLipsZ} holds, by applying Theorem \ref{thHJBszih}, in particular from (\ref{ResultValue}), for every $ \xi=(\xi_p,\xi_\tau) \in \Pi(\overline{C}) \cup \Pi(D)$ the value function is $\J^*(\xi)
    x= \xi_p^\top P(\xi_\tau)\xi_p$. From (\ref{kHJBeqzsihc}) and (\ref{kBeqzsihc}),  when $P_1$ plays $u_1^*$ defined by $\kappa_1=(\kappa_{C1},\kappa_{D1})$ with values as in (\ref{NashkCLQzsih}) and (\ref{NashkDLQzsih}), and $P_2$ plays any disturbance $u_2$ such that solutions to $\HS$ with data as in (\ref{eq: CDFG priodic}) are complete, then the cost is upper bounded by $\J(\xi,u^*)$,  
satisfying (\ref{SaddlePointIneq}).
\end{proof}
\NotAutomss{Notice that the saddle-point equilibrium $\kappa=(\kappa_1, \kappa_2)$ is composed by $P_1$ playing the {minimizing} strategy $\kappa_1$, and $P_2$ playing the maximizing disturbance $\kappa_2$ with values as in (\ref{NashkCLQzsih}) and (\ref{NashkDLQzsih}).}

\begin{figure}[t]
  \includegraphics[width=8.4cm]{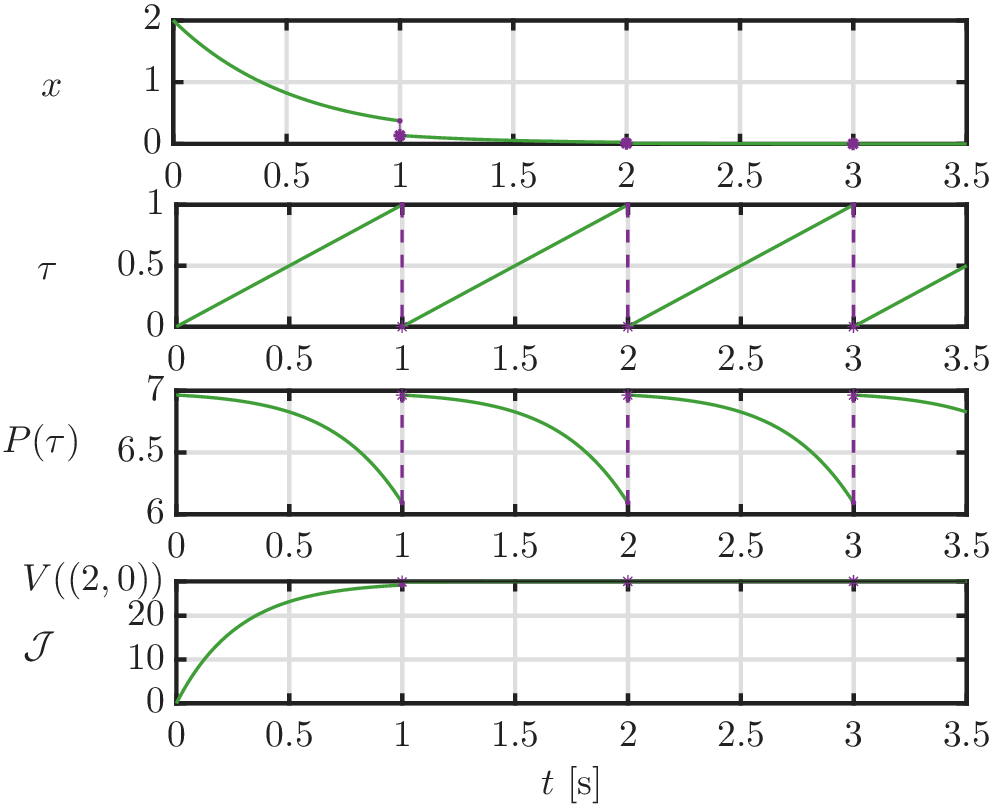}
  \vspace{-0.5cm}
  \caption{1D Robust hybrid LQR with periodic jumps. Dynamics as in \eqref{eq: CDFG priodic} with $A_C=1.8, B_C=[1, 1], A_D=2, B_D=[1, 1], Q_C=0.1,$ $R_C=\textup{diag}(1.304, -4), Q_D=1, R_D=\textup{diag}(1.304, -8),$ $P(0)=6.9653, T_1=T_2=1$.}
  \label{PeriodicGame}
\end{figure}

\IfAutomss{N}{Furthermore, n}otice that when $T_1 < T_2$ are finite, the jumps are not necessarily periodic, since they can occur when $\tau=T_1$ or when $\tau=T_2$. 
When $T_1 = T_2 =0$ we recover the discrete-time LQR robust problem, when $T_1 = T_2 =\infty$ we recover the
continuous-time LQR robust problem, and when $T_1 = T_2$ are finite, we have a hybrid game with periodic jumps as in in Figure \ref{PeriodicGame}.
}{
\begin{equation}\label{eq: CDFG RoHyLQR}
\begin{array}{rll}
C&\subset&\mathbb{R}^n  \times \mathbb{R}^{m_C} 
\\
F(x,u_C)&:=&A_C x+B_C u_C \\ &=&\left[\begin{smallmatrix}A_{C1} & 0 \\ 0 & A_{C2}\end{smallmatrix}\right]\left[ \begin{smallmatrix}x_1\\ x_2\end{smallmatrix}\right] +\left[ \begin{smallmatrix}B_{C1} \>B_{C2}\end{smallmatrix}\right]\left[ \begin{smallmatrix}u_{C1}\\ u_{C2}\end{smallmatrix}\right]
\\
D&\subset&\mathbb{R}^n \times  \mathbb{R}^{m_D} 
\\
G(x,u_D)&:=&A_D x+B_D u_D\\ &=&\left[\begin{smallmatrix}A_{D1} & 0 \\ 0 & A_{D2}\end{smallmatrix}\right]\left[ \begin{smallmatrix}x_1\\ x_2\end{smallmatrix}\right] +\left[ \begin{smallmatrix}B_{D1} \>B_{D2}\end{smallmatrix}\right]\left[ \begin{smallmatrix}u_{D1}\\ u_{D2}\end{smallmatrix}\right]
\end{array}
\end{equation}
where 
$C\cup D$ is nonempty. 
\sj{Here,} the input $u_1:=(u_{C1},u_{D1})$ plays the role of the control and $u_2:=(u_{C2},u_{D2})$ is the disturbance input. 
\pno{The problem of 
 upper bounding the effect of the disturbance $u_2$ in the cost of complete solutions to $\HS$ is formulated as the two-player zero-sum game as in Section \ref{Sec: Game form}.} Thus, \textcolor{black}{by solving Problem $(\diamond)$ for every $\xi \in \Pi(C)\cup \Pi(D)$,} 
the control objective is achieved. 
\IfPers{Note that the control objective does not include attenuation of $u_2$.} 

The following result presents a tool for the solution of the optimal control problem for hybrid systems with linear maps under the presence of disturbances.
\begin{corollary}{Hybrid Riccati equation for disturbance rejection} 
Given {a hybrid system $\HS$ as in \eqref{Heq} defined by $(C,F,D,G)$ as in \eqref{eq: CDFG RoHyLQR}, }
and,
{with the aim of pursuing minimum energy and distance to the origin, consider the cost functions 
 $L_C(x,u_C):=x^\top Q_C x+ u_{C1}^\top R_{C1} u_{C1}+ u_{C2}^\top R_{C2} u_{C2}$, $L_D(x,u_D):=x^\top Q_D x+  u_{D1}^\top R_{D1} u_{D1}+ u_{D2}^\top R_{D2} u_{D2} $, and terminal cost $q(x):=x^\top P x$, {defining $\mathcal{J}$ as in (\ref{defJTNCinc}),}
{with} $Q_C,Q_D\in \mathbb{S}^n_+$, $R_{C1} \in \mathbb{S}^{m_{C_1}}_+$, $-R_{C2} \in \mathbb{S}^{m_{C_2}}_+$, $R_{D1} \in \mathbb{S}^{m_{D_1}}_+$, $-R_{D2} \in \mathbb{S}^{m_{D_2}}_+$   and  $P\in \mathbb{S}^n_+ $. }
Suppose there exists a matrix $P\in \mathbb{S}^n_+ $ such that
%
\begin{multline}
0 =-P(  B_{C2} R_{C2}^{-1}B_{C2}^\top+B_{C1} R_{C1}^{-1}B_{C1}^\top)P\\+Q_C+P A_C+A_C^\top P,
\label{DiffRiccatizsihrlqr} 
\end{multline}
\begin{equation}
\begin{array}{rrl}
-R_{D2} -B_{D2}^\top P B_{D2},\>&
R_{D1} +B_{D1}^\top P B_{D1}&
\in \mathbb{S}_{0+}^{m_D},
\end{array}
\label{zlqeqinvrlqr}
\end{equation}
{the matrix $R_v=\left[\begin{smallmatrix}
R_{D1}+B_{D1}^\top PB_{D1}
&
B_{D1}^\top P B_{D2}
 \\
B_{D2}^\top P B_{D1}
&
R_{D2}+B_{D2}^\top P B_{D2}
\end{smallmatrix}\right] $
is invertible,} and
\begin{equation}
\begin{split}
0 =-P+Q_D+A_D^\top P A_D \hspace{3.1cm}\IfConf{\\}{}
-
\begin{bmatrix}
A_D^\top P B_{D1}
 &
A_D^\top P B_{D2}
\end{bmatrix}
R_v
^{-1}
\begin{bmatrix}
B_{D1}^\top P A_D
\\
B_{D2}^\top P A_D
\end{bmatrix}
\end{split}
\label{Riccatizsihrlqr} 
\end{equation}
Then, {for} the feedback law $\kappa_1:=(\kappa_{C1},\kappa_{D1})$ with values\footnote{The notation $R_v^{-1}(p,q)$ denotes the $(p,q)$ entry of the matrix $R_v^{-1}$.}
\begin{equation}
\kappa_{C1}(x)=-R_{C1}^{-1}B_{C1}^\top P x
\hspace{1cm} 
\forall x \in \Pi(C),
\label{NashkCLQzsihrlqr}
\end{equation}
{\small
\begin{equation}
\kappa_{D1}(x)=-[R_v^{-1}(1,1) \> R_v^{-1}(1,2)]
\left[\begin{smallmatrix}
B_{D1}^\top P A_D
\\
B_{D2}^\top P A_D
\end{smallmatrix}\right] x \hspace{0.4cm} \forall x \in \Pi(D),
\label{NashkDLQzsihrlqr}
\end{equation}}
\noindent the cost of complete solutions to $\HS$ from $\xi$ in the presence of any disturbance $u_2$ is upper bounded by $\xi^\top P\xi$. In addition, for each $x \in \Pi(\overline{C}) \cup \Pi(D)$, the value function is equal to $V(x):=x^\top P x$ and the maximizing disturbance is  given by $\kappa_2:=(\kappa_{C2},\kappa_{D2})$, with values
\begin{equation}
\kappa_{C2}(x)=-R_{C2}^{-1}B_{C2}^\top P x 
\hspace{1cm} 
\forall x \in \Pi(C),
\label{NashkC2LQzsihrlqr}
\end{equation}
{\small
\begin{equation}
\kappa_{D2}(x)=-[R_v^{-1}(2,1) \> R_v^{-1}(2,2)]
\left[\begin{smallmatrix}
B_{D1}^\top P A_D
\\
B_{D2}^\top P A_D
\end{smallmatrix}\right] x \hspace{0.4cm} \forall x \in \Pi(D).
\label{NashkD2LQzsihrlqr}
\end{equation}}
\label{HREqRobust}
\end{corollary}
{\begin{proof}
We show that when conditions (\ref{DiffRiccatizsihsec})-(\ref{Riccatizsihsec}) hold, by using the result in Theorem \ref{thHJBszih}, the value function is equal to the function $V$ and under 
 the feedback law as in (\ref{NashkCLQzsihrlqr}) and (\ref{NashkDLQzsihrlqr}) such a cost is attained in the presence of the maximizing disturbance given by (\ref{NashkC2LQzsihrlqr}) and (\ref{NashkD2LQzsihrlqr}). 
We can write (\ref{HJBzsih}) in Theorem \ref{thHJBszih} as
\begin{eqnarray}
&&0=
\underset{u_C=(u_{C1},u_{C2}) \in \Pi_u^C(x)}
{\min_{u_{C1}} \max_{u_{C2}}}  \mathcal{L}_C(x,u_C),
 \nonumber
\\ 
&&\mathcal{L}_C(x,u_C)= 
 x^\top Q_C x+ u_{C1}^\top R_{C1} u_{C1}+ u_{C2}^\top R_{C2} u_{C2}
 \nonumber \\
 &&\hspace{3.7cm}+2x^\top P (A_C x+B_C u_C)  
\label{HJBzsihlqrlqr}
\end{eqnarray}
First, 
given that (\ref{DiffRiccatizsihrlqr}) holds, and $x^\top(P A_C+A_C^\top P)x=x^\top (2PA_C)x$ for every $x\in \Pi(C)$, one has
\begin{eqnarray*}
&&\mathcal{L}_C(x,u_C)=
 x^\top P(  B_{C2} R_{C2}^{-1}B_{C2}^\top +B_{C1} R_{C1}^{-1}B_{C1}^\top) P x \hspace{0.8cm}
 \nonumber \\
  &&\hspace{1.6cm}
 + u_{C1}^\top R_{C1} u_{C1}+ u_{C2}^\top R_{C2} u_{C2}
+ 2 x^\top  P B_C u_C
\end{eqnarray*}
The first order necessary conditions for optimality
\NotAutom{\begin{equation*}
\frac{\partial}{\partial u_{C1}} \mathcal{L}_C(x,u_C)
\left.\right|_{u_{C}^*} 
=0,
\quad
\frac{\partial}{\partial u_{C2} } \mathcal{L}_C(x,u_C)
\left.\right|_{u_{C}^*} 
=0
\end{equation*}}
are satisfied by the point $u_C^*=(u_{C1}^*,u_{C2}^*)$, with values at any $x \in \Pi(C)$
\begin{equation}
u_{C1}^*=-R_{C1}^{-1}B_{C1}^\top P x,
\quad
u_{C2}^*=-R_{C2}^{-1}B_{C2}^\top P x.
\label{zLQHPu2*rlqr}
\end{equation}
Given that $R_{C1},-R_{C2} \in \mathbb{S}^{m_D}_+$, the second-order sufficient conditions for optimality 
\NotAutom{
\begin{equation*}
\frac{\partial^2}{\partial u_{C1}^2}  \mathcal{L}_C(x,u_C)
\left.\right|_{u_{C}^*}
 \succeq 0,
 \quad
\frac{\partial^2}{\partial u_{C2}^2 } \mathcal{L}_C(x,u_C)
\left.\right|_{u_{C}^*} 
\preceq 0
\end{equation*}}
hold, rendering $u_{C}^*$ with values as in 
(\ref{zLQHPu2*rlqr}) for each $x\in \Pi(C)$ as an optimizer of the min-max problem in (\ref{HJBzsihlqrlqr}).
In addition, 
it satisfies $\mathcal{L}_C(x,u_C^*)=0$, making $V(x)=x^\top P x$ a solution to (\ref{HJBzsih}) in Theorem \ref{thHJBszih}. 

On the other hand, we can write (\ref{Bellmanzsih}) in Theorem \ref{thHJBszih} as
\begin{eqnarray}
&&x^\top P x =\underset{u_D=(u_{D1},u_{D2}) \in \Pi_u^D(x)}
{\min_{u_{D1}} \max_{u_{D2}}}  \mathcal{L}_D(x,u_D),
 \nonumber
\\ 
&&\mathcal{L}_D(x,u_D)= 
x^\top Q_D x+ u_{D1}^\top R_{D1} u_{D1}+ u_{D2}^\top R_{D2} u_{D2}\nonumber\\
&&\hspace{1.7cm}
+ (A_D x+B_D u_D)^\top P (A_D x+B_D u_D) 
\label{Bellmanzsihlqrlqr}
\end{eqnarray}
\NotAutom{which can be expanded as 
{\scriptsize
\begin{multline}
\mathcal{L}_D(x,u_D)=  
 x^\top (Q_D+A_D^\top P A_D) x + 2 x^\top A_D^\top  P B_D u_D \\
 + u_{D1}^\top (R_{D1} +B_{D1}^\top P B_{D1}) u_{D1}
 + u_{D2}^\top (R_{D2} +B_{D2}^\top P B_{D2}) u_{D2}\\
 + u_{D1}^\top (B_{D1}^\top P B_{D2}) u_{D2}
 + u_{D2}^\top (B_{D2}^\top P B_{D1}) u_{D1}
\end{multline}}
}
The first order necessary conditions for optimality
\NotAutom{\begin{equation*}
\frac{\partial}{\partial u_{D1}}\left.
\mathcal{L}_D(x,u_D)  \right|_{u_{D}^*} 
=0,
\quad
\frac{\partial}{\partial u_{D2}}\left.
\mathcal{L}_D(x,u_D) \right|_{u_{D}^*} 
=0
\end{equation*}
}
are satisfied by the point $u_D^*=(u_{D1}^*,u_{D2}^*)$, with values for each $x \in \Pi(D)$
\begin{equation}
{\tiny
u_{D}^*=-\begin{bmatrix}
R_{D1}+B_{D1}^\top P B_{D1}
&
B_{D1}^\top P B_{D2}
 \\
B_{D2}^\top P B_{D1}
&
R_{D2}+B_{D2}^\top P B_{D2}
\end{bmatrix} 
^{-1}
\begin{bmatrix}
B_{D1}^\top P A_D
\\
B_{D2}^\top P A_D
\end{bmatrix} x_p
}
\label{zLQHPud*rlqr}
\end{equation}

Given that (\ref{zlqeqinvrlqr}) holds, the second-order sufficient conditions for optimality of $u_D^*$\NotAutom{, namely
\begin{equation*}
\frac{\partial^2}{\partial u_{D1}^2}\left. 
 \mathcal{L}_D(x,u_D) \right|_{u_D^*} \succeq 0,
\quad
\frac{\partial^2}{\partial u_{D2}^2}\left.
 \mathcal{L}_D(x,u_D)\right|_{u_D^*} \preceq 0,
\end{equation*}}
are satisfied, rendering $u_D^*$ in (\ref{zLQHPud*rlqr}) as an optimizer of the min-max problem in (\ref{Bellmanzsihlqrlqr}).
In addition, $u_D^*$
satisfies $\mathcal{L}_D(x,u_D^*)=x^\top P x$, with $P$ as in (\ref{Riccatizsihrlqr}), making $V(x)=x^\top P x$ a solution to (\ref{Bellmanzsih}) in Theorem \ref{thHJBszih}. 

Thus, given that $V$ is continuously differentiable in $\reals^n$  and Assumption \ref{AssLipsZ} holds,
 by applying Theorem \ref{thHJBszih}, in particular from (\ref{ResultValue}), {for every $ \xi \in \Pi(\overline{C}) \cup \Pi(D)$ the value function is $\J^*(\xi)=\J(\xi,((u_{C1}^*,u_{D1}^*),(u_{C2}^*,u_{D2}^*))= \xi^\top P \xi$. From (\ref{kHJBeqzsihc}) and (\ref{kBeqzsihc}), when $P_1$ plays $u_1^*$ defined by $\kappa_1=(\kappa_{C1},\kappa_{D1})$ with values as in (\ref{NashkCLQzsihrlqr}), (\ref{NashkDLQzsihrlqr}), and $P_2$ plays any disturbance $u_2$ such that solutions to $\HS$ with data as in (\ref{eq: CDFG RoHyLQR}) are complete, then the cost is upper bounded by $\J(\xi,u^*)$,  
satisfying (\ref{SaddlePointIneq}).
 }
\end{proof}}
\NotAutomss{Notice that the saddle-point equilibrium $\kappa=(\kappa_1, \kappa_2)$ is composed by $P_1$ playing the \sj{minimizing} strategy $\kappa_1$ with values as in (\ref{NashkCLQzsihrlqr}) and (\ref{NashkDLQzsihrlqr}), and $P_2$ playing the maximizing disturbance $\kappa_2$ with values as in (\ref{NashkC2LQzsihrlqr}) and (\ref{NashkD2LQzsihrlqr}).}
}
%
\NotAutomss{
\subsection{{
Robust control with flows-actuated nonunique solutions}}\label{NumericalEx}
{
  As illustrated next, there are useful families of hybrid systems for which a pure strategy saddle-point equilibrium exists. 
 The following problem {which has nonunique solutions to $\HS$ for a given feedback law} characterizes both the pure strategy saddle-point equilibrium and the value function in a two-player zero-sum game with a one-dimensional state, that is associated to player $P_1$, i.e., $n_1=1, n_2=0$. 

Consider a {hybrid} system $\HS$ with state $x \in \reals$, input $u_C:=(u_{C1},u_{C2}) \in \reals^2$, and dynamics 
\begin{equation}
\begin{array}{rll}
\dot{x}=&F(x,u_C):=a x+\langle B, u_C \rangle &  x \in [0,\delta]
\\
x^+=&G(x):=\sigma & x = \mu
\end{array}
\label{nuHGamez}
\end{equation}
where $a<0,B=(b_1, b_2) \in \reals^2$  and {$ \mu>\delta >\sigma>0$\footnote{Given that $ \mu>\delta$, flow from $\mu$ is not possible.}}. 
 Consider the cost functions $L_C(x,u_C):=x^2Q_C+u_C^\top R_C u_C$, $L_D(x):=P(x^2-\sigma^2)$, and terminal cost $q(x):=P x^2$, {defining $\mathcal{J}$ as in (\ref{defJTNCinc}),} with $R_C :=\left[ \begin{smallmatrix}R_{C1} & 0 \\ 0 & R_{C2}\end{smallmatrix}\right]$, {$Q_C,$ $R_{C1},$ $-R_{C2},$ $P\in \reals_{>0}$, such that} 
 \begin{equation}
 Q_C+2Pa-P^2(b_1^2 R_{C1}^{-1}+b_2^2 R_{C2}^{-1})=0.
 \label{PQ}  
 \end{equation}
\pn{Setting \pn{$X =\emptyset$,} the input $u_1:=(u_{C1},u_{D1})$ designed by {player} $P_1$ plays the role of the control and $u_2:=(u_{C2},u_{D2})$ is the disturbance input assigned by {player} $P_2$.} 
  This is formulated as a two-player zero-sum hybrid game \textcolor{black}{via solving Problem ($\diamond$)} {in Section \ref{subsec: probstat}}. The function $V(x):=P x^2$ is such that
\IfConf{
 {
\begin{equation}
\begin{split}
&\underset{u_C=(u_{C1},u_{C2}) \in \reals^2}
{\underset{u_{C1}}{\min}\> \underset{u_{C2}}{\max}}  
\mathcal{L}_C(x,u_C)
\\ 
&=
\underset{u_{C1} \in \mathbb{R}}{\min} \>\underset{u_{C2} \in \mathbb{R}}{\max} 
\left\{ (Q_C+ 2P a) x^2+  R_{C1}u_{C1}^2  \right.
\\ 
&\hspace{1.2cm}
+R_{C2}u_{C2}^2+  2xP ( b_1u_{C1}+ b_2u_{C2})  \big\} 
 = 0
\end{split}
 \label{HJBLQGnusCz}
\end{equation}
}
}{
\begin{eqnarray}
&&
\underset{u_C=(u_{C1},u_{C2}) \in \reals^2}
{\min_{u_{C1}} \max_{u_{C2}}}  
\left\{L_C(x,u_C) +\left\langle \nabla V(x),F(x,u_C) \right\rangle \right\}
\nonumber\\ 
&=&
\IfPers{%
{\min_{u_{C1}\in \mathbb{R}} \max_{u_{C2}\in \mathbb{R}}}  
\left( x^\top (Q_C+ 2P a) x+ u_C^\top R_C u_C
+  2x^\top P B u_C \right)
 \nonumber \\ &=&
{\min_{u_{C1}\in \mathbb{R}} \max_{u_{C2}\in \mathbb{R}}}  
\left\{ x^2 0.142+ 1.304u_{C1}^2-2916u_{C2}^2 
+  0.858x ( b_1u_{C1}+ b_2u_{C2}) \right\}\nonumber \\&=& }%
{\min_{u_{C1}\in \mathbb{R}} \max_{u_{C2}\in \mathbb{R}}}  
\left\{ (Q_C+ 2P a) x^2+  R_{C1}u_{C1}^2
+R_{C2}u_{C2}^2+  2xP ( b_1u_{C1}+ b_2u_{C2})  \right\} \nonumber
 \\ &=&
0 \label{HJBLQGnusCz}
\end{eqnarray}}
holds for all $x \in [0,\delta]$. In fact, the min-max in (\ref{HJBLQGnusCz}) is attained by 
$\kappa_C(x)=(-R_{C1}^{-1} b_1P x,-R_{C2}^{-1}b_2P x)$.
In particular, thanks to \eqref{PQ}, {we have} \IfAutom{{$\mathcal{L}_C(x,\kappa_C(x))=0.$}}
{
\begin{equation*}\hspace{-0.2cm}
\begin{array}{l}
L_C(x,\kappa_C(x)) +\left\langle \nabla V(x),F(x,\kappa_C(x)) \right\rangle 
  \nonumber\\ \hspace{2cm}
=
 (Q_C+2Pa) x^2+  R_{C1}\kappa^2_{C1}(x)
+R_{C2}\kappa_{C2}^2(x)\nonumber
\\ \hspace{4.3cm}
+  2xP ( b_1\kappa_{C1}(x)+ b_2\kappa_{C2}(x)) \nonumber
 \\\hspace{2cm}
=
 P^2(b_1^2R_{C1}^{-1}+b_2^2 R_{C2}^{-1}) x^2+  R_{C1}\kappa^2_{C1}(x)
 \nonumber
\\ \hspace{2.3cm} 
+R_{C2}\kappa_{C2}^2(x)+  2xP ( b_1\kappa_{C1}(x)+ b_2\kappa_{C2}(x))
 \nonumber
 \\ \hspace{2cm} =
\left[P^2(b_1^2 R_{C1}^{-1}+b_2^2 R_{C2}^{-1}) +  R_{C1}^{-1} b_1^2 P^2
+R_{C2}^{-1} b_2 P^2 \right.\nonumber
\\  \hspace{3.7cm}
\left.-  2P ( b_1R_{C1}^{-1} b_1 P+ b_2R_{C2}^{-1} b_2 P)\right]x^2 \nonumber
 \\ \hspace{2cm}
 =0 
\end{array}
\label{HJBLQGnusCzk}
\end{equation*}}
Then,  $V(x)= Px^2$ is a solution to (\ref{HJBzsih}).
 In addition, the function $V$ is such that
\begin{equation}
\begin{array}{r@{}l}
\underset{(u_{D1},u_{D2}) \in \reals^2}
{\underset{u_{D1}}\min\> \underset{u_{D2}}\max}  
\left\{
 L_D(x)
+ V(G(x)) \right\} =&
P x^2
\end{array}
\label{HJBLQGnusDz}
\end{equation}
at $x=\mu$\pno{,} which makes $V(x)= Px^2$ a solution to (\ref{Bellmanzsih}) \textcolor{black}{with saddle-point equilibrium $\kappa_C$}.
Given that $V$ is continuously differentiable on $\reals$, and that (\ref{HJBzsih}) and (\ref{Bellmanzsih}) hold thanks to (\ref{HJBLQGnusCz}) and (\ref{HJBLQGnusDz}), from Theorem \ref{thHJBszih} {we have that} the value function is
$
\mathcal{J}^*(\xi) :=
 P \xi^2
\label{HyLCost}
$
for any $\xi \in [0,\delta] \cup \{\mu\}$.

{To investigate the case of nonunique solutions {yielded by the feedback law $\kappa_C$}, now let $\delta {\geq}\mu>\sigma>0$ and notice that solutions can potentially flow or jump at $x=\mu$. The set of all maximal {solutions} from $\xi=\delta$ is denoted $\mathcal{R}_\kappa(\xi)=\{\phi_\kappa,\phi_h\}$. {The continuous solution} {$\phi_\kappa$} is such that $ \dom \phi_\kappa=\realsgeq \times \{0\}$, and is given by
$ \phi_\kappa(t,0)=  \delta \exp((a-R_{C1}^{-1} b_1P-R_{C2}^{-1}b_2P) t)  $ for all $t\in[0,\infty)$.
In simple words, $\phi_\kappa $ flows from $\delta${, and converges (exponentially fast)} to $0$. 
The maximal {solution} $\phi_h$ has domain $ \dom \phi_h=([0,t^h]\times \{0\}) \cup ([t^h, \infty)\times \{1\})$, and is given by
$ \phi_h(t,0)= \delta \exp((a-R_{C1}^{-1} b_1P-R_{C2}^{-1}b_2P) t) ,\hspace{0.5cm} \phi_h(t,1)=
\sigma \exp((a-R_{C1}^{-1} b_1P-R_{C2}^{-1}b_2P) (t-t^h)).$
In simple words, $\phi_h$ flows from $\delta$ to $\mu$ in $t^h$ seconds, then it jumps to $\sigma$, and flows converging (exponentially fast) to {zero}. Figure \ref{fig:my_label} illustrates this behavior.
%
%
By denoting the corresponding input signals as $u_\kappa=\kappa(\phi_\kappa)$ and $u_h=\kappa(\phi_h)$, 
we show in the bottom of Figure \ref{fig:my_label} that the cost of {the solutions $\phi_\kappa$ and $\phi_h$, yielded by $\kappa_C$,} equal $P\delta^2$. 
\begin{figure}[h]
    \centering
    \includegraphics[width=8.5cm]{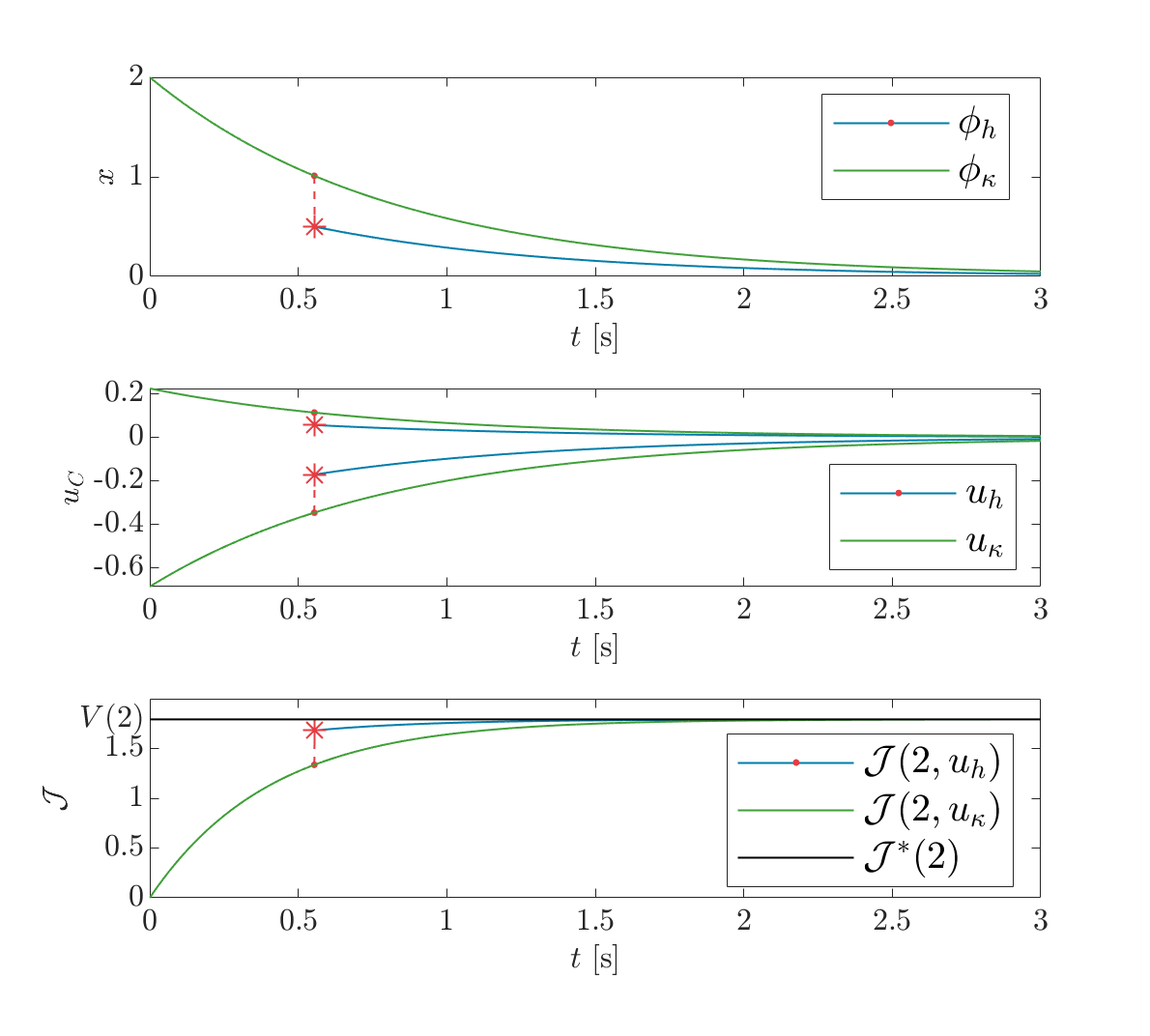}
    \vspace{-0.8cm}
    \caption{Nonunique solutions attaining min-max optimal cost for $a=-1, b_1=b_2=1, \delta=\xi=2, \mu=1, \sigma=0.5, Q_C=1,$ $R_{C1}=1.304, R_{C2}=-4$, and  $P =0.4481$. Continuous solution (green). Hybrid solution (blue and red).}
    \label{fig:my_label}
\end{figure}
This corresponds to {the optimal value with every maximal solution rendered by the equilibrium $\kappa_C$} from $\xi=2$ attaining {it}.}
\IfPers{
\begin{figure}[h]
    \centering
    \includegraphics[width=8cm]{Figures/NonuniqueZsCost.png}
    \caption{Saddle point behavior in the cost of continuous solutions from $\xi=2$ when varying the feedback gains around the optimal value. The cost is evaluated on solutions $(\phi,u)\in \pn{\mathcal{S}^X_\HS} (\xi) $ with feedback law variations specified by  $\epsilon_u$ and $\epsilon_u$ in $u=(\epsilon_u \kappa_1(\phi),\epsilon_w \kappa_2(\phi))$.}
    \label{Costs}
\end{figure}
This corresponds to the saddle-point equilibrium in Definition \ref{NashEqNonCoop}, and Remark \ref{SolExistenceZ}, with every maximal solution rendered by $\kappa$ from $\xi=2$ attaining the optimal cost.}

The next example illustrates Theorem \ref{zCorStability} and shows that our results, 
in the spirit of the Lyapunov theorem, only require that the conditions in Corollary \ref{zCorStability} hold. 
{
\begin{example}{Hybrid game with nonunique solutions}
 Let $\A =\{0\}$ and 
given that $L_C \in \mathcal{PD}_{\kappa_C }(\A 
)$, 
%
 (\ref{kHJBeqihza})-(\ref{kBeqihzf}) hold, and the function $s \mapsto \eta(s)=:P \frac{s^2}{2}$ is such that $L_D(x,\kappa_D(x)) \geq \eta(|x|_\A)$ for all $x \in D_\kappa$, by setting $\alpha_1(|x|_\A)=\underline{\lambda}(P)|x|^2$ and $\alpha_2(|x|_\A)=\overline{\lambda}(P)|x|^2$, 
from Corollary \ref{zCorStability} we have that $\kappa_C$ is the saddle-point equilibrium and renders $\A$ uniformly globally asymptotically stable for $\HS$ as in (\ref{nuHGamez}).
 \end{example}
 }
 }
}
%
%
\subsection{
Security jumps-actuated hybrid game}
Consider a hybrid system with state $x\in \reals^{n} $, input $ u_D=(u_{D1},u_{D2}) \in  \reals^{m_D}$, and dynamics $\HS$ as in (\ref{Heq})\IfTp{,}{ with $N=2$ and} described by
\begin{equation}
\begin{array}{rl@{}r@{}l}
\dot{x}=&F(x)&  x\>\>\> &\in C
\\
x^+=&A_D x +\left[ \begin{matrix}B_{D1} \>B_{D2}\end{matrix}\right]
\IfAutomss{{u_D}}{\left[ \begin{matrix}u_{D1}\\ u_{D2}\end{matrix} \right]}&\hspace{0.5cm} (x,u_D)& \in D
\end{array}
\label{eq: CDFG security}
\end{equation}
with 
$\sj{F: \mathbb{R}^{n}\rightarrow \mathbb{R}^{n}},
A_D 
\in \mathbb{R}^{n\times n}$,  and  $C \subset \mathbb{R}^{n}, D\subset\mathbb{R}^n \times  \mathbb{R}^{m_D}$, such that $C\cup \Pi(D)$ is nonempty. 
The input $u_{D1}$ plays the role of the control and $u_{D2}$ the disturbance input\sj{\footnote{\sj{``Jumps-actuated'' makes reference to the lack of inputs during flows. Since this example is general, any {condition involving $x$ and $u_D$} can be specified in $D$ to trigger jumps. }}}. 
{Here}, the problem of 
minimizing a cost \NotAutomss{functional }$\mathcal{J}$ in the presence of the maximizing attack {$u_{D2}$} is formulated as a two-player zero-sum game 
as in Problem $(\diamond)$. 

The following result presents a tool for the solution of optimal control problems for jumps-actuated hybrid systems and state-affine flow maps under a malicious input attack designed to cause as much damage as possible.
\begin{corollary}{Hybrid Riccati equation for security}
Given {a hybrid system $\HS$ as in \eqref{Heq} defined by $(C,F,D,G)$ as in \eqref{eq: CDFG security}, }
and,
{with the aim of pursuing minimum energy and distance to the origin, consider the cost functions 
$L_C(x,u_C):=0
$, $L_D(x,u_D):=x^\top Q_D x+  u_{D1}^\top R_{D1} u_{D1}+ u_{D2}^\top R_{D2} u_{D2} $, and  terminal cost $q(x):=x^\top P x$, {defining $\mathcal{J}$ as in (\ref{defJTNCinc}),} {with} $Q_D\in \mathbb{S}^n_+$, 
$R_{D1} \in \mathbb{S}^{m_{D_1}}_+$, $-R_{D2} \in \mathbb{S}^{m_{D_2}}_+$  and  $P\in \mathbb{S}^n_+ $. }
Suppose there exists a matrix $P\in \mathbb{S}^n_+ $ such that
\begin{equation}
0=
2x^\top P F(x)\hspace{1cm} \forall x \in \Pi(C),
\label{DiffRiccatizsihsec} 
\end{equation}
%
%
\begin{equation}
\begin{array}{rl}
-R_{D2} -B_{D2}^\top P B_{D2},R_{D1} +B_{D1}^\top P B_{D1} &\in \mathbb{S}_{0+}^{m_D},
\end{array}
\label{zlqeqinvsec}
\end{equation}
the matrix \IfAutomss{{$R_v:=R_D+[B_{D1}\>B_{D2}]^\top P [B_{D1}\>B_{D2}]$, with $R_D:=\textup{diag}(R_{D1},R_{D2})$,}}{$R_v=\left[\begin{smallmatrix}
R_{D1}+B_{D1}^\top PB_{D1}
&
B_{D1}^\top P B_{D2}
 \\
B_{D2}^\top P B_{D1}
&
R_{D2}+B_{D2}^\top P B_{D2}
\end{smallmatrix}\right] $} is invertible, and 
\begin{equation}
\begin{split}
0=-P+Q_D+A_D^\top P A_D \hspace{3.1cm}\IfConf{\\}{}
-\IfAutomss{{A_D^\top P [B_{D1}\>B_{D2}]}}{
  \begin{bmatrix}
A_D^\top P B_{D1}
 &
A_D^\top P B_{D2}
\end{bmatrix}}
  R_v
  ^{-1} \IfAutomss{{[B_{D1}\>B_{D2}]^\top P A_D}}{
\begin{bmatrix}
B_{D1}^\top P A_D
\\
B_{D2}^\top P A_D
\end{bmatrix}}
\end{split}
\label{Riccatizsihsec} 
\end{equation}
Then, the feedback law 
\IfAutomss{
\begin{equation}
\begin{split}
&\kappa_{D1}(x)=-{[R_{v_{11}}^{-1}\> R_{v_{12}}^{-1}]
[B_{D1}\>B_{D2}]^\top P A_D}
x  \\& \hspace{5.5cm}
 \forall x \in \Pi(D)
\end{split}
\label{NashkDLQzsihsec}
\end{equation}
}{\small
\begin{equation}
\kappa_{D1}(x)=-[R_v^{-1}(1,1) \> R_v^{-1}(1,2)]
\begin{bmatrix}
B_{D1}^\top P A_D
\\
B_{D2}^\top P A_D
\end{bmatrix} x  \hspace{0.4cm} \forall x \in \Pi(D)
\label{NashkDLQzsihsec}
\end{equation}}
minimizes the cost functional $\mathcal{J}$ in the presence of the maximizing attack $u_2$, given by
\IfAutomss{
\begin{equation}
\begin{split}
&\kappa_{D2}(x)=-{[R_{v_{21}}^{-1}\> R_{v_{22}}^{-1}]
[B_{D1}\>B_{D2}]^\top P A_D}
x  \\& \hspace{5.5cm}
 \forall x \in \Pi(D)
\end{split}
\label{NashkD2LQzsihsec}
\end{equation}
}{\begin{equation}
\kappa_{D2}(x)=-[R_{v_{21}}^{-1} \> R_{v_{22}}^{-1}]
\begin{bmatrix}
B_{D1}^\top P A_D
\\
B_{D2}^\top P A_D
\end{bmatrix} x  \hspace{0.4cm} \forall x \in \Pi(D)
\label{NashkD2LQzsihsec}
\end{equation}}
In addition, for each $x \in \Pi(\overline{C}) \cup \Pi(D)$, the value function is equal to $V(x):=x^\top P x$. 
\label{HREqSec}
\end{corollary}
{\begin{proof}
We show that when conditions (\ref{DiffRiccatizsihsec})-(\ref{Riccatizsihsec}) hold, by using the result in Theorem \ref{thHJBszih} \pn{with $X =\emptyset$,},  the value function is equal to the function $V$ and under 
 the feedback law as in (\ref{NashkDLQzsihsec}) such a cost is attained in the presence of the maximizing attack given by (\ref{NashkD2LQzsihsec}). 
\sj{We can write (\ref{HJBzsih}) 
as $0=2x^\top P F(x)$ for all $x \in \Pi(C)$, 
which is satisfied thanks to (\ref{DiffRiccatizsihsec}).
Likewise, we can write (\ref{Bellmanzsih}) 
as
\begin{eqnarray}
&&x^\top P x =\underset{u_D=(u_{D1},u_{D2}) \in \Pi_u^D(x)}
{\min_{u_{D1}} \max_{u_{D2}}}  \mathcal{L}_D(x,u_D),
 \nonumber
\\ 
&&\mathcal{L}_D(x,u_D)= 
x^\top Q_D x+ u_{D1}^\top R_{D1} u_{D1}+ u_{D2}^\top R_{D2} u_{D2}\nonumber\\
&&\hspace{1.7cm}+ (A_D x+B_D u_D)^\top P (A_D x+B_D u_D) 
\label{Bellmanzsihlqsec}
\end{eqnarray}
}
\NotAutom{which can be expanded as
\begin{multline}
\mathcal{L}_D(x,u_D)=  
 x^\top (Q_D+A_D^\top P A_D) x + 2 x^\top A_D^\top  P B_D u_D \\
 + u_{D1}^\top (R_{D1} +B_{D1}^\top P B_{D1}) u_{D1}
 + u_{D2}^\top (R_{D2} +B_{D2}^\top P B_{D2}) u_{D2}\\
 + u_{D1}^\top (B_{D1}^\top P B_{D2}) u_{D2}
 + u_{D2}^\top (B_{D2}^\top P B_{D1}) u_{D1}
\end{multline}}
The first order necessary conditions for optimality\NotAutom{
\begin{equation*}
\frac{\partial}{\partial u_{D1}}\left.
\mathcal{L}_D(x,u_D)  \right|_{u_{D}^*} =0,
\quad
\frac{\partial}{\partial u_{D2}}\left.
\mathcal{L}_D(x,u_D) \right|_{u_{D}^*} =0
\end{equation*}}
are satisfied by $u_D^*\hspace{-0.1cm}=\hspace{-0.1cm}(u_{D1}^*,u_{D2}^*)$, defined for each $x\in \Pi(D)$ as
\begin{equation}
u_{D}^*=-
{\tiny
\IfAutomss{
{R_v^{-1} [B_{D1}\>B_{D2}]^\top P A_D
}}{\begin{bmatrix}
R_{D1}+B_{D1}^\top P B_{D1}
&
B_{D1}^\top P B_{D2}
 \\
B_{D2}^\top P B_{D1}
&
R_{D2}+B_{D2}^\top P B_{D2}
\end{bmatrix} 
^{-1}
\begin{bmatrix}
B_{D1}^\top P A_D
\\
B_{D2}^\top P A_D
\end{bmatrix}}} x
\label{zLQHPud*sec}
\end{equation}

Given that (\ref{zlqeqinvsec}) holds, the second-order sufficient conditions for optimality
\NotAutom{
\begin{equation*}
\frac{\partial^2}{\partial u_{D1}^2}\left. 
 \mathcal{L}_D(x,u_D) \right|_{u_D^*} \succeq 0,
\quad
\frac{\partial^2}{\partial u_{D2}^2}\left.
 \mathcal{L}_D(x,u_D)\right|_{u_D^*} \preceq 0,
\end{equation*}}
are satisfied, rendering $u_D^*$ as in (\ref{zLQHPud*sec}) as an optimizer of the \NotAutomss{min-max }problem in (\ref{Bellmanzsihlqsec}).
In addition, $u_D^*$
satisfies $\mathcal{L}_D(x,u_D^*)=x^\top P x$, with $P$ as in (\ref{Riccatizsihsec}), leading $V(x)=x^\top P x$ as a solution of (\ref{Bellmanzsih})\NotAutomss{\>in Theorem \ref{thHJBszih}}. 

Thus, given that $V$ is continuously differentiable in $\reals^n$ and Assumption (\ref{AssLipsZ}) holds, by applying Theorem \ref{thHJBszih}, in particular from (\ref{ResultValue}), 
{for every $ \xi \in \Pi(\overline{C}) \cup \Pi(D)$ the value function is $\J^*(\xi)=\NotAutomss{\J(\xi,(u_{D1}^*,u_{D2}^*))=} \xi^\top P \xi$. From (\ref{kHJBeqzsihc}) and (\ref{kBeqzsihc}), when $P_2$ plays $u_2^*$ defined by $\kappa_{D2}$  as in (\ref{NashkD2LQzsihsec}), $P_1$ minimizes the cost of complete solutions to $\HS$ 
by playing $u_1^*$ defined by $\kappa_{D1}$ as in (\ref{NashkDLQzsihsec}), attaining $ \J(\xi,u^*)$, and satisfying (\ref{SaddlePointIneq}).}
\end{proof}}
\textcolor{black}{Notice that the saddle-point equilibrium $\kappa_D:=(\kappa_{D1}, \kappa_{D2})$ is composed by $P_1$ playing the minimizer strategy $\kappa_{D1}$ as in (\ref{NashkDLQzsihsec}), and $P_2$ playing the maximizing attack $\kappa_{D2}$ as in (\ref{NashkD2LQzsihsec}).}
\sj{Given that the flow map $F$ does not have inputs, as long as it satisfies \eqref{DiffRiccatizsihsec}, the two-player discrete time Riccati algebraic equation \eqref{Riccatizsihsec} allows to charazterize the optimal at-jumps-only strategy.}
\IfPers{
\begin{example}{Linear quadratic differential games}
Consider a two-player zero-sum hybrid game with state $x \in \reals^n$, input $u=(u_C, u_D) \in \reals^{m_C}\times \reals^{m_D}$, and dynamics $\HS$ as in (\ref{Heq}) with \IfTp{}{$N=2$ and}
\begin{eqnarray}
C&=&\mathbb{R}^n \times \mathbb{R}^{m_C} \nonumber
\\
F(x,u_C)&=&A_C x+B_C u_C \nonumber
\label{zsLinearGame}
\end{eqnarray}
and jump set $D=\emptyset$, where $A_C:=\textup{blkdiag}_{i \in \IfTp{\{1,2\}}{\mathcal{V}}}\{A_i\} \in \mathbb{R}^{n\times n}$ and $B_C:=\textup{blkdiag}_{i \in \IfTp{\{1,2\}}{\mathcal{V}}} \{B_{Ci}\} \in \mathbb{R}^{n\times m_C}$. 
Notice that $B_C u_C = B_{C1} u_{C1} + B_{C2} u_{C2}$, where $B_{C1}=(B_{C1},\textbf{0}_{n_2\times m_{C1}})$ and $B_{C2}=(\textbf{0}_{n_1\times m_{C2}},B_{C2})$. The jump map is arbitrary, as $D$ is empty.
Consider the stage costs $L_C(x,u_C):=x^\top Q_C x+ u_{C1}^\top R_{C1} u_{C1}+ u_{C2}^\top R_{C2} u_{C2}$ and $L_C(x,u_D):=0$, the terminal cost $q(x):=x^\top P x$ and the function $V(x):=x^\top P x$, where $Q_C \in \mathbb{S}^n_+$, $R_{C1} \in \mathbb{S}^{m_{C_1}}_+$, $-R_{C2} \in \mathbb{S}^{m_{C_2}}_+$  and  $P\in \mathbb{S}^n_+$. We refer to this game as a zero-sum infinite horizon linear quadratic differential game $(\mathcal{ZLQC})$. The function $V$ is continuously differentiable in $\reals^n$, and
\begin{eqnarray}
 &&\nonumber 
\underset{u_C=(u_{C1},u_{C2}) \in \Pi_u^C(x)}
{ \min_{u_{C1}} \max_{u_{C2}}}
\left\{L_C(x,u) +\left\langle \nabla V(x),F(x,u_C) \right\rangle \right\}
\\&=&
{  \min_{u_{C1} \in \reals^{m_{C_1}}} \max_{u_{C2} \in \reals^{m_{C_2}}}  }
\left\{x^\top Q_C x+ u_{C1}^\top R_{C1} u_{C1}+ u_{C2}^\top R_{C2} u_{C2}+2x^\top P (A_C x+B_C u_C ) \right\}
\nonumber
\\
&=&
{  \min_{u_{C1} \in \reals^{m_{C_1}}} \max_{u_{C2} \in \reals^{m_{C_2}}}  }
\left\{ x^\top (Q_C+2 P A_C) x   
+u_{C1}^\top R_{C1} u_{C1} +  2 x^\top  P B_{C1} u_{C1} 
\right. 
\nonumber \\  &&
\hspace{9cm} 
\left. + u_{C2}^\top R_{C2} u_{C2} + 2 x^\top  P  B_{C2} u_{C2} \right\}
\nonumber
\\  
&=& 0
\label{zsHJBLQGC}
\end{eqnarray}
is satisfied $ \forall x \in \reals^n$, and attained by $\kappa_C(x)=(-R_{C1}^{-1}B_{C1}^\top P x, -R_{C2}^{-1}B_{C2}^\top P x)$, if the continuous time $\HS_\infty$-like Riccati equation 
\begin{equation}
0
=Q_C+PA_C+A_C^\top P-P(  B_{C2} R_{C2}^{-1}B_{C2}^\top+B_{C1} R_{C1}^{-1}B_{C1}^\top) P
\label{zsRiccati}
\end{equation}
is satisfied. 
Given that $V$ is continuously differentiable, and that (\ref{Bellmanzsih}) and (\ref{HJBzsih}) are satisfied thanks to  (\ref{zsHJBLQGC}) and $D=\emptyset$, from Theorem \ref{thHJBszih}, the value function is
\begin{eqnarray}
\mathcal{J}^*(\xi) &:=& 
 \xi^\top P \xi,
\label{zsContLCost}
\end{eqnarray}
for any $\xi \in \reals^n$ 
%
{Given that this is a purely continuous system, the calculation of the cost for when applying the feedback law $\kappa_C(x)=-R_C^{-1}B_C^\top P x$, which renders a solution $(\phi,u_C^*)$ with $\dom \phi \ni (t,j) \mapsto u_C^*(t,j)= \kappa_C(\phi(t,j))$, is as follows
\begin{equation*}
\mathcal{J}(\xi,u_C^*) =
 \int_{t_0}^{\infty} L_C(\phi(t,j),u_{C}^*(t,j))dt  + \underset{(t,j) \in \textup{dom}\phi}{\limsup_{t+j\rightarrow \infty} } q(\phi(t,j))
\end{equation*}
and, when (CARE) is satisfied, one has
\begin{eqnarray*}
 \int_{t_0}^{\infty} L_C(\phi(t,j),u_{C}^*(t,j))dt  &=& 
 \int_{0}^{\infty} \phi^\top(t,j) Q_C \phi(t,j) +{u_{C}^*}^\top(t,j) R_C u_{C}^*(t,j)dt  
 \\&=&  \int_{0}^{\infty} \phi^\top(t,j) (Q_C +P B_C R_C^{-1}B_C ^\top P)\phi(t,j) dt  
  \\&=&  \int_{0}^{\infty} \phi^\top(t,j) (2P B_C R_C^{-1}B_C ^\top P-2PA_C)\phi(t,j) dt  
  \\&=&-\int_{0}^{\infty} \frac{d V(\phi(t,j))}{d t} dt
    \\&=&V(\phi(0,0))-\limsup_{t\rightarrow \infty}V(\phi(t,0))
  \\&=& 
  \phi^\top(0,0) P \phi(0,0)- \limsup_{t\rightarrow \infty} \phi^\top(t,0) P \phi(t,0)
\end{eqnarray*}
which, given that $\phi(0,0)=\xi$ and $q(\phi(t,j))=\phi^\top(t,0) P \phi(t,0)$, implies
$	\mathcal{J}(\xi,u_C^*) =\xi^\top P \xi$, showing that the feedback law $\kappa_C$ attains the optimal cost.} 
and the feedback law $\kappa_C$ is a pure strategy saddle-point equilibrium for the $(\mathcal{ZLQC})$ game.
\end{example}
\begin{example}{Linear quadratic difference games}
Consider a two-player zero-sum hybrid game with state $x \in \reals^n$, input $u=(u_C, u_D) \in \reals^{m_C}\times \reals^{m_D}$, and dynamics $\HS$ as in (\ref{Heq}) with \IfTp{}{$N=2$ and}
\begin{eqnarray}
D&=&\mathbb{R}^n \times \mathbb{R}^{m_D} \nonumber
\\
G(x,u_D)&=&A_D x+B_D u_D
\label{zsLinearDGame}
\end{eqnarray}
and flow set $C=\emptyset$, where $A_D:=\textup{blkdiag}_{i \in \IfTp{\{1,2\}}{\mathcal{V}}}\{A_i\}  \in \mathbb{R}^{n\times n}$ and $B_D:=\textup{blkdiag}_{i \in \IfTp{\{1,2\}}{\mathcal{V}}} \{B_{Di}\} \in \mathbb{R}^{n\times m_D}$. 
Notice that $B_D u_D = B_{D1} u_{D1} + B_{D2} u_{D2}$, where $B_{D1}=(B_{D1},\textbf{0}_{n_2\times m_{D1}})$ and $B_{D2}=(\textbf{0}_{n_1\times m_{D2}},B_{D2})$. The flow map is arbitrary, as $C$ is empty.
Consider the stage costs $L_C(x,u_C):=0$ and $L_D(x,u_D):=x^\top Q_D x+u_{D1}^\top R_{D1} u_{D1}+ u_{D2}^\top R_{D2} u_{D2}$, the terminal cost $q(x):=x^\top P x$ and the function $V(x):=x^\top P x$, where $Q_D\in \mathbb{S}^n_+$, $R_{D1}\in \mathbb{S}^{m_{D_1}}_+, -R_{D2}\in \mathbb{S}^{m_{D_2}}_+$ and  $P \in \mathbb{S}^n_+$. We refer to this game as an infinite horizon one player linear quadratic difference game $(\mathcal{ZLQD})$. The function $V$ is such that
\begin{eqnarray}
&&
\underset{u_D=(u_{D1},u_{D2}) \in \Pi_u^D(x)}
{\min_{u_{D1}} \max_{u_{D2}}}  
\left\{ 
 L_D(x,u_D)
+ V(G(x,u_D)) \right\}  \nonumber
\\&=&
{  \min_{u_{D1} \in \reals^{m_{D_1}}} \max_{u_{D2} \in \reals^{m_{D_2}}}  }
\left\{x^\top Q_D x+ u_{D1}^\top R_{D1} u_{D1}+ u_{D2}^\top R_{D2} u_{D2} \right. \nonumber \\
&&\hspace{8cm}\left. +(A_D x+B_D u_D )^\top P (A_D x+B_D u_D ) \right\}
\nonumber
\\&=&
{  \min_{u_{D1} \in \reals^{m_{D_1}}} \max_{u_{D2} \in \reals^{m_{D_2}}}  }
\left\{  x^\top (Q_D+A_D^\top P A_D) x 
+ u_{D1}^\top (R_{D1}+B_{D1}^\top P B_{D1}) u_{D1} 
\right. 
 \nonumber \\ && 
 \left. \hspace{2cm}
+ u_{D2}^\top (R_{D2}+B_{D2}^\top P B_{D2}) u_{D2} 
+2u_{D1}^\top B_{D1}^\top P B_{D2} u_{D2} 
+ 2 x^\top A_D^\top  P B_D u_D\right\} 
\nonumber
\\
&=&x^\top P x,
\label{zsHJBLQGD}
\end{eqnarray}
is satisfied $ \forall x \in \reals^n$, and attained by $\kappa_D(x)=-R^{-1}
(B_{D1}^\top P A_D, 
B_{D2}^\top P A_D)x$, with $R=\begin{bmatrix}
R_{D1}+B_{D1}^\top P B_{D1}
&
B_{D1}^\top P B_{D2}
 \\
B_{D2}^\top P B_{D1}
&
R_{D2}+B_{D2}^\top P B_{D2}
\end{bmatrix}$,
when the discrete-time $\HS_\infty$-like algebraic Riccati equation
\begin{multline}
P =Q_D+A_D^\top P A_D \IfConf{\\}{}
-
\begin{bmatrix}
A_D^\top P B_{D1}
 &
A_D^\top P B_{D2}
\end{bmatrix}
R 
^{-1}
\begin{bmatrix}
B_{D1}^\top P A_D
\\
B_{D2}^\top P A_D
\end{bmatrix}
\label{zsDARE} 
\end{multline}
is satisfied. 
Given that $V$ is continuously differentiable, and that (\ref{Bellmanzsih}) and (\ref{HJBzsih}) are satisfied thanks to $C=\emptyset$ and (\ref{zsHJBLQGC}), from Theorem \ref{thHJBszih}, the value function is
\begin{eqnarray}
\mathcal{J}^*(\xi) &:=& 
 \xi^\top P \xi,
\label{zsContLCost}
\end{eqnarray}
{Given that this is a purely discrete system, the calculation of the cost for when applying the feedback law $\kappa_D(x)=-(R_D+B_D^\top PB_D)^{-1}B_D^\top P A_D x$, which renders a solution $(\phi,u_D^*)$ with $\dom \phi \ni (t,j) \mapsto u_D^*(t,j)= \kappa_D(\phi(t,j))$, is as follows
\begin{equation*}
\mathcal{J}(\xi,u_D^*) =
  \sum_{j=0}^{\infty} L_D(\phi(t_{j+1},j),u_{D}^*(t_{j+1},j))+ \underset{(t,j) \in \textup{dom}\phi}{\limsup_{t+j\rightarrow \infty} } q(\phi(t,j))
\end{equation*}
and, when (DARE) is satisfied, one has
\begin{eqnarray*}
 \sum_{j=0}^{\infty} L_D(\phi(t_{j+1},j),u_{D}^*(t_{j+1},j)) &=& 
 \sum_{j=0}^{\infty} \phi^\top(t,j)Q_D \phi(t,j) +{u_{D}^*}^\top(t,j) R_D u_{D}^*(t,j)  
 \\&=&  \sum_{j=0}^{\infty} \phi^\top(t,j) (Q_D +A_D^\top P B_D (R_D+B_D^\top P B_D)^{-1}B_D^\top P A_D)\phi(t,j)   
  \\&=&  \sum_{j=0}^{\infty} \phi^\top(t,j) (2A_D^\top P B_D (R_D+B_D^\top P B_D)^{-1}B_D^\top P A_D+P-A_D^\top PA_D)\phi(t,j)  
  \\&=&\sum_{j=0}^{\infty} V(\phi(t,j))-V(G(\phi(t,j),u_D^*(t,j))
  \\&=&V(\phi(0,0))-\limsup_{j\rightarrow \infty}V(\phi(0,j))
  \\&=& 
  \phi^\top(0,0) P \phi(0,0)- \limsup_{j\rightarrow \infty} \phi^\top(0,j) P \phi(0,j)
\end{eqnarray*}
which, given that $\phi(0,0)=\xi$ and $q(\phi(t,j))=\phi^\top(0,j) P \phi(0,j)$, implies
$	\mathcal{J}(\xi,u_D^*) =\xi^\top P \xi$, showing that the feedback law $\kappa_D$ attains the optimal cost.} 
and the feedback law $\kappa_D$ is a pure strategy saddle-point equilibrium for the $(\mathcal{ZLQD})$ game.
\end{example}}
\NotAutom{
\pn{
\section{Application: Variable terminal-time games}
In this section, we address solutions to a class of games with hybrid dynamics. Consider the following pursuit-evasion case that motivates the application of game theoretical tools with variable terminal time for hybrid systems.

\subsection{Continuous Pursuit-Evasion Game with Periodic State Jumps}

A slower-defender target-defense problem is formulated with a defender $P_1$ that moves in the plane with velocity $V_D$ and an attacker $P_2$ that also moves in the plane with velocity $V_A$ such that $\frac{V_A}{V_D}:=v>1$. The goal of the attacker is to reach the origin before it is intercepted by the defender while the defender aims to stop him. The dynamics of the players are given by
\begin{equation}
\begin{split}
\dot{x}_D=V_D\cos u_{CD},
\\ 
\dot{y}_D=V_D\sin u_{CD},
\\ 
\dot{x}_A=V_A\cos u_{CA},
\\ 
\dot{y}_A=V_A\sin u_{CA},
\end{split}
\label{PuEvC}
\end{equation}
Consider the case in which the defender $P_1$ is allowed to jump but only at certain instants of time with the following dynamics
\begin{equation}
  \begin{split}
x_D^+=&x_D+r\cos u_{DA} \\
y_D^+=&y_D+r\sin u_{D2}
  \end{split}
\end{equation}
for a given $r>0$.
%
Inspired by \cite{123,150} we model 
the ability to jump 
by defining a timer state $\tau \in [0,\bar T]$. 
This state decreases with ordinary time and every time it reaches $0$, it gets reset to  $\bar T$. The dynamics of the timer state are given by
\begin{equation}
\begin{split}
&\dot{\tau}=-1, \hspace{1.5cm} \tau \in [0,\bar T]\\ 
&\tau^+=\bar T, \hspace{1.5cm} \tau=0. 
\end{split}
\label{timer}
\end{equation}

Thus, this game has dynamics $\HS$ as in (\ref{Heq}) with state $x=[x_D,y_D,x_A,y_A,\tau]\in \reals^4 \times [0,\bar T]$, input $u=(u_1,u_2)=((u_{CD},u_{DD}),u_{CA})$ and described by
\begin{equation}
\begin{split}
&C=\{(x,u) \in\reals^4\times  [0,\bar T]: \tau>0
\}   \\ 
&F(x,u)=\begin{bmatrix}V_D\cos u_{CD}
  \\ 
  V_D\sin u_{CD}
  \\ 
  V_A\cos u_{CA}
  \\ 
  V_A\sin u_{CA}\\-1 \end{bmatrix}\\
&D=\{(x,u) \in \reals^4\times  [0,T_2]
: \tau=0
\}\\
&G(x,u)=\begin{bmatrix}x_D+r\cos u_{DD}\\y_D+r \sin u_{DD} \\ x_A\\y_A\\ \bar T \end{bmatrix}
\end{split}
\end{equation}
The game is over when the attacker reaches the origin or when defender intercepts the attacker, i.e., $x_D=x_A, y_D=y_A$. In this work, we are interested in studying the latter case and the terminal condition is defined in terms of the capture set $X:=\{x\in \reals^4 \times [0,T_2]:x_D=x_A, y_D=y_A\}$.

We encode the objective of each of the players by setting $L_C(x,u):=0$, $L_D(x,u):=0$, and $q(x)=x_A^2+y_A^2$. Thus, the defender seeks to maximize the capture distance to the origin, while the attacker wants to minimize it. 
}
}

\color{black}
\begin{example}{Bouncing ball with terminal set} 
Inspired by the problem in \cite{79c}, consider a simplified model of a juggling system as in \cite{69c},  
  with state $x=(\sj{x_p,x_v}) \in \reals^2$, input $ u_D:=(u_{D1},u_{D2}) \in \reals^2$, and dynamics $\HS$ as in (\ref{Heq}), with data
  \pno{\begin{eqnarray}
  \nonumber
  C&=&\realsgeq \times \reals, 
  \quad \>\> F(x)=(\sj{x_v},-1) \quad \forall x \in C
  \\
  D&=&\{0\} \times \reals_{\leq 0} \times \reals^2,  
  \\G(x,u_D)&=&(0, -\lambda  \sj{x_v}+u_{D1}+u_{D2})
  \quad \forall (x, u_D) \in D
  \nonumber
  \label{BouncingBallDynamics}
  \end{eqnarray}}
  where $u_{D1}$ is the control input, $u_{D2}$ is the action of an attacker, and $\lambda\in (0,1)$ is the coefficient of restitution of the ball. 
 The scenario in which $u_{D1}$ is designed to minimize a cost functional $\J$ until the game ends, which occurs when the state enters a set $X$, under the presence of the worst-case disturbance $u_{D2}$ {is formulated as a two-player zero-sum game.}
  With the aim of pursuing minimum velocity and control effort at jumps, consider the cost functions $L_C(x,u_C):=0$, $L_D(x,u_D):=\sj{x_v}^2Q_D+u_D^\top R_D u_D$, and terminal cost $q(x):=\frac{1}{2} \sj{x_v}^2+\sj{x_p}$ defining $\J$ as in (\ref{defJTNCinc}), with $R_D :=\left[ \begin{smallmatrix}R_{D1} & 0 \\ 0 & R_{D2}\end{smallmatrix}\right]$ and $Q_D,$ $R_{D1},$ $-R_{D2}>0$.  Here, $u_{D1}$ is designed by player $P_1$, which aims to minimize $\mathcal{J}$, while player $P_2$ seeks to maximize it by choosing $u_{D2}$.
  \pno{A game of kind \cite[Section 5.2]{basar1999dynamic}} arises and its solution characterizes a division of the state space into two dominance regions, $\mathcal{M},\Psi \subset \Pi(C) \cup \Pi(D)$, in which, under optimal play, it can be determined whether the terminal set $X$ is reached or not as a function of the initial condition. 
  If the initial state satisfies $\xi \in \mathcal{M}$ \pno{(the feasible set)}, then, under optimal play, the ball reaches the terminal set $X$ at some time $(T,J)$ and the game ends. On the other hand, if $\xi \in \Psi$, under optimal play, we have an infinite horizon game \pn{(if maximal solutions are complete \NotAutomss{after the inputs are assigned}). 
  }
  \IfPers{
  Given that $P=[\begin{smallmatrix}
    1 && 0 \\ 0 &&1
  \end{smallmatrix}]$ satisfies the conditions in Corollary \ref{HREqSec} when 
  \begin{equation}
  Q_D=\frac{- 2R_{D1} R_{D2} \lambda^2 + R_{D1} + R_{D2} + 2 R_{D1}R_{D2}}{2 R_{D1} + 2 R_{D2} + 4 R_{D1} R_{D2}},
  \label{QdBB}
  \end{equation} the feedback law  
    $\gamma_{D1}(z)=\frac{R_{D2}\lambda }{R_{D1} + R_{D2} + 2 R_{D1} R_{D2}}z_2$ 
  minimizes the cost functional $\mathcal{J}$ in the presence of the worst-case disturbance given by  $\gamma_{D2}(z)=\frac{R_{D1}\lambda }{R_{D1} + R_{D2} + 2 R_{D1} R_{D2}}z_2$.
  Thus,} 

  The function $V(x):=\sj{x_p}+\frac{1}{2}\sj{x_v}^2$  is such that $\left\langle \nabla V(x),F(x) \right\rangle=0$
for all $x \in C$, making $V$ a solution to (\ref{HJBzsih}). In addition, the function $V$ is such that
%
\begin{equation}
  \begin{split}
 &\underset{u_D=(u_{D1},u_{D2}) \in \reals^2}
{  \underset{u_{D1}}{\min}\>  \underset{u_{D2}}{\max}  }  
\IfAutom{ \mathcal{L}_D(x,u_D)}
{
\left\{
 L_D(x,u_D)
+ V(G(x,u_D)) \right\} }
\NotAutom{
\\
&  
= \underset{u_{D1}\in \mathbb{R}}{\min} \> \underset{u_{D2}\in \mathbb{R}}{\max}  
\Big\{ \sj{x_v}^2Q_D+u_D^\top R_D u_D
\\ & \hspace{2cm} +\frac{(-\lambda \sj{x_v} + u_{D1}+u_{D2})^2}{2} \Big\} 
}
= \frac{1}{2} \sj{x_v}^2
\end{split}
\label{IBB}
\end{equation}
for all $(x,u_D)\in D$. Equality \eqref{IBB} is attained by $\kappa_D(x)=(
\kappa_{D1}(x),\kappa_{D2}(x))$ with 
$\kappa_{D1}(x)=\frac{R_{D2}\lambda }{R_{D1} + R_{D2} + 2 R_{D1} R_{D2}}\sj{x_v}$ and $\kappa_{D2}(x)=\frac{R_{D1}\lambda }{R_{D1} + R_{D2} + 2 R_{D1} R_{D2}}\sj{x_v}$ when 
\begin{equation}
Q_D=\frac{- 2R_{D1} R_{D2} \lambda^2 + R_{D1} + R_{D2} + 2 R_{D1}R_{D2}}{2 R_{D1} + 2 R_{D2} + 4 R_{D1} R_{D2}},
\label{QdBB}
\end{equation}
which 
makes $V$ a solution to (\ref{Bellmanzsih}) \textcolor{black}{with saddle-point equilibrium $\kappa_D$}.
Thus, given that $V$ is continuously differentiable on $\reals^2$, and that (\ref{HJBzsih}) and (\ref{Bellmanzsih}) hold thanks to (\ref{IBB}) and (\ref{QdBB}), from Theorem \ref{thHJBszih}\pno{,} the value function is
$
\mathcal{J}^*(\xi) =
  \frac{\sj{\xi_v}^2}{2}+\sj{\xi_p}.
\label{HyLCostBB}
$
  Figure \ref{fig:bb1} displays this behavior \pn{with $\xi \in \mathcal{M}$ and both players playing the saddle point equilibrium. The terminal set $X$ is reached at $t=8$s and the cost of the displayed solution is $V(\xi)$.} 
  \begin{figure}[h]
      \hspace{-0.2cm}
      \includegraphics[width = 8.4cm]{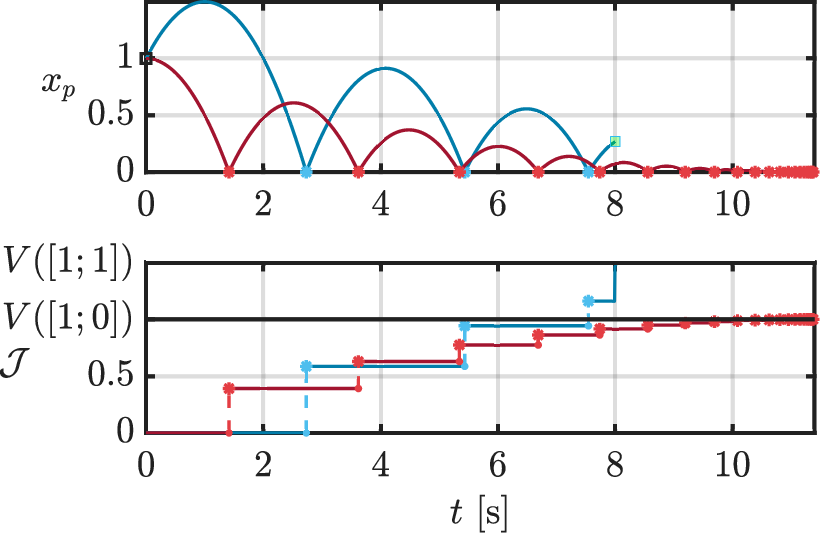}
      \vspace{-0.3cm}
      \caption{Bouncing ball solutions attaining minimum cost under worst-case $u_{2}$, with 
      $X = \{x \in \reals^2: 0 \leq \sj{x_p} \leq 0.3, -0.37 \leq \sj{x_v} \leq 0.37\}$, 
      $\lambda=0.8$, $R_{D1}=10, R_{D2}=-20$, and $Q_D=0.189$. 
      \sj{Solution in reachable set (blue). Complete solution (red). Value function (black). Initial conditions (squares). Terminal set (green).}
      }
      \label{fig:bb1}
  \end{figure}
  \begin{figure}[h]
    \vspace{-0.3cm}
    \hspace{-0.2cm}
    \includegraphics[width = 8.4cm]{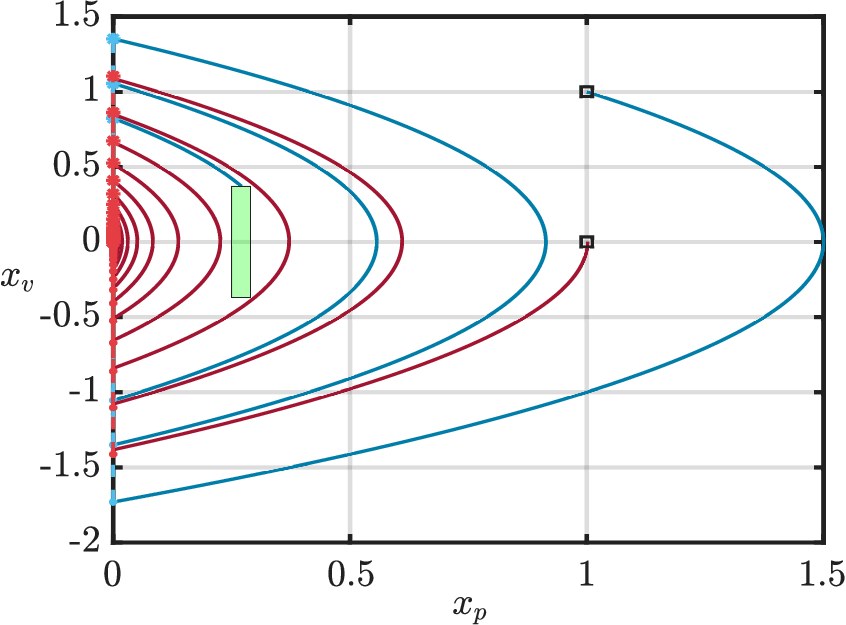}
    \vspace{-0.3cm}
    \caption{Bouncing ball phase portrait. Terminal set (green) and initial condition (square).
    }
    \label{fig:bbpp}
\end{figure}
\NotAutomss{
  Figure \ref{fig:bbppc} displays a solution that does not enter the terminal set, the cost associated to it over time and the value function. \pn{Notice that the cost of such solution from  $\xi \notin \mathcal{M}$ under both players playing the saddle point equilibrium, is equal to $V(\xi)$.}
\begin{figure}[h]
  \vspace{-0.4cm}
  \hspace{-0.2cm}
  \includegraphics[width = 9cm]{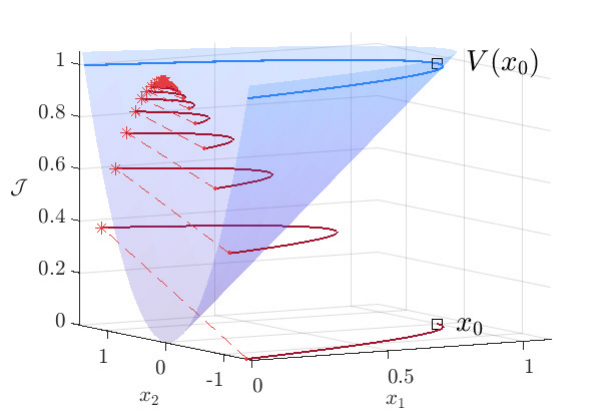}
  \vspace{-0.6cm}
  \caption{Bouncing ball cost. Initial condition (square). Value function and saddle-point equilibrium trajectory attaining evaluated cost at initial condition. 
  }
  \label{fig:bbppc}
\end{figure}
}
  
 \NotAutom{ In Figure \ref{CostsBB}, we let the players select feedback laws close with the Nash equilibrium and calculate the cost associated to the new laws. The variation of the cost along the changes in the feedback laws makes evident the saddle-point geometry. 
  
  \begin{figure}
      \centering
      \includegraphics[width=7.5cm]{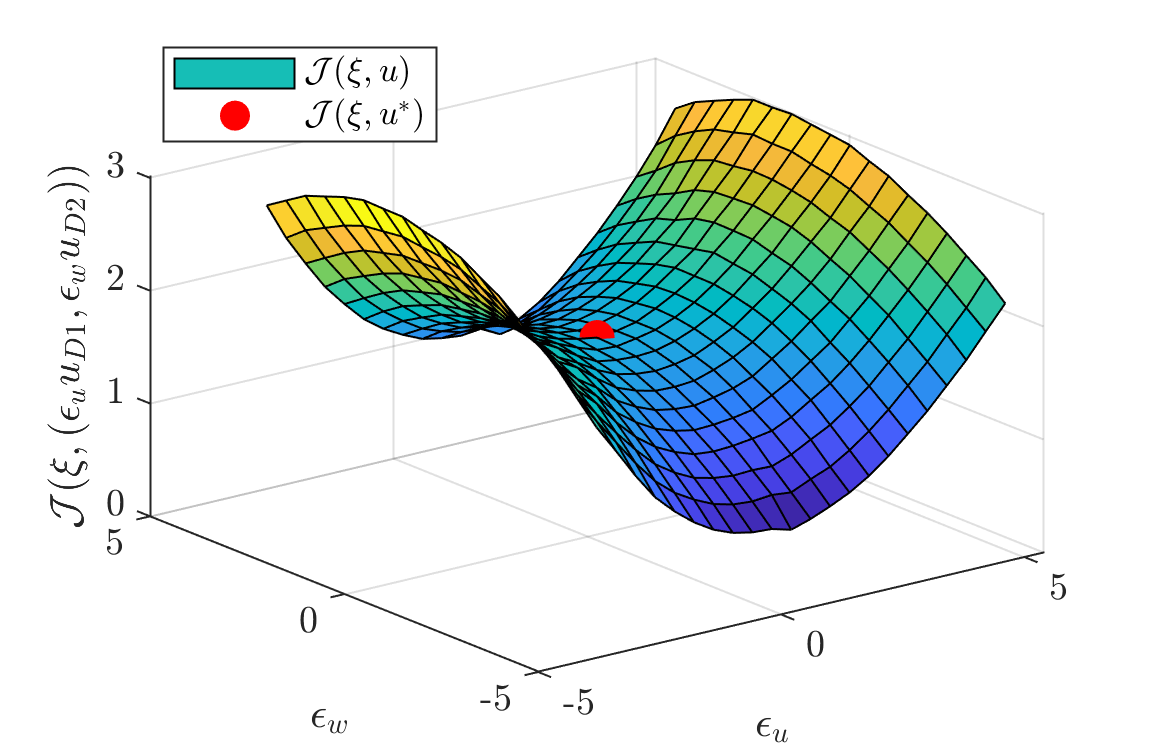}
      \vspace{-0.8cm}
      \caption{Saddle point behavior in the cost of solutions to bouncing ball from $\xi=(1,1)$ when the feedback gains vary around the optimal value. The cost is evaluated on solutions $(\phi,u)\in \mathcal{S}^\mathcal{T}_\HS (\xi) $ with feedback law variations specified by  $\epsilon_u$ and $\epsilon_u$ in $u=(\epsilon_u \kappa_1(t,j,\phi),\epsilon_w \kappa_2(t,j,\phi))$.}
      \label{CostsBB}
  \end{figure}}
%
\pno{As an alternative version of this game, consider the case in which $X =\emptyset$ and let $\A=\{0\}$, encoding the goal of stabilizing the ball to {rest} under the effect of an attacker. This implies that $\mathcal{M} = \emptyset$ and $\Psi = \Pi(C) \cup \Pi(D)$.} 
This is formulated as a two-player zero-sum \pno{infinite horizon hybrid game via solving Problem ($\diamond$) over the set of complete input actions}. 
The function $V$ is a solution to (\ref{Bellmanzsih}) \textcolor{black}{with saddle-point equilibrium $\kappa_D$}.
Similarly, given that $V$ is continuously differentiable on $\reals^2$, and that (\ref{HJBzsih}) and (\ref{Bellmanzsih}) hold thanks to (\ref{IBB}) and (\ref{QdBB}), from Theorem \ref{thHJBszih}, the value function is
$
\mathcal{J}^*(\xi) =
  \frac{\sj{\xi_v}^2}{2}+\sj{\xi_p}.
\label{HyLCostBB}
$
Figure \ref{fig:bb1} displays this behavior. 

{Furthermore, given that 
$L_D \in \mathcal{PD}_{\kappa_D}(\A 
)$, and (\ref{kHJBeqihza})-(\ref{kBeqihzf}) hold, by setting $\alpha_1(s)=\min \left\{\frac{1}{2} \left(\frac{s}{ \sqrt 2}\right)^2, \frac{s}{ \sqrt 2}\right\}$ and $\alpha_2 (s)=\frac{1}{2}s^2+s$, 
from Corollary \ref{zCorStability}, we have that $\kappa_D$ is the saddle-point equilibrium and renders $\A=\{0\}$ uniformly globally asymptotically stable for $\HS$}.  

\NotAutomss{In Figure \ref{CostsBB}, we let the players select feedback laws close to the Nash equilibrium and calculate the cost associated to the new laws. The variation of the cost along the changes in the feedback laws makes evident the saddle-point geometry. This example illustrates how our results apply to Zeno systems.
\begin{figure}[h]
    \centering
    \includegraphics[width=9cm]{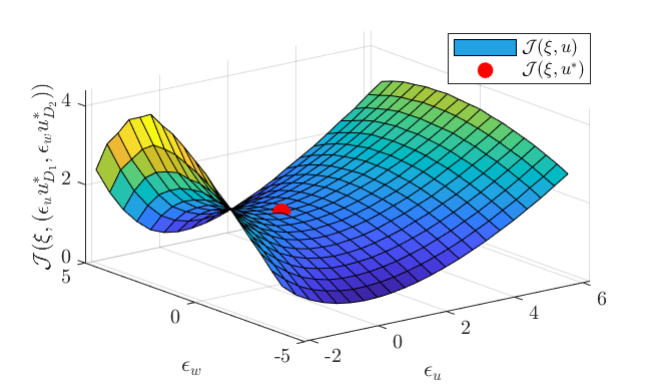}
    \caption{Saddle point behavior in the cost of solutions to bouncing ball from $\xi=(1,1)$ when varying the feedback gains around the optimal value. The cost is evaluated on solutions $(\phi,u)\in \mathcal{S}^\infty_\HS (\xi) $ with feedback law variations specified by  $\epsilon_u$ and $\epsilon_w$ in $u=(\epsilon_u \kappa_1(\phi),\epsilon_w \kappa_2(\phi))$.}
    \label{CostsBB}
\end{figure}}
\label{BouncingBallEx}
\vspace{-0.3cm}
\end{example}
{A special \NotAutomss{case of a }hybrid game emerges under a \emph{capture-the-flag} setting, in which two teams compete to grab the opponents flag and return it to their \NotAutomss{own }base without getting tagged.} This is formulated as in Problem $(\diamond)$ in \cite{leudohycflag}. 
\IfPers{
\begin{example}{Linear quadratic differential games}
Taking up the $(\mathcal{ZLQC})$ game, recall that $\kappa_C:=(-R_{C1}^{-1}B_{C1}^\top P x,-R_{C2}^{-1}B_{C2}^\top P x)$ and consider the case in which $A_C-B_{C1}R_{C1}^{-1}B_{C1}^\top P-B_{C2}R_{C2}^{-1}B_{C2}^\top P$ is Hurwitz and $\A=\{0\}$. Given that $L_C \in \mathcal{PD}_{\kappa_C}(\A)$, $V$ is continuously differentiable, and (\ref{kHJBeqihza}), (\ref{kHJBeqihzb}), and (\ref{kHJBeqihzc}) hold, by setting $\alpha_1(|x|_\A)=x^\top (P-I) x$ and $\alpha_1(|x|_\A)=x^\top (P+I) x$, from Corollary \ref{zCorStability}, we have that $\kappa:=(\kappa_C, \cdot )$ is the pure strategy saddle-point equilibrium and renders $\A$ uniformly globally asymptotically stable for $\HS_{\kappa}$.
\end{example}

\begin{example}{Linear quadratic difference games}
Taking up the $(\mathcal{ZLQD})$ game, recall that $\kappa_D:=-
(R^{-1}B_{D1}^\top P A_D x, 
R^{-1}B_{D2}^\top P A_D x)$ and consider the case in which $A_D-[B_{D1},B_{D2}] R^{-1}(B_{D1}^\top P(0) A_D,
B_{D2}^\top P(0) A_D)$ is Schur and $\A=\{0\}$. Given that $L_D \in \mathcal{PD}_{\kappa_D}(\A)$, and  (\ref{kBeqihzd}), (\ref{kBeqihze}), and (\ref{kBeqihzf}) hold, by setting $\alpha_1(|x|_\A)=x^\top (P-I) x$ and $\alpha_1(|x|_\A)=x^\top (P+I) x$, from Corollary \ref{zCorStability}, we have that $\kappa:=(\cdot, \kappa_D)$ is the pure strategy saddle-point equilibrium and renders $\A$ uniformly globally asymptotically stable for $\HS_{\kappa}$. 
\end{example}}

\sj{\section{Further Connections with the Literature}}
Some results provided in this {paper} have direct counterparts in the continuous-time and discrete-time game {theory} literature. The definition of a game in terms of its elements can be\NotAutomss{\>directly }traced back to \cite{basar1999dynamic}, {as explained below.}

  \pno{ Given a discrete-time \IfTp{two-player}{} zero-sum game with final time{\footnote{This corresponds to the hybrid time $(0,J)$ for $\HS\IfIncd{_s}{}$.}} $``J$'', $f_k$ and $X$ defining the {single-valued} jump map and jump set, respectively, as in \cite{basar1999dynamic}, setting the data of $\HS\IfIncd{_s}{}$ as $C = \emptyset$, $G=f_k$ for $k \in \nats_{\leq J}$, and $D=X$ reduces Definition \ref{elements} to  \cite[Def. 5.1]{basar1999dynamic} for the case in which
  the output of each player is equal to its state and there is a feedback information structure as in \cite[Def. 5.2]{basar1999dynamic}. 
  Thus, items $(vi)-(ix)$ in \cite[Def. 5.1]{basar1999dynamic} are omitted in the formulation  herein and items $(i)-(v)$ and $(x)-(xi)$ are covered by Definition \ref{elements}, the definition of the hybrid time domain with final time $(0,J)$, and the set $\mathcal{S}_{\HS\IfIncd{_s}{}}$.
  \\
  \indent Given a continuous-time \IfTp{two-player}{} zero-sum game with final time{\footnote{This corresponds to the hybrid time $(T,0)$ for $\HS\IfIncd{_s}{}$.}} $``T$'', $f$ and $\mathcal{S}^0$ defining the {single-valued} flow map and flow set, respectively, as in \cite{basar1999dynamic}, setting the data of $\HS\IfIncd{_s}{}$ as $D = \emptyset$, $F=f$, and $C=\mathcal{S}^0$ reduces Definition \ref{elements} to \cite[Def. 5.5]{basar1999dynamic} for the case in which
  the output of each player is equal to its state and there is a feedback information structure as in \cite[Def. 5.6]{basar1999dynamic}. 
  Thus, items $(vi)-(vii)$ in \cite[Def. 5.5]{basar1999dynamic} are omitted in the formulation herein and items $(i)-(v)$ and $(viii),(ix)$ are covered by Definition \ref{elements}, the definition of the hybrid time domain with final time $(0,T)$, and the set $\mathcal{S}_{\HS\IfIncd{_s}{}}$.
}
\pno{
\begin{remark}{Equivalent costs}
  Given $\xi \in \Pi(\overline{C}) \cup \Pi(D)$ and a strategy \IfIh{$\kappa^*=(\kappa_1^*\IfTp{\kappa_2^*}{, \dots, \kappa_N^*})\in \mathcal{K}$,}{$\gamma^*=(\gamma_1^*, \dots, \gamma_N^*)\in \mathcal{K}$,} denote by \IfIh{$\mathcal{U}^*(\xi, \kappa^*)$}{$\mathcal{U}^*(\xi, \gamma^*)$} the set of joint actions  $u=(u_1, \IfTp{u_2}{\dots, u_N
  })$ rendering a maximal trajectory $\phi$ to $\HS\IfIncd{_s}{}$ from $\xi$ with components defined as  $\dom \phi \ni (t,j) \mapsto u_i(t,j)= \IfIh{\kappa}{\gamma}_i^*(\IfIh{}{t,j,}\phi(t,j))$ for each $i\in \IfTp{\{1,2\}}{\mathcal{V}}$. 
  By expressing the largest cost associated to the solutions to $\HS\IfIncd{_s}{}$ from $\xi$ under the strategy $\IfIh{\kappa}{\gamma}^*$ 
   as $\hat{\mathcal{J}}(\xi,\IfIh{\kappa}{\gamma}^*):=\sup_{u \in \mathcal{U}^*(\xi, \IfIh{\kappa}{\gamma}^*)} \mathcal{J}(\xi, u)$, 
    an equivalent condition to 
   (\ref{SaddlePointIneq}) for when
  {$\mathcal{J}(\xi, u)=\hat{\mathcal{J}}(\xi, \IfIh{\kappa}{\gamma}^*)$ for every $u \in \mathcal{U}^*(\xi, \IfIh{\kappa}{\gamma}^*)$ }
  is
  \IfPers{\begin{eqnarray*}
  \hat{\mathcal{J}}_1(\xi,\IfIh{\kappa}{\gamma}^*) 
  &\leq&  \hat{\mathcal{J}}_1(\xi,\IfIh{(\kappa_{1}, \kappa_{2}^*\IfTp{}{, \dots, \kappa_{N}^*})}{(\gamma_{1}, \gamma_{2}^*, \dots, \gamma_{N}^*)})
   \\
  \hat{\mathcal{J}}_2(\xi,\IfIh{\kappa}{\gamma}^*) 
  &\leq& \hat{\mathcal{J}}_2(\xi,\IfIh{(\kappa_{1}^*, \kappa_{2}\IfTp{}{, \dots, \kappa_{N}^*})}{(\gamma_{1}^*, \gamma_{2}, \dots, \gamma_{N}^*)})
  \IfTp{}{, \\
  &\vdots&
  \\
  \hat{\mathcal{J}}(\xi,\IfIh{\kappa}{\gamma}^*) 
  &\leq&  \hat{\mathcal{J}}(\xi,\IfIh{(\kappa_{ 1}^*, \kappa_{2}^*, \dots, \kappa_{N})}{(\gamma_{ 1}^*, \gamma_{2}^*, \dots, \gamma_{N})})}.
  \end{eqnarray*}}
  {$\hat{\mathcal{J}}(\xi,\IfIh{(\kappa_{1}^*, \kappa_{2}\IfTp{}{, \dots, \kappa_{N}^*})}{(\gamma_{1}^*, \gamma_{2}, \dots, \gamma_{N}^*)})\leq\hat{\mathcal{J}}(\xi,\IfIh{\kappa}{\gamma}^*) 
  \leq  \hat{\mathcal{J}}_1(\xi,\IfIh{(\kappa_{1}, \kappa_{2}^*\IfTp{}{, \dots, \kappa_{N}^*})}{(\gamma_{1}, \gamma_{2}^*, \dots, \gamma_{N}^*)})
  $}
   for all $\IfIh{\kappa}{\gamma}_i \in \mathcal{K}_i$, $i \in \IfTp{\{1,2\}}{\mathcal{V}}$.
  \label{EqCosts} 
  \end{remark}
  }
  \NotAutomss{
  \begin{remark}{Relation of definition of solution to literature} 
  By considering a discrete-time system with the single-valued function $G$ or by considering a continuous-time system with $F$ \IfIncd{single valued and}{} Lipschitz continuous in $\overline{C}$, and by removing the initial condition as an argument of the cost functionals and specifying it in the state equation, Remark \ref{EqCosts} presents equivalent conditions to those in \cite[(6.3)]{basar1999dynamic}. 
  Thus, Definition \ref{NashEqNonCoop} covers the definitions of a pure strategy Nash equilibrium in \cite[Sec. 6.2, 6.5]{basar1999dynamic}  for the zero-sum case.
  \end{remark}
  }
{Conditions for computing value functions for linear quadratic problems have been widely studied, concerning solving differential and algebraic Riccati equations. The computation of value functions for systems with nonlinear dynamics is an open research problem and has seen interesting learning-based contributions in the last years, e.g., \cite{MontegroLeudoDataDrivenAS}.
The computation of value functions for DAEs is discussed in \cite{GARDNER1978355}, \cite{voigt2015linear}, for the case of linear differential games under algebraic constraints. Such value functions have a similar structure to the ones provided herein for hybrid systems with linear jump and flow maps and algebraic constraints encoded by the flow set $C$. }

\pno{The design of value functions for switched DAEs imposes additional challenges that follow the discussion in \cite{LIBERZON2012954} on the existence of Lyapunov functions and asymptotic stability. In some cases, a common Lyapunov function for all the subsystems of a switched DAE does not exist and even when it exists, it is not enough to guarantee asymptotic stability due to arbitrary switching. To solve this, conditions over switching are provided in \cite[Theorem 4.1]{LIBERZON2012954}, and for the optimality of hybrid systems, such conditions are resembled by the point-wise conditions on the change of V along jumps. 
In \cite{tanwani2019feedback}, there are coupled value functions associated to each subsystem of a switched DAE in a zero-sum game, which result in coupled Riccati differential equations with optimal feedback strategies described by linear-time-varying functions of the state.
Note that both scenarios are accounted for in the design of a value function for hybrid games based on optimality pointwise conditions provided in this work.
}
\section{Conclusion and Future Work}
%
%
%

In this paper, {we formulate a two-player zero-sum game under dynamic constraints given in terms of hybrid dynamical systems, as in \cite{65}. }
Scenarios in which the control action is selected by a player $P_1$ to accomplish an objective and counteract the damage caused by an adversarial player $P_2$ are studied. By encoding the objectives of the players in the optimization of a cost functional, sufficient conditions in Hamilton–Jacobi–Bellman-Isaacs form are provided to upper bound the cost for any disturbance. The main result allows the {optimal strategy of $P_1$} to minimize the cost under the maximizing {adversarial action}. Additional conditions are proposed to allow the saddle-point strategy to render a set of interest asymptotically stable by \IfAutomss{using}{letting\>} the value function \IfAutomss{as}{take the role of\>} a Lyapunov function. 

Future work includes \sj{generalizing results to \NotAutomss{the space of }mixed strategies. 
Structural} conditions on the system that do not involve $V$ and guarantee the existence of a solution to Problem $(\diamond)$ based on the smoothness and regularity of the data of the system, \IfAutomss{as }{similar to those }in \cite{GOEBEL2019153}, will be studied.

\bibliographystyle{plain}
\bibliography{GT,RGSOnline}

\NotAutomss{
\appendix
\section{Appendix}
\NotAutom{
  \par\noindent\textbf{Proof of Proposition \ref{Pp:ProofDerivation}.}
Given a 
$(\phi,u)\in \mathcal{S}^\infty_\HS (\xi)$, where $\{t_j\}_{j=0}^{\sup_j \dom \phi}$ is a nondecreasing sequence associated with the hybrid time domain of $(\phi,u)$ as in \IfAutomss{\cite[Definition 2.3]{65}}{Definition \ref{htd}},  
{
for each $j \in \nats$ such that $I_{\phi}^j=[t_j, t_{j+1}]$ has a nonempty interior int$I_{\phi}^j$, 
}
by integrating \eqref{eq:upperboundhamiltonian} over 
$I_{\phi}^j$, 
we obtain
\begin{equation*}
\begin{split}
0\geq  \int_{t_{j}}^{t_{j+1}} \left( L_C(\phi(t,j),u_{ C}(t,j))
+ \frac{d}{d t}V(\phi(t,j)) \right )dt
\end{split}
\end{equation*}
from where we have
\begin{equation*}
\begin{split}
0\geq  \int_{t_{j}}^{t_{j+1}} L_C(\phi(t,j),u_{ C}(t,j)) dt
\hspace{3cm}
\\
+  V(\phi(t_{j+1},j))- V(\phi(t_j,j))
\end{split}
\end{equation*}
{Pick $(t^*,j^*) \in \dom (\phi,u)$.} Summing from $j=0$ to $j={j^*}$ 
we obtain 
\begin{equation*}
\begin{split}
0
\geq\sum_{j=0}^{{j^*}}  \int_{t_{j}}^{t_{j+1}} L_C(\phi(t,j),u_{ C}(t,j)) dt
\hspace{3cm}
\\+\sum_{j=0}^{{j^*}} \left( V(\phi(t_{j+1},j))-V(\phi(t_j,j)) \right)
\end{split}
\end{equation*}
Then, solving for $V$ at the initial condition $\phi(0,0)$, we obtain
\begin{equation} \label{ContinCost2goihzol}
\begin{split} 
&  \hspace{-0.5cm}V(\phi(0,0))
\geq\sum_{j=0}^{{j^*}}  \int_{t_{j}}^{t_{j+1}} L_C(\phi(t,j),u_{ C}(t,j)) dt
\\  +& V(\phi(t_1,0))
+\sum_{j=1}^{{j^*}} \left( V(\phi(t_{j+1},j))-V(\phi(t_j,j)) \right) 
\end{split}
\end{equation}
%
In addition, {if $j^*>0$,} 
adding \eqref{upperbounddiscrete} from $j=0$ to $j={j^*}-1$, 
we obtain
\begin{equation}
\begin{split} 
\sum_{j=0}^{{j^*}-1}V(\phi(t_{j+1},j)) 
\geq 
\sum_{j=0}^{{j^*}-1}L_D(\phi(t_{j+1},j),u_{ D}(t_{j+1},j))\nonumber
\\
+\sum_{j=0}^{{j^*}-1} V(\phi(t_{j+1},j+1))
\nonumber
\end{split}
\end{equation}
Then, solving for $V$ at the first jump time, we obtain
\begin{equation}
\begin{split} 
  & V(\phi(t_{1},0))
\geq V( \phi(t_{1},1))
\hspace{3cm}
\\
& \hspace{0.5cm}
+ \sum_{j=0}^{{j^*}-1}L_D(\phi(t_{j+1},j),u_{ D}(t_{j+1},j)) 
\\
&\hspace{0.5cm}+\sum_{j=1}^{{j^*}-1}\left( V(\phi(t_{j+1},j+1))-V(\phi(t_{j+1},j)) \right) 
\end{split}
\label{DiscreteCost2goNOihzwincol}
\end{equation}
In addition, given that $\phi(0,0)=\xi$, lower bounding $V(\phi(t_{1},0))$ in (\ref{ContinCost2goihzol}) by the right-hand side of (\ref{DiscreteCost2goNOihzwincol}), we obtain 
\begin{equation*}
\begin{split} 
V(\xi)
\geq& 
  \sum_{j=0}^{{j^*}} \int_{t_{j}}^{t_{j+1}} L_C(\phi(t,j),u_{C}(t,j))dt+ V(\phi(t_1,0)) 
\\&
+\sum_{j=1}^{{j^*}} \left( V(\phi(t_{j+1},j))-V(\phi(t_j,j)) \right)  
\\
\geq&
\sum_{j=0}^{{j^*}} \int_{t_{j}}^{t_{j+1}} L_C(\phi(t,j),u_{C}(t,j))dt
\\&
+ \sum_{j=0}^{{j^*}-1}L_D(\phi(t_{j+1},j),u_{D}(t_{j+1},j))
 \\
&+\sum_{j=1}^{{j^*}-1}\left( V(\phi(t_{j+1},j+1))-V(\phi(t_{j+1},j)) \right) 
\\&
+ V(\phi(t_{1},1))+\sum_{j=1}^{{j^*}} \left( V(\phi(t_{j+1},j))-V(\phi(t_j,j)) \right) 
\end{split}
\end{equation*}
Since 
\begin{eqnarray*} \nonumber
&& V(\phi(t_{1},1))
+\sum_{j=1}^{{j^*}-1}\left( V(\phi(t_{j+1},j+1))-V(\phi(t_{j+1},j)) \right) \nonumber 
\\&&
+\sum_{j=1}^{{j^*}} \left( V(\phi(t_{j+1},j))-V(\phi(t_j,j)) \right) \nonumber
  \\
 &=& V(\phi(t_{{{j^*}}+1},{j^*})) +V(\phi(t_{1},1))\nonumber
\\&&
+\sum_{j=1}^{{j^*}-1}\left( V(\phi(t_{j+1},j+1)) \right) \nonumber 
-\sum_{j=1}^{{j^*}} \left(V(\phi(t_j,j)) \right) \nonumber
\\ &=&
V(\phi(t_{{{j^*}}+1},{j^*})) 
\end{eqnarray*}
then we have
{
\begin{eqnarray*}\nonumber
V(\xi)
&\geq&\sum_{j=0}^{{j^*}} \int_{t_{j}}^{t_{j+1}} L_C(\phi(t,j),u_{C}(t,j))dt
\\&&
+ \sum_{j=0}^{{j^*}-1}L_D(\phi(t_{j+1},j),u_{D}(t_{j+1},j))
\\&&+V(\phi(t_{{j^*}+1},{j^*}))
\nonumber
\end{eqnarray*}}
By taking the limit when {$t_{{j^*}+1}+{j^*} \rightarrow \infty$, 
we establish \eqref{eq:costboundVxi}.
Notice that if $j^*=0$, the solution $(\phi,u)$ is continuous and \eqref{eq:costboundVxi} reduces to 
{
\begin{equation*}\nonumber
V(\xi)
\geq \underset{t^* \rightarrow \infty}{\limsup}\int_{t_{0}}^{t^*} L_C(\phi(t,0),u_{C}(t,0))dt
+V(\phi(t^*,0)).
\nonumber
\end{equation*}}
On the other hand, if $t_{j^*+1}=0$ for all $j^*$, {the solution $(\phi,u)$ is discrete} and 
\eqref{eq:costboundVxi} reduces to 
{
\begin{equation*}\nonumber
V(\xi)
\geq
\underset{j^* \rightarrow \infty}{\limsup}
 \sum_{j=0}^{{j^*}-1}L_D(\phi(0,j),u_{D}(0,j))
+V(\phi(0,{j^*})).
\nonumber
\end{equation*}}}
%
\eop\smallskip\vskip 3 pt
}{}
\IfIncd{
\begin{lemma}{Solutions to $\HS_{\textup{max}}$}\label{Pp:sol2Hmaxsol2Hk}
  Any solution to $\HS_{\textup{max}}$ as in \eqref{OptimalSystProof} is a solution to $\HS\IfIncd{_s}{} $ as in \IfIncd{\eqref{Hinc}}{\eqref{Heq}}.
\end{lemma}
}{}
\NotConf{\par\noindent\textbf{Proof of Corollary \ref{HREqPJ}.}
We show that when conditions (\ref{DiffRiccatizsih})-(\ref{Riccatizsih}) hold, by using the result in Theorem \ref{thHJBszih}, the value function is equal to the function $V$ and with 
the feedback law with values as in (\ref{NashkCLQzsih}) and (\ref{NashkDLQzsih}), such a cost is attained. 
We can write (\ref{HJBzsih}) in Theorem \ref{thHJBszih} as
\begin{eqnarray}
&&0=
\underset{u_C=(u_{C1},u_{C2}) \in \Pi_u^C(x)}
{\min_{u_{C1}} \max_{u_{C2}}}  \mathcal{L}_C(x,u_C),
\nonumber
\\ 
&&\mathcal{L}_C(x,u_C)= 
x_p^\top Q_C x_p+ u_{C1}^\top R_{C1} u_{C1}+ u_{C2}^\top R_{C2} u_{C2}
+2x_p^\top P(\tau) (A_C x_p+B_C u_C)  
+ x_p^\top \frac{{d} }{{d} \tau}P(\tau) x_p\nonumber \\
\label{HJBzsihlq}
\end{eqnarray}
First, 
given that (\ref{DiffRiccatizsih}) holds, and $x_p^\top(P(\tau) A_C+A_C^\top P(\tau))x_p=x_p^\top (2P(\tau)A_C)x_p$ for every $x\in C$, one has
\begin{equation}
\mathcal{L}_C(x,u_C)=
x_p^\top P(\tau)(  B_{C2} R_{C2}^{-1}B_{C2}^\top+B_{C1} R_{C1}^{-1}B_{C1}^\top)P(\tau)  x_p+ u_{C1}^\top R_{C1} u_{C1}+ u_{C2}^\top R_{C2} u_{C2}
+ 2 x_p^\top  P(\tau) B_C u_C 
\end{equation}
The first order necessary conditions for optimality
\begin{equation*}
\frac{\partial}{\partial u_{C1}} \left( 
x_p^\top P(\tau)(  B_{C2} R_{C2}^{-1}B_{C2}^\top+B_{C1} R_{C1}^{-1}B_{C1}^\top)P(\tau)  x_p\left.+ u_{C1}^\top R_{C1} u_{C1}+ u_{C2}^\top R_{C2} u_{C2}
+ 2 x_p^\top  P(\tau)  (B_{C1} u_{C1} + B_{C2} u_{C2}) \right)\right|_{u_{C}^*} =0
\end{equation*}
\begin{equation*}
\frac{\partial}{\partial u_{C2} }\left( 
x_p^\top P(\tau)(  B_{C2} R_{C2}^{-1}B_{C2}^\top+B_{C1} R_{C1}^{-1}B_{C1}^\top)P(\tau)  x_p\left.+ u_{C1}^\top R_{C1} u_{C1}+ u_{C2}^\top R_{C2} u_{C2}
+ 2 x_p^\top  P(\tau) (B_{C1} u_{C1} + B_{C2} u_{C2}) \right)\right|_{u_{C}^*} =0
\end{equation*}
are satisfied by the point $u_C^*=(u_{C1}^*,u_{C2}^*)$, with values for each $x_p\in \Pi(C)$
\begin{equation}
u_{C1}^*=-R_{C1}^{-1}B_{C1}^\top P(\tau) x_p
\label{zLQHPu1*}
\end{equation}
\begin{equation}
u_{C2}^*=-R_{C2}^{-1}B_{C2}^\top P(\tau) x_p
\label{zLQHPu2*}
\end{equation}
Given that $R_{C1},-R_{C2} \in \mathbb{S}^{m_D}_+$, the second-order sufficient conditions for optimality 
\begin{equation*}
\frac{\partial^2}{\partial u_{C1}^2} \left( 
x_p^\top P(\tau)(  B_{C2} R_{C2}^{-1}B_{C2}^\top+B_{C1} R_{C1}^{-1}B_{C1}^\top)P(\tau)  x_p\left.+ u_{C1}^\top R_{C1} u_{C1}+ u_{C2}^\top R_{C2} u_{C2}
+ 2 x_p^\top  P(\tau)  (B_{C1} u_{C1} + B_{C2} u_{C2}) \right)\right|_{u_{C1}^*} \succeq 0
\end{equation*}
\begin{equation*}
\frac{\partial^2}{\partial u_{C2}^2 }\left( 
x_p^\top P(\tau)(  B_{C2} R_{C2}^{-1}B_{C2}^\top+B_{C1} R_{C1}^{-1}B_{C1}^\top)P(\tau)  x_p\left.+ u_{C1}^\top R_{C1} u_{C1}+ u_{C2}^\top R_{C2} u_{C2}
+ 2 x_p^\top  P(\tau) (B_{C1} u_{C1} + B_{C2} u_{C2}) \right)\right|_{u_{C2}^*} \preceq 0
\end{equation*}
hold, rendering $u_{C}^*$ as in (\ref{zLQHPu1*}) and (\ref{zLQHPu2*}) as an optimizer of the min-max problem in (\ref{HJBzsihlq}).
In addition, 
it satisfies $\mathcal{L}_C(x,u_C^*)=0$, making $V(x)=x_p^\top P(\tau)x_p$ a solution of (\ref{HJBzsih}) in Theorem \ref{thHJBszih}. 
On the other hand, we can write (\ref{Bellmanzsih}) in Theorem \ref{thHJBszih} as
\begin{eqnarray}
&&x_p^\top P(\bar{T}) x_p =\underset{u_D=(u_{D1},u_{D2}) \in \Pi_u^D(x)}
{\min_{u_{D1}} \max_{u_{D2}}}  \mathcal{L}_D(x,u_D),
\nonumber
\\ 
&&\mathcal{L}_D(x,u_D)= 
x_p^\top Q_D x_p+ u_{D1}^\top R_{D1} u_{D1}+ u_{D2}^\top R_{D2} u_{D2}+ (A_D x_p+B_D u_D)^\top P(0) (A_D x_p+B_D u_D) 
\label{Bellmanzsihlq}
\end{eqnarray}
which can be expanded as
\begin{multline}
\mathcal{L}_D(x,u_D)=  
x_p^\top (Q_D+A_D^\top P(0) A_D) x_p + 2 x_p^\top A_D^\top  P(0) B_D u_D \\
+ u_{D1}^\top (R_{D1} +B_{D1}^\top P(0) B_{D1}) u_{D1}
+ u_{D2}^\top (R_{D2} +B_{D2}^\top P(0) B_{D2}) u_{D2}
+ u_{D1}^\top (B_{D1}^\top P(0) B_{D2}) u_{D2}
+ u_{D2}^\top (B_{D2}^\top P(0) B_{D1}) u_{D1}
\end{multline}
The first order necessary conditions for optimality
\begin{equation*}
\frac{\partial}{\partial u_{D1}}\left.
\mathcal{L}_D(x,u_D)  \right|_{u_{D}^*} =0
\end{equation*}
\begin{equation*}
\frac{\partial}{\partial u_{D2}}\left.
\mathcal{L}_D(x,u_D) \right|_{u_{D}^*} =0
\end{equation*}
are satisfied by the point $u_D^*=(u_{D1}^*,u_{D2}^*)$, such that for each $x_p \in \Pi(D)$
%
\begin{equation}
u_{D}^*=-\begin{bmatrix}
R_{D1}+B_{D1}^\top P(0) B_{D1}
&
B_{D1}^\top P(0) B_{D2}
\\
B_{D2}^\top P(0) B_{D1}
&
R_{D2}+B_{D2}^\top P(0) B_{D2}
\end{bmatrix} 
^{-1}
\begin{bmatrix}
B_{D1}^\top P(0) A_D
\\
B_{D2}^\top P(0) A_D
\end{bmatrix} x_p
\label{zLQHPud*}
\end{equation}
Given that (\ref{zlqeqinv}) holds, the second-order sufficient conditions for optimality 
\begin{equation*}
\frac{\partial^2}{\partial u_{D1}^2}\left. 
\mathcal{L}_D(x,u_D) \right|_{u_D^*} \succeq 0,
\end{equation*}
\begin{equation*}
\frac{\partial^2}{\partial u_{D2}^2}\left.
\mathcal{L}_D(x,u_D)\right|_{u_D^*} \preceq 0,
\end{equation*}
are satisfied, rendering $u_D^*$ as in (\ref{zLQHPud*}) as an optimizer of the min-max problem in (\ref{Bellmanzsihlq}).
In addition, $u_D^*$
satisfies $\mathcal{L}_D(x,u_D^*)=x_p^\top P(\bar{T})x_p$, with $P(\bar{T})$ as in (\ref{Riccatizsih}), making $V(x)=x_p^\top P(\tau)x_p$ a solution of (\ref{Bellmanzsih}) in Theorem \ref{thHJBszih}. 
Then, given that $V$ is continuously differentiable on a neighborhood of $\Pi(C)$ and Assumption \ref{AssLipsZ} holds, by applying Theorem \ref{thHJBszih}, in particular from (\ref{ResultValue}), for every $ \xi=(\xi_p,\xi_\tau) \in \Pi(\overline{C}) \cup \Pi(D)$ the value function is $\J^*(\xi)=\J(\xi,((u_{C1}^*,u_{D1}^*),(u_{C2}^*,u_{D2}^*))= \xi_p^\top P(\xi_\tau)\xi_p$. From (\ref{kHJBeqzsihc}) and (\ref{kBeqzsihc}) the feedback law $\kappa=(\kappa_C,\kappa_D)$ with values as in (\ref{NashkCLQzsih}) and (\ref{NashkDLQzsih}) is a pure strategy saddle-point equilibrium.}
%
%
%
    \par\noindent\textbf{Proof of Lemma \ref{Lemma:EquivCond}.}
$(\Rightarrow)$
From (\ref{kHJBeqzsihc}) and (\ref{kBeqzsihc}) we have
{
\begin{equation}
\begin{split} 
{ \underset{u_{C1}}{\min}\> \underset{u_{C2}}{\max}}_
{\mathclap{\substack{\\ \\ {u_C=(u_{C1},u_{C2}) \in \Pi_u^C(x)}} }}
\
\IfAutom{\mathcal{L}_C(x,u_C)}{\left\{ 
L_C(x,u_C) \right.
+  \left.\left\langle \nabla V(x),F(x,u_C) \right\rangle \right\} 
=
\hspace{0.6cm}
\\
\hspace{0.6cm}}
\IfAutom{=\mathcal{L}_C(x,\kappa_C(x))}{
L_C(x,\kappa_C(x))
+\left\langle \nabla V(x) ,F(x,\kappa_C(x)) \right\rangle}  \hspace{0.6cm} \forall x \in \Pi(C)
\end{split} 
\label{kHJBeqihzproof}
\end{equation}}
and
{
\begin{equation}
\begin{split}
{ \underset{u_{D1}}{\min}\> \underset{u_{D2}}{\max}}_ 
{\mathclap{\substack{\\ \\ {u_D=(u_{D1},u_{D2}) \in \Pi_u^D(x)}} }}\
\IfAutom{\mathcal{L}_D(x,u_D)}{\left\{ 
L_D(x,u_D)
\right. + \left. V(G(x,\kappa_D(x))) \right\}
=
\hspace{0.4cm}
\\
\hspace{0.6cm}}
\IfAutom{=\mathcal{L}_D(x,\kappa_D(x))}{L_D(x,\kappa_D(x))
+ V(G(x,\kappa_D(x)))} \hspace{0.6cm} \forall x \in \Pi(D)
\end{split} 
\label{kBeqihzproof}
\end{equation}}
Thus, (\ref{HJBzsih}) and (\ref{kHJBeqihzproof}) imply 
\begin{equation}
  \IfAutom{\mathcal{L}_C(x,\kappa_C(x))}{
    L_C(x,\kappa_C(x))
    +\left\langle \nabla V(x) ,F(x,\kappa_C(x)) \right\rangle}=0  \hspace{1cm} \forall x \in \Pi(C), 
\label{kHJBeqihzproofa}
\end{equation}
while (\ref{Bellmanzsih})  and (\ref{kBeqihzproof}) imply  
\begin{equation}
  \IfAutom{\mathcal{L}_D(x,\kappa_D(x))}{L_D(x,\kappa_D(x))
+ V(G(x,\kappa_D(x)))}
 =V(x) \hspace{1cm} \forall x \in \Pi(D).
\label{kBeqihzproofa}
\end{equation}
From (\ref{kHJBeqihzproofa}) and (\ref{HJBzsih}), we have 
{
\begin{equation*}
\hspace{0.5cm}
\begin{split}
{ \underset{u_{C1}}{\min} 
}_
{\mathclap{\substack{\\ \\ {u_{C1}:(u_{C1},\kappa_{C2}(x)) \in \Pi_u^C(x)}} }}\
 \IfAutom{\mathcal{L}_C(x,(u_{C1},\kappa_{C2}(x)))}{\left\{ 
L_C(x,(u_{C1},\kappa_{C2}(x))) \right.
\hspace{0.4cm}
\\
\hspace{0.35cm}
+  \left.\left\langle \nabla V(x),F(x,(u_{C1},\kappa_{C2}(x))) \right\rangle \right\}
}
\geq 0 \hspace{0.4cm}\forall x \in \Pi(C)
\end{split} 
\label{kHJBeqihzproofb}
\end{equation*}}
and
{
\begin{equation*}
\hspace{0.5cm}
\begin{split}
{ 
\underset{u_{C2}}{\max}
}_
{\mathclap{\substack{\\ \\ {u_{C2}:(\kappa_{C1}(x), u_{C2}) \in \Pi_u^C(x)}} }}\
\IfAutom{\mathcal{L}_C(x,(\kappa_{C1}(x), u_{C2}))}{
\left\{ 
L_C(x,(\kappa_{C1}(x), u_{C2})) \right.
\hspace{0.4cm}
\\
\hspace{0.35cm}
+  \left.\left\langle \nabla V(x),F(x,(\kappa_{C1}(x), u_{C2})) \right\rangle \right\}
}
 \leq 0  \hspace{0.4cm} \forall x \in \Pi(C)
\end{split} 
\label{kHJBeqihzproofc}
\end{equation*}
}
which imply (\ref{kHJBeqihzb}) and (\ref{kHJBeqihzc}), respectively. Likewise, 
{from} (\ref{kBeqihzproofa}) and (\ref{Bellmanzsih}), we have 
{
\begin{equation*}
\hspace{0.6cm}
\begin{split}
&
{\underset{u_{D1}}{\min}}_
{\mathclap{\substack{\\ \\ {u_{D1}:(u_{D1},\kappa_{D2}(x)) \in \Pi_u^D(x)}} }}\
 \IfAutom{\mathcal{L}_D(x,(u_{D1},\kappa_{D2}(x)))
 }{
 \left\{ 
L_D(x,(u_{D1},\kappa_{D2}(x))) \right.
\hspace{0.4cm}
\\
&\hspace{0.6cm}
+  \left.    V(G(x,(u_{D1},\kappa_{D2}(x))))   \right\}
 }
\geq V(x)
\hspace{0.4cm}
\forall x \in \Pi(D)
\end{split} 
\label{kBeqihzproofb}
\end{equation*}}
and
\begin{equation*}
\hspace{0.6cm}
\begin{split}
&
{ 
\underset{u_{D2}}{\max}}_
{\mathclap{\substack{\\ \\ {u_{D2}:(\kappa_{D1}(x), u_{D2}) \in \Pi_u^D(x)}} }}\
\IfAutom{\mathcal{L}_D(x,(\kappa_{D1}(x), u_{D2}))  
}{
\left\{ 
L_D(x,(\kappa_{D1}(x), u_{D2})) \right.
\hspace{0.4cm}
\\
&\hspace{0.6cm}
+ \left.  V(G(x,(\kappa_{D1}(x), u_{D2}))) \right\} 
}
\leq V(x)
\hspace{0.4cm} \forall x \in \Pi(D)
\end{split} 
\label{kBeqihzproofc}
\end{equation*}
which imply (\ref{kBeqihze}) and (\ref{kBeqihzf}), respectively.

\noindent $(\Leftarrow)$ Given $V$ and $\kappa:=(\kappa_C,\kappa_D)=((\kappa_{C1},\kappa_{C2}),(\kappa_{D1},\kappa_{D2}))$ such that (\ref{kHJBeqihza})-(\ref{kBeqihzf}) are satisfied, and such that $C_\kappa = \Pi(C)$ {and} $D_\kappa = \Pi(D)$, {we show} that $V$ and $\kappa$ satisfy (\ref{HJBzsih}), (\ref{Bellmanzsih}), (\ref{kHJBeqzsihc}), and (\ref{kBeqzsihc}). From (\ref{kHJBeqihza}) and (\ref{kHJBeqihzb}) we have 
{
\begin{equation}
\begin{split}
&\hspace{0.8cm}
{ \underset{u_{C1}}{\min} 
}_
{\mathclap{\substack{\\ \\ {u_{C1}:(u_{C1},\kappa_{C2}(x)) \in \Pi_u^C(x)}} }}\
\IfAutom{\mathcal{L}_C(x,(u_{C1},\kappa_{C2}(x)))  \\ &\hspace{1.5cm}
=\mathcal{L}_C(x,\kappa_C(x))
}{
\left\{ 
L_C(x,(u_{C1},\kappa_{C2}(x))) \right.
\\&\hspace{4cm}+  \left.\left\langle \nabla V(x),F(x,(u_{C1},\kappa_{C2}(x))) \right\rangle \right\}  
\\ &=L_C(x,\kappa_C(x))
+\left\langle \nabla V(x) ,F(x,\kappa_C(x)) \right\rangle} =0 \hspace{0.4cm} \forall x \in \Pi(C)
\end{split}
\label{kHJBeqihzproof1left}
\end{equation}}
and from (\ref{kHJBeqihza}) and (\ref{kHJBeqihzc}) we have 
{
\begin{equation}
\begin{split}
&\underset{u_{C2}:(\kappa_{C1}(x), u_{C2}) \in \Pi_u^C(x)}
{ 
\underset{u_{C2}}{\max}
} 
\IfAutom{\mathcal{L}_C(x,(\kappa_{C1}(x), u_{C2}))  \\ &\hspace{1.5cm}
=\mathcal{L}_C(x,\kappa_C(x))
}{
\left\{ 
L_C(x,(\kappa_{C1}(x), u_{C2})) \right.
\\&\hspace{4cm}+  \left.\left\langle \nabla V(x),F(x,(\kappa_{C1}(x), u_{C2})) \right\rangle \right\}  
\\&=L_C(x,\kappa_C(x))
+\left\langle \nabla V(x) ,F(x,\kappa_C(x)) \right\rangle
}
=0 \hspace{0.4cm} \forall x \in \Pi(C)
\end{split}
\label{kHJBeqihzproof2left}
\end{equation}}
Thus, (\ref{kHJBeqihzproof1left}) and (\ref{kHJBeqihzproof2left}) imply (\ref{HJBzsih}) and  (\ref{kHJBeqzsihc}).
Similarly, from (\ref{kBeqihzd}) and (\ref{kBeqihze}) we have 
{
\begin{equation}
\begin{split}
&\underset{u_{D1}:(u_{D1},\kappa_{D2}(x)) \in \Pi_u^D(x)}
{ \underset{u_{D1}}{\min} 
} 
\IfAutom{\mathcal{L}_D(x,(u_{D1},\kappa_{D2}(x)))  \\ &\hspace{1cm}
=\mathcal{L}_D(x,\kappa_D(x))
}{
\left\{ 
L_D(x,(u_{D1},\kappa_{D2}(x))) \right.
\\&\hspace{4.7cm}+  \left.    V(G(x,(u_{D1},\kappa_{D2}(x))))   \right\} 
\\&=L_D(x,\kappa_D(x))
+ V(G(x,\kappa_D(x)))
}
=V(x) \hspace{0.4cm} \forall x \in \Pi(D)
\end{split}
\label{kBeqihzproof1left}
\end{equation}}
and from (\ref{kBeqihzd}) and (\ref{kBeqihzf}) we have 
{
\begin{equation}
\begin{split}
&\underset{u_{D2}:(\kappa_{D1}(x), u_{D2}) \in \Pi_u^D(x)}
{ 
\underset{u_{D2}}{\max}
}
\IfAutom{\mathcal{L}_D(x,(\kappa_{D1}(x), u_{D2}))  \\ &\hspace{1cm}
=\mathcal{L}_D(x,\kappa_D(x))
}{
\left\{ 
L_D(x,(\kappa_{D1}(x), u_{D2})) \right.
\\&\hspace{4.8cm}+ \left.  V(G(x,(\kappa_{D1}(x), u_{D2}))) \right\}
\\& =L_D(x,\kappa_D(x))
+ V(G(x,\kappa_D(x)))
}
=V(x) \hspace{0.4cm} \forall x \in \Pi(D)
\end{split}
\label{kBeqihzproof2left}
\end{equation}}
Thus, (\ref{kBeqihzproof1left}) and (\ref{kBeqihzproof2left}) imply (\ref{Bellmanzsih}) and  (\ref{kBeqzsihc}).

\eop\smallskip\vskip 3 pt
%
\NotConf{\par\noindent\textbf{Proof of Corollary \ref{zCorStability}.}
Since by assumption we have that $C_\kappa = \Pi(C), D_\kappa = \Pi(D)$, and $V,\kappa:=(\kappa_C,\kappa_D)=((\kappa_{C1},\kappa_{C2}),(\kappa_{D1},\kappa_{D2}))$ are such that (\ref{kHJBeqihza})-(\ref{kBeqihzf}) hold, then, thanks to Lemma \ref{Lemma:EquivCond}, $V$ and $\kappa$ satisfy (\ref{HJBzsih}), (\ref{Bellmanzsih}), (\ref{kHJBeqzsihc}), and (\ref{kBeqzsihc}). 
Since in addition, for each $\xi \in \overline{C_\kappa} \cup D_\kappa$, each $\phi \in \mathcal{S}_{\HS_{\kappa}  }^\infty(\xi)$ satisfies (\ref{TerminalCondCG}), we have from Theorem \ref{thHJBszih} that $V$ is the value function as in (\ref{cost2gozsih}) for $\HS_{\kappa}$ at $\overline{C_\kappa} \cup D_\kappa$ 
and the feedback law
$\kappa$ with values (\ref{kHJBeqzsihc}), (\ref{kBeqzsihc}) is the pure strategy Nash equilibrium for this game. 
\pnn{Given that maximal solutions to $\HS_\kappa$ are complete by assumption, $G(D_\kappa)\subset \overline{C_\kappa} \cup D_\kappa$. Then, $V$ is a Lyapunov candidate for $\HS_{\kappa}$ \cite[Definition 3.16]{65} since $\overline{C_\kappa} \cup D_\kappa \subset \dom V=\reals^n$ and $V$ is continuously differentiable on an open set containing $\overline{C_\kappa}$}. \pnn{From (\ref{kHJBeqihza}), (\ref{kBeqihzd}), we have
\begin{eqnarray}
\left\langle \nabla V(x),F(x,\kappa_C(x)) \right\rangle&\leq&- L_C(x,\kappa_C(x)) \hspace{0.5cm} \forall x \in C_\kappa, \label{kHJBeqihzabound}\\
V(G(x,\kappa_D(x)))-V(x)&\leq&-L_D(x,\kappa_D(x)) \hspace{0.5cm} \forall x \in D_\kappa. 
\label{kBeqihzbound}
\end{eqnarray}
Moreover, one of the following holds
\begin{enumerate}[label=\arabic*)]
\item 
$L_C \in \mathcal{PD}_{\kappa_C}(\A)$ and $L_D \in \mathcal{PD}_{\kappa_D}(\A )$;

\item 
$L_D \in \mathcal{PD}_{\kappa_C}(\A )$ and given that $L_C(x,\kappa_C(x)) \geq \eta(|x|_\A)$ for all $x \in C_\kappa$, and from (\ref{kHJBeqihzabound}) we have
\begin{eqnarray*}
\left\langle \nabla V(x),F(x,\kappa_C(x)) \right\rangle&\leq&- \eta(|x|_\A) \hspace{1.5cm} \forall x \in C_\kappa; 
\end{eqnarray*}

\item 
$L_C \in \mathcal{PD}_{\kappa_C}(\A )$ and given that $L_D(x,\kappa_D(x)) \geq \eta(|x|_\A)$ for all $x \in D_\kappa$, and from (\ref{kBeqihzbound}) we have
\begin{eqnarray*}
V(G(x,\kappa_D(x)))-V(x)&\leq&-\eta(|x|_\A) \hspace{1.5cm} \forall x \in D_\kappa; 
\end{eqnarray*}
\end{enumerate}}
\noindent Thus,  given the functions $\alpha_1, \alpha_2$ satisfying (\ref{zalphabound}), from \cite[Theorem 3.18]{65} we have that $\A$ is uniformly globally asymptotically stable for $\HS_{\kappa}$. \pnn{Notice that the role of $\rho$ therein is played by the stage costs $L_C \in \mathcal{PD}_{\kappa_C}(\A
)$, $L_D \in \mathcal{PD}_{\kappa_D}(\A 
)$ or the function $\eta \in \mathcal{PD}$ according to the corresponding case.}
Furthermore, via Assumption \ref{AssLipsZ}, maximal solutions to $\HS$ are complete, which allows arguing uniform non-pre-asymptotic stability of $\A$.}
\NotConf{\par\noindent\textbf{Proof of Corollary \ref{HREqRobust}.}
We show that when conditions (\ref{DiffRiccatizsihsec})-(\ref{Riccatizsihsec}) hold, by using the result in Theorem \ref{thHJBszih},  the value function is equal to the function $V$ and under 
the feedback law as in (\ref{NashkCLQzsihrlqr}) and (\ref{NashkDLQzsihrlqr}) such a cost is attained in the presence of the maximizing disturbance given by (\ref{NashkC2LQzsihrlqr}) and (\ref{NashkD2LQzsihrlqr}). 
We can write (\ref{HJBzsih}) in Theorem \ref{thHJBszih} as
\begin{eqnarray}
&&0=
\underset{u_C=(u_{C1},u_{C2}) \in \Pi_u^C(x)}
{\min_{u_{C1}} \max_{u_{C2}}}  \mathcal{L}_C(x,u_C),
\nonumber
\\ 
&&\mathcal{L}_C(x,u_C)= 
x^\top Q_C x+ u_{C1}^\top R_{C1} u_{C1}+ u_{C2}^\top R_{C2} u_{C2}
+2x^\top P (A_C x+B_C u_C)  
\nonumber \\
\label{HJBzsihlqrlqr}
\end{eqnarray}
First, 
given that (\ref{DiffRiccatizsihrlqr}) holds, and $x^\top(P A_C+A_C^\top P)x=x^\top (2PA_C)x$ for every $x\in \Pi(C)$, one has
\begin{equation}
\mathcal{L}_C(x,u_C)=
x^\top P(  B_{C2} R_{C2}^{-1}B_{C2}^\top+B_{C1} R_{C1}^{-1}B_{C1}^\top) P x+ u_{C1}^\top R_{C1} u_{C1}+ u_{C2}^\top R_{C2} u_{C2}
+ 2 x^\top  P B_C u_C
\end{equation}
The first order necessary conditions for optimality
\begin{equation*}
\frac{\partial}{\partial u_{C1}} \mathcal{L}_C(x,u_C)
\left.\right|_{u_{C}^*} 
=0
\end{equation*}
\begin{equation*}
\frac{\partial}{\partial u_{C2} } \mathcal{L}_C(x,u_C)
\left.\right|_{u_{C}^*} 
=0
\end{equation*}
are satisfied by the point $u_C^*=(u_{C1}^*,u_{C2}^*)$, with values at any $x \in \Pi(C)$
\begin{equation}
u_{C1}^*=-R_{C1}^{-1}B_{C1}^\top P x
\label{zLQHPu1*rlqr}
\end{equation}
\begin{equation}
u_{C2}^*=-R_{C2}^{-1}B_{C2}^\top P x
\label{zLQHPu2*rlqr}
\end{equation}
Given that $R_{C1},-R_{C2} \in \mathbb{S}^{m_D}_+$, the second-order sufficient conditions for optimality 

\begin{equation*}
\frac{\partial^2}{\partial u_{C1}^2}  \mathcal{L}_C(x,u_C)
\left.\right|_{u_{C}^*}
\succeq 0
\end{equation*}
\begin{equation*}
\frac{\partial^2}{\partial u_{C2}^2 } \mathcal{L}_C(x,u_C)
\left.\right|_{u_{C}^*} 
\preceq 0
\end{equation*}
hold, rendering $u_{C}^*$ with values as in (\ref{zLQHPu1*rlqr}) and (\ref{zLQHPu2*rlqr}) for each $x\in \Pi(C)$ as an optimizer of the min-max problem in (\ref{HJBzsihlqrlqr}).
In addition, 
it satisfies $\mathcal{L}_C(x,u_C^*)=0$, making $V(x)=x^\top P x$ a solution to (\ref{HJBzsih}) in Theorem \ref{thHJBszih}. 

On the other hand, we can write (\ref{Bellmanzsih}) in Theorem \ref{thHJBszih} as
\begin{eqnarray}
&&x^\top P x =\underset{u_D=(u_{D1},u_{D2}) \in \Pi_u^D(x)}
{\min_{u_{D1}} \max_{u_{D2}}}  \mathcal{L}_D(x,u_D),
\nonumber
\\ 
&&\mathcal{L}_D(x,u_D)= 
x^\top Q_D x+ u_{D1}^\top R_{D1} u_{D1}+ u_{D2}^\top R_{D2} u_{D2}+ (A_D x+B_D u_D)^\top P (A_D x+B_D u_D) 
\label{Bellmanzsihlqrlqr}
\end{eqnarray}
which can be expanded as 
\begin{multline}
\mathcal{L}_D(x,u_D)=  
x^\top (Q_D+A_D^\top P A_D) x + 2 x^\top A_D^\top  P B_D u_D \\
+ u_{D1}^\top (R_{D1} +B_{D1}^\top P B_{D1}) u_{D1}
+ u_{D2}^\top (R_{D2} +B_{D2}^\top P B_{D2}) u_{D2}
+ u_{D1}^\top (B_{D1}^\top P B_{D2}) u_{D2}
+ u_{D2}^\top (B_{D2}^\top P B_{D1}) u_{D1}
\end{multline}
The first order necessary conditions for optimality
\begin{equation*}
\frac{\partial}{\partial u_{D1}}\left.
\mathcal{L}_D(x,u_D)  \right|_{u_{D}^*} 
=0
\end{equation*}
\begin{equation*}
\frac{\partial}{\partial u_{D2}}\left.
\mathcal{L}_D(x,u_D) \right|_{u_{D}^*} 
=0
\end{equation*}
are satisfied by the point $u_D^*=(u_{D1}^*,u_{D2}^*)$, with values for each $x \in \Pi(D)$
\begin{equation}
u_{D}^*=-\begin{bmatrix}
R_{D1}+B_{D1}^\top P B_{D1}
&
B_{D1}^\top P B_{D2}
\\
B_{D2}^\top P B_{D1}
&
R_{D2}+B_{D2}^\top P B_{D2}
\end{bmatrix} 
^{-1}
\begin{bmatrix}
B_{D1}^\top P A_D
\\
B_{D2}^\top P A_D
\end{bmatrix} x_p
\label{zLQHPud*rlqr}
\end{equation}
Given that (\ref{zlqeqinvrlqr}) holds, the second-order sufficient conditions for optimality of $u_D^*$, namely
\begin{equation*}
\frac{\partial^2}{\partial u_{D1}^2}\left. 
\mathcal{L}_D(x,u_D) \right|_{u_D^*} \succeq 0,
\end{equation*}
\begin{equation*}
\frac{\partial^2}{\partial u_{D2}^2}\left.
\mathcal{L}_D(x,u_D)\right|_{u_D^*} \preceq 0,
\end{equation*}
are satisfied, rendering $u_D^*$ in (\ref{zLQHPud*rlqr}) as an optimizer of the min-max problem in (\ref{Bellmanzsihlqrlqr}).
In addition, $u_D^*$
satisfies $\mathcal{L}_D(x,u_D^*)=x^\top P x$, with $P$ as in (\ref{Riccatizsihrlqr}), making $V(x)=x^\top P x$ a solution to (\ref{Bellmanzsih}) in Theorem \ref{thHJBszih}. 

Thus, given that $V$ is continuously differentiable in $\reals^n$  and Assumption \ref{AssLipsZ} holds,
by applying Theorem \ref{thHJBszih}, in particular from (\ref{ResultValue}), \pnn{for every $ \xi \in \Pi(\overline{C}) \cup \Pi(D)$ the value function is $\J^*(\xi)=\J(\xi,((u_{C1}^*,u_{D1}^*),(u_{C2}^*,u_{D2}^*))= \xi^\top P \xi$. From (\ref{kHJBeqzsihc}) and (\ref{kBeqzsihc}), when $P_1$ plays $u_1^*$ defined by $\kappa_1=(\kappa_{C1},\kappa_{D1})$ with values as in (\ref{NashkCLQzsihrlqr}), (\ref{NashkDLQzsihrlqr}), and $P_2$ plays any disturbance $u_2$ such that solutions to $\HS$ with data as in (\ref{lrlqr}) are complete, then the cost is upper bounded by $\J(\xi,u^*)$,  
satisfying (\ref{SaddlePointIneq}).
}
}
\NotConf{\par\noindent\textbf{Proof of Corollary \ref{HREqSec}.}
We show that when conditions (\ref{DiffRiccatizsihsec})-(\ref{Riccatizsihsec}) hold, by using the result in Theorem \ref{thHJBszih},  the value function is equal to the function $V$ and under 
the feedback law as in (\ref{NashkDLQzsihsec}) such a cost is attained in the presence of the maximizing attack given by (\ref{NashkD2LQzsihsec}). 
We can write (\ref{HJBzsih}) in Theorem \ref{thHJBszih} as
\begin{equation*}
0=
2x^\top P F(x) \hspace{1cm} \forall x \in \Pi(C)
\label{HJBzsihlqsec}
\end{equation*}
which is satisfied thanks to (\ref{DiffRiccatizsihsec}).
On the other hand, we can write (\ref{Bellmanzsih}) in Theorem \ref{thHJBszih} as
\begin{eqnarray}
&&x^\top P x =\underset{u_D=(u_{D1},u_{D2}) \in \Pi_u^D(x)}
{\min_{u_{D1}} \max_{u_{D2}}}  \mathcal{L}_D(x,u_D),
\nonumber
\\ 
&&\mathcal{L}_D(x,u_D)= 
x^\top Q_D x+ u_{D1}^\top R_{D1} u_{D1}+ u_{D2}^\top R_{D2} u_{D2}+ (A_D x+B_D u_D)^\top P (A_D x+B_D u_D) 
\label{Bellmanzsihlqsec}
\end{eqnarray}
which can be expanded as
\begin{multline}
\mathcal{L}_D(x,u_D)=  
x^\top (Q_D+A_D^\top P A_D) x + 2 x^\top A_D^\top  P B_D u_D \\
+ u_{D1}^\top (R_{D1} +B_{D1}^\top P B_{D1}) u_{D1}
+ u_{D2}^\top (R_{D2} +B_{D2}^\top P B_{D2}) u_{D2}
+ u_{D1}^\top (B_{D1}^\top P B_{D2}) u_{D2}
+ u_{D2}^\top (B_{D2}^\top P B_{D1}) u_{D1}
\end{multline}
The first order necessary conditions for optimality
\begin{equation*}
\frac{\partial}{\partial u_{D1}}\left.
\mathcal{L}_D(x,u_D)  \right|_{u_{D}^*} =0
\end{equation*}
\begin{equation*}
\frac{\partial}{\partial u_{D2}}\left.
\mathcal{L}_D(x,u_D) \right|_{u_{D}^*} =0
\end{equation*}
are satisfied by the point $u_D^*=(u_{D1}^*,u_{D2}^*)$, such that for each $x \in \Pi(D)$,
\begin{equation}
u_{D}^*=-\begin{bmatrix}
R_{D1}+B_{D1}^\top P B_{D1}
&
B_{D1}^\top P B_{D2}
\\
B_{D2}^\top P B_{D1}
&
R_{D2}+B_{D2}^\top P B_{D2}
\end{bmatrix} 
^{-1}
\begin{bmatrix}
B_{D1}^\top P A_D
\\
B_{D2}^\top P A_D
\end{bmatrix} x
\label{zLQHPud*sec}
\end{equation}
Given that (\ref{zlqeqinvsec}) holds, the second-order sufficient conditions for optimality 
\begin{equation*}
\frac{\partial^2}{\partial u_{D1}^2}\left. 
\mathcal{L}_D(x,u_D) \right|_{u_D^*} \succeq 0,
\end{equation*}
\begin{equation*}
\frac{\partial^2}{\partial u_{D2}^2}\left.
\mathcal{L}_D(x,u_D)\right|_{u_D^*} \preceq 0,
\end{equation*}
are satisfied, rendering $u_D^*$ as in (\ref{zLQHPud*sec}) as an optimizer of the min-max problem in (\ref{Bellmanzsihlqsec}).
In addition, $u_D^*$
satisfies $\mathcal{L}_D(x,u_D^*)=x^\top P x$, with $P$ as in (\ref{Riccatizsihsec}), leading $V(x)=x^\top P x$ as a solution of (\ref{Bellmanzsih}) in Theorem \ref{thHJBszih}. 

Thus, given that $V$ is continuously differentiable in $\reals^n$ and Assumption (\ref{AssLipsZ}) holds, by applying Theorem \ref{thHJBszih}, in particular from (\ref{ResultValue}), 
\pnn{for every $ \xi \in \Pi(\overline{C}) \cup \Pi(D)$ the value function is $\J^*(\xi)=\J(\xi,(u_{D1}^*,u_{D2}^*))= \xi^\top P \xi$. From (\ref{kHJBeqzsihc}) and (\ref{kBeqzsihc}), when $P_2$ plays $u_2^*$ defined by $\kappa_{D2}$  as in (\ref{NashkD2LQzsihsec}), $P_1$ minimizes the cost of complete solutions to $\HS$ 
by playing $u_1^*$ defined by $\kappa_{D1}$ as in (\ref{NashkDLQzsihsec}), attaining $ \J(\xi,u^*)$, and satisfying (\ref{SaddlePointIneq}).}
}
}
\IfAutomss{
\begin{wrapfigure}{l}{0.14\textwidth}
\includegraphics[height=1.5in,width=1in,clip,keepaspectratio]{Figures/Santiagos}
\end{wrapfigure}
\noindent {\bf Santiago J. Leudo }\
received the B.S. degrees in Electronic Engineering and Mechanical Engineering in 2014, and the M.S. degree in Electronic and Computer Engineering in 2017,  from University of Los Andes, Colombia.
Santiago was a fellow of the Control Research Group GIAP and an Instructor Professor at University of Los Andes in 2017-2018. In 2018, Santiago received the Fulbright Ph.D. Scholarship for high-level skills in Information Technologies - Robotics. 
He received the M.S. and Ph.D degrees in Electrical and Computer Engineering - Robotics and Control from the University of California, Santa Cruz, in 2021 and 2024, respectively. %
His research interests are in connections between control theory and game theory for nonlinear, hybrid, and multi-agent systems, and their applications to robotic networks, UAVs, security, and aerospace systems. 

\begin{wrapfigure}{l}{0.14\textwidth}
  \includegraphics[height=1.5in,width=1in,clip,keepaspectratio]{Figures/Ricardos}
  \end{wrapfigure}
  \noindent {\bf Ricardo G. Sanfelice}\
   received the B.S. degree in Electronics Engineering from the Universidad de Mar del Plata, Buenos Aires, Argentina, in 2001, and the M.S. and Ph.D. degrees in Electrical and Computer Engineering from the University of California, Santa Barbara, in 2004 and 2007, respectively. In 2007 and 2008, he held postdo.positions at the Laboratory for Information and Decision Systems at \IfAutomss{MIT}{the Massachusetts Institute of Technology} and at the Centre Automatique et Systèmes at the École de Mines de Paris. In 2009, he joined the faculty of the Department of Aerospace and Mechanical Engineering at the University of Arizona, Tucson, where he was an Assistant Professor. In 2014, he joined the University of California, Santa Cruz, where he is currently Professor and Chair in the Department of Electrical and Computer Engineering. Prof. Sanfelice is the recipient of the 2013 SIAM Control and Systems Theory Prize, the National Science Foundation CAREER award, the Air Force Young Investigator Research Award, the 2010 IEEE Control Systems Magazine Outstanding Paper Award, and the 2020 Test-of-Time Award from the Hybrid Systems: Computation and Control Conference. He is Associate Editor for Automatica, Communicating Editor for the Journal of Nonlinear Science, and a Fellow of the IEEE. His research interests are in modeling, stability, robust control, observer design, and simulation of nonlinear and hybrid systems with applications to power systems, aerospace, and biology.
}{}
\end{document}